\documentclass{article}
\usepackage{amsfonts}
\usepackage{amsmath, amssymb, amsthm, amscd}

\def\Box{\hbox{$\sqcap \unskip \kern -6.5pt 
\sqcup$}}
\title{The Ending Laminations Theorem direct
from Teichm\"uller Geodesics}
\author{Mary Rees}

\begin{document}

\maketitle

\begin{abstract}
A proof of the Ending Laminations Theorem is given which 
uses Teichm\"uller geodesics directly.
\end{abstract}
\newenvironment{ulemma}{\par\noindent\textbf{Lemma}\,\,\em}{\rm}
\newenvironment{utheorem}{\par\noindent\textbf{Theorem}\,\,\em}{\rm}
\newenvironment{ucorollary}{\par\noindent\textbf{Corollary}\,\,\em}{\rm}
\renewcommand{\theequation}{\thesubsection.\arabic{equation}}
\newcommand{\ssubsection}[{1}]{\subsection {#1}
\setcounter{equation}{0} }
\newtheorem{theorem}[subsection]{Theorem}
\newtheorem{lemma}[subsection]{Lemma}
\newtheorem{corollary}[subsection]{Corollary}
\newtheorem{sslemma}[subsubsection]{Lemma}

\section{Introduction}\label{1}

 Over 30 years ago, Scott proved the 
remarkable result \cite{Sco1}, \cite{Sco2} that any 
finitely-generated fundamental group of any three-manifold $N$ is 
finitely presented, and, further, that the manifold is homotopy 
equivalent to a compact three-dimensional submanifold with boundary. 
This submanifold, known as the {\em{Scott core}}, was shown 
\cite{McC-M-S} 
to be unique up to homeomorphism --- but not, in general, up to 
isotopy. 
If, in addition, $N$ is a hyperbolic manifold, 
then it is useful to consider
 the  manifold with boundary 
$N_{d}$ obtained by removing from $N$ components of the $\varepsilon 
_{0}$-thin part round cusps, for a fixed Margulis constant 
$\varepsilon _{0}$. Sullivan \cite{Sul2} proved that there are only 
finitely many cusps. We shall call $N_{d}$ 
{\em{the horoball deletion of $N$.}} Topologically, it does not 
depend on the choice of Margulis constant $\varepsilon _{0}$, 
although metrically it does. There is a relative version of Scott's 
result \cite{McC} which says that $(N_{d},\partial N_{d})$ is 
homotopy equivalent to $(N_{c},N_{c}\cap \partial N_{d})$ where 
$N_{c}\subset N_{d}$ is 
the {\em{relative Scott core}}. The relative Scott core is, again, 
unique 
up to homeomorphism \cite{McC-Mil}. The closure of each  component of 
$\partial N_{c}\setminus \partial N_{d}$ is an orientable compact
surface $S_{d}$ with boundary, and is also the boundary in $N_{d}$ of 
$U(S_{d})$, for a unique component $U(S_{d})$ of 
$N_{d}\setminus N_{c}$. The surface $S_{d}$ can be any compact 
orientable surface with boundary apart from the closed disc, annulus 
or torus. The set $U(S_{d})=U(e)$ is a neighbourhood of a 
unique end $e$ of $N_{d}$ \cite{Bon}. 
The correspondence between end $e$, and the 
surface $S_{d}(e)$ bounding 
its neighbourhood, is thus  one-to-one. We shall also write $S(e)$ 
for the 
surface without boundary obtained up to homeomorphism by attaching a 
punctured disc to each boundary component of $S_{d}(e)$.

 In 2004, following some close approaches (e.g. 
\cite{Sou}), two proofs of the Tameness conjecture were 
announced \cite{Ag}, \cite{C-G}, that $N_{c}$ can be chosen so 
that each closed end neighbouhood 
$S_{d}(e)\cup U(e)$ is 
homeomorphic 
to 
$S_{d}(e)\times [0,\infty )$, that is, each end of $N_{d}$ is 
{\em{topologically tame}}. (See also \cite{Som}.)
The topology of $N$ is thus uniquely 
determined by the topology of the pair $(N_{c},N_{c}\cap 
\partial N_{d})$. 

The history of the geometry of $N$ runs parallel to this. The 
concept of 
geometric tameness was originally developed by Thurston \cite{T}. An 
end of $N_{d}$ is {\em{geometrically finite}} if it has a 
neighbourhood disjoint from the the set of closed geodesics in $N$. 
(There are a number of equivalent definitions.)   
An 
end of $N_{d}$ is {\em{simply degenerate}} if each 
neighbourhood of 
the end contains a {\em{simple}} closed geodesic. These two 
possibilities are mutually exclusive, and if one or the other holds, 
the end is said to be {\em{geometrically tame}}. It is not clear a 
priori that one or the other possibilities must hold, but in 1986 
Bonahon \cite{Bon} 
published a proof 
that an end $e$  of $N_{d}$ is 
geometrically tame if the inclusion of $S_{d}(e)$ in $N_{c}$ is 
injective on 
fundamental groups, that is, if $S_{d}(e)$ is {\em{incompressible}} 
in 
$N_{d}$. Importantly for future developments, this made rigorous the 
end invariant suggested by Thurston. In the geometrically finite 
case, if $N=H^{3}/\pi _{1}(N) $, and $\Omega \subset \partial 
H^{3}$ is the domain of discontinuity of $\pi _{1}(N) $, then the 
closure 
of $U(e)$ in 
$(H^{3}\cup \Omega /\pi _{1}(N) $  intersects a unique 
punctured hyperbolic surface  
homeomorphic to $S(e)$, and contains a subsurface  which is isotopic 
in 
$(H^{3}\cup \Omega )/\pi _{1}(N) $ to $S_{d}(e)$. Thus, this 
surface in $\Omega /\pi _{1}(N) $ 
determines a point $\mu (e)$ in the Teichm\" uller space ${\cal 
T}(S(e))$. Bonahon's work 
dealt with the simply degenerate ends. He showed that if $S(e)$ was 
endowed with any complete hyperbolic structure,
then the Hausdorff limit of any sequence of simple closed geodesics 
exiting 
$e$ was the same arational geodesic lamination, independent, up 
to homeomorphism, of the hyperbolic metric chosen. A {\em{geodesic 
lamination}} on $S$ is a closed set of nonintersecting geodesics on 
$S$, and a lamination is {\em{arational}} if every simple closed loop 
intersects a recurrent leaf in the lamination transversally. So in 
both cases, geometrically finite and infinite, 
an invariant $\mu (e)$ is 
obtained, and if there are $n$ ends $e_{i}$, $1\leq i\leq n$, this 
gives an invariant $(\mu (e_{1}),\cdots \mu (e_{n}))$. Bonahon
also showed that geometric tameness implies topological tameness. 
Later, Canary 
\cite{Can} extended Bonahon's result about geometric tameness to any 
end of any topologically 
tame $N$. The resolutions of the Tameness Conjecture mentioned above 
means that any end of any three-dimensional hyperbolic manifold with 
finitely generated fundamental group is 
topologically tame, and hence, by the work of Bonahon and Canary, 
geometrically tame.  

Bonahon's result on an invariant $(\mu (e_{1},\cdots \mu (e_{i}))$ 
associated to $N$ for which $\partial N_{c}\setminus \partial N_{d}$ 
is incompressible, generalises to any tame hyperbolic manifold $N$. 
Let $G$ denote the group of 
orientation-preserving homeomorphisms $\varphi $ of $N_{c}$, modulo 
isotopy,  and are homotopic in $N_{c}$ to the identity. Then 
$G$ acts on $\prod _{e}{\cal T}(S(e))\cup {\cal {GL}}(S(e))$, where 
the 
product is over ends $e$ of $N_{d}$ and  ${\cal{GL}}(S(e))$ is the 
space  of geodesic laminations on $S(e)$. We write $[\mu _{1},\cdots 
\mu _{n}]$ to denote an element of the quotient space. Note that $G$ 
is trivial if each $S_{d}(e)$ is incompressible in 
$N_{d}$. 
We also write ${\cal {GL}}_{a}(S)$ for the space of arational 
laminations 
in ${\cal{GL}}(S)$. The ending invariant of a tame hyperbolic 
manifold is now a point
$$[\mu (e_{1}),\cdots \mu (e_{n})]\in \prod 
_{i=1}^{n}({\cal{T}}(S(e_{i}))\cup {\cal{O}}_{a}(S(e_{i}),N))/G$$
Here, ${\cal O}(S,N)\subset {\cal {GL}}(S)$ is the {\em{Masur 
domain}} 
as defined by Otal 
\cite{Ot1},  since the concept arose in the case of handlebodies in 
work by H. Masur \cite{Mas}. The precise definition will be given 
later, 
but ${\cal O}(S,N)$ is open 
in
${\cal {GL}}(S)$, invariant under the action of the group 
$G$ (as above), contains no closed geodesics in $S$ which are trivial 
in $N$, and ${\cal 
O}(S,N)=\cal {GL}(S)$ if $(S_{d},\partial S_{d})$ is incompressible 
in 
$(N_{d},\partial N_{d})$. As for ${\cal{GL}}(S)$, 
${\cal{O}}_{a}(S,N)$ 
denotes the set of arational laminations in ${\cal{O}}(S,N)$.

The obvious questions to ask about the invariants $[\mu 
(e_{1}),\cdots \mu (e_{n})]$ of hyperbolic manifolds of a fixed 
topological type are: does the invariant $[\mu (e_{1},\cdots \mu 
(e_{n})]$ of a hyperbolic manifold uniquely determine that manifold 
up to isometry, and which invariants can occur. 

There is a natural compactification of ${\cal T}(S)$, or, more 
generally, of \\ $\prod _{e}{\cal T}(S(e))/G$, in which $\prod 
_{e}{\cal 
T}(S(e))/G$ is the 
interior and $\prod _{e}{\cal O}_{a}(S(e),N)/G$ is contained in the 
boundary. The 
topology 
at the boundary will be specified later. The boundary is actually 
larger than $\prod _{e}{\cal O}_{a}(S(e),N)/G$. It also includes 
countably many {\em{split boundary 
pieces}}, one for each multicurve  $\Gamma $  in the Masur domain. 
(Multicurves are defined in \ref{2.0}, and the Masur domain in 
\ref{3.7}.) 
Given a multicurve $\Gamma _{i}$ on
$S(e_{i})$ for each end $e_{i}$ of $N$, we define a 
topological manifold 
with boundary $N_{d}(\Gamma _{1},\cdots \Gamma _{n})$ which 
contains $N_{c}$ and is contained in $N$, as follows. We can assume 
that $\gamma \subset S_{d}(e_{i})$ for all $\gamma \in \Gamma _{i}$. 
Fix a 
homeomorphism $\Phi _{i}$ from $S(e_{i})\times [0,\infty )$ to the 
closed neighbourhood of $e_{i}$ in $N$ bounded by $S(e_{i})$. Choose 
a collection of disjoint open
annulus neighbourhoods $A(\gamma )\subset S_{d}(e_{i})$ of the loops 
$\gamma \in \Gamma _{i}$. Then
$$N_{d}(\Gamma _{1},\cdots \Gamma _{n})=N_{d}\setminus \cup 
_{i=1}^{n}\Phi _{i}(\cup _{\gamma \in \Gamma _{i}}A(\gamma )\times 
(0,\infty )),$$
$$N_{c}(\Gamma _{1},\cdots \Gamma _{n})=
N_{c}\cap N_{d}(\Gamma _{1},\cdots \Gamma _{n}).$$
 The closures of components 
of $\partial N_{c}(\Gamma _{1},\cdots \Gamma _{n})\setminus 
\partial N_{d}(\Gamma _{1},\cdots \Gamma _{n})$ are 
sets $S_{d}(e_{i},\alpha )$, with interior homeomorphic to $\alpha 
$, for each component $\alpha $ of $S(e_{i})\setminus (\cup \Gamma 
_{i})$.
We write $\Sigma (\Gamma _{i})$ for the set of such components. We 
let $S(e_{i},\alpha )$ be the 
topological surface obtained by adding a punctured disc round each 
boundary component.  Write
$${\cal{W}}(S,\Gamma )=\prod _{\alpha \in \Sigma 
(\Gamma }({\cal 
{T}}(S(\alpha ))\cup {\cal{O}}_{a}(S(\alpha ),N).$$
 The boundary 
of $\prod _{e}{\cal{T}}(S(e))/G$ consists of 
$$\prod _{e}({\cal{O}}_{a}(S(e),N)\cup  \cup _{\Gamma 
}{\cal{W}}(S(e),\Gamma))/G,$$
where the union is over all nonempty multicurves $\Gamma $ as above.
We continue to denote an element of this space by $\mu $, so that 
points in ${\cal{T}}(S)/G$ and its boundary are denoted by $\mu $.

The Ending Laminations Theorem can then be formulated as follows, 
given the proofs of the Tameness conjecture. It says that the ending 
lamination invariants are unique. This is the first main result of 
this 
paper.

\begin{theorem}\label{1.1} Let $N$ be any three-dimensional 
hyperbolic 
manifold with finitely generated fundamental group such that $N_{d}$ 
has ends $e_{i}$, $1\leq i\leq n$. Then $N$ is uniquely determined up 
to isometry by its topological type and the ending lamination data 
$[\mu (e_{1}),\cdots ,\mu (e_{n})]$. 

\end{theorem}

As for existence, we have the following, which has extensive overlap 
with 
earlier results in the literature.

\begin{theorem}\label{1.2} Let $N$ be any geometrically finite 
three-dimensional 
hyperbolic 
manifold with finitely generated fundamental group such that the 
horoball deletion $N_{d}$ 
has ends $e_{i}$, $1\leq i\leq n$. Let $N(\mu _{1}',\cdots \mu 
_{n}')$ 
be the manifold in the quasi-isometric deformations space 
$Q(N)$ with ending data $[\mu _{1}',\cdots \mu 
_{n}']\in \prod _{i}{\cal{T}}(S(e_{i}))/G$. Then the map
$$[\mu _{1}',\cdots \mu _{n}']\to N([\mu _{1}',\cdots \mu _{n}'])$$
extends continuously, with respect to both algebraic and geometric 
convergence, to map any point $[\mu _{1},\cdots \mu _{n}]$ of the 
boundary
$$(\prod _{i=1}^{n}{\cal{O}}_{a}(S(e_{i}),N)\cup \cup _{\Gamma 
_{i}}
{\cal{W}}(S(e_{i}),\Gamma _{i}))/G$$
to a hyperbolic manifold $N([\mu _{1},\cdots \mu _{n}])=N')$, with 
$(N_{c}',\partial N_{c}'\setminus \partial N_{d}')$ homeomorphic to 
$(N_{c},S_{d}(\Gamma 
_{1},\cdots \Gamma _{n})$, if $\mu _{i}\in {\cal{W}}(S(e_{i}),\Gamma 
_{i})$ for $1\leq i\leq n$
except in the following cases.
\begin{description}
    \item[1.] $(N_{c},\partial N_{d}\cap N_{c})=(S_{d}\times 
    [0,1],\partial S_{d}\times [0,1] )$ up to homeomorphism, so that, 
    if $n=2$, 
we identify $S\times \{ j
    \} $ with $S$ under $(x,j)\mapsto x$, and $(\mu _{1},\mu 
_{2})=(\mu 
    ,\mu )$ for some $\mu \in \cal{GL}(S)$
    \item[2.] For some $i\neq j$, and nonempty $\Gamma _{i}$, $\Gamma 
    _{j}$, on $S(e_{i})$, $S(e_{j})$, $\mu _{i}\in \cal{W}(S,\Gamma 
    _{i})$, $\mu _{j}\in {\cal{W}}(S,\Gamma _{j})$ and some loops 
    $\gamma \in \Gamma _{i}$ and $\gamma '\in \Gamma _{j}$ are 
    isotopic in $N$.
    \end{description}

   \end{theorem}
    
This overlaps with a substantial and longstanding literature of 
existence results, 
dating back to Thurston's proof of algebraic convergence of 
subsequences in the deformation space with converging ending 
laminatins data, initially for acylindrical manifolds in \cite{T1}. 
The deformation space is 
the space 
of discrete faithful representations of $\pi _{1}(N)$ in 
${\rm{Isom}}(H^{3})$, modulo conjugation in ${\rm{Isom}}(H^{3})$, 
that is, the topolology of algebraic convergence. 
Other examples are \cite{Oh1}, \cite{Oh2}, and, recently, 
\cite{K-L-O}, \cite{Oh4}. This well-developed approach, which often  
uses a 
compactification 
of the deformation space by a space of ${\mathbb R}$-trees, is 
by-passed in \ref{1.2}. Another approach, for Masur domain 
laminations  for handlebodies, is given in Namazi's thesis \cite{Nam} 
in collaboration with Souto, using the Brock-Canary-Minsky Ending 
Lamination Theorems of \cite{B-C-M}.

Thurston's Geometrization Theorem (\cite{Ot2}, \cite{Mor}, \cite{McM} 
says that any compact  
three-dimensional {\em{pared}} manifold $(M,P)$ with boundary is 
homeomorphic to $(H^{3}\cup \Omega 
)/\Delta \setminus P',\partial P')$ for some discrete geometrically 
finite group of isometries 
$\Delta \cong \pi _{1}(M)$ with domain of discontinuity $\Omega $ and 
horoball neighbourhoods of cusps $P'$. Let $N=H^{3}/\Delta $ for such 
a $\Delta $. It 
is then classical, a consequence of the Measurable Riemann-Mapping 
Theorem 
\cite{A-B}, that the space $Q(N)$ of manifolds quasi-isometric to 
$N$, using 
the topology of algebraic convergence, is homeomorphic to $\prod 
_{i=1}^{n}{\cal T}(S(e_{i}))$, where $e_{i}$  ($1\leq i\leq n$) are 
the ends of $N$. The map from   $\prod 
_{i=1}^{n}{\cal T}(S(e_{i}))$ to $Q(N)$ is also continuous with 
respect to the geometric topology on $Q(N)$. Theorem \ref{1.2} then 
implies that $Q(N)$ is dense, in the algebraic topology, in the 
space 
of hyperbolic manifolds of the topological type of $N$. This, then, 
suggests an alternative proof of the Bers-Sullivan-Thurston Density 
Conjecture, that is, density of geometrically finite groups. This was 
proved by Bromberg \cite{Brom} and Brock-Bromberg \cite{B-B}
in the case of incompressible 
boundary and without parabolic elements, derived the Ending 
Laminations Theorem and work of Thurston and Ohshika (\cite{T1}, 
\cite{T2}, \cite{Oh1}, 
\cite{Oh2}, \cite{Oh4}) by 
Brock-Canary-Minsky \cite{B-C-M} in 
their proof of the Ending Laminations Theorem, and can apparently be 
derived from 
the work of Inkang Kim, Lecuire and Ohshika  (\cite{K-L-O}, 
\cite{Oh4}) and the 
Brock-Canary-Minsky proof of the full version of the Ending 
Laminations Theorem. Another proof is in preparation by Bromberg and 
Souto \cite{B-S}

A proof of \ref{1.1}, for both geometrically finite and infinite 
ends was announced in the Kleinian surface case --- 
when $N_{c}$ is homeomorphic to $S_{d}\times [0,1]$, possibly with 
more 
than two ends --- in 2002 by Y. Minsky, J. Brock and D. Canary, the 
first part  being due to Minsky alone \cite{Min2}. They 
also announced the result in general the following year. Their full 
proof in the Kleinian surface case \cite{B-C-M} became publicly 
available at the end of 2004, and includes a number of other results 
including a result related to \ref{1.2}. The general case is in 
preparation.

The Brock-Canary-Minsky proof of the Ending Laminations Conjecture
is the culmination of a number of papers of Minsky in 
which the Ending Laminations Theorem was successively proved in 
important special cases, especially the  Kleinian surface 
once-punctured 
torus case \cite{Min0} and the Kleinian surface bounded geometry case 
\cite{Min1}, \cite{Min3}. The proof of the once-punctured torus case 
was striking because of the strategy of proof, which was then used 
in other cases of the result. The Teichm\" uller space ${\cal T}(T)$ 
of 
the once-punctured torus $T$ identifies with the unit disc $\{ 
z:\vert z\vert <1\} $, and the set of geodesic laminations with the 
boundary $\{z:\vert z\vert =1\} $.
So  the ending laminations data 
identifies an element $(\mu _{+}(N),\mu _{-}(N))$ of $\{ (z,w):\vert 
z\vert ,\vert w\vert \leq 1\} \setminus \{ (z,z):\vert z\vert =1\} $, 
using a slightly nonstandard identification between the two 
components of the relative Scott core. Minsky 
constructs a geometric model for a Kleinian surface hyperbolic 
manifold $N$ homeomorphic to $T\times {\mathbb R}$ with 
given end invariants in terms of the paths in the Farey graph to 
$z$, $w$, at least in the case when $\vert z\vert =\vert w\vert =1$. 
 He was then able to show that 
the geometric model manifold was biLipschitz equivalent to $N$, in 
fact boundedly so, with bounds independent of  $(\mu _{+}(N),\mu 
_{-}(N))$ . From this, the fact that $N$ was the unique manifold in 
this homeomorphism class with invariant  $(\mu _{+}(N),\mu _{-}(N))$ 
was deduced. 

The bounded geometry case of the Ending Laminations Theorem has an 
interesting history, 
\cite{Min1} being a return to the problem of {\em{bounded geometry 
ending invariants}} (as it is reasonable to call them)
8  years after the results of \cite{Min3}, when the Ending 
Laminations 
Theorem was proved for hyperbolic manifolds with bounded geometry, 
which is a stronger assumption than bounded geometry of the end 
invariants . Some of the techniques of 
\cite{Min1} have wider application, and they are of fundamental 
importance to the current work. The other big new input to 
\cite{Min1} was the theory of the curve complex
developed by Minsky and H. Masur (\cite{M-M1}, \cite{M-M2}) which is 
replaced by a 
different theory in the current work.
The general strategy which Minsky 
developed in the punctured torus and bounded geometry cases was then 
used in his resolution, with Brock and Canary for the last part, for 
the proof of the general case.  In 
summary, the Brock-Canary-Minsky proof of the 
Ending Laminations Theorem  can be considered 
to 
consist of three stages: 

1. the construction of a geometric model $M(\mu _{1},\cdots \mu 
_{mn})$
up to quasi-isometry, for given ending lamination data $(\mu 
_{1},\cdots \mu _{n})$ and fixed topological type;

2. given any hyperbolic manifold $N$ with ending laminations
the construction of a map from the 
geometric model $M(\mu _{1},\cdots \mu _{n})$ 
to $N$ which is Lipschitz;

3. a proof (with Brock and Canary) that this map is, in fact, 
biLipschitz.

Minsky's construction of the geometric model and the Lipschitz map, 
and the final proof of biLipschitz, all
 use the curve complex and deep and extensive work of Masur and 
Minsky (\cite{M-M1}, \cite{M-M2}) on 
 hierarchies of tight geodesics in the curve complex. 
One purpose of 
 the current  paper is to show that it is possible to carry out the 
 programme working  directly with 
 Teichm\"uller geodesics. The general strategy is that developed by 
Minsky, but the detail is quite different. The work is built on two 
planks. One, as already mentioned, is some results of \cite{Min3} 
about pleated surfaces, which apply in a wider context than used 
there. The other is results about Teichm\"uller geodesics which 
were developed in \cite{R1} for quite another purpose. These Teichm\" 
uller geodesic results 
are used both to define a geometric model $M=M(\mu _{1},\mu 
_{2}\cdots 
\mu _{n})$ up to Lipschitz equivalence for a hyperbolic manifold in a 
given homeomorphism class with end invariant $(\mu _{1},\cdots \mu 
_{n})$, and to construct a Lipschitz map from the model manifold to a 
hyperbolic manifold with this ending data, and to show that this map 
is coarse biLipschitz.

  We say that a map 
$\Phi :M\to N$ between complete Riemannian manifolds $M$ and $N$ is 
{\em{coarse biLipschitz}} if  $d_{1}$, $d_{2}$ are the 
lifted metrics on the universal covers ${\tilde M}$, ${\tilde N}$ of  
$M$, $N$ and
 there  is $K$ such that for 
all  $x$, $y\in {\tilde M}$, 
$$K^{-1}d_{1}(x,y)-K\leq d_{2}({\tilde \Phi }(x),{\tilde \Phi 
}(y))\leq 
Kd_{1}(x,y)+K,$$
where ${\tilde \Phi }$ is any fixed lift of $\Phi $.

We refer to data in the compactificaton  of $\prod 
_{i}{\cal{T}}(S(e_{i}))$
as {\em{excluded}} if it is as in cases 1 or 2 of \ref{1.2} and 
otherwise it is {\em{permissible}}. The exclusion of 2 of \ref{1.2} 
is 
trivially necessary, if one examines the condition, because if 
$\gamma 
\in \Gamma _{i}$ and $\gamma '\in \Gamma _{j}$ are isotopic in 
$N=N(\mu _{1},\cdots \mu _{n})$ then they must represent the same 
parabolic element. It may be that hyperbolic manifolds in 
different  homeomorphism classes are in the boundary of $Q(N)$ but we 
do not pursue this. The exclusion of 1 of \ref{1.2} is certainly 
necessary, but this is nontrivial, being tied up with the deepest 
arguments in \cite{Bon}. 
The theorem which will imply both \ref{1.1} and \ref{1.2} is 
as follows. 
 
\begin{theorem}\label{1.3} Let $N_{0}$ be a complete hyperbolic 
$3$-manifold with finitely generated fundamental group with horoball 
deletion $N_{0,d}$ and ends $e_{i}$ of $N_{0,d}$, $1\leq i\leq n$.

Then for   any quotiented $n$-tuple of 
permissble 
invariants for the ends $e_{i}$,  $[\mu _{1},\cdots \mu _{n}]$, 
 there is a 
manifold with boundary $M=M([\mu _{1},\cdots \mu _{n}])$, with a 
basepoint $x_{0}(M)$ with interior 
homeomorphic to $N_{0}$ under a map $\Psi _{M}:M\to N_{0}$ and a 
Riemannian 
metric $\sigma _{M}$  with the 
following properties.
\begin{description}

    \item[1.] There is 
a constant $K_{1}$ which depends only on the topological type,
such that the following holds. Any geometric 
limit 
 of the 
structures $(M',x_{0}(M'))$, for  $M'=M([\mu _{1}',\cdots \mu _{n}'])$
converges to $(M,x_{0}(M))$ up to $K_{1}$-Lipschitz equivalence, with 
an isometry between any boundary components,
as\\ $(\mu _{1}',\cdots \mu 
_{n}')\to (\mu _{1},\cdots \mu _{n})$, for 
$M=M(\mu _{1},\cdots \mu _{n})$, and geometrically finite invariants 
$[\mu _{1}',\cdots \mu _{n}']$, and  under a limit of maps homotopic 
to 
the maps $\Psi _{M}^{-1}\circ \Psi _{M'}$.

\item[2.] There is a constant $K_{2}=K_{2}([\mu _{1},\cdots \mu 
_{n}])$ for any permissible 
$[\mu 
_{1},\cdots \mu _{n}]$ which varies continuously with $[\mu 
_{1},\cdots \mu _{n}]$,  such that the following holds. 
Let $N$ be any hyperbolic manifold with $N_{d}$ 
homeomorphic to $N_{0,d}$. Let $\overline{N}$ 
denote the closure of $N$ in $(H^{3}\cup \Omega (\Gamma ))/\Gamma $, 
where $M=H^{3}/\Gamma $ and $\Omega (\Gamma )$ is the domain of 
discontinuity. Let $CH(N)$ denote the convex hull of $N$. 
Then there is  a map 
$\Phi :M(\mu _{1},\cdots \mu _{n})\to \overline{N}$ 
which 
is $K_{2}$-coarse-biLipshitz with respect to the 
hyperbolic metric on $CH(N)$ on the preimage of $CH(N)$, a 
homeomorphism on the preimage of any component of 
$\overline{N}\setminus 
CH(N)$ and an isometry between any boundary components, using the 
Poincar\'e metric on $\Omega (\Gamma )/\Gamma $.

If $N_{c}$ has incompressible boundary, then $K_{2}$ can be 
chosen independent of $(\mu _{1},\cdots \mu _{n})$.
\end{description}
\end{theorem}

A version of this for Kleinian surface groups (or part 2, 
at least) occurs in \cite{B-C-M}. It is stated there that a version 
holds in the general case, but with $K_{2}$ depending on $N$ 
in the case of compressible boundary.
 
 Theorem \ref{1.3} will imply that, given hyperbolic manifolds 
$N_{1}$ and $N_{2}$ with the same ending data $[\mu _{1},\cdots \mu 
_{n}]$, and maps $\Phi _{1}$ and $\Phi _{2}$ as in \ref{1.3}, 
the set-valued map $\tilde {\Phi _{2}}\circ \tilde{\Phi _{1}}^{-1}$ 
 extends to a 
quasiconformal map of $\partial H^{3}$, invariant with respect to the 
actions of the covering groups $\pi _{1}(N_{1})$ and $\pi 
_{1}(N_{2})$ 
on $\partial H^{3}$ which is conformal on the (possibly empty)
domain of discontuity. This is fairly standard in the case of no 
geometrically finite ends. There is slightly more to do if there are 
geometrically finite ends. Details will be given 
later. Such 
a quasiconfomal map must be conformal by Sullivan's result 
\cite{Sull} that there are no $\pi _{1}(N_{1})$-invariant line fields 
on 
$\partial H^{3}$, and hence $N_{1}$ and $N_{2}$ are isometric. So 
Theorem \ref{1.1} follows directly from \ref{1.3}. The deduction of 
\ref{1.2} from \ref{1.3} is almost as direct, once we have defined 
the topology on the compactification of ${\cal T}(S)$ and the model 
manifolds.
 
The paper proceeds as follows. 

\begin{description}
    
    \item[Section 2.] Teichm\"uller space.
    
    \item[Section 3.] Pleated surfaces and  geodesic laminations.
    \item[Section 4.] More on pleated surfaces.
     \item[Section 5.] Teichm\"uller geodesics: long thick and 
dominant 
    definitions.
    
    \item[Section 6.] Long thick and dominant ideas.
    \item[Section 7.] Geometric model manifolds.
    \item[Section 8.] Model-adapted families of pleated surfaces.
  
    \item[Section 9.] Proof of \ref{1.3} in the combinatorial 
bounded geometry Kleinian surface case.

    \item[Section 10.] Lipschitz bounds.
    \item[Section 11.] BiLipschitz bounds.
    \item[Section 12.] Proofs of Theorems \ref{1.1}, \ref{1.2}, 
\ref{1.3}.
           
\end{description}
Sections \ref{2}, \ref{5} and \ref{6} are concerned with the Teichm\" 
uller space
${\cal T}(S)$ of a finite type surface $S$, with no reference to 
$3$-dimensional hyperbolic geometry. In Section \ref{3}, pleated 
surfaces are 
introduced, Bonahon's work \cite{Bon} recalled and its use to define 
the ending laminations, and the work of Canary in the general tame 
case.
Minsky's rather astonishing result  
\cite{Min1}, 
that there is a bounded homotopy between pleated 
surfaces  in 
$N$ whose pleating loci are related by an elementary move, at least 
in the thick part of $N$, is discussed, reinterpreted and extended 
 in Section 4. This 
extension of 
Minsky's pleated surface result is one of two main inputs into the 
current paper. The other main input  is 
the theory of Techm\"uller geodesics from 
\cite{R1}, of which relevant features are described in Section 
\ref{5} and \ref{6}. The most important result (although not the most 
difficult)
is probably \ref{5.5}, which 
states 
how to decompose $\ell \times S$, for any geodesic segment $\ell $ in 
the Teichm\" uller space ${\cal T}(S)$, into disjoint product 
setswhich are either {\em{bounded}} or 
{\em{long, thick and dominant}} (ltd). Definitions are given in 
Section \ref{5}. For the moment, suffice it to say that this 
decomposition is vital to all three parts of the strategy.  From 
Section \ref{7} onwards we work explicitly towards proving the main 
theorems. 
I am following a fairly direct suggestion of the second 
reader 
of an earlier version in the Kleinian surface case, and shall 
highlight the proof 
of the Kleinian surface combinatorial bounded geometry case at each 
stage of 
the proof. Some specialisation to the case of combinatorial bounded 
geometry 
occurs even in Section \ref{6}. Subsections are devoted to this case 
in sections \ref{7} and \ref{8} and as the title of the section given 
above 
indicates, section 9 proves Theorem \ref{1.3} in this case. 
Section \ref{7} goes into considerable detail about the 
model manifolds, especially about the model of the Scott core. The 
geometric convergence of model manifolds --- the first part of 
Theorem \ref{1.3} --- is proved at the end of Section \ref{7}. 
Section 
\ref{8} includes estimates on the geometry of pleated surfaces given 
certain purely combinatorial information about their pleating loci. 
This is used to obtain information about the geometry of the Scott 
core. 
In the case of incompressible boundary, the geometry of the 
non-interval-bundle part of the Scott core is eventually shown to be 
independent of the ending laminations data, up to bounded Lipschitz 
equivalence. As their titles indicate, explicit work on proving the 
main 
theorems is carried out in Sections \ref{10} to \ref{12}, but a lot 
of 
groundwork is done before this. 

Other approaches to the Ending Laminations Theorem are available, or 
are in preparation. One, by Bowditch, has been emerging over the 
last few years (\cite{Bow1}, \cite{Bow2}, \cite{Bow3}). A 
``prehistoric approach'' is in preparation by Bromberg and Souto 
(to which \cite{B-S} is relevant), with contributions from Evans, 
including his paper \cite{Ev}, and 
Brock.

The idea of the approach in this paper emerged during the course of 
other work. 
I first contacted Minsky directly about it in 2001. I thank Yair 
Minsky 
for some helpful discussion of these ideas over 
the last five years. I also thank his collaborators, and, in 
particular, Dick Canary, for his generosity as organiser/secretary of 
the 
Ahlfors Bers Colloquium in 2005.  I am indebted to the 
referee and second reader of the earlier paper on the Kleinian 
surface case 
for their careful reading and (extensive) detailed and useful 
criticisms, and 
the editors of the Newton Institute Proceedings (Kleinian Groups 
2003) for their 
fair and tactful handling of the earlier submission. In particular, 
I should also like to thank Caroline Series for her very helpful and 
pertinent questioning, and for facilitating discussion of this work. 
I also thank Marc Lackenby in this respect, and Brian Bowditch for 
recent comments. I thank Kasra Rafi --- and 
also Misha Kapovich --- for their contributions, and for
fruitful interchange on further 
developments. This is not the end of the story.

\section{Teichmuller space.}\label{2}
\ssubsection{Very basic objects in surfaces}\label{2.0}

Unless otherwise stated, in this work, $S$ always denotes an oriented 
finite type surface without boundary, that is, obtained from a 
compact oriented surface by removing finitely many points, called 
{\em{punctures}}. One does not of course need an explicit realisation 
of $S$ as a compact minus finitely many points. One can define a 
puncture simply to be an end of $S$. A {\em{multicurve}} $\Gamma $ on 
$S$ is a union of simple close nontrivial nonperipheral loops on $S$, 
which are isotopically distinct, and disjoint. A multicurve is 
{\em{maximal}} if it is not properly contained in any other 
multicurve. Of course, this simply means that the number of loops in 
the multicurve is $3g-3+b$, where $g$ is the genus of $S$ and $b$ the 
number of punctures. A {\em{gap}} is a subsurface $\alpha $ of a 
given surface $S$ such that the topological boundary $\partial \alpha 
$ of $\alpha $ in $S$ is a multicurve. A {\em{gap}} of a multicurve 
on a surface $S$ is simply a component of $S\setminus (\cup \Gamma 
)$. If $\alpha $ is any gap, $\Gamma $ is a  {\em{multicurve in 
$\alpha $}} if it satisfies all the above conditions for a closed 
surface, and, in addition, $\cup \Gamma \subset \alpha $ and no loops 
in $\Gamma $ are homotopic to components of $\partial \alpha $. A 
positively oriented Dehn twist round a loop $\gamma $ on an oriented 
surface $S$ will always be denoted by $\tau _{\gamma }$.

\ssubsection{Teichm\"uller space}\label{2.1} 
We consider Teichm\" uller space 
${\cal T}(S)$ of a surface $S$. If $\varphi _{i}:S\to 
S_{i}=\varphi _{i}(S)$ is an orientation preserving homeomorphism, 
and $S_{i}$ is a complete hyperbolic surface with constant curvature 
$-1$, then we define the equivalence relation $\varphi _{1}\sim 
\varphi _{2}$ if and only if there is an orientation-preserving 
isometry $\sigma :S_{1}\to S_{2}$ such that $\sigma \circ \varphi 
_{1}$ is isotopic to $\varphi _{2}$. We define $[\varphi ]$  to be 
the equivalence class of $\varphi $, and ${\cal T}(S)$ to be the set 
of all such $[\varphi ]$, this being regarded as sufficient since 
definition of a function includes definition of its domain. We shall 
often fix   a complete hyperbolic 
metric of constant curvature $-1$ on $S$ itself, which we shall also 
refer to as ``the'' Poincar\'e metric on $S$. 

Complete hyperbolic structure in dimension two is equivalent to 
complex 
structure, for any orientable surface $S$ of finite topological type 
and  
negative Euler characteristic, by the Riemann mapping theorem.  
So endowing such a surface $S$ with a complex 
structure defines an element of the Teichm\"uller space $\cal{T}(S)$. 
More generally, the Measurable Riemann Mapping Theorem \cite{A-B} 
means that supplying a bounded measurable conformal structure for 
$S$ 
is 
enough to define an element of $\cal{T}(S)$, and indeed $\cal{T}(S)$ 
is 
often (perhaps usually) defined in this way.

\ssubsection{Teichm\"uller distance}\label{2.2}
We shall use $d$ to denote Teichm\"uller distance, so long as the 
Teichm\"uller space ${\cal T}(S)$ under consideration is regarded as 
clear. Moreover a metric $d$ will always be Teichm\"uller metric 
unless otherwise specified.
If more than one space is under consideration, we shall use 
$d_{S}$ to denote Teuchm\"uller distance on ${\cal T}(S)$. The 
distance is defined as 
$$d([\varphi _{1}],[\varphi _{2}])={\rm{inf}}\{
{1\over 2} \log \Vert \chi \Vert 
_{qc}:[\chi \circ \varphi _{1}]=[\varphi _{2}]\}$$
where
$$\Vert \chi \Vert _{qc}=\Vert K(\chi )\Vert 
_{\infty }\vert ,\ \ K(\chi )(z)=\lambda (z)/\mu 
(z),$$
where $\lambda (z)^{2}\geq \mu 
(z)^{2}\geq 0$ are the eigenvalues of $D\chi 
_{z}^{T}D\chi _{z}$, and $D\chi _{z}$ is the 
derivative of $\chi $ at $z$ (considered as a 
$2\times 2$ matrix) and $D\chi _{z}^{T}$ is its transpose. The 
infimum is achieved uniquely at a map $\chi $ which is given locally 
in terms of a unique quadratic mass 1 differential $q(z)dz^{2}$ on 
$\varphi 
_{1}(S)$, and its {\em{stretch}} $p(z)dz^{2}$ on $\varphi _{2}(S)$. 
The local coordinates are 
$$\zeta =x+iy=\int _{z_{0}}^{z}\sqrt{q(t)}dt,$$
$$\zeta '=\int _{z_{0}'}^{z'}\sqrt{p(t)}dt.$$
With respect to these local coordinates, 
$$\chi (x+iy)=\lambda x+i{y\over \lambda }.$$
So the distortion $K(\chi )(x+iy)=\lambda $ is constant. The singular 
foliations $x={\rm{constant}}$ and $y={\rm{constant}}$ on $\varphi 
(S)$ 
given locally by the coordinate $x+iy$ for $q(z)dz^{2}$ are known as 
the {\em{stable and unstable foliations for $q(z)dz^{2}$.}} We also 
say that $q(z)dz^{2}$ is the {\em{quadratic differential at $[\varphi 
_{1}]$ for $d([\varphi _{1}],[\varphi _{2}])$}}, and $p(z)dz^{2}$ is 
its {\em{stretch}} at $[\varphi _{2}]$.

\ssubsection{Thick and thin parts}\label{2.4}

Let $\varepsilon $ be any fixed Margulis constant for dimension two, 
that is, for any hyperbolic surface $S$, if $S_{<\varepsilon }$ is 
the 
set of points of $S$ through which there is a nontrivial closed loop 
of length 
$<\varepsilon $, then $S_{<\varepsilon }$ is a (possibly empty) union 
of cylinders with disjoint closures. Then 
$({\cal{T}}(S))_{<\varepsilon }$ is the set of $[\varphi ]$ for 
which $(\varphi (S))_{<\varepsilon }$ contains an least one 
nonperipheral cylnder. The complement of 
$({\cal{T}}(S))_{<\varepsilon }$ is $({\cal{T}}(S))_{\geq \varepsilon 
}$. 
We shall sometimes write simply ${\cal{T}}_{<\varepsilon }$ and 
${\cal{T}}_{\geq \varepsilon }$ if it is clear from the context what 
is 
meant. We shall also write ${\cal{T}}(\gamma ,\varepsilon )$ for the 
set 
of $[\varphi ]$ such that $(\varphi (S))_{<\varepsilon }$ contains a 
loop homotopic to $\varphi (\gamma )$. If $\Gamma $ is a set of 
loops, we write 
$${\cal{T}}(\Gamma ,\varepsilon )=\cup \{ {\cal{T}}(\gamma 
,\varepsilon ):\gamma \in \Gamma \} .$$

\ssubsection{Length and the interpretation of Teichm\"uller 
distance}\label{2.3}
We fix a surface $S$. It will sometimes be convenient to fix a 
hyperbolic metric on $S$, in which case we shall use $\vert \gamma 
\vert $ to denote length of a geodesic path $\gamma $ with respect to 
this metric. With abuse of notation, for 
$[\varphi ]\in {\cal{T}}(S)$ and a
nontrivial nonperipheral closed 
loop $\gamma $ on $S$, we write $\vert \varphi (\gamma )\vert $ for 
the length, with respect to the Poincar\'e metric on the hyperbolic 
surface $\varphi (S)$, of the geodesic homotopic to $\varphi 
(\gamma )$. We write $\vert \varphi (\gamma )\vert '$ for a 
modification of this, obtained as follows. We change the metric in 
$\varepsilon _{0}$-Margulis tube of $\varphi (S)$, for some fixed 
Margulis constant $\varepsilon _{0}$, to the Euclidean metric for 
this 
complex structure in the $\varepsilon _{0}/2$-Margulis tube, so that 
the loop round the annulus is 
length $\sqrt{\vert \varphi (\gamma )\vert }$, and a convex-linear 
combination with the Poincar\'e metric between the $\varepsilon 
_{0}$-Margulis tubes and $\varepsilon _{0}/2$-Margulis tubes. Then we 
take $\vert \varphi (\gamma )\vert '$ 
to be 
the length of the geodesic isotopic to $\varphi (\gamma)$ with 
respect 
to 
this modified metric. If the geodesic homotopic to $\varphi 
(\gamma )$ does not intersect any Margulis tube, then, of course, 
$\vert \varphi (\gamma )\vert =\vert \varphi (\gamma )\vert '$. 
Then for a constant $C$ depending only on $S$ and $\varepsilon _{0}$.
\begin{equation}\label{2.3.1}\vert {\rm{Max}}\{ \vert \log  \vert 
\varphi 
_{2}(\gamma )\vert '
-\log \vert \varphi _{1}(\gamma )\vert '\vert :\gamma \in \Gamma \} 
-d([\varphi _{1}],[\varphi _{2}])\vert \leq C.\end{equation}
Here, $\Gamma $ can be taken to be the set of all simple closed 
nonperipheral 
nontrivial closed loops on $S$. This estimate on Teichm\"uller 
distance derives from the fact that $\vert \varphi (\gamma )\vert '$ 
is inversely proportional to the largest possible square root of 
modulus of an embedded annulus in 
$S$ 
homotopic to $\varphi (\gamma )$. See also 14.3, 14.4 and 14.7 of 
\cite{R1} 
(although the square root of modulus was mistakenly left out of 
\cite{R1}) 
but this estimate appears in other places, for example \cite{M-M2}.
We can simply take $\Gamma $ to be any set of simple 
closed nontrivial nonperipheral loops on $S$ such that   
 that every component of $S\setminus (\cup \Gamma )$ is a 
topological disc with at most one puncture. We shall call such a loop 
set {\em{cell-cutting}}

Fix a Margulis constant $\varepsilon _{0}$. We define
$$\vert \varphi (\gamma )\vert ''=\vert \varphi (\gamma )\cap 
(\varphi (S))_{\geq \varepsilon _{0}}\vert +\sum \{ \vert n_{\zeta 
,\gamma }([\varphi ])\vert :\vert \varphi (\zeta )<\varepsilon _{0}\} 
,$$
where $n_{\zeta 
,\gamma }([\varphi ])$ is such that $\vert \tau _{\zeta 
}^{m}(\gamma )\vert $ is minimised when $m=n$, for $\tau _{\zeta }$ 
as 
in \ref{2.0}. (See 15.9 of \cite{R1} and \ref{2.5}.)

There is $L(\varepsilon _{0})$ such that a cell-cutting loop set 
$\Gamma _{i}$ can 
always be chosen with $\vert \varphi _{i}(\Gamma _{i})\vert ''\leq 
L(\varepsilon _{0})$. Having fixed such loop sets, there is a 
constant $C(\varepsilon _{0})$ such that 
\begin{equation}\label{2.3.2}\vert {\rm{Max}}\{ \log (\# (\gamma 
_{1}\cap \gamma _{2}):\gamma _{1}\in \Gamma _{1},\gamma _{2}\in 
\Gamma _{2}\} 
-d([\varphi _{1}],[\varphi _{2}])\vert \leq C(\varepsilon 
_{0}).\end{equation}

\ssubsection{ Projections to subsurface Teichm\" uller 
spaces}\label{2.5}

For any gap $\alpha \subset S$, we define a 
topological surface $S(\alpha )$ without boundary and a continuous 
map $\pi 
_{\alpha }:{\cal{T}}(S)\to {\cal{T}}(S(\alpha ))$. We define $\varphi 
_{\alpha }(S(\alpha ))$ by defining its conformal structure. After 
isotopy of $\varphi $, we can assume that all the components of 
$\varphi (\partial \alpha )$ are geodesic. We now write 
$\overline{\varphi (\alpha )}$ for the compactification of $\varphi 
(\alpha )$ obtained by cutting along $\varphi (\partial \alpha )$ and 
adding boundary components, each one isometric to some component of 
$\varphi (\partial \alpha)$. Then we form the Riemann 
surface $\varphi 
_{\alpha }(S(\alpha ))$ by attaching a once-punctured disc 
$\{ z:0<\vert z\vert \leq 1\} $ 
to $\overline{\varphi (\alpha )}$ along each component of $\varphi 
(\partial 
\alpha )$, taking the 
attaching map to have constant derivative with respect to length on 
the geodesics $\varphi (\partial \alpha )$ and length on the unit 
circle. Then we define $\varphi _{\alpha }=\varphi $ on $\alpha $ and 
then extend the map homeomorphically across each of the punctured 
discs. 
Then $[\varphi _{\alpha }]$ is a well-defined element of ${\cal 
T}(S(\alpha ))$.

Now let $\alpha $ be a nontrivial nonperipheral simple closed loop. 
Fix an 
orientation on $\alpha $ Then we define
$$S(\alpha )=\overline{\mathbb C}\setminus\{ \pm 2,\pm {1\over 2}\} 
.$$
Now we define an element $[\varphi _{\alpha }]=\pi _{\alpha 
}([\varphi ])\in \cal{T}(S(\alpha ))$, 
for each $[\varphi ]\in \cal{T}(S)$, as follows. Fix a Margulis 
constant $\varepsilon $. If $\varphi (\alpha )\leq \varepsilon $, 
let $A$ be the closed $\varepsilon $-Margulis tube in $\varphi (S)$ 
homotopic 
to $\varphi (\alpha )$. If $\vert \varphi (\alpha )\vert 
>\varepsilon $, let $A$ be the closed $\eta $-neighbourhood of the 
geodesic homotopic to $\varphi (\alpha )$ where $\eta $ is chosen so 
that $A$ is an embedded annulus, and thus can be chosen bounded from 
$0$ if $\vert \varphi (\alpha )\vert $ is bounded above. Fix a simple 
closed geodesic $\beta (\alpha )$ which intersects $\alpha $ at most 
twice and at least once, depending on whether or not $\alpha $ 
separates $S$. We can assume after isotopy that
$\varphi (\alpha )$ and $\varphi (\beta (\alpha ))$ are both 
geodesic, and we fix a  point $x_{1}(\alpha )\in \varphi (\alpha \cap 
\beta (\alpha))$. We make a Riemann surface $S_{1}$ homeomorphic to 
the sphere, by attaching a unit disc to each component of $\partial 
A$, taking the attaching maps to have constant derivative with 
respect to length.  Then we define $\varphi _{\alpha }$ to map 
$\overline{\mathbb C}$ to $S_{1}$ by mapping $\{ z:\vert z\vert =1\} 
$ 
to $\varphi (\alpha )$, $1$ to $x_{1}(\alpha )$, $\{ z:{1\over 
2}\leq \vert z\vert \leq 2\} $ to $A$ and $\{ z\in \mathbb R:{1\over 
2}\leq z\leq 2\} $ to the component of $\varphi (\beta (\alpha 
))\cap A$ containing $\alpha $. Then $\varphi _{\alpha }(S(\alpha ))$ 
is a four-times punctured sphere and so we have an element 
$[\varphi _{\alpha }]\in \cal{T}(S(\alpha ))$. Now the Teichm\"uller 
space $\cal{T}(S(\alpha ))$ is isometric to the upper half plane 
$H^{2}$ with metric ${1\over 2}d_{P}$, where $d_{P}$ deonotes the 
Poincar\'e metric $(1/y)(dx^{2}+dy^{2})$. This is well-known. We now 
give an identification. 
Let  $n_{\alpha }([\varphi ])=n_{\alpha ,\beta (\alpha )}([\varphi 
])$ be the integer assigning 
the minimum value to 
$$m\to \vert \varphi (\tau _{\alpha }^{m}(\beta (\alpha ))\vert .$$
If there is more than one such integer then we take the smallest one. 
There is a bound on the number of such integers of at most two 
consecutive ones. We 
see this as follows. Let $\ell $ be a geodesic in the hyperbolic 
plane 
and let $g$ be a M\"obius involution such that $g.\ell $ is disjoint 
from, and not asymptotic, to $\ell $, and such that the common 
perpendicular 
 geodesic 
segment from $\ell $ to $g.\ell $ meets them in points $x_{0}$, 
$g.x_{0}$, for some $x_{0}\in \ell $. Then the complete
geodesics  meeting both $\ell $ and 
$g.\ell $ and crossing them both at the same angle, are precisely 
those that pass through points $x$ and $g.x$ for some $x\in \ell $, 
and the hyperbolic length of the segment joining $x$ and $g.x$ 
increases strictly with the length between $x_{0}$ and $x$. This 
implies the essential uniqueness of $n$, as follows. We take $\ell $ 
to be a lift of $\varphi (\alpha)$ to the universal cover, and let 
$\ell _{1}$ be another lift of $\varphi (\alpha )$, such that some 
lift of $\varphi (\beta (\alpha ))$ has endpoints on $\ell $ and 
$\ell _{1}$. Then $g$ is determined by making $\ell _{1}=g.\ell $ for 
$g$ as above. But also $\ell _{1}=g_{2}.\ell $, where $g_{2}$ is the 
element of the covering group corresponding to $\varphi (\beta 
(\alpha ))$. We also  have an element $g_{1}$ of the covering group 
corresponding to $\varphi(\alpha )$, which preserves $\ell $ and 
orientation on $\ell $ Then $\vert \varphi (\tau _{\alpha }^{m}(\beta 
(\alpha )))\vert $ is the distance between $x$ and $g.x$ for the 
unique $x$ such that some lift of a loop freely homotopic to  
$\varphi (\tau _{\alpha }^{m}(\beta (\alpha )))$ has endpoints at $x$ 
and $g.x$. The endpoints are also $g_{1}^{-m}.y$ and $g_{2}.y$ for 
$y$ such that $x=g_{1}^{-m}.y$. So $x$ is determined by the $y=y_{m}$ 
such that $g.x=g_{2}.g_{1}^{m}.x$. So then $d(x,x_{0})={1\over 
2}d(x_{0},g^{-1}g_{2}g_{1}^{m}.x_{0})$, which takes its minimum at 
either one, or two adjacent, values of $m$.

Then the isometric identification 
with $H^{2}$ can be chosen 
so that, if we use the identification to regard $\pi _{\alpha }$ as a 
map to $H^{2}$, 
\begin{equation}\label{2.5.1}\pi _{\alpha }([\varphi ])=n_{\alpha 
}([\varphi ])+i\vert \varphi 
(\alpha 
)\vert ^{-1}+O(1).\end{equation}

If $\alpha $ is either a gap or a loop we now define a semimetric 
$d_{\alpha }$ by
$$d_{\alpha }([\varphi _{1}],[\varphi _{2}])=d_{S(\alpha )}(\pi 
_{\alpha }([\varphi _{1}]),\pi _{\alpha }([\varphi _{2}])).$$

\ssubsection{Use of the semimetrics to bound metric 
distance}\label{2.6}

Using \ref{2.3.1}, we see that for a constant $C=C(L_{0})$, 
if $\alpha $ is a gap and $\vert \varphi (\partial \alpha )\vert \leq 
L_{0}$ then
for all $[\varphi _{1}]$, $[\varphi _{2}]\in \cal{T}(S)$,
\begin{equation}\label{2.6.1}d_{\alpha }([\varphi _{1}],[\varphi 
_{2}])\leq 
d([\varphi _{1}],[\varphi _{2}])+C(L_{0}).\end{equation}
This is simply because, for any simple closed nontrivial 
nonperipheral loop $\gamma \subset \alpha $, an annulus homotopic to 
$\varphi (\gamma )$ with modulus boundedly proportional to the 
maximum possible is contained in the surface homotopic to $\varphi 
(\alpha )$, and bounded by the geodesics homotopic to $\varphi 
(\partial \alpha )$. But (\ref{2.6.1}) also holds if $\alpha $ is a 
loop with $\vert \varphi (\alpha )\vert \leq L_{0}$, by considering 
(\ref{2.3.1}) applied to the 
loops $\alpha $ and $\tau _{\alpha }^{n}(\beta (\alpha ))$ ($n\in 
\mathbb Z$).

There is a converse to (\ref{2.6.1}), again using (\ref{2.3.1}), 
which as noted works for a restricted set of loops. 
Suppose that we have a  set 
$\Gamma \subset S$ of simple closed nontrivial nonperipheral loops 
which are all isotopically distinct and disjoint and
such that 
$$\vert \varphi (\Gamma )\vert \leq L_{0}.$$
 Let 
$\Sigma (\Gamma )$ denote the set of gaps of $\Gamma $. 
 Then for a constant $C(L_{0})$,
\begin{equation}\label{2.6.2} d_{S}([\varphi _{1}],[\varphi 
_{2}])\leq 
{\rm{Max}}\{ d_{\alpha }([\varphi _{1}],[\varphi _{2}]):\alpha \in 
\Gamma \cup \Sigma (\Gamma )\} +C(L_{0})\} .\end{equation}

\ssubsection{}\label{2.7}
In 11.1 of \cite{R1}, a projection $\pi _{\alpha }$ was defined 
differently
in the case of marked (equivalently punctured) 
spheres, the projection being done by simply deleting some of the 
punctures. So in those cases the condition $d_{\alpha }\leq d$ was 
automatic. The identification of the image of the projection with the 
$H^{2}$ in the case of $\alpha $ a loop was done in 15.8 of 
\cite{R1}.
If we denote the projection above by $\pi _{1,\alpha }$ 
and the projection of \cite{R1} by $\pi _{2,\alpha }$, then  for 
$L_{0}'=L_{0}'(L_{0})$, if $\vert \varphi _{1}(\partial \alpha )\vert 
\leq 
L_{0}$, 
$$d_{S(\alpha )}(\pi _{1,\alpha }([\varphi _{1}]),
\pi _{2,\alpha }([\varphi _{1}]))\leq L_{0}'.$$
This is proved simply by constructing a bounded distortion 
homeomorphism between the surfaces given by $\pi _{1,\alpha 
}([\varphi _{1}])$ and $\pi _{2,\alpha 
}([\varphi _{1}])$.

\ssubsection{$d_{\alpha }'$, $d_{\alpha _{1},\alpha 
_{2}}'$.}\label{2.8} 
The quantity $d_{\alpha }'$ was defined 14.10 in \cite{R1}, and is an 
extension of the definition of $d_{\alpha }$. Here, $\alpha $ is 
either a nontrivial nonperipheral simple closed loop, or is a 
subsurface of $S$ bounded by such loops, all isotopically disjoint 
and distinct. We use this concept when $\vert \varphi _{1}(\alpha 
)\vert \leq L_{0}$, or $\vert \varphi _{1}(\partial \alpha )\vert 
\leq 
L_{0}$ for some fixed constant $L_{0}$. We fix a Margulis constant 
$\varepsilon _{0}$.  If $\alpha $ is a loop, we take $[\varphi _{3}]$ 
to be the first point on the geodesic segment $[[\varphi 
_{1}],[\varphi _{2}]]$ for which $\varphi _{3}(\alpha )\vert \geq 
\varepsilon _{0}$, and 
$$d_{\alpha }'([\varphi _{1},[\varphi _{2}])=d_{\alpha }([\varphi 
_{1},\varphi _{3}])+\vert \log \vert \varphi _{2}(\alpha )\vert 
'\vert .$$
If $\alpha $ is a subsurface, then
$$d_{\alpha }'([\varphi _{1}],[\varphi _{2}])={\rm{Max}}\{ \vert \log 
\vert \varphi _{1}(\gamma)\vert '-\log \vert \varphi _{2}(\gamma 
)\vert ' \vert \} ,$$
where the maximum is taken over multicurves in $\alpha $ which are 
{\em{cell-cutting in $\alpha $}}, that is,  every component of 
$\alpha \setminus (\cup 
\Gamma )$ is either a topological disc with at most one puncture, or 
an annulus parallel to the boundary.
As in \ref{2.3}, we can take the maximum over a restricted set of 
 multicurves $\Gamma _{1}$ which are cell-cutting in $\alpha $, with 
$\vert \varphi _{1}(\Gamma _{1})\vert ''\leq 
L(\varepsilon _{0})$, at the expense of changing $d_{\alpha 
}'([\varphi 
_{1}],[\varphi 
_{2}])$ by an additive constant. We can 
then also take a multicurve $\Gamma _{2}$ which is cell-cutting in 
$\alpha $ with $\vert 
\varphi _{2}(\Gamma _{2})\vert \leq L(\varepsilon _{0})$, and we then 
have an analogue of \ref{2.3.2}, which is more symmetric in $[\varphi 
_{1}]$ and $[\varphi _{2}]$:
\begin{equation}\label{2.8.2}\vert {\rm{Max}}\{ \log (\# (\gamma 
_{1}\cap \gamma _{2}):\gamma _{1}\in \Gamma _{1},\gamma _{2}\in 
\Gamma _{2}\} 
-d_{\alpha }'([\varphi _{1}],[\varphi _{2}])\vert \leq C(\varepsilon 
_{0},L_{0}).\end{equation}

It therefore makes sense to define $d_{\alpha _{1}\alpha 
_{2}}'([\varphi _{1}],[\varphi _{2}])$ if $\vert \varphi 
_{i}(\partial 
\alpha _{i})\vert \leq L_{0}$, by taking loop sets $\Gamma _{i}$ 
relative to $\alpha _{i}$ like 
$\Gamma _{1}$ relative to $\alpha _{1}$ above, with $\varphi 
_{i}(\Gamma _{i})\vert ''\leq L(\varepsilon _{0})$, and then defining
$$d_{\alpha _{1},\alpha _{2}}'([\varphi _{1}],[\varphi _{2}])=
 {\rm{Max}}\{ \log (\# (\gamma 
_{1}\cap \gamma _{2}):\gamma _{1}\in \Gamma _{1},\gamma _{2}\in 
\Gamma _{2}\} .$$
This is symmetric in $\alpha _{1}$ and $\alpha _{2}$, and changing 
the 
loop sets $\Gamma _{1}$, $\Gamma _{2}$ only changes the quantity 
$d_{\alpha _{1},\alpha _{2}}'([\varphi _{1}],[\varphi _{2}])$ by an 
additive constant. Also, if we 
write $S$ as a disjoint union of gaps and loops $\alpha _{2}$, then 
$d_{\alpha _{1}}'([\varphi _{1}],[\varphi _{2}])$ is the maximum of 
all 
$d_{\alpha _{1},\alpha _{2}}'([\varphi _{1}],[\varphi _{2}])$ up to 
an 
additive constant.

If $\vert \varphi _{3}(\partial \alpha )\vert \leq L_{0}$ then 
it is clear from \ref{2.6} that 
$$\vert d_{\alpha }'([\varphi _{1}],[\varphi _{3}])-d_{\alpha 
}([\varphi _{1}],[\varphi _{2})\vert \leq C(L_{0}).$$
If in addition $[\varphi _{3}]$ is on the geodesics segment 
$[[\varphi _{1}],[\varphi _{2}]]$ then we have from the definitions, 
and from \ref{2.6}:
$$d_{\alpha }'([\varphi _{1}],[\varphi _{2}])\leq d_{\alpha 
}([\varphi _{1}],[\varphi _{3}])+d_{\alpha }'([\varphi _{3}],[\varphi 
_{2}])\leq C'(L_{0}).$$
Actually, a converse inequality holds, and will be discussed in 
section 
\ref{6}.

\ssubsection{Projection is a single point.}\label{2.9}

Suppose that $S=\alpha _{1}\cup \alpha _{2}$ is a union of two closed
subsurfaces, with disjoint interiors with $\alpha _{2}$ not 
necessarily connected, and suppose that the common boundary consists 
of 
nontrivial nonperipheral loops. Fix a homeomorphism $\varphi 
_{1}:\alpha _{1}\to \alpha _{1}'=\varphi _{1}(\alpha _{1})$, where 
$\alpha _{1}$ is a complete hyperbolic surface with geodesic 
boundary $\varphi _{1}(\alpha _{1})$. Consider the set $\cal {X}$ of 
$[\varphi ]$ 
in $\cal{T}(S)$ such that $\varphi (S)=S'$ is the union of the 
hyperbolic surface $\alpha 
_{1}'$ and another subsurface joined along the geodesic boundary, 
and $\varphi =\varphi _{1}$ on $\alpha _{1}$. Then the definition 
of $\pi _{\alpha _{1}}$ in \ref{2.5} is such that $\pi _{\alpha 
}({\cal 
X})$ is a single point in ${\cal T}(S(\alpha _{1}))$.

\section{Pleated surfaces and geodesic laminations.}\label{3}

Throughout this section, $S$ is a finite type surface and $N$ is a 
complete hyperbolic $3$-manifold  with 
finitely generated fundamental group. We fix a hyperbolic metric on 
$S$ and use the length conventions described in \ref{2.3}. 
 In later subsections, we shall use the notation for subsets of $N$ 
as 
described in section \ref{1}: $N_{d}$ for the horoball deletion, 
$N_{c}$ for the relative Scott core, $U(e)$ for the component of 
$N_{d}\setminus N_{c}$ which is a neighbourhood of the end $e$, and 
so on. We shall assume the result of \cite{Ag} and \cite{C-G} where 
necessary, that is, that $N$ is tame.

\ssubsection{}\label{3.1}
The powerful tool of pleated surfaces was introduced by Thurston 
\cite{T}. 
A basic reference is \cite{C-E-G}. A {\em{pleated surface}} is a 
continuous map 
$f:S\to N$ such that peripheral loops are mapped to cusps,
and there is a geodesic lamination $\mu $ on $S$ with respect to some 
hyperbolic structure on $S$,  such that each component of $f(S
\setminus \mu )$ is totally geodesic in $N$ with  boundary consisting 
of 
complete geodesics in $N$. A {\em{geodesic lamination}} is a closed 
set of nonintersecting geodesics on $S$. We then 
call $\mu $ the {\em{bending locus of }} $f$. One could quibble about 
this because $\mu $ is then not quite uniquely defined given  a map 
$f$: 
there may be no bending along some leaves of $\mu $. But one can at 
least be sure that there is no bending along any geodesic which 
intersects $\mu $ transversally, and it will be convenient in the 
present work to assume that the structure of a pleated surface 
includes a lamination $\mu $ such that any bending takes place 
inside $\mu $ and none outside it.

A pleated surface $f$ defines an element of ${\cal T}(S)$ which we 
call $[f]$. This is done as follows, basically just pulling back the 
hyperbolic structure from $f(S)$. Of course, $f$ is not an embedding 
in general, but it is a local embedding restricted to each component 
of 
$S\setminus \mu $, whose image is a complete geodesic triangle in 
$N$, and we take the new hyperbolic structure on $S$ so that $f$ is 
an 
isometry restricted to each such component. This actually defines the 
hyperbolic structure uniquely, in the given homotopy class. The 
transverse length of the geodesic lamination $\mu $ in the new 
hyperbolic structure is $0$, as it was before: bounded length arc 
intersecting a geodesic lamination has zero one-dimensional Lebesgue 
measure intersection with that lamination.

We shall sometimes write $S(f)$ for the abstract Riemann surface with 
hyperbolic structure induced by $f$, if we want to make clear that 
the 
metric on it is that induced from $f(S)$ locally but not globally. 
The map 
from $S$ to $S(f)$ is homotopic to a homeomorphism, and we denote the 
corresponding point of $\cal{T}(S)$ by $[f]$.
There is then a map defined using $f$ defined pointwise from $S(f)$ 
to 
$f(S)$, which we shall  sometimes call ${\rm{Imp}}_{1}(f)$, and which 
is distance-decreasing from $S(f)$ to $N$,
a map  up to homotopy from $S$ to $f(S)$, 
which we shall  sometimes call ${\rm{Imp}}(f)$. By 
abuse 
of notation we shall often denote all three of these maps by $f$, 
although 
${\rm{Imp}}(f))={\rm{Imp}}_{1}(f)\circ f$. We shall sometimes call 
either 
${\rm{Imp}}(f)$ or ${\rm{Imp}}_{1}(f)$ the {\em{impression}} of $f$.

If $\gamma $ is a closed loop in the pleating locus of $f$ then 
$f(\gamma )=\gamma _{*}$, as we shall always denote 
the closed geodesic in $N$ which is freely 
homotopic to $f(\gamma )$, and so, of course, $\vert f(\gamma )\vert 
=\vert \gamma _{*}\vert $. Similarly suppose that $f_{1}$ and $f_{2}$ 
are two pleated surfaces with the same pleating loci restricted to a 
subsurface $\alpha $, where $\partial \alpha $ is also in the 
pleating locus of both. Then the subsets of $S(f_{i})$ bounded by 
$f_{i}(\partial \alpha )$, and isotopic to $f_{i}(\alpha )$ are 
isometric under a homeomorpism $\sigma $ with $[\sigma \circ \varphi 
_{1}]=[\varphi _{2}]$. This, combined with \ref{2.9}, will be of 
crucial importance later.

\ssubsection{Making examples.}\label{3.2}

Fix a homotopy class $f:S\to 
N$ which maps peripheral loops to peripheral loops.
The easiest way to make a pleated surface in this homotopy class 
is to choose a lamination 
on 
$S$ in which every leaf is either a simple closed geodesic or has 
each end asymptotic to either a closed geodesic of the lamination, or 
to a puncture, and so that any complementary component of the 
lamination is a triangle. There are many such examples.
A pleated surface is obtained by 
{\em{spinning round closed loops}} if the nonclosed leaves of the 
pleating locus are obtained as follows. Take a maximal multicurve  
$\Gamma $ 
 on $S$. If $f$ is 
injective on $\pi _{1}$, no further conditions on $\Gamma $ are 
needed.
But if $f$ is not injective on $\pi _{1}$, let $f_{*}(\gamma )$ 
be
homotopically nontrivial and 
nonperipheral for each $\gamma \in 
\Gamma $.  This ensures that each  component $P$ of $S\setminus (\cup 
\Gamma )$ 
is homotopic in 
$S$ to a closed incompressible pair of pants  
$P'$ in $N$. 
 We shall say that such a loop set 
$\Gamma $ is {\em{noncollapsing}} (for $f$).

Fix a decomposition of $S$ into hexagons, possibly 
with some sides of zero length, by taking a set of arcs in each pair 
of pants $P$ with 
endpoints on boundary components or punctures. Fix a homeomorphism 
$\varphi $ 
which is a nontrivial Dehn twist round each loop of $\Gamma $. If 
$\alpha $ is any of the arcs then as $n\to +\infty $, $\varphi 
^{n}(\alpha )$ 
converges to an infinite arc which is asymptotic to a loop of $\Gamma 
$ 
or a puncture of $S$
at each end and the hexagons converge to ideal triangles. 
The union of $\Gamma $ and the limits of the arcs is 
the 
pleating locus of a pleated surface {\em{obtained by spinning round 
$\Gamma $.}} These are, in fact, the basic examples mentioned in 
Thurston's notes \cite{T} Chapter 8.

\ssubsection{Bounded Diameter and Injectivity Radius.}\label{3.3}

Bonahon made the following important observation (essentially from  
1.8 of 
\cite{Bon}, although finite simplicial surfaces are considered there).

\newenvironment{dlemma}
{\par\noindent\textbf{Bounded Diameter Lemma}\,\,\em}{\rm}

\begin{dlemma} For a constant $L(\varepsilon _{0})$ and any pleated 
surface $f:S\to N$, the image in $N$ of each component of the thick 
part 
$(S(f))_{\geq 
\varepsilon _{0}}$ of $S(f)$ has diameter $\leq L(\varepsilon _{0})$. 
\end{dlemma} 
  
As for the thin part, there is a simple but important estimate on 
injectivity radius for 
pleated surfaces, which is actually part of the Bounded Diameter 
Lemma 
as stated by Bonahon \cite{Bon} or the more general form in 2.1 of 
\cite{Can}, where the injective-on-$\pi _{1}$ condition was also 
relaxed.

\newenvironment{ilemma}
{\par\noindent\textbf{Radius of Injectivity  Lemma}\,\,\em}{\rm}

\begin{ilemma} Fix a Margulis constant $\varepsilon _{0}$ for both 
two 
and three dimensions. There are constants $D_{0}$ and 
$C=C(\varepsilon _{0})>1$ 
such that the 
following holds, with $D_{0}$ depending only on the topological type 
of $S$ and $C(\varepsilon _{0})$ depending only on $\varepsilon 
_{0}$. Let 
$f:S\to N$ be a 
pleated surface. Let $\zeta $ be any simple closed geodesic on $S$ 
such that $f(\zeta )$ is homotopically nontrivial and nonperipheral 
in 
$N$. If $\vert f(\zeta )\vert <\varepsilon\leq \varepsilon _{0}$, let 
$T_{S}(f(\zeta ),\varepsilon )$ denote the $\varepsilon $-Margulis 
tube in $S(f)$, if this is nonempy and let $T_{N}(\zeta 
_{*},\varepsilon )$ 
denote the 
Margulis tube in $N$, where $\zeta _{*}$ is the closed geodesic 
freely 
homotopic to $f(\zeta )$. Then
\begin{equation}\label{3.3.1}
    T_{S}(f(\zeta ),\varepsilon )\subset T_{N}(\zeta 
    _{*},\varepsilon ).\end{equation}
    Suppose in addition that $f(\zeta )$ is nontrivial in $N$ 
    whenever $\zeta $ is nontrivial in $S$ and $\vert f(\zeta 
    )\vert \leq D_{0}$. Then
\begin{equation}\label{3.3.2}(S(f)\setminus 
T_{S}(\zeta ,\varepsilon ))\cap T_{N}(\zeta _{*},\varepsilon 
/C)=\emptyset .\end{equation}

\end{ilemma}

(\ref{3.3.1}) is simply because the pointwise map $f:S(f)\to 
f(S)$ preserves length on paths in $S$ and one can apply this to 
closed loops freely homotopic to $f(\zeta )$. (\ref{3.3.2}) is a 
little more 
involved, and uses the extra hypothesis. The extra hypothesis, 
together with the Bounded Diameter Lemma, shows 
that $f((S(f))_{\geq \varepsilon _{0}})$ cannot intersect 
$N_{<\varepsilon _{0}/C}$ for a suitable $C$, which gives 
(\ref{3.3.2}) for $\varepsilon =\varepsilon _{0}$. Then to get the 
result for a general $\varepsilon $ we use the fact that $f:S(f)\to 
N$ 
decreases length of paths joining $\partial T(f(\zeta ),\varepsilon 
_{0})$ and 
$\partial T(f(\zeta ),\varepsilon )$.\Box

The following suggests why this hypothesis of ``no bounded trivial 
loops'' arises.
 
 \begin{lemma}\label{3.3.3} Given $L$, there is $\Delta $ depending 
only on 
$L$ 
 and $N$
 such that the following holds. Let $f:S(e)\to N$ be any pleated 
surface 
 such that the restriction to $S_{d}(e)$ is homotopic to the natural 
inclusion of $S_{d}(e)$ in $N_{c}$, 
 and such that  
with $f(S(e))$ contains at least one point in $U(e)$ distance 
$\geq \Delta $ from $N_{c}$. Suppose that $\zeta \subset f(S(e))$ 
has 
length $\leq L$ and is 
homotopic to a nontrivial nonperipheral loop on $S(e)$. Then $\zeta $ 
is 
nontrivial nonperipheral in $N$.\end{lemma}

\noindent {\em{Proof}} Choose $\varepsilon _{L}$ such that every 
component of $N_{<\varepsilon _{L}}$is distance $\geq L$ from 
$N_{\geq \varepsilon _{0}}$ for some fixed Margulis constant 
$\varepsilon _{0}$. Then let $N_{d}$ be the complement of cuspoidal 
components of $N_{<\varepsilon _{L}}$. We can take our original Scott 
core and extend it to a core for this new $N_{d}$. We continue to 
call 
the core $N_{c}$. Fix $x\subset f(S)$ distance $\geq \Delta $ 
from $N_{c}$.
By the 
Bounded Diameter and Injectivity Radius Lemmas, any point $x'\in 
f(S)$ can be joined to $x$ 
by a path whose intersection with $N_{\geq \varepsilon _{0}}$ has 
length $\leq \Delta _{1}(\varepsilon _{0})$, where 
$\Delta _{1}(\varepsilon _{0})$ depends only on $\varepsilon _{0}$. 
If $\Delta $ is sufficiently large given $L$ and $N$, we can 
assume that none of these thin parts of $N$ is within $2L +\Delta 
_{1}(\varepsilon _{0})$ of $N_{c}$. First suppose that $\zeta $ is 
trivial in $N$. Then $\zeta $ lifts to a closed 
loop $\tilde{\zeta }$ in the universal cover $H^{3}$, and since it 
has length 
$\leq L$ it also 
has diameter $\leq L$, and we can find a continuous map of the disc 
$D$
into $H^{3}$ with $\partial D$ mapped to $\tilde{\zeta }$, by joining 
all points on $\tilde{\zeta }$ to some fixed point on 
$\tilde{\zeta }$ by geodesic segments. The image of this disc in 
$H^{3}$ then has diameter $\leq 2L$, as does its projection in $N$. 
Then by the choice of $\Delta $, this disc does not intersect 
$N_{c}$, which is impossible. If $\zeta $ is peripheral in $N$, then 
argument is similar. This time, we have a bound in terms of $L$ on 
the diameter of the image of the homotopy between $\zeta $ and the 
corresponding element of $\partial N_{d}$. So we again deduce that 
the homotopy cannot intersect $N_{c}$, if $\Delta $ is sufficiently 
large.
\Box 

\ssubsection{Short Bridge Arcs.}\label{3.3.4}

Here is another result which weakens a common  hypothesis
of doubly incompressible for pleated surfaces, but {\em{strengthens another 
hypothesis}}  
to obtain a result in the not-injective-on-$\pi _{1}$ case. The 
result 
which is being generalised is the Short Bridge Arc Lemma 2.2 of 
\cite{Min1} or 
5.5, Uniform 
Injectivity, of \cite{T1}. There is a Uniform Injectivity result in a 
somewhat different direction in 
 Otal's thesis \cite{Ot1}, and  others in \cite{Nam} and \cite{Bow3}. 
We need a notion of {\em{badly bent 
annuli}} for $f_{0}$ and for a fixed Margulis constant $\varepsilon 
_{0}$. 
A {\em{badly bent annulus}} (for $f_{0}$ and $\varepsilon _{0}$) is a 
Margulis tube
$T(f_{0}(\gamma ),\varepsilon (\gamma ))\subset S(f_{0})$ for which 
$\gamma $ is in the pleating locus of $f_{0}$,
$\vert f(\gamma )\vert <\varepsilon _{0}$, and  $\varepsilon (\gamma 
)$ 
is the largest number $ \leq \varepsilon _{0}$ 
such that either $\varepsilon (\gamma 
)=\vert f_{0}(\gamma )\vert $ or  
the images under $f_{0}$ two components of 
$\partial T(f_{0}(\gamma ),\varepsilon 
(\gamma ))$ are distance $\leq 1$ apart in $N$.

\begin{ulemma} The following holds for $L_{0}'$, $L_{1}$ sufficiently 
large, 
depending only on the topological type of $S$ and a given  constant 
$L_{0}$.  Let $N$ be a hyperbolic 
manifold. Fix a Margulis constant $\varepsilon _{0}$.
Let $f_{0}:S(f_{0})\to N$  be a pleated surface.  
For $i=1$, $2$ let $\ell _{i}=\{ x_{i,t}:0\leq t\leq T\} $ 
be geodesic segments in 
$S(f_{0})$, such that $f_{0}(\ell _{i})$ is a geodesic segment in 
$N$, and $t$ is the length parameter. Let $\tau \subset S(f_{0})$ be 
a 
geodesic segment, with respect to the hyperbolic structure on 
$S(f_{0})$, joining $x_{1,0}$ and $x_{2,0}$, and 
let $\tau _{t}$ be the continuously varying geodesic segment joining 
$x_{1,t}$ 
and $x_{2,t}$. 

Suppose that $f_{0}(\tau _{t})$ is homotopic in $N$, 
via homotopy fixing endpoints, to 
a geodesic segment in $N$ of length $\leq L_{0}$, for all $t$.

If $f_{0}$ is not injective on 
$\pi _{1}$, let $A\subset (S(f_{0}))_{<\varepsilon _{0}}$ be the 
union of 
badly bent annuli for $f_{0}$ and $\varepsilon _{0}$ which intersect 
$\tau $, and make two further assumptions.
\begin{description}
\item[1.] The length of $\tau $ in $S(f_{0})\setminus A$ is $\leq 
L_{0}$. 
\item[2.] For any nontrivial loop $\gamma \subset S(f_{0})$ 
for which $f_{0}(\gamma )$ is trivial in $N$, 
$\vert 
f_{0}(\gamma )\setminus A\vert \geq L_{1}$. 
\end{description}

Then, after translating  the length parameter on one of $\ell _{1}$, 
$\ell _{2}$ by $\leq 2\log L_{0}$ if necessary, 
either $\tau _{t}$ has length $\leq 1$ in $S(f_{0})$ for all 
$L_{0}'\leq t\leq 
T-L_{0}'$, or $x_{1,t}$ and $x_{2,t}$ are in the same badly bent 
annulus, on opposite sides of the core loop, either for all 
$L_{0}'\leq t\leq T$, or for all $0\leq t\leq T-L_{0}'$. \end{ulemma}

This has been stated a bit differently from 
2.2 of \cite{Min1}, even leaving aside the different hypotheses in 
the not-injective-on-$\pi _{1}$ case. It is probably worth saying at 
this point that the proof does not use geometric limits, in contrast 
to the proofs of similar results that I am aware of.

\noindent {\em{Proof.}} 
Write $d_{2}$ and $d_{3}$ for 
the hyperbolic distances in the universal covers $H^{2}$  and 
$H^{3}$ of $S(f_{0})$ and $N$ respectively, 
where $H^{2}$ projects to the lift  of 
$f_{0}(S(f_{0}))$ in $H^{3}$. Then $\pi _{1}(S(f_{0}))$ acts on 
$H^{2}$. Identify $x_{i,t}$, $\tau _{t}$, $\ell _{i}$ with lifts to 
$H^{2}$. Now 
there is a lift $\tilde{f_{0}}:H^{2}\to H^{3}$ of $f_{0}$. By abuse 
of 
notation, we also write $x_{i,t}$ for $\tilde{f_{0}}(x_{i,t})$. Note 
that 
$$d_{3}(x_{i,t},x_{i,s})=\vert t-s\vert $$
for all $s$, $t$, but that in general
$$d_{3}(x_{1,t},x_{2,t})\leq d_{2}(x_{1,t},x_{2,t}).$$
Assume without loss of generality that $T\geq 3L_{0}$. Also assume 
without loss of generality that $x_{1,t}$ and 
$x_{2,t}$ are equidistant 
from the ends of the common perpendicular between $\ell _{1}$ and 
$\ell _{2}$ in $H^{3}$, or from the intersection point, if there is 
one, or are the closest points on their respective geodesics, if 
these geodesics are asymptotic in $H^{3}$ . This can be done by 
removing length $\leq \log L_{0}+O(1)$ from the $\ell _{i}$, and 
translating  the length parameters by $\leq $ that amount. 
We can also assume that $d_{3}(x_{1,0},x_{2,0})\leq 1$.
This is because two 
sufficiently 
long geodesics in $H^{3}$ which are distance $\leq L_{0}$ apart at 
the 
two pairs of 
endpoints are distance $\leq 1$ apart in the interiors, apart from 
within $\log L_{0} +O(1)$ of the endpoints. If $f_{0}$ 
is not injective on $\pi _{1}$, since $\tau $ changes as a result 
this 
translation, 
we may replace the hypothesis  $\vert f_{0}(\tau )\setminus A\vert 
\leq L_{0}$ by  $\vert f_{0}(\tau )\setminus A\vert 
\leq 2L_{0}$.

Now fix the greatest $t_{1}\leq T$ such that 
$\vert \tau _{t_{1}}\setminus A\vert 
\leq 2L_{0}$, if such a 
$t_{1}$ 
exists. If it does not 
exist, choose $t_{1}$ so that $d_{2}(x_{1,t},x_{2,t})=\vert \tau 
_{t}\vert $ is minimised 
at $t=t_{1}$. In both cases, 
$$d_{2}(x_{1,t},x_{2,t})=\vert \tau _{t}\vert \geq \vert \tau 
_{t}\setminus A\vert 
\geq 
2L_{0}$$
for all $t_{1}<t\leq T$. If  $f_{0}$ is not injective on $\pi _{1}$, 
we also have a bound of $3L_{0}$ on $\vert \tau _{t}\setminus A\vert 
$ for 
$0\leq t\leq t_{1}$, assuming $L_{0}$ is large enough, depending only 
on a universal constant. We see this as follows. The geodesic segment 
$[x_{i,0},x_{i,T}]$ can intersect at most 
two components of $A$, 
one at the start and one at the end, because the core loops of $A$ 
are 
geodesic and cannot be crossed by other geodesics. Also, because of 
the properties of polygons of geodesics in $H^{2}$, $\tau _{t}$ is in 
a $C_{0}$-neighbourhood of $\tau \cup [x_{1,0},x_{1,t}]\cup 
[x_{2,0},x_{2,t}]$, for a universal constant $C_{0}$. To within a 
universal constant, $\vert \tau _{t}\setminus A\vert $ is the sum 
of the length of $3$ segments, on $\tau \setminus A$, 
$[x_{1,0},x_{1,t}]\setminus A$ and on $[x_{2,0},x_{2,t}]\setminus A$, 
where these 3 segments are maximal with respect to the property that 
for each segment, no point is distance $\leq C_{0}$ in $H^{2}$
from any point on either of the other two segments. The length of 
such 
maximal segments on $[x_{1,0},x_{1,t}]\setminus A$
and $[x_{2,0},x_{2,t}]\setminus A$ is essentially increasing with 
$t$: note that 
$d_{2}(x_{1,s},x_{2,u})\geq d_{3}(x_{1,s},x_{2,u})\geq \vert 
s-u\vert -1$. So $\vert \tau _{t}\setminus A\vert $ is bounded to 
within a universal constant by the 
maximum of $\vert \tau _{t_{1}}\setminus A\vert $ and $\vert \tau 
\setminus 
A\vert $, for $t\leq t_{1}$.  It $t_{1}=T$, $\tau _{T}\subset A$, 
and $\tau _{T}$ is contained in a single component of $A\setminus 
f_{0}(\gamma )$, where $f_{0}(\gamma )$ is the core loop in this 
component of $A$, then the bound on $d_{3}(x_{1,T},x_{2,T})$ shows 
that, as elements of $N$, $x_{j,T}$,  is in a component of $\partial 
N_{\varepsilon _{j}}$ with $e^{-2L_{0}}\leq \varepsilon 
_{1}/\varepsilon _{2}\leq 
e^{2L_{0}}$. Then by the Radius of Injectivity Lemma, 
as elements of $S(f_{0})$, $x_{j,T}$ is in a 
component of $\partial (S(f_{0}))_\geq \varepsilon _{j}'$ with 
$C^{-2}e^{-2L_{0}}\leq \varepsilon _{1}'/\varepsilon _{2}'\leq 
C^{2}e^{2L_{0}}$. Since the components of $\partial (S(f_{0}))_\geq 
\varepsilon 
_{j}'$ are not separated by a core loop $f_{0}(\gamma )$, we have a 
bound on $d_{2}(x_{1,T},x_{2,T})$ in this case.  So we have a bound 
on $d_{2}(x_{1,t},x_{2,t})$ for $t\leq t_{1}$ in terms of the bound 
on $\vert \tau \setminus A\vert $, in the non-injective case, and 
the proof is completed if we can bound $T$ above in terms of $L_{0}$. 
In 
the case when $f_{0}$ is injective on $\pi _{1}$, we need to do a 
similar procedure for decreasing $t$. So in that case, we
similarly define $t_{0}$  be the least $t\geq 0$ such that 
$\vert \tau _{t_{0}}\setminus A\vert 
\leq 2L_{0}$, if such a $t_{0}$ exists, and if not, define 
$t_{0}=t_{1}$. 
We then need to bound $t_{0}$ above in terms of $L_{0}$.

We return to $t\geq t_{1}$ and assume that $t_{1}<T$. 
We also claim, in both injective and noninjective cases, that the 
minimum of
$\vert \tau _{t}\vert =d_{2}(x_{1,t},x_{2,t})$ occurs
at some $t\leq t_{1}+C_{0}'$ for a universal constant $C_{0}'$,  
even if we are assuming that $t=t_{1}$ is a minimum 
of  $\vert \tau _{t}\setminus A\vert $, rather 
than  of $\vert \tau _{t}\vert $. This is the same argument as 
before. To within a universal constant, 
$d_{2}(x_{1,t},x_{2,t})=\vert \tau _{t}\vert $ is the sum of 
the lengths of $3$ segments, on 
$\tau $, $[x_{1,0},x_{1,t}]$ and on $[x_{2,0},x_{2,t}]$, 
where these 3 segments are maximal with respect to the property that 
for each segment, no point is distance $\leq C_{0}$ in $H^{2}$
from any point on either of the other two segments. Since we know 
that $\vert \tau _{t}\setminus A$ is not decreasing for $t>t_{1}$ 
near $t_{1}$ $[x_{1,t_{1}},x_{1,t}]$ and $[x_{2,t_{1}},x_{2,t}]$ are 
bounded from $\tau $ for $t\geq t_{1}$, and also bounded from each 
other, because otherwise $\vert \tau _{t_{1}}\setminus A\vert 
<2L_{0}$.
 It follows that, for a universal 
constant $C_{0}''$, for all $t\geq t_{1}$,
$$d_{2}(x_{1,t},x_{2,t})\geq d_{2}(x_{1,0},x_{2,0})+2t-C_{0}''.$$

Now, to simplify the notation, write $t_{1}=0$.
Define $\varepsilon =\varepsilon _{0}/4C$, for $C$ the constant in 
the Radius of Injectivity Lemma. Fix a constant 
$L_{2}$. If $L_{1}$ is sufficiently large given $L_{0}$ and $L_{2}$, 
and $T\geq L_{2}$,
we can assume that $x_{1,t}\notin (S(f_{0}))_{<\varepsilon }$ so 
long as $0\leq t\leq  L_{2}$. For otherwise, since 
$d_{3}(x_{1,t},x_{2,t})<1$, $x_{1,t}$ and 
$x_{2,t}$ are in the same  lift of the same component of 
$N_{<\varepsilon _{0}/C}$ in $H^{3}$. Assuming that $L_{1}$ is large 
enough for the Radius of Injectivity Lemma to hold, both 
$x_{1,t}$ and $x_{2,t}$ are in $(S(f_{0}))_{<\varepsilon _{0}}$.
Let $\gamma $ be the core loop or 
parabolic of 
the component of $(S(f_{0}))_{<\varepsilon _{0}}$ up to homotopy,  
whose lift 
contains $x_{1,t}$. If $x_{1,t}$ and $x_{2,t}$ are in the same lift 
in 
$H^{2}$ of the same component of $(S(f_{0}))_{<\varepsilon _{0}}$ in 
$H^{2}$, then $\tau _{t}$ is in this component of 
$(S(f_{0}))_{<\varepsilon _{0}}$, which means we have 
$d_{2}(x_{1,t},x_{2,t})\leq Cd_{3}(x_{1,t},x_{2,t})+O(1)$, for $C$ as 
in the Radius of Injectivity Lemma, and either $T=0$, in which case 
the proof is finished,  or we have a contradiction to our assumption 
that 
$d_{2}(x_{1,t},x_{2,t})\geq 2L_{0}$ for all $0\leq t\leq T$, 
assuming only that $L_{0}$ is sufficiently large given $\varepsilon 
_{0}$. 
If $x_{1,t}$ and $x_{2,t}$ are not in the same lift in 
$H^{2}$ of the same component of $(S(f_{0}))_{<\varepsilon _{0}}$, 
then we consider the lift to $H^{2}$ of 
$\tau _{t}*\beta *\tau _{t}^{-1}$, where 
 $\beta $ is 
the core loop of the component of $(S(f_{0}))_{<\varepsilon _{0}}$ 
containing $x_{2,t}$. Replacing $\beta $ by $\beta ^{-1}$ if 
necessary, $\gamma $ and $\tau _{t}*\beta *\tau _{t}^{-1}$ represent 
the same short loop in $N$ but different loops in $S(f_{0})$.
Then $\tau _{t}*\beta 
*\tau _{t}^{-1}*\gamma ^{-1}$ then gives a nontrivial 
loop  in $\pi _{1}(S(f_{0}))$ with length $\leq L_{1}$ in 
$S(f_{0})\setminus A$, which is trivial 
in $\pi _{1}(N)$, if $L_{1}$ is 
sufficiently large given $L_{0}$ and $L_{2}$, which contradicts our 
hypothesis. 

So now we 
assume that for $t\in [0,L_{2}]$, $x_{1,t}$ 
and $x_{2,t}$ lie in $(S(f_{0}))_{\geq \varepsilon }$, and that 
$T=d(x_{1,0},x_{1,T})\geq L_{2}$. We shall obtain a contradiction 
for 
$T$ and $L_{2}$ sufficiently large given a universal constant --- 
$L_{1}$ 
sufficiently 
large given $L_{0}$, in the case when $f_{0}$ is not injective on 
$\pi 
_{1}$. 
So now we can assume that $[x_{1,0},x_{1,T}]$ lies in a component of 
$(S(f_{0}))_{\geq \varepsilon }$, and similarly for 
$[x_{1,0},x_{1,T}]$ --- possibly for a different component of 
$(S(f_{0}))_{\geq \varepsilon }$.  So for $C_{1}$ depending 
only on the topological 
type of $S$, and  on $\varepsilon $, we can find  sequences $g_{j}$ 
and $g_{j}'$ 
in $\pi _{1}(S(f_{0}))$ with $g_{j+1}=g_{j}h_{j}$, 
$g_{j+1}'=g_{j}'h_{j}'$, with 
$$d_{2}(x_{1,0},h_{j}.x_{1,0})\leq C_{1},\ \ 
d_{2}(x_{2,0},h_{j}'.x_{2,0})\leq C_{1},$$
$$d_{2}(x_{1,t_{j}},g_{j}.x_{1,0})\leq C_{1},\ \ 
 d_{2}(x_{2,t_{j}},g_{j}'.x_{2,0})\leq C_{1},$$ 
$$5C_{1}\leq t_{j+1}-t_{j}\leq 6C_{1}.$$
Then
$$d_{3}(x_{1,0},g_{j}^{-1}g_{j}'.x_{2,0})=d_{3}(g_{j}.x_{1,0},g_{j}'.x_{2,0})$$

$$\leq 
d_{2}(g_{j}.x_{1,0},x_{1,t_{j}})+d_{3}(x_{1,t_{j}},x_{2,t_{j}})
+d_{2}(x_{2,t_{j}},g_{j}'.x_{2,0})\leq 
2C_{1}+1,$$
while
$$d_{2}(x_{1,0},g_{j}^{-1}g_{j}'.x_{2,0})=d_{2}(g_{j}.x_{1,0},g_{j}'.x_{2,0})$$

$$\geq 
d_{2}(x_{1,t_{j}},x_{2,t_{j}})-2C_{1}\geq 
2t_{j}+d_{2}(x_{1,0},x_{2,0})-C_{0}''-2C_{1}$$
$$>2t_{j-1}+d_{2}(x_{1,0},x_{2,0})+3C_{1}-C_{0}''\geq 
d_{2}(x_{1,t_{j-1}},x_{2,t_{j-1}})+3C_{1}-C_{0}''$$
$$>d_{2}(g_{j-1}.x_{1,0},g_{j-1}'.x_{2,0})+C_{1}-C_{0}''$$
$$>d_{2}(g_{j-1}.x_{1,0},g_{j-1}'.x_{2,0})=
d_{2}(x_{1},g_{j-1}^{-1}g_{j-1}'.x_{2,0}).$$
So, as points in $H^{3}$, all the points 
$g_{j}^{-1}g_{j}'.x_{2}$
lie in a ball in $H^{3}$ centred on $x_{1}$ of radius $2C_{1}+1$, but 
as points in $H^{2}$, they are distinct. So assuming that $T$ is 
sufficiently large given $C_{1}$,  we can choose 
$j$ 
and $k$ so that $d_{2}(g_{j}.x_{1,0},g_{k}.x_{1,0})$ and
$d_{2}(g_{j}'.x_{2,0},g_{k}'.x_{2,0})$ are bounded in terms of 
$C_{1}$ 
and 
$d_{3}(g_{j}^{-1}g_{j}'.x_{2},g_{k}^{-1}g_{k}'.x_{2})=0$ but 
$d_{2}(g_{j}^{-1}g_{j}'.x_{2},g_{k}^{-1}g_{k}'.x_{2})\neq 0$. Because 
of the assumptions on $\tau $, both $g_{j}g_{k}^{-1}$ and 
$g_{j}'(g_{k}')^{-1}$ can be represented, up to free homotopy, by 
closed loops at both $x_{1}$ and $x_{2}$ with length outside $A$ 
bounded in terms of $L_{0}$ and $C_{1}$, in the case of not injective 
on $\pi 
_{1}$. So, 
 if $L_{2}$ is sufficiently large given 
$C_{1}$, and $L_{1}$ sufficiently large given $L_{0}$, we obtain a 
loop given by the free homotopy class of 
$g_{k}g_{j}^{-1}g_{j}'(g_{k}')^{-1}$ 
which is nontrivial in $S(f_{0})$, with length $\leq 
L_{1}$ in $S(f_{0})\setminus A$ in the case of $f_{0}$ not being 
injective 
on $\pi _{1}$, 
and trivial in $N$. This gives the required 
contradiction.\Box

\newenvironment{elemma}
{\par\noindent\textbf{Efficiency of Pleated surfaces}\,\,\em}{\rm}
\ssubsection{Efficiency of Pleated Surfaces.}\label{3.3.5}

Following Thurston in \cite{T1}, Minsky used his Short Bridge Arc 
Lemma to prove the following. Actually, he proved more, but we only 
state what we need.
\begin{ulemma} The following holds for a suitable $L_{1}$ and 
constant 
$C_{1}(n)$ given $n$. Let $f_{0}:S\to N$ be a pleated surface with 
pleating locus including a maximal multicurve $\Gamma $.  Let $\gamma 
\subset S$ 
have $\leq n$ 
intersections with closed loops in the pleating locus of $f_{0}$. Let 
$A\subset (S(f_{0}))_{<\varepsilon _{0}}$ be the union  of badly bent 
annuli for 
$f_{0}$ for a fixed Margulis constant $\varepsilon _{0}$ 
(\ref{3.3.4}) which 
intersect $f_{0}(\gamma )\subset S(f_{0})$. Suppose also that
 for any 
loop $\gamma '$ for which $f_{0}(\gamma ')$ is trivial in $N$, $\vert 
f_{0}(\gamma ')\setminus A\vert \geq L_{1}$.   Then 
$$\vert f_{0}(\gamma \setminus A)\vert \leq \vert \gamma _{*}\vert 
+C_{1}(n),$$
and $f_{0}(\gamma \cap A)$ is a union of $\leq C_{1}(n)$ segments, 
each of which is homotopic in $N$, via homotopy fixing endpoints, to 
a 
segment of length $\leq 1$.
\end{ulemma}

\noindent {\em{Proof.}} For some $n'\leq 2n$, which depends on how 
the arcs of $\gamma \setminus \Gamma $ compare with the inifinite 
geodesics outwide $\Gamma $ which are in the pleating locus of 
$f_{0}$, 
we choose a connected union of $4n'$ geodesic 
segments $\ell _{i}$, $1\leq i\leq 4n'$ in the universal cover 
$H^{2}$ of 
$S(f_{0})$ which projects to a closed loop homotopic to $f_{0}(\gamma 
)$. 
The segments $\ell _{2i+1}$ project to short segments in $S(f_{0})$, 
of length $<1$, say. It is convenient to extend this sequence, using 
the action of the covering group to a 
bi-infinite sequence which is homotopic to the lift of 
$f_{0}(\gamma )$ in $H^{2}$. So $\ell _{i}$ and $\ell _{i+4kn'}$ 
project to the same segment in $S(f_{0})$ for all integers $k$. The 
segments $\ell _{4i}$ project to loops in $\Gamma $, probably not 
injectively. The geodesic containing $\ell _{4i+2}$ is asymptotic in 
$H^{3}$ to the geodesic containing $\ell _{4i}$ 
at one end, and to the geodesic containing $\ell _{4i+4}$ at the 
other, 
and, again, $\ell _{4i+2}$ projects to the pleating locus of $f_{0}$. 
We apply the Short 
Bridge Arc Lemma \ref{3.3.4} to each pair of geodesics containing 
$\ell _{4i}$ and 
$\ell _{4i+2}$ and to each pair containing $\ell _{4i+2}$ and $\ell 
_{4i+4}$. The role of $\tau $ is played by $\ell _{4i+1}$ and $\ell 
_{4i+3}$ respectively. By abuse of notation, we write $A$ for 
the preimage in $H^{2}$ of $A\subset 
S(f_{0})$. Then \ref{3.3.4} implies that $\ell _{4i}$ 
and $\ell _{4i+2}$ spread apart in $H^{2}$ only when they spread 
apart 
in $H^{3}$. We can then reduce the lengths of even-indexed 
segments, possibly increasing the lengths of odd-numbered segments by 
a bounded amount, 
possibly removing some segments altogether and renumbering, but 
keeping  the segments preserved by the action of the element of the 
covering group determined by $f_{0}(\gamma )$. Carrying out this 
procedure a bounded number of times, bounded in terms of $n'$ and 
$n$, 
we reach the stage 
where, for some $p\leq 4n'$, $\cup _{0\leq i<p}$ projects to 
$f_{0}(\gamma )$ up to homotopy, and 
 each segment $\ell _{i}$ is the union of at most two end 
segments in $A$, one or two segments of length bounded in terms of 
$i$, and a segment which is bounded from $\cup _{0<\leq j< i}\ell 
_{j}$ in $H^{2}$. So then either $\cup _{0\leq i<p}\ell _{i}\setminus 
A$ has 
bounded length in terms of $n$, or  every point in 
$\cup _{0\leq i }\ell _{i}\setminus A$
is a bounded distance from some  geodesic in $H^{3}$, which projects 
to 
the 
closed geodesic $\gamma _{*}$ in $N$ which is homotopic to 
$f_{0}(\gamma )$, with components of $\cup _{0\leq i<p}\ell _{i}\cap 
A$ 
having endpoints a bounded distance apart. This gives the required 
results on $f_{0}(\gamma )$.
\Box 

\ssubsection{Recurrent and Minimal Laminations}\label{3.4}
A point 
$x_{0}$ in a leaf of a lamination 
 is {\em{recurrent}} if, for every open set $U$ 
containing $x_{0}$, and every $L>0$, there are points $x$, 
$y\in U$ in the same leaf of the lamination and distance $>L$ apart 
along that leaf. The set of recurrent points of a lamination is 
nonempty and closed, as is the case for any compact dynamical system. 
In fact, all nonisolated points of a geodesic lamination on a finite 
type hyperbolic surface are recurrent.

A geodesic lamination $\mu $ on $S$ is {\em{minimal}} if given 
$\varepsilon >0$ 
there is $L>0$ such that for every pair of points $x$, $y$ in the 
lamination, there is a point $z$ on the same leaf as $y$, and 
distance 
$\leq 
L$ along the leaf from $y$, such that the hyperbolic distance between 
$x$ and $z$ is $<\varepsilon $. Geodesic laminations are  
exceptional dynamical systems in that a geodesic lamination is 
necessarily minimal if it is recurrent and intersects any simple 
closed geodesic transversally. This is shown in \cite{F-L-P}, where 
the language is of 
measured foliations, but measured foliations and geodesic laminations 
are basically equivalent concepts. In fact, for measured foliations, 
minimality and arationality are exactly the same, if minimality is 
defined in 
the right way.   The method of proof is to endow a recurrent 
lamination with a finite transverse invariant measure, apply 
Poincar\'e 
recurrence, and then analyse the ways in which return can occur. For 
a geodesic lamination or measured foliations, the ways in which a 
return can occur are pretty restricted.  A related exceptional 
property of geodesic 
laminations is that any recurrent lamination is a union of finitely 
closed 
geodesics and minimal laminations on subsurfaces with boundary

The only recurrent leaves in the 
lamination in \ref{3.2} are the closed geodesics in $\Gamma $.

\ssubsection{Intersection number.}\label{3.5}

 Any 
lamination can be given a finite transverse invariant measure, which 
is then automatically supported on the recurrent set. If the 
lamination is minimal, then the support is full. Even if the 
lamination 
is minimal, there may be more than one transverse invariant measure 
up 
to scalar multiple, but the space of measures is finite dimensional. 
The simplest examples of measured  laminations are simple 
closed geodesics, or disjoint unions of simple closed geodesics. If a 
geodesic 
lamination $\mu $ is a single closed geodesic then one can assign a 
finite transverse measure to $\mu $ by saying that the measure of a 
travsersal $I$ is $\# (I\cap \mu )$. One can of course also do this 
if $\mu $ is a finite disjoint union of closed simple geodesics $\mu 
_{i}$, $1\leq i\leq n$. More 
generally one can a assign weight $m_{i}$ to each $\mu _{i}$, and can 
define 
a 
finite transverse invariant measure by defining the measure of a 
transversal $I$ to be
$$\sum _{i=1}^{n}m_{i}\# (I\cap \mu _{i}).$$

A geodesic lamination with a transverse invariant measure is called a 
{\em{measured geodesic lamination}}. An {\em{intersection number}} 
$i(\mu _{1},\mu _{2})$ is defined for each pair of measured geodesic 
laminations $\mu _{1}$, $\mu _{2}$. If $\mu _{1}$ and $\mu _{2}$ are 
simple closed geodesics with tranverse measure assigning weight one 
to 
each intersection $\mu _{i}\cap I$, then $i(\mu _{1},\mu _{2})=\# 
(\mu _{1}\cap \mu _{2})$ or $0$ if $\mu _{1}=\mu _{2}$. We extend 
this linearly to the case when $\mu _{j}$ is a disjoint union of 
simple closed geodesics for $j=1$, $2$. If a disjoint union of simple 
closed geodesics is being considered as a measured geodesic 
lamination, we shall always take each of the geodesics to have weight 
one, unless otherwise stated.
More generally, if $\mu _{1}$ is any 
measured geodesic lamination on $S$ and $ \mu _{2}$ is still a simple 
closed geodesic, then $\mu _{1}$ and $\mu _{2}$ are either disjoint, 
or $\mu _{2}$ is contained in $\mu _{1}$, or $\mu _{2}$ is transverse 
to $\mu _{1}$. In the first two cases, $i(\mu _{1},\mu _{2})=0$. In 
the last case, $i(\mu _{1},\mu _{2})$ is the measure of the 
transversal $\mu _{2}$ with respect to $\mu _{1}$. This then 
generalises easily to the case when $\mu _{2}$ is a finite disjoint 
union of weighted simple closed geodesics. In the general 
case, $i(\mu _{1},\mu _{2})$ can be defined using a partition into 
product rectangles, at least for transverse minimal laminations. But 
one can also note that the set of measured 
geodesic laminations has the structure of a piecewise linear 
manifold \cite{F-L-P}, using a natural correspondence between the 
space of 
measured geodesic laminations ${\cal{MGL}}(S)$ and the space of  
measured foliations 
${\cal {MF}}(S)$ of 
\cite{F-L-P}, for which the transverse invariant masures are 
equivalent to Lebesgue measure. The map $i(\mu _{1},\mu _{2})$ is 
uniformly Lipschitz for $\mu _{2}$ a finite disjoint 
union of weighted simple closed geodesics, and so has a unique 
continuous Lipschitz extension $\mu _{2}$ being any measured geodesic 
lamination (1.10 of \cite{R3}).

The natural projection from measured geodesic laminations with 
nonzero measure to 
recurrent geodesic 
laminations is continuous, with respect to the piecewise linear 
manifold structure on measured geodesic laminations mentioned above 
and the Hausdorff 
topology on geodesic laminations. It is not a bijection, trivially so 
since any transverse measure can be scaled. The natural map from 
projective measured laminations to geodesic laminations is also to a 
bijection, even on the inverse image of minimal laminations. 
Nevertheless, restricted to minimal measured geodesic laminations, 
the 
relation described by : $\mu _{1}\sim \mu _{2}$ and only if 
$i(\mu _{1},\mu _{2})=0$, 
is an equivalence relation, and all elements of the equivalence class 
are the same minimal lamination.

\ssubsection{Geodesic laminations, measured foliations and the 
compactification of Teichm\" uller space.}\label{3.6}

There is a topology on the union of Teichm\"uller space ${\cal T}(S)$ 
and the  space ${\cal{PMF}}(S)$  of projective measured foliations 
--- 
or, equivalently, the space \\ ${\cal{PMGL}}(S)$ of projective 
measured 
geodesic laminations --- 
which makes this space compact, homeomorphic to a closed ball, such 
that the subspace topology on ${\cal T}(S)$ is the usual topology on 
${\cal T}(S)$, and ${\cal{PMF}}(S)$ has the piecewise linear topology 
referred to above \cite{F-L-P}. With respect to this topolology 
${\cal{PMF}}(S)$ is known as the {\em{Thurston boundary}}. The 
topology is obtained from projectivising an embedding of 
${\cal{T}}(S)\cup 
{\cal{MF}}(S)$ in $\mathbb R_{+}^{\cal{S}}$ where $\cal{S}$ is the 
set 
of simple closed nontrivial nonperipheral loops on $S$, and the 
embedding is
$$[\varphi ]\mapsto (\vert \varphi (\gamma )\vert ){\rm{\ (}}[\varphi 
] 
\in \cal{T}(S){\rm{ )}},$$
$$\mu \mapsto i(\mu ,\gamma ){\rm{\ (}}\mu \in {\cal{MF}}(S).$$
We shall 
never make direct use of this topology, but it has the property that 
if $[\varphi _{n}]$ is a sequence in $\cal{T}(S)$ converging to a 
projective class of an arational measured lamination $\mu $
and $\gamma _{n}$ is a sequence of simple closed loops such that 
$\vert \varphi _{n}(\gamma _{n})\vert \leq L$ for all $n$,
then the recurrent part 
of any Hausdorff limit of $\gamma _{n}$ is $\mu $, and $i(\gamma 
_{n},\mu)/\vert \gamma _{n}\vert \to 0$ as $n\to \infty $.

There is a related 
boundary of ${\cal {T}}(S)$, which is of more direct 
relevance, 
in 
which the boundary points are either arational geodesic laminations - 
with no measure - or in a set ${\cal{X}}(\Gamma )$ for some  
multicurve 
$\Gamma $ on $S$. Let $\Sigma $ 
be the set of components of $S\setminus \Gamma $ which are not 
topologically 3-holed spheres. If $\Sigma =\emptyset $ then 
${\cal{X}}(\Gamma )$ is a single point. Otherwise, the elements of 
${\cal{X}}(\Gamma )$ are of the form $(\mu _{1},\cdots \mu _{n})$, 
where 
the elements of $\Sigma $ are numbered, $\Sigma =\{ \alpha _{i}:1\leq 
i\leq m\} $ and each $\mu _{i}$ is either in ${\cal{T}}(S(\alpha 
_{i}))$ 
or an arational geodesic lamination on $S(\alpha _{i})$. Here, 
$S(\alpha 
_{i})$ is as in \ref{2.5}. Convergence of a 
sequence $[\varphi _{n}]$ in ${\cal{T}}(S)$ to 
$(\mu _{1},\cdots \mu _{m}) \in {\cal{X}}(\Gamma )$ 
(including $\Gamma =\emptyset $) is then 
defined as follows. We must have $\vert \varphi _{n}(\Gamma )\vert 
\to 
0$ as $n\to \infty $, which of course is an empty condition if 
$\Gamma =\emptyset $.
 We also have $\pi _{\alpha _{i}}([\varphi _{n}])\to \mu _{i}$ as 
 $n\to \infty $ if $\mu 
 _{i}\in {\cal{T}}(S(\alpha _{i}))$, and, if $\mu _{i}$ is an 
 arational 
 lamination on $S(\alpha _{i})$, any limit of $\pi _{\alpha 
 _{i}}([\varphi _{n}])$ in the Thurston compactification 
 ${\cal{T}}(S(\alpha _{i}))\cup {\cal{PMGL}}(S(\alpha _{i}))$ is $\mu 
 _{i}$, endowed with some transverse measure. 

\ssubsection{The Masur Domain.}\label{3.7}
Let $S\subset N$ be an embedded surface. The {\em{Masur domain}} 
${\cal O}(S,N)$ according to the original 
definition 
(\cite{Ot1}, \cite{Can}, \cite{Le1}) is a set of measured geodesic 
laminations, 
but in fact the property of being in the Masur domain is independent 
of the transverse invariant measure chosen, simply because a 
statement $i(\mu ,\mu ')=0$ 
about measured geodesic laminations is  purely topological, provided 
each minimal component of $\mu $, $\mu '$ is in the support of the 
measure. 

The definition we shall use is:
$\mu \in {\cal{O}}(S,N)$ if and only if there is a constant $c=c(\mu 
)>0$ 
such that the 
following holds, for any transverse invariant measure on $\mu $. For 
each  simple nontrivial 
$\gamma \subset S$ 
which bounds a disc in $N$, and any geodesic laminations $\mu '$ 
and $\mu ''$ 
with 
normalised transverse invariant measures, 
\begin{equation}\label{3.7.1}i(\mu ,\mu ')\geq c{\rm{\ 
whenever\ }}i(\mu '',\mu ')=0{\rm{\ and\ }}i(\gamma ,\mu 
'')\leq c\vert \gamma \vert .\end{equation} 

In particular, if we take $\mu ''=\mu '=\gamma /\vert \gamma \vert $, 
we 
obtain
\begin{equation}\label{3.7.2}i(\mu ,\gamma )\geq c(\mu )\vert \gamma 
\vert 
.\end{equation}
For large compression bodies, as they are called, (\ref{3.7.2}) 
actually implies (\ref{3.7.1}), for a different constant $c(\mu )$ in 
(\ref{3.7.1}) from that in (\ref{3.7.2}). 
It also does so for arational 
geodesic laminations in all cases. This is simply because, if 
$i(\mu ,\mu ')=0$ and $\mu $ is arational, $\mu '=\mu $, and because 
intersection number is uniformly Lipschitz in each variable.
According to \cite{Le1}, (\ref{3.7.2}) can be weakened further for 
arational 
laminations, and 
the term $\vert \gamma \vert $ can be omitted. 

\ssubsection{Bonahon's far-out pleated surfaces and the 
invariants.}\label{3.8} 

 Let $e$ be any end of $N_{d}$, and $S=S(e)$ the corresponding 
surface 
as in 
Section \ref{1}. It was shown by Bonahon \cite{Bon} in the case of 
$S_{d}$ incompressible, and by Canary \cite{Can} in the compressible 
(tame) case,  that there is a sequence $\zeta _{n}$ of simple closed 
geodesics on $S$ with $\zeta _{n}$ nontrivial nonperipheral in $N$
such that if $(\zeta _{n})_{*}$ is the geodesic in $N$ 
freely homotopic to $\zeta _{n}$ then 
 $\zeta _{n}^{*}$ converges to $e$ in $N$.  By extending 
 $\zeta _{n}$ to a maximal multicurve, possibly with some boundary 
components of zero length,
 there is a
 pleated surface $g_{n}$ in the homotopy class of  
 inclusion of $S$ in $N$ with pleating locus including $\zeta 
_{n}$. 
 There are only finitely many Margulis tubes 
 intersecting any compact set, and given any Margulis tube $T$, one 
can 
 find a neighbourhood of $e$ disjoint from $T$ --- the same proof as 
 in \ref{3.3.3}. So by the Bounded 
 Diameter and Injectivity Radius Lemmas in \ref{3.3}, $g_{n}(S)\to e$.

 In the compressible case, $g_{n}$ still exists with $g_{n}(S)\to e$. 
 The following lemma shows the existence of a pleated surface with 
 pleating locus including $\zeta _{n}$. Every short loop on this 
 pleated surface which is nontrivial nonperipheral in $S(g_{n})$ is 
 also nontrivial nonperipheral in $N$ by \ref{3.3.3}, and then 
 $g_{n}(S)\to e$ as in the incompressible case.
 
 \begin{ulemma} Any nontrivial nonperipheral simple loop $\gamma $ on 
$S$ which is 
 nontrivial in $N$ is contained 
 in a noncollapsing maximal multicurve.
\end{ulemma}
 
\par \noindent{\em{Proof.}}
  If not, we have a multicurve $\Gamma $ 
 containing $\gamma $, and such that some complementary component 
 $Q$ of $S\setminus (\cup \Gamma )$ is not a pair of pants, and such 
 that every 
 simple closed loop in the interior of $Q$, which is not homotopic to 
a 
 boundary component, is trivial in $N$, but the loops of 
$\Gamma $ are nontrivial in $N$. Take any pair $\gamma 
 _{1}$, $\gamma _{2}$ of 
 boundary components of $Q$, positively oriented with as elements 
 of $\partial Q$. Then $\gamma _{1}*\gamma _{2}$ is trivial. If this 
 is true for all pairs then there are at most two boundary 
components, 
 and $Q$ must have genus at least one. Except in the case when $Q$ 
 has genus one and one boundary component, we can find two simple 
 loops $\zeta _{1}$, $\zeta _{2}$ in $Q$ which are not homotopic to 
 the boundary but such that $\zeta _{1}*\zeta _{2}$ is homotopic to a 
 boundary component, which must then be trivial, a contradiction. If 
 $Q$ is a torus with one boundary component we can again find $\zeta 
 _{1}$, $\zeta _{2}$ on $S$ such that $\zeta _{1}*\zeta _{2}*\zeta 
 _{1}^{-1}*\zeta _{2}^{-1}$ is trivial, again a contradiction.\Box
  
Bonahon and Canary  also showed, in the incompressible and 
compressible 
cases respectively, that 
 any such sequence $\zeta _{n}$ converging to $e$ converged in the 
 Hausdorff topology to a geodesic lamination $\mu $. It was shown 
 further (Section 5  of \cite{Bon} and 10.1 of \cite{Can}) that $\mu 
$ 
 was minimal, and,Canary showed  
 that $\mu $ is in the Masur domain (10.2 of 
 \cite{Can}). His proof was said to be for the non-cusp case only but 
 does in fact work in general.
 The proofs in \cite{Bon}, which are also part of the 
 argument in \cite{Can}, are inextricably linked with the 
 proof of the main result, the existence of the sequence $\zeta 
 _{n}$, and thus very delicate. 
 
 We shall deal with the case of geometrically finite ends in Section 
 \ref{4}.

\section{More on pleated surfaces.}\label{4}
\ssubsection{Removing badly bent annuli.}\label{4.1}

The following lemma shows that it is possible to avoid badly bent 
annuli.

\begin{ulemma} The following holds for a suitable constant $D_{0}$ 
 and, given any integer $k$ and Margulis 
constant $\varepsilon _{0}$,  an integer $r(\varepsilon _{0},k)$.
    Let $f:S\to N$ be continuous. Let $\gamma $ be a simple 
    nontrivial nonperipheral loop on 
    $S$ such that $f(\gamma )$ is nontrivial 
with $\vert \gamma _{*}\vert <\varepsilon _{0}$. We allow $f(\gamma 
)$ to be a cusp.
Let $\zeta $ be a simple loop with $\leq k$ intersections with 
$\gamma $. 
Suppose that 
the loop between each two consecutive returns of $\zeta $ to $\gamma 
$ is not 
homotopic in $N$ a multiple of 
$\gamma$.
Let $\zeta _{n}=\tau _{\gamma }^{n}(\zeta )$
  
Then  $\vert (\zeta _{n})_{*}\vert <\varepsilon _{0}$ 
for at most $r$ integers $n$, which can be taken to be consecutive. 
Indeed, $r$ can be chosen so that $\zeta _{n}$ intersects $T(\gamma 
_{*},\varepsilon _{0})$ for all but an interval of $r$ consectuive 
$n$, and $\zeta _{n})_{*}\neq \gamma _{*}$. 

A similar result holds if $\gamma $ is a multicurve $(\gamma 
_{1},\cdots \gamma _{s})$, $\vert (\gamma _{i})_{*}\vert <\varepsilon 
_{0}$ for all $i$, all the $(\gamma _{i})_{*}$ are distinct, and 
the loop between each two consecutive returns of $\zeta $ to $\gamma 
_{i}$ is not 
homotopic in $N$ a multiple of 
$\gamma _{i}$. In this case, write $\tau _{\underline{n}}=\tau 
_{\gamma 
_{1}}^{n_{1}}\circ \cdots \circ \tau _{\gamma _{s}}^{n_{s}}$, 
$\underline{n}=(n_{1},\cdots n_{s})$. Then the lower bound on 
$\vert (\tau _{\underline{n}}(\zeta ))_{*}\vert $, and the 
intersection with the $T((\gamma _{i})_{*},\varepsilon _{0})$, hold 
for each $n_{i}$ excluded from an interval of length $\leq 
r(\varepsilon _{0},k)$.   \end{ulemma}

\noindent{\em{Proof.}} We use the following fundamental fact. There 
are constants $C_{d}$ and $\Delta _{d}$ such that the following 
holds. Let 
$\gamma _{i}'$ be any sequence of geodesics in $H^{d}$ such that 
positive end on $\partial H^{d}$ of $\gamma _{i}'$ coincides with the 
negative end of $\gamma _{i+1}'$, and such that there is a segment 
on $\gamma _{i}'$ of length $\geq \Delta _{d}$ which is distance 
$\geq 1$ from both $\gamma _{i+1}'$ and $\gamma _{i-1}'$. Then the 
union 
of the $\gamma _{i}'$ comes within  distance $\leq C_{d}$ of a unique 
geodesic, 
on each segment of each $\gamma _{i}'$ which is distance $\geq 1$ 
from 
both $\gamma _{i-1}'$ and $\gamma _{i+1}'$.
Also $\gamma _{i}'$ comes within within $\varepsilon _{d}$
of the geodesic on this segment, if the segment has length $\geq 
\Delta _{d}'$, where 
$\varepsilon _{d}\to 0$ as $\Delta _{d}'\to \infty $

 Up to free  homotopy in $N$, we 
can make $\zeta $ out of $k$ infinite geodesics, such that the 
positive and negative ends of each one are asymptotic to $\gamma $. 
This gives a bi-infinite sequence of geodesics in $H^{3}$ which we 
call $\gamma _{i}$ such that the positive end of $\gamma _{i}$ is 
asymptotic to the negative end of $\gamma _{i+1}$, and this common 
endpoint 
is an endpoint of a lift of $\gamma $. 
The endpoint is fixed by 
$g_{i}$, whose conjugacy class represents $\gamma $. The sequence for 
a lift of 
$\zeta _{n}$ is obtained by taking as an adjacent pair 
$h_{i,n}.\gamma 
_{i}'$ and $h_{i,n}.g_{i}^{n}.\gamma _{i+1}'$, where 
$h_{i,n}=g_{1}^{n}.\cdots g_{i-1}^{n}$ if $i\geq 1$, with 
modifications if $i<1$.  

Distinct 
Margulis tubes, and components of $N_{<\varepsilon _{0}}$ round 
cusps, 
are a definite distance apart. So the hypothesis on 
consecutive returns means the following. Let $d$ denote hyperbolic 
distance in $H^{3}$. Let $X_{i}$ be the set of 
points $x$ on $\gamma _{i}$ with $d(g_{i}.x,x)\leq \varepsilon _{0}$, 
and $X_{i}'$ the set of points on $\gamma _{i+1}$ with 
$d(g_{i}.x,x)\leq \varepsilon _{0}$. Then the distance between $\cup 
_{n}g_{i}^{n}.(X_{i}\cup X_{i}')$ and $\cup _{j\neq i}\cup 
_{n}g_{i}^{n}.(X_{j}\cup X_{j}')$ is $\geq C(\varepsilon _{0})>0$, 
where 
$C(\varepsilon _{0})$ can be taken arbitrarily large by choice of 
sufficiently small 
$\varepsilon _{0}$. So, from considering $j=i+1$, the biinfinite 
sequence for $\zeta _{n}$ is a 
distance $\leq C_{3}$ from a lift of the geodesic representing it 
provided that the sets $\gamma _{i}\setminus X_{i}$ and 
$g_{i}^{n}(\gamma _{i+1}\setminus X_{i}')$ are a distance $\geq 
(2\Delta _{3}-C_{3})$ apart. Since points of 
$g_{i}^{n}(\gamma _{i+1}'\setminus X_{i}')$ and 
$g_{i}^{m}(\gamma _{i+1}'\setminus X_{i}')$ are a distance $\geq 
\vert 
m-n\vert \varepsilon _{0}$ apart, this is true for all but a bounded 
interval of $n$. For the statement about $\zeta _{n}$ entering 
$T(\gamma _{*},\varepsilon _{0})$, we simply need 
$\gamma _{i}\setminus X_{i}$ and 
$g_{i}^{n}(\gamma _{i+1}'\setminus X_{i}')$ somewhat further apart.

The statement with $\gamma $ replaced by a multicurve is proved in 
exactly the same way.\Box 

\begin{corollary}\label{4.2} The following holds for a suitable 
constant $L_{0}$, $L_{1}$, and 
integers $k_{0}$ and $r$. Here, $k_{0}$, $r$ and $L_{1}$ depend on 
$L_{0}$.
Suppose that $f:S\to N$ is a pleated surface, 
homotopic to an embedding. Let 
$\Gamma \subset S$ be a maximal multicurve such that for any simple 
nontrivial loop $\gamma '\subset S$ with $\leq D_{0}$ intersections 
with 
$\Gamma $, $f(\gamma )$ is nontrivial in $N$. 
Let $\vert \varphi (\Gamma )\vert \leq L_{0}$ 
for some $[\varphi ]\in {\cal{T}}(S)$. 

Then we can choose $\Gamma '$ 
with $\vert \varphi (\Gamma ')\vert \leq L_{1}$
 and $\leq r$ intersections with $\Gamma $ such that 
$\vert \zeta _{*}\vert \geq \varepsilon _{0}$ for all $\zeta \in 
\Gamma '$.\end{corollary}

\noindent{\em{Proof.}} Let $\Gamma ''=\{ \gamma _{i}:1\leq i\leq 
s\} $ be the set of loops in $\Gamma $ for which the 
corrsponding geodesics in $N$ which are of length $<\varepsilon 
_{0}$. If $\Gamma ''=\emptyset $, there is nothing to prove.
Otherwise $\Gamma '$ be a multicurve such that 
$\Gamma ''$ does not intersect $\Gamma \setminus \Gamma '$, and such 
that each loop of $\Gamma '$ is intersected by a loop of $\Gamma ''$, 
and $\Gamma ''\cup (\Gamma \setminus \Gamma ')$ is a maximal 
multicurve. 
We can choose $\Gamma '$ so that $\vert \varphi (\Gamma ')\vert \leq 
L_{0}'$, where $L_{0}'$ depends only on $\L_{0}$. Then $\# (\Gamma 
'\cap \Gamma '')$ is bounded in terms of $L_{0}'$. 
We only need to show that the other hypotheses of 
\ref{4.1}, about distinct Margulis tubes and nontrivial loops, 
hold for the multicurve 
$(\gamma _{1},\cdots \gamma _{s})$ and each $\zeta \in \Gamma ''$, 
if $k_{0}$ is sufficiently large. 
If Margulis tubes round $(\gamma _{i})_{*}$ and $(\gamma _{j})_{*}$ 
coincide, where $\gamma _{i}$ and $\gamma _{j}$ are both intersected 
by $\zeta $, with an arc $\zeta _{1}$ from $\gamma _{i}$ to $\gamma 
_{j}$, then consider the loop $[\gamma _{i}*\zeta _{1}*\gamma 
_{j}*\overline{\zeta _{1}}]$. This loop is nontrivial in $S$ but 
trivial in $N$, and has a bounded number of intersections with 
$\Gamma $ with bound in terms of $L_{0}$. By the Loop 
Lemma, we can then find a simple loop with these properties. The same 
argument works if $i=j$: a closed loop formed from $\zeta _{1}$ by 
adding an arc along $\gamma _{i}$ is a multiple of $\gamma _{i}$.
\Box
 
\ssubsection{Bounded homotopy distance between pleated 
surfaces.}\label{4.3}

In this subsection we rework Minsky's remarkable estimate in 4.1 and 
4.2 of \cite{Min1} on the 
distance between two pleated surfaces related by elementary moves.
Note that we do not require our  pleating loci to be related by 
elementary moves, just to have bounded intersections. We also drop 
the condition of injective-on-$\pi 
_{1}$, although, of course, we do need something to replace it. By 
\ref{4.2}, the assumption of no badly bent annuli on one of the 
surfaces is not much of a restriction. Furthermore, we refine the 
concept of badly bent annuli that was used in Section \ref{3}. Let 
$C=C(\varepsilon _{0})$ 
be the constant of the Radius of Injectivity Lemma of \ref{3.3}. 
Given 
homotopic pleated surfaces $f_{0}$, $f_{1}:S\to N$, a {\em{badly bent 
annulus for $(f_{0},f_{1})$}} is a Margulis tube $T(\gamma 
,\varepsilon (\gamma ))\subset S(f_{0})$ round a loop $\gamma $ which 
is in the pleating locus of $f_{0}$, such that either $\varepsilon 
(\gamma )=\varepsilon _{0}$, or both components of 
$\partial T(\gamma ,\varepsilon (\gamma ))$ are mapped by $f_{0}$ to 
within $2C$ of the same point on a single  geodesic in the pleating 
locus of $f_{1}$. Moreover, the latter is not true for 
$\varepsilon (\gamma )$ replaced by $\varepsilon (\gamma )/C$.

\begin{ulemma} The following holds given an integer $r$ for any 
Margulis constant 
$\varepsilon _{0}>0$, $C>1$ as in \ref{3.3.3}, and some  constants  
$D_{0}$ and 
$L(\varepsilon _{0},r)$. Let $N$ be any complete hyperbolic 
$3$-manifold. Let $\Gamma _{0}$, 
 $\Gamma _{1}$ be noncollapsing maximal multicurves on $S$ with 
$\# (\Gamma 
 _{0}\cap \Gamma _{1})\leq r$. Let $f_{0}$, $f_{1}:S\to N$ be 
homotopic 
 pleated surfaces with pleating loci $\Gamma _{0}$ and $\Gamma _{1}$ 
 respectively. Suppose that there are no badly bent annuli for one of 
 $(f_{0},f_{1})$, $(f_{1},f_{0})$. 
 Let $A\subset S(f_{1})$ (or $A\subset S(f_{0})$)
 be the union of badly bent annuli for $(f_{1},f_{0})$ (or 
 $(f_{0},f_{1})$).
Suppose that for any 
homotopically nontrivial loop 
$\gamma \subset S$ with $f_{0}(\gamma )$ trivial in $N$, $\vert 
f_{0}(\gamma )\vert \geq D_{0}$ (or $\vert f_{0}(\gamma )\setminus 
A\vert \geq D_{0}$).

Then composing $f_{1}$ on the right with a homeomorphism isotopic to 
the identity 
if necessary, there is a 
homotopy $(x,t)\mapsto f_{t}(x):S\times [0,1]\to N$
between $f_{0}$ and $f_{1}$ whose homotopy 
tracks have length $\leq L(\varepsilon _{0},r)$ outside  
$f_{1}^{-1}(A)\times [0,1]$ (or $f_{0}^{-1}(A)\times [0,1]$). 
Moreover, the image under the homotopy of
$f_{j}^{-1}(A)\times [0,1]$ is contained in the union of the 
corresponding Margulis tubes, and the image under $f_{0}$ of each 
component of 
$f_{1}^{-1}(A)$ (or under $f_{1}$ of each component of 
$f_{0}^{-1}(A)$) has bounded diameter.

\end{ulemma}

This lemma is true without any restriction on the length of loops in 
$(\Gamma _{0})_{*}$, $(\Gamma _{1})_{*}$. The assumption on no badly 
bent annuli for one of $(f_{0},f_{1})$, $(f_{1},f_{0})$ is not really 
necessary.
It simplifies the proof 
slightly -- and even the statement. A general statement and proof can 
be deduced by applying the statement as given, twice.

\noindent{\em{Proof.}}

We start by assuming that there are no badly bent annuli for 
$(f_{0},f_{1})$. 

Using  the version of Efficiency of pleated surfaces   which has 
been reformulated from 
 Minsky's version, and proved, in \ref{3.3.5}, we see that, for a 
suitable $C_{0}=C_{0}(r,\varepsilon _{0})$ and  
 $C_{1}(r,\varepsilon _{0})$ as in \ref{3.3.5},  for 
any 
 $\gamma \in \Gamma _{1}$ such that $\vert f_{1}(\gamma )\vert \geq 
\varepsilon _{0}/C$,  $f_{0}(\gamma )$ is, up to homotopy in 
 $S(f_{0})$, a distance $\leq C_{0}$ from $\gamma 
 _{*}=f_{1}(\gamma )$ and
 $$\vert \vert f_{0}(\gamma )\vert -
\vert \gamma _{*}\vert \leq C_{1}.$$
If $\vert \gamma _{*}\vert \leq \varepsilon _{0}/C$, we first deduce 
that $\vert f_{0}(\gamma )\vert \leq C_{0}$ and then, by considering 
a loop $\zeta $ with nonzero intersection with $\gamma $, and $\leq 
r$ intersections with $\Gamma _{0}\cup \Gamma _{1}$, if
$A(\gamma ,f_{1})= \emptyset $,
$$\vert f_{0}(\gamma )\vert \leq C_{1}\vert \gamma 
_{*}\vert .$$
Let $A(\gamma )$ denote the badly bent annulus for $(f_{1},f_{0})$ 
round 
$\gamma $, if this exists. 
If  the modulus in $S(f_{1})$ of either 
component of $T(f_{1}(\gamma ),\varepsilon _{0}\setminus A$ is $L$, 
and 
$\varepsilon $ is ${\rm{Min}}(1/L,1)$, then
$$\vert f_{0}(\gamma )\vert \leq C_{1}\varepsilon .$$
Regarding the domain of $f_{1}$ as being $S(f_{0})$ and composing on 
the right with a homeomorphism isotopic to the identity if 
necessary, we can assume that $f_{1}(\gamma )=\gamma _{*}$ for all 
$\gamma \in \Gamma _{1}$. So then we have a bound $C_{2}(r)$ on the 
distance between $f_{0}(\gamma )$ and $f_{1}(\gamma )$ for 
all 
$\gamma \in \Gamma _{1}$ with $\vert \gamma _{*}\vert \geq 
\varepsilon 
_{0}$. If $\vert \gamma _{*}\vert \leq \varepsilon 
_{0}$, then we have a bound on $\vert f_{0}(\gamma )\vert $, and on 
the distance between $f_{0}(\gamma )$ and $f_{1}(A(\gamma ))$, and 
between $f_{0}(\gamma )$ and $f_{1}(\gamma )$ if there is no badly 
bent annulus round $\gamma $.  
Now we are going to  regard the homotopy as having domain 
$S(f_{0})\times [0,1]$. So far, the homotopy has been constructed on 
$f_{0}(\Gamma _{1})$. We can extend it to the boundary of a bounded 
modulus annulus 
round $f_{0}(\gamma )$ for each $\gamma \in \Gamma _{1}$ for which 
$A(\gamma )$ exists. Let $A'$ be the union of these annuli in 
$S(f_{0})$. The homotopy tracks on $f_{0}^{-1}((f_{0}(\Gamma 
_{1})\setminus 
A')\cup \partial A')\subset S$ are bounded, composing $f_{1}$ on the 
right with a homeomorphism isotopic to the identity if necessary, so 
that the image under $f_{1}$ in $S(f_{1})$ is $(f_{1}(\Gamma 
_{1})\setminus A)\cup \partial A$. Let $A''$ be the union of 
components of $(S(f_{0}))_{<\varepsilon _{0}}\setminus A'$ homotopic 
to 
loops  $f_{0}(\gamma )$ for which $\gamma \in \Gamma _{1}$, $\vert 
f_{1}(\gamma )\vert =\vert \gamma _{*}\vert <\varepsilon _{0}/2C$. 
Here, $C$ 
is, again, the constant of \ref{3.3.3}. The corresponding components 
of 
$(S(f_{1}))_{<\varepsilon _{0}}\setminus A$ have moduli which are 
boundedly proportional. The sets  
$f_{0}(S(f_{0}))$ and $f_{1}(S(f_{1}))$ are a bounded 
distance apart in $N$. So the homotopy can be extended, with bounded 
tracks, to $f_{0}^{-1}((f_{0}(\Gamma _{1})\setminus 
A')\cup A'')\subset S$, composing $f_{1}$ on the right with a 
homeomorphism isotopic to the identity if necessary, so that the 
image 
under $f_{1}$ in $S(f_{1})$ is 
$$f_{1}(\Gamma _{1})\setminus \cup \{ T(f_{1}(\gamma ),\varepsilon 
_{0}):\gamma \in \Gamma _{1},\vert \gamma _{*}\vert <\varepsilon 
_{0}/2C\} .$$

So now, after suitably right-composing $f_{1}$ the homotopy with 
bounded tracks is defined on a set whose 
complement of this set consists of the 
annuli $A'$ and a union of pairs of pants whose images under $f_{0}$ 
in $S(f_{0})$, and under $f_{1}$ in $S(f_{1})$, have boundaries which 
have length $\geq \varepsilon _{0}/2C$, and similarly for the images 
in $N$. So we now need to define the homotopy on these pairs of 
pants. For this, we again follow Minsky \cite{Min1}, but care is 
needed, because we can only apply the Short Bridge Arc Lemma 
\ref{3.3.4} 
to $f_{0}$, not to $f_{1}$. So, for $j=0$, $1$, we foliate the pairs 
of pants in 
$S(f_{j})$ by arcs between boundary components, of bounded length, 
apart from a tripod, of which the arms are of bounded length, one arm 
ending on each boundary component, the other arms meeting at the 
centre of the tripod, as we shall call it. Map the foliation plus 
tripod into 
$N$ by $f_{j}:S(f_{j})\to N$. Lift to the universal 
cover. The images are a bounded distance 
from geodesics in $N$ with lifts which are geodesics in $H^{3}$. 
These three geodesics in $H^{3}$ come within a bounded distance of 
the 
lift of the image  of the centre of the tripod defined using 
$S(f_{0})$, and similarly for the tripod defined using $S(f_{1})$. We 
claim that for any 
$C_{2}>0$ the set 
of points which are distance $\leq C_{2}$ from all three of these 
geodesics has diameter bounded in terms of $C_{2}$, if $D_{0}$ is 
large enough. This will imply that the images of the tripod are a 
bounded distance apart, and the homotopy can be extended to match up 
the images of the foliations. To prove the claim, apply the short 
Bridge Arc Lemma \ref{3.3.4} to 
each of the pairs of geodesics or components of 
$\partial (S(f_{0}))_{<\varepsilon _{0}}$ in $S(f_{0})$. The set of 
points on 
each pair whose lifts are a bounded distance $\leq C_{2}$ apart in 
$H^{3}$
corresponds, up to bounded distance, to the set of points whose lifts 
in the universal cover $H^{2}$ of $S(f_{0})$ are a distance $\leq 
C_{2}$ apart, for $C_{2}$ sufficiently large (but 
universally bounded). In $H^{2}$,  the set of points which are 
distance $\leq 
C_{2}$ from all three geodesics is nonempty and of diameter bounded 
in 
terms of $C_{2}$. So the set in $H^{3}$ is bounded also, as required. 
So the images of the tripods are a bounded distance apart, and the 
homotopy can be extended to match up the images of the foliations, 
after suitabel right-composition of $f_{1}$.

The proof when $f_{0}$, rather than $f_{1}$, is allowed to have badly 
bent annuli, is essentially exactly the same. The original homotopy 
is, after right-composition with a homeomorphism, between $f_{0}$ on 
$f_{0}^{-1}(f_{0}(\Gamma _{1})\setminus A)$ and $f_{1}$ on 
$f_{1}^{-1}(f_{1}(\Gamma _{1}\setminus A')$, where $A$ is the union 
of badly bent annuli round loops of $f_{0}(\Gamma _{0})$ in 
$S(f_{0})$ 
and $A'$ is the corresponding set in $S(f_{1})$. The set 
$(f_{0}(\Gamma _{1})\setminus A)\cup A''$ is then a union of 
pairs of pants minus rectangles in $A$ between bounndary components, 
with sides between boundary components being of bounded length. The 
foliations between boundary components can then be taken to include 
$\partial A$ and to foliate $A$, but the homotopy does not extend 
with bounded tracks across $f_{0}^{-1}(A)\times [0,1]$.  \Box

\ssubsection{Bounded Teichm\" uller distance between pleated 
surfaces.}\label{4.4}
 
We have the following, which is an extension of \ref{4.3} and 
 essentially a reinterpretation of 4.1 and 4.2 of \cite{Min1}. 
 
 \begin{ulemma} Fix a Margulis constant $\varepsilon _{0}$. 
     The following holds for a suitable constant $L_{0}$. Let $N$, 
     $S$, $\Gamma _{0}$, $\Gamma _{1}$, $f_{0}$ and $f_{1}$ satisfy 
     the hypotheses of \ref{4.3}.
     
     Let $\alpha $ be the surface which is the complement in $S$ of 
     loops whose images under $f_{0}$ (or $f_{1}$) are cores of $A$. 
Then 
     $$\vert f_{j}(\partial \alpha )\vert \leq L_{0}$$
     for $j=0$, $1$ and
     and
     $$d_{\alpha }([f_{0}],[f_{1}])\leq L_{0}.$$
 \end{ulemma}

 \noindent {\em{Proof.}}  
 There is a constant $L'(\varepsilon _{0})$ such that the following 
 holds. We use the homotopy whose existence is given by \ref{4.3} 
 to relate the maps $f_{0}$ 
 and $f_{1}$ pointwise (not just up to homotopy). For simplicity, we 
 assume that there are no badly bent annuli for $(f_{0},f_{1})$. As a 
result of 
 \ref{4.3}, we can, in any case, interchange $f_{0}$ and $f_{1}$ and 
 get the first hypotheses of \ref{4.3}.
For any path $\zeta \subset \alpha $ and such that 
$f_{1}(\zeta )$ is in 
the 
$\varepsilon _{0}$-thick part of $S(f_{1})$, assuming that 
$f_{1}(\zeta )$ is geodesic,
\begin{equation}\label{4.4.1}\vert f_{0}(\zeta )\vert \leq 
L'(\varepsilon _{0})(\vert 
f_{1}(\zeta )\vert +1).\end{equation}

  We see this as follows.
 The image under $f_{1}$ of 
$(S(f_{1}))_{\geq \varepsilon _{0}}$ cannot intersect 
$N_{<\varepsilon 
_{0}/C'}$, for a suitably large $C'$ depending only on the 
topological 
type of $S$ and $L(\varepsilon _{0})$, that is, only on 
the topological type of $S(f_{0})$. For if it does so intersect, the 
image 
under the homotopy to $f_{0}$ means that there is a map $f_{0}:U\to 
N_{\varepsilon _{0}/C}$, where $U\subset S(f_{0})$ is connected and 
bounded and 
carries a nonabelian 
subgroup of $\pi _{1}(S)$, contradicting the assumptions
on $f_{0}$. So $f_{1}(\zeta )$ is in $N_{\geq 
\varepsilon _{0}/C'}$. Then  $f_{0}(\zeta )$ is in $N_{\geq 
\varepsilon _{0}/C''}$ for a suitable $C''$. Then  by the Injectivity 
Radius Lemma, $f_{0}(\zeta )$ is contained in $(S(f_{0}))_{\geq 
\varepsilon _{0}/C_{3}}$, where, again, both $C''$ and $C_{3}$ depend 
only on the topological type of $S$.  Then it suffices to show that 
if $\vert 
f_{1}(\zeta 
)\vert $ is bounded, then $\vert f_{0}(\zeta )\vert $ is also,
with a less good bound. Suppose 
it is not so. Then we can find a large number of points in a 
single orbit under the covering group of $S(f_{0})$ within a bounded 
distance of the lift of $f_{1}(\zeta)$. In fact these points can be 
put in a sequence $\{ x_{n}\} $ such that the distance between 
$x_{n}$ 
and $x_{n+1}$ is $\leq C_{1}$, where $C_{1}$ depends only on the 
Margulis 
constant. Then, regarding this as a 
subset of the universal cover of $N$, we have a large number of 
points in a single orbit of the covering group of $N$ within a 
bounded distance of the lift of $f_{0}(\zeta )$ --- which is a 
bounded 
set in the thick part of $N_{\geq \varepsilon _{0}/C_{3}}$.
Then for some $n$ depending only on 
the Margulis constant, and the various constants involved, 
that is, just on the Margulis constant and the topological type of 
$S$, 
two of the orbit 
points $x_{m}$ and $x_{p}$, 
for $m<p\leq m+n$, must be identified in the universal cover of $N$, 
giving a nontrivial loop in $f_{0}(S)$ which is 
trivial in $N$ and of 
length $\leq D_{0}$, assuming $D_{0}$ is large enough 
given the Margulis constant and the topological type of $S$, which 
contradicts our hypothesis (stated at the start of \ref{4.3}). 

 Now 
if $\vert f_{1}(\gamma )\vert \leq \varepsilon /2$, $f_{1}(\gamma )$ 
can be 
realised up to free homotopy by the union of two geodesic segments of 
bounded length in $(S(f_{1}))_{\geq \varepsilon _{0}/2}$, by taking 
segments with endponts in $\partial (S(f_{1}))_{\geq \varepsilon 
_{0}}$. So then $f_{0}(\gamma )$ has length bounded above, and $\leq 
D_{0}$ for $D_{0}$ large enough. So then $f_{0}(\gamma )$ is 
nontrivial.  So we 
now assume that the Injectivity Radius Lemma of \ref{3.3} holds for 
both $f_{0}$ and $f_{1}$. Then the bounded homotopy given by 
\ref{4.3} 
implies that for some $C_{4}$, if $\vert f_{0}(\gamma )\vert $, 
$\vert f_{1}(\gamma )\vert <\varepsilon _{0}/C_{4}$, then both have 
length $<\varepsilon _{0}/2$, and the images under $f_{0}$, $f_{1}$, 
of the boundaries of a Margulis tube round $\gamma $ in $S(f_{0})$, 
$S(f_{1})$ are a bounded distance apart. This includes loops of 
$\partial \alpha $. It follows 
that, for a suitable constant $C_{5}$, if 
$\gamma $ is not homotopic to $\partial \alpha $, the shortest paths 
in $(S(f_{j}))_{<\varepsilon _{0}}$ between boundary components 
differ by a constant, and hence 
$${1\over C_{5}}\vert f_{0}(\gamma )\vert \leq \vert f_{1}(\gamma 
)\vert 
\leq C_{5} \vert f_{0}(\gamma )\vert .$$
It then follows that we can extend (\ref{4.4.1}) to all closed paths 
in $S(f_{1})$, that is, for a constant $C_{6}$,
whenever $f_{1}(\zeta )$ is a closed geodesic,
\begin{equation}\label{4.4.2}\vert f_{0}(\zeta )\vert \leq C_{6}\vert 
f_{1}(\zeta )\vert 
.\end{equation}
It now 
follows that if $D_{0}$ is large enough, given $D_{1}$ depending only 
on the topological type of $S$,  for every nontrivial loop $\gamma $ 
in  $\alpha \subset S(f_{1})$ 
of length $\leq D_{1}$, $f_{1}(\gamma )$ is nontrivial in $N$. Now 
that we know that $\vert f_{0}(\partial \alpha )\vert \leq 
C_{5}\varepsilon 
_{0}$, we can apply the above arguments with $f_{1}$ and $f_{0}$ 
reversed, and we obtain (\ref{4.4.2}) with $f_{0}$ and $f_{1}$ 
interchanged. 

Now to bound $d_{\alpha }([f_{0}],[f_{1}])$ we want to use \ref{2.6}. 
This means that we need to bound
$\vert f_{1}(\zeta )\vert '/\vert f_{0}(\zeta )\vert '$ for a 
cell-cutting set of loops $\zeta $ for which $\vert f_{0}(\zeta 
)\vert '' $ is bounded. Here, $\vert .\vert '$ and $\vert .\vert ''$ 
are as in \ref{2.6}. So $\vert f_{0}(\zeta )\vert '=\vert 
f_{0}(\zeta )\vert ''=\vert f_{0}(\zeta )\vert $ unless $\vert 
f_{0}(\zeta )\vert <\varepsilon _{0}$ or $\zeta $ is transverse to a 
loop $\gamma $ for which $\vert f_{0}(\gamma )\vert <\varepsilon 
_{0}$. If $\vert f_{0}(\zeta )\vert <\varepsilon _{0}$ then 
$\vert f_{0}(\zeta )\vert '$ is boundedly proportional to 
$\sqrt {\vert f_{0}(\zeta )\vert }$ and we already have the bound on
$\vert f_{1}(\zeta )\vert /\vert f_{0}(\zeta )\vert $, which 
suffices. 
So it remains to obtain an estimate for one $\zeta $ crossing each 
loop $\gamma $ which is not the core of a badly bent annulus for 
$(f_{1},f_{0})$, with $\vert f_{0}(\gamma )\vert <\varepsilon 
_{0}$, $\gamma $  and with $\vert f_{0}(\zeta )\vert ''$ bounded. For 
such 
a loop, by the 
definition of $\vert .\vert '$,
$\vert f_{0}(\zeta 
)\vert '$ is boundedly proportional to $\exp \vert f_{0}(\zeta 
)\vert /2$ and to $1/\sqrt{\vert f_{0}(\gamma )\vert }$. If we can 
show that $\exp \vert f_{1}(\zeta 
)\vert $ is boundedly proportional to $1/\vert f_{1}(\gamma )\vert $, 
then $\vert f_{1}(\zeta )\vert '$ is also 
boundedly proportional to $1/\vert f_{1}(\gamma )\vert $ and the 
result follows. So fix $\zeta $. We write $f_{0}(\zeta )$ as a union 
of 
finitely many 
components: one or two components in $(S(f_{0}))_{\geq \varepsilon 
_{0}}$, not necessarily in the pleating locus of $f_{0}$, but of 
bounded length, and four or two components --- which are in the 
pleating locus --- in the Margulis tube 
$T(f_{0}(\gamma ),\varepsilon _{0})$, depending on whether or not 
$\gamma $ is in the pleating locus of $f_{0}$. We make a similar 
decomposition of  $f_{1}(\zeta )$. The components in 
$(S(f_{1}))_{\geq \varepsilon 
_{0}}$  are again bounded, by \ref{4.4.1} (or \ref{4.4.2}). 

First, 
suppose that for at least one of $j=0$, $1$, $\gamma $ is not in the 
pleating locus of $f_{j}$. Then these unions of geodesic 
segments are a bounded distance from the geodesic $\zeta _{*}$ in 
$N$, and from the geodesic representatives in $S(f_{0})$, $S(f_{1})$ 
respectively, assuming that $\varepsilon _{0}$ is small enough. This 
uses the fact mentioned in \ref{4.1}. The long geodesic segments used 
here are bounded apart along most of their length, in the lift of $N$ 
to $H^{3}$,  and the lifts of $S(f_{0})$ 
and  $S(f_{1})$ to $H^{2}$, because the different lifts of Margulis 
tubes are 
distinct. So this means we have a bound on $\vert f_{j}(\zeta )\vert 
-\vert \zeta _{*}\vert $ for $j=0$, $1$, and hence a bound on 
$\exp \vert f_{1}(\zeta )\vert / \exp \vert f_{0}(\zeta )\vert $. 
Since we already know that $\vert f_{1}(\gamma )\vert /\vert 
f_{0}(\gamma )\vert $ is bounded, the result follows. Now suppose 
that $\gamma $ is in the pleating locus of both $f_{0}$, and $f_{1}$. 
We can no longer deduce that the four or two components of 
$f_{j}(\zeta )$ in  
$T(f_{j}(\gamma ),\varepsilon _{0})$ are within a bounded distance of 
$\zeta _{*}$. But the Short Bridge Arc Lemma \ref{3.3.3} implies that 
a point $x_{0}$ in $f_{0}(\zeta )$ is within a bounded distance of 
$f_{0}(\gamma )=\gamma _{*}$ precisely when the same is true for a 
point $x_{1}\in f_{1}(\zeta )$ for which $d_{3}(x_{0},x_{1})$ is 
bounded. 
So in this case, also, we have a bound on $\vert \vert f_{0}(\zeta 
)\vert -\vert f_{1}(\zeta )\vert \vert $, and, hence, upper and lower 
bounds on $\vert f_{0}(\zeta )\vert '/\vert f_{1}(\zeta )\vert '$.

\Box
 
\ssubsection {}\label{4.7}

We have the following extension of the ideas of \ref{4.4}, \ref{4.1}. 
For simplicity, we remove as many badly bent annuli as possible. 
\begin{ulemma} The following holds for a constant $L_{0}$. Let 
$\Gamma _{0}$ be a maximal multicurve. Let $\zeta \in \Gamma 
_{0}$ 
 and let $\gamma $ be a simple loop intersecting $\zeta $ finitely 
 many times  
  but no 
 other loop in $\Gamma _{0}$. Let $f:S\to N$ have pleating locus 
including 
$\Gamma 
 =(\Gamma \setminus \{ \zeta \} )\cup \{ \gamma \} $. Let $\vert 
\gamma 
 _{*}\vert <\varepsilon _{0}$ but $\vert \beta _{*}\vert \geq 
 \varepsilon _{0}$ for all $\beta \in \Gamma _{0}\setminus \{ \zeta 
 \} $. Let $f_{n}:S\to N$ be a pleated surface 
 with pleating locus including $\Gamma _{n}$ and homotopic to $f$.
 Let $\zeta _{n}=\tau _{\gamma 
 }^{n}(\zeta )$ for all $n$, or $\tau _{\gamma 
 }^{-n}(\zeta )$ for all $n$. Suppose that both $f_{0}(\zeta )$ and 
 $f_{0}(\zeta _{1})$ are nontrivial in $N$.
 Let $\Gamma _{n}=(\Gamma \setminus \{ \gamma \} 
 )\cup \{ \zeta _{n}\} $.
 
 Suppose that either $f$ satisfies the hypothesis of $f_{0}$ in 
 \ref{4.3} (and \ref{4.4}), or that $f_{j}$ does, for all $0\leq 
j\leq n$.
 Then so long as $\vert (\zeta _{j})_{*}\vert \geq \varepsilon _{0}$
  for $j=0$ and $j=n$,
 $$d([f_{0}],[f_{n}])\leq L_{0}(\log \vert n\vert +1),$$
 with similar bounds on the length of homotopy tracks of some 
homotopy 
 between the impressions of $f_{0}$ and $f_{n}$,
 and if $\alpha =S\setminus \gamma $, for all $n$,
 $$d_{\alpha }([f_{0}],[f_{n}])\leq L_{0},$$
 $$d_{\alpha }([f_{0}],[f])\leq L_{0}.$$

\end{ulemma}
\noindent{\em{Proof.}} 

From \ref{4.4}, since the hypotheses are satisfied, we have, for 
suitable $L_{0}$, 
$$d_{\alpha }([f],[f_{j}])\leq L_{0}/2{\rm{\ \ whenever\ \ }}\vert 
f_{j}(\zeta 
_{j})\vert \geq \varepsilon _{0}.$$ This then gives $d_{\alpha 
}([f_{0}],[f_{n}])\leq L_{0}$.

 Let $\tilde{T}$ be a 
lift of $T(\gamma _{*},\varepsilon _{0})$ to the universal cover 
$H^{3}$ of $N$, and let $\tilde{\zeta _{i}}$ be a lift of $(\zeta 
_{i})_{*}$ such that  any intersection $\tilde {\zeta _{i}}\cap 
\tilde{T}$ corresponds to the intersection between $f_{i}(\zeta 
_{i})$ 
and $f_{i}(\gamma )$ in $S(f_{i})$. Then, by the proof of \ref{4.1}, 
we can assume that
$\tilde{\zeta _{i}}\cap \tilde{T}\neq \emptyset $  and 
$\vert \tilde{\zeta _{i}}\cap 
\tilde{T}\vert $ is increasing in $i$. 
Moreover, for a constant $L_{1}$, for $i\geq 0$, for some constant 
$s_{0}\geq 0$,
\begin{equation}\label{4.7.1}\vert \vert \tilde{\zeta _{i}}\cap 
\tilde{T}\vert -\log (\vert i+s_{0}\vert +1)\vert \leq 
L_{1}.\end{equation}
This can be seen by considering $\tilde{\zeta 
_{i}}$, , translating one component of $\tilde{\zeta _{i}}\cap 
\tilde{T}$ using the covering element of $\gamma $ and putting in a 
new connecting segment.  Since the Radius of Injectivity Lemma of 
\ref{3.3} holds for $f_{j}$ for all $j$, the image under $f_{
i}$ of the Margulis tube in $S(f_{i})$ round 
$f_{i}(\gamma )$ is, to within bounded distance, the 
intersection with $T(\gamma _{*},\varepsilon _{0})$. So it follows 
from (\ref{4.7.1}) that, for a constant $L_{2}$, for $i$ between $0$ 
and $n$,
\begin{equation}\label{4.7.2}
  \vert f_{i}(\gamma )\vert \leq {L_{2}\over  i+1}\end{equation}

Now suffices to show 
that 
for a suitable $L_{3}$, for  all 
$0<i\leq  n\vert $, and assuming that the  pleating lamination 
for $f_{i-1}$ maps under $\tau _{\gamma }$ to the pleating locus 
for $f_{i}$ (as it can be chosen to do)
\begin{equation}\label{4.7.3}d([f_{i-1}],[f_{i}])\leq 
    L_{3}\vert f_{i}(\gamma )\vert .\end{equation}
For then
$$d([f_{0}],[f_{n}])\leq L_{4}\sum _{i=1}^{n}{1\over 
\vert i\vert }\leq L_{4}(\log \vert n\vert +1),$$
as required. So it remains to prove (\ref{4.7.3}). We do this by 
constructing $\chi 
_{i}:S(f_{i-1})\to S(f_{i})$ with $[\chi _{i}\circ 
f_{i-1}]=[f_{i}]$ and with $\chi _{i}$ of distortion $1+O(\vert 
f_{i}(\gamma )\vert )$. We already have the bound by \ref{4.4}, 
if $\vert f_{i}(\gamma )\vert $ is bounded from $0$. So now we 
assume 
that $\vert f_{i}(\gamma )\vert <\varepsilon _{0}/2$, and similarly 
for $i-1$ replacing $i$. Then we construct $\chi _{i}$ to map 
$\partial T(f_{i-1}(\gamma ),\varepsilon _{0})$ to 
$\partial T(f_{i}(\gamma ),\varepsilon _{0})$. Outside 
these sets the pleating loci are homeomorphic and we simply map them 
across. Inside the Margulis tubes, we choose $\chi _{i}$ to map the 
geodesic segments in 
$\partial T(f_{i-1)}(\gamma ),\varepsilon _{0})$, 
geodesic with respect to the hyperbolic metric on $S(f_{ 
i-1)})$, whose 
images
under $f_{i -1}$ are homotopic to the segments of 
$f_{i}(\zeta _{i})$ in 
$\partial T(f_{i}(\gamma ),\varepsilon _{0})$, to those 
segments. In all cases, 
geodesic side lengths differ by 
$O(\vert f_{i}(\gamma )\vert )$. So we can take constant derivative 
with respect to length on each geodesic segment , with derivative 
$1+O(\vert f_{i}(\gamma )\vert )$. The union of the polygons has 
full measure. Then we extend across polygons and  we get the required 
bound on 
distortion.

\Box 
\ssubsection{Pleated Surfaces near the domain of 
discontinuity.}\label{3.9}
 
We now return to a topic left open at the end of Section \ref{3}. 
\ref{3.8} is concerned with geometrically infinite ends. We shall 
also need a corresponding result for geometrically finite ends. So 
now, let $e$ be a geometrically finite end of $(N_{d},\partial 
N_{d})$. Write $S(e)=S$ and $S_{d}(e)=S_{d}$. Then there is a 
component $S_{1}$ of the boundary of the convex hull 
bounding a neighbourhood of $e$, whose intersection with $N_{d}$ is 
homeomorphic to $S_{d}$, under a homeomorphism of $N$ which is 
isotopic to the identity. The neighbourhood of $e$ can also be 
compactified by adding a component $S_{2}$ of $\Omega /\pi _{1}(N)$, 
where 
$\Omega $ is the domain of discontinuity of $\pi _{1}(N)$, where we 
are identifying $\pi _{1}(N)$ with the covering group of hyperbolic 
isometries of $N$. These two surfaces $S_{1}$ and $S_{2}$ both have 
hyperbolic structures, or equivalently 
complex structures, and thus give points of $\cal{T}(S)$, which 
for the moment we call respectively $[f_{1}]$ and $[f_{2}]$. 
It is proved in \cite{E-M} that there is a natural bounded 
distortion 
map between the two surfaces, and that the two points in Teichm\" 
uller space are a bounded distance apart. Actually, there is a fairly 
direct proof of the following. But note the assumption that $\Gamma $ 
is noncollapsing.

\begin{ulemma} Fix a $2$- and $3$-dimensional Margulis constant 
$\varepsilon _{0}$. Given $L_{1}$ there is $L_{2}$ such that the 
following 
hold. Let $\Gamma \subset S$ be a  noncollapsing maximal multicurve
with $\vert f_{2}(\gamma )\vert \leq 
L_{1}$ for all $\gamma \in \Gamma $.

Then $\vert \gamma _{*}\vert \leq 
L_{2}$ for all $\gamma \in \Gamma $, where $\gamma _{*}$ denotes the 
geodesic in $N$ homotopic to $f_{2}(\gamma )$. 

There is  an embedded surface  $S_{3}$ in $N$ whose preimage in 
$H^{3}$ 
bounds a convex subset of $H^{3}$, such that inclusion 
$f_{3}:S\to N$ is homotopic to inclusion $f_{2}:S\to \overline{N}$, 
such that, using the metric on $S_{3}$ induced by the hyperbolic 
metric on $H^{3}$,
$$L_{2}^{-1}\leq {\vert f_{3}(\gamma )\vert \over \vert 
f_{2}(\gamma )\vert }\leq L_{2}$$
for all $\gamma \in \Gamma $.

There is  a pleated surface 
$f_{4}$ with pleating locus a maximal multicurve $\Gamma '$,
such that 
\begin{description}
    \item[.] $\# (\Gamma \cap \Gamma ')\leq L_{2}$, 
    \item[.] $n_{\gamma ,\beta }([f_{2}])\leq L_{2}\vert f_{2}(\gamma 
    )\vert ^{-1}$ for all $\beta \in \Gamma '$ transverse to $\gamma 
    \in \Gamma $ (where $n_{\gamma ,\beta }(.)$ is as in \ref{2.5}),
    \item[.] $f_{4}$ has no 
badly bent annuli, 
\item[.] there is a homotopy between $f_{3}$ and $f_{4}$
 with homotopy tracks of  hyperbolic length $\leq 
L_{2}$,
\item[.]
$$d([f_{2}],[f_{4}])\leq L_{2}.$$
\end{description}

\end{ulemma}

\ssubsection{Proof of \ref{3.9} in the case of bounded geometry on 
$S_{2}$.}\label{4.8}
 To start with, we assume that $[f_{2}]\in 
(\cal{T}(S))_{\geq \nu }$ for some $\nu >0$, and obtain estimates in 
terms of $\nu $. We can cover $S_{2}$ by $\leq n(\nu )$ topological 
discs, which are round discs up to bounded distortion in 
the hyperbolic metric on $S_{2}$ (independent of $\nu $), all with 
boundedly proportionate radii in terms 
of $\nu $,  and lift to Euclidean discs in $\Omega $.  

We can 
also assume that Euclidean discs of half the Euclidean radii still 
cover $\Omega $. The discs cover $\Omega $ 
with index $\leq n(\nu )$, assuming the discs have radius less than 
half the radius of injectivity. Each Euclidean disc is the base of a 
Euclidean half-ball in $H^{3}$, using the half-space model for 
$H^{3}$. The intersection of all the complementary half-balls is a 
convex set. Fix one boundary component, $U$. Then the stabilisier of 
$U$ in $\pi _{1}(N)$ is, up to conjugacy in $\pi _{1}(N)$,
$i_{*}(\pi _{1}(S))$, where 
$i:S\to N$ denotes inclusion. One component $\Omega _{2}$ of 
$\Omega $ is separated by $U$ from the convex hull of $\pi _{1}(N)$, 
is stablised by  $i_{*}(\pi _{1}(S))$, and covers $S_{2}$.
Because Euclidean discs of half the radius 
still cover,  each hemisphere boundary of a half-ball which 
intersects 
$U$ 
does so in a set of bounded hyperbolic diameter, and each bounded 
loop on $S_{2}$ has a lift to the upper boundary of the hemispheres 
such 
that the projection to $U/i_{*}(\pi _{1}(S))$ is of bounded 
hyperbolic length. This is our surface $S_{3}$. We write $f_{3}:S\to 
S_{3}$ for the incusion map.

Now let $\Gamma $ be given. For each $\gamma \in \Gamma $, if 
we write $\gamma $ also for the corresponding element of $i_{*}(\pi 
_{1}(S))$, 
there 
is $x=x(\gamma )\in U$ such that the hyperbolic distance between $x$ 
and $\gamma .x$ in $U$ is $\leq L_{2}=L_{2}(L_{1},\nu )$. Hence, 
$\vert \gamma _{*}\vert \leq L_{2}$. It could be that $\vert \gamma 
_{*}\vert <\varepsilon _{0}$, where $\varepsilon _{0}$ is a fixed 
Margulis constant. But if so, we can find a loop $\beta (\gamma )$ 
such 
that $\beta (\gamma )$ intersects $\gamma $ at most twice, but no 
other loop in $\Gamma $, and $\vert f_{2}(\beta (\gamma ))\vert \leq 
L_{1}'$ 
where $L_{1}'$ depends only on $L_{1}$. If our first choice of $\beta 
(\gamma )$ is trivial in $N$, then we can make it nontrivial, just by 
composing with a single twist round $\gamma $. This must be 
nontrivial, because otherwise $\gamma $ is trivial. Then $\beta 
(\gamma )$ also 
has a lift to $U$ of length $\leq L_{2}$, assuming that 
$L_{2}$ large enough given $L_{1}$. Also, the points $x(\gamma )$ and 
$x(\beta (\gamma ))$ can be chosen a bounded hyperbolic distace 
apart. Then $\vert (\beta (\gamma 
))_{*}\vert \leq L_{2}$. If both $\vert \gamma _{*}\vert $ and 
$\vert (\beta (\gamma ))_{*}\vert $ are small, then their 
$\varepsilon _{0}$-Margulis tubes are disjoint, and both come within 
a 
bounded distance of the lifts of $\gamma $, $\beta (\gamma )$ 
respectively on $U$. So they cannot both be above $U$. But they do 
not 
intersect $U$. But the region above $U$ is convex, and its closure 
contains the limit set. So both $\gamma _{*}$ and $(\beta 
(\gamma ))_{*}$ must be above $U$. The loop $(\beta (\gamma ))_{*}$ 
can also not 
be separated from $U$ by the Margulis tube of $\gamma _{*}$, if there 
is one, because it is a bounded distance from $U$. So if $\gamma $ is 
short, $\beta (\gamma )$ is not, and is not separated from $U$ by the 
orbit of $\gamma $. Now we change $\Gamma $ to $\Gamma '$ through a 
sequence of maximal multicurves $\Gamma _{i}$, $0\leq i\leq r$, with 
$\Gamma 
_{0}=\Gamma $ and $\Gamma _{r}=\Gamma '$. At each stage, we replace 
some short loop $\gamma \in '\Gamma _{i}\cap \Gamma $ by a loop 
$\beta (\gamma )$ disjoint from all loops of $\Gamma _{i}\setminus 
\{ \gamma \} $, such that $\vert f_{2}(\beta (\gamma )\vert $ is 
bounded 
and bounded from $0$. We can choose the $\beta (\gamma )$ at each 
stage so that $\# (\Gamma _{i}\cap \Gamma _{j})\leq L_{2}$ for all 
$i$, $j$. 
In particular, $\# (\Gamma \cap \Gamma ')\leq L_{2}$. 
We then take 
$f_{4}$ to be a pleated surface with pleating locus containing 
$\Gamma '$.
We have bounds on $\vert f_{4}(\gamma )\vert $ for all $\gamma \in 
\cup _{i}\Gamma _{i}$  because
either $\gamma \in \Gamma '$ or $\# (\gamma 
\cap \Gamma ')$ is bounded, in which case we can apply \ref{3.3.5}. 
We 
also have bounds on $\vert f_{4}(\beta (\gamma ))\vert $ for all 
$\gamma \in \Gamma \cap \Gamma '$, again by applying \ref{3.3.5}. So 
we have
$$d([f_{2}],[f_{4}])\leq L_{2}.$$

In order to construct a homotopy with bounded tracks between $f_{3}$ 
and $f_{4}$, use the same 
procedure as in \ref{4.3}, which is itself derived from \cite{Min1}. 
The loop set $\Gamma '$ has the property that $f_{3}(\gamma )$ and 
$f_{4}(\gamma )$ are a bounded hyprbolic distance apart, and 
$\Gamma '$ is a maximal multicurve. We are assuming, here, that 
$f_{3}(\gamma )$ and $f_{4}(\gamma)$ are geodesics in the respective 
Poincar\'e metrics.  So we take the homotopy to 
homotope $f_{3}(\gamma )$ to $f_{4}(\gamma )$, for all $\gamma \in 
\Gamma '$. Then, as in \ref{4.3}, we foliate each complementary pairs 
of pants by arcs, except for two tripods. Because of bound on 
$d([f_{3}],[f_{4}])$, bounded arcs on $f_{3}(S)=S_{3}$ are homotopic 
to bounded arcs on $f_{4}(S)$. Then we homotope tripods to tripods, 
and 
arcs in between as dictated by the endpoints of arcs. 

\ssubsection{Proof of \ref{3.9} in general.}\label{4.10}

Now suppose that $\vert f_{2}(\gamma )\vert $ is small for some loop 
$\gamma \in 
\Gamma $. Then we construct $U$ over the preimage in $\Omega 
_{2}$ 
of $(S_{2})_{\geq \varepsilon _{1}}$ in the same way as before. But 
over   $(S_{2})_{<\varepsilon _{1} }$, for a sufficiently small 
$\varepsilon _{1}$, we change the construction. As 
before, let $\gamma $ denote both a loop with $\vert f_{2}(\gamma 
)\vert <\varepsilon _{1}$, and a corresponding element of the 
covering 
group, which leaves invariant a component of  the preimage of 
$T(f_{2}(\gamma ),\varepsilon _{0})\subset S_{2}$ in $\partial 
H^{3}$. 
We identify 
$\partial H^{3}$ with $\overline{\mathbb C}$ in the usual way. Let 
$\lambda $ be the complex length of $\gamma $. Now $\gamma $ embeds 
in 
a $1$-parameter subgroup $\gamma _{t}$ of $PSL(2,\mathbb C)$
which acts on $\overline{\mathbb C}$. The eigenvalues of $\gamma 
_{t}$, considered as an element of $PSL(2,\mathbb C)$, are 
$e^{\pm t\lambda /2}$. The action of the group $\gamma _{t}$ has two 
fixed points $x_{1}$, $x_{2}$, the endpoints of the geodesic 
$\gamma _{*}$ in $H^{3}$ lifting  
the loop $\gamma $. Let $V$ be the connected component of the 
preimage of $T(f_{2}(\gamma ),\varepsilon _{0})\subset S_{2}$ in 
$\overline{\mathbb C}$ whose closure contains the points $x_{1}$ and 
$x_{2}$. The closure of $V$ in $\overline{\mathbb C}$ is then the 
union of the lift of $\overline{T(f_{2}(\gamma ),\varepsilon 
_{0})}\subset S_{2}$ 
and $x_{1}$, $x_{2}$. 
Normalise $\overline{\mathbb C}$ so that $\overline{V}\subset \mathbb 
C$. In fact, we can, and shall, 
normalise so that the diameters are both boundedly proportional to 
$1$.
The orbits of the $\gamma _{t}$ action are spirals connecting $x_{1}$ 
and $x_{2}$. 

If $C_{1}$ is sufficiently large, a 
ball of Poincar\'e radius $\varepsilon _{0}$ centred on a point $x$ 
of 
$(S_{2})_{\leq \varepsilon _{0}/C_{1}}$ lifts to 
$V$ in such a way that the disc $D(\tilde x)$ 
of Euclidean radius 
$2\vert \gamma .\tilde x-\tilde x\vert $ centred on any lift $\tilde 
x\in V$ of $x$ is contained in the lift. Assuming $C_{1}$ is 
sufficiently large, the section of spiral between $\tilde x$ and 
$\gamma .\tilde x$ is contained in this disc. By invariance, the 
whole 
spiral must be contained in $V$, as claimed.  One interesting feature 
of this is that the spiral only intersects $D(\tilde x)$ in one arc. 
If it intersected in more than one arc, then one endpoint of the 
geodesic $\gamma $ would be an isolated point in the closure of $V$, 
which is impossible, because the limit set is a perfect set. So now 
let $V_{1}$ be the connected union of spirals.  Then the closure 
$\overline{V_{1}}$, apart from the endpoints $x_{1}$, $x_{2}$ is 
contained in $V$.
We can further normalise $\overline{\mathbb C}$ so that the two 
components of 
$\partial V_{1}\setminus \{ x_{1},x_{2}\} $, as 
for $\partial V$, have 
proportional Euclidean diameters, and so that the Euclidean diameter 
of $V_{1}$ is $1$.  The modulus of the annulus $V_{1}/<\gamma >$ is
boundedly proportional to the modulus of $T(f_{2}(\gamma ),\varepsilon
_{0})$, and hence inversely proportional, for a multiplicative
constant depending on $\varepsilon _{0}$, to $\vert f_{2}(\gamma
)\vert $.  Now let $\tilde \gamma $ be the spiral $\{ \gamma
_{t}.\tilde x:t\in \mathbb R\} $ for some $\tilde x\in V_{1}$ which is
in the lift of $f_{2}(\gamma )$ in $V$.  The Euclidean distance
between spirals for different choices of $\tilde x$ is bounded by
$C_{0}\vert x_{1}-x_{2}\vert $, for a universal constant $C_{0}$.  The
modulus of $V_{1}/<\gamma >$ is also boundedly proportional to
${\rm{Max}}(1/d_{1},1/d_{2})$, where $d_{1}$ is the maximum Euclidean
distance between $\tilde{x}$ and $\gamma .\tilde x$, for $\tilde x\in
\partial V_{1}$, and $d_{2}$ is similarly defined for $\tilde x\in
\tilde \gamma $.  In fact, ${\rm{Max}}(1/d_{1},1/d_{2})$ is always
boundedly proportional to $1/d_{2}$, and $d_{1}$ and $d_{2}$ are
boundedly proportional if and only if $\vert x_{1}-x_{2}\vert $ is
bounded from $0$.

 Now we claim that there is a $\gamma 
$-invariant set 
$V_{1}'$ with $V_{1}'\subset V_{1}$, 
such that the Poincar\'e distance between $V_{1}'$ and $V_{1}$ is 
bounded, and $V_{1}'$ is the $\gamma $-orbit of 
between one and five Euclidean discs, where two of the discs have 
Euclidean radius which are bounded and bounded to $0$, and the 
smallest 
has Euclidean radius boundedly proportional to $d_{2}$.
We see this as follows. We can transfer back under a 
M\"obius transformation $\sigma $  for the 
moment, to the situation when $x_{1}=0$ and $x_{2}=\infty $. So under 
this transformation, $V_{1}$ transforms to  
$$V_{0}=\sigma (V_{1})=\{ e^{\lambda t+i\theta }:\theta \in [0,\theta 
_{1}]\} $$
for some $0<\theta _{1}<2\pi $. We can assume without loss of 
generality that ${\rm{Re}}(\lambda )<0$ (interchanging $0$ and 
$\infty $ if necessary). The element $\gamma $ 
transfers to multiplication by $e^{\lambda }$, and $\lambda $ is 
small. 
We only need to produce an orbit of up to five for this set,
because our normalisation of $V_{1}$ 
has $V_{1}$ bounded, and so discs will transfer to discs under 
$\sigma $. If we can do so, we choose just one disc in $\sigma 
(V_{1})$ 
which is of Euclidean diameter which is bounded and bounded from $0$, 
and tangent to both components of $\sigma (\partial V_{1})\setminus 
\{ 
0,\infty \} $. This is possible if $\theta _{1}<\pi /3$, and also for 
any $\theta _{1}<2\pi $, if ${\rm{Im}}(\lambda 
)/{\rm{Re}}(\lambda )$ is bounded from $0$. In general, we can always 
choose a connected union $V_{0}''$ $\leq 5$ discs of 
Euclidean diameters which are bounded, and 
the outer two bounded from $0$, and each tangent to one of the 
components of $\sigma (\partial V_{1})\setminus \{ 
0,\infty \} $, at points which are bounded from $0$. We can choose 
$V_{0}''$ so that $V_{0}\setminus (V_{0}''\cup e^{\lambda }V_{0}'')$ 
has at most two components whose closure does not include $0$, 
$\infty $, at most one intersecting each component of $\partial 
V_{0}\setminus \{ 0,\infty \} $. If more than two discs are needed, 
we can also choose all but the outer two so that each one intersects 
the 
two nearest discs in two sets, each of diameter proportional to its 
own diameter.  Let $V_{0}'=\sigma 
(V_{1}')$ be the orbit of the union of these discs under 
multiplication by $e^{\lambda }$. Then $V_{1}\setminus V_{1}'$ is 
contained in a neighbourhood of $\partial V_{1}$ of bounded 
Poincar\'e diameter, with respect to the Poincar\'e metric on $\Omega 
_{2}$.

Then we take $U$ to be the surface bounded 
by the union of the hemispheres over the lift of $(S_{2})_{\geq 
\varepsilon _{1}}$
and over $V_{1}'$, for a suitable $\varepsilon _{1}$ with 
$\varepsilon _{0}/\varepsilon _{1}$ bounded, so that the sets 
$V_{1}'(\tilde \gamma )$, for varying $\tilde \gamma $, and the lifts 
of $(S_{2})_{\geq \varepsilon _{1}}$, cover $\Omega _{2}$.  Then 
$S_{3}=U/i_{*}(\pi _{1}(S))$ is an embedded surface, with full 
preimage in $H^{3}$ bounding a convex subset of $H^{3}$, with 
corresponding map $f_{3}:S\to S_{3}$ homotopic, as a map from $S$ to 
$\overline{N}$, to $f_{2}:S\to S_{2}$. To find the shortest loop in 
$S_{3}$
homotopic to $f_{3}(\gamma )$, we can draw paths on the tops of the 
hemispheres. As before, we can estimate hyperbolic length by the 
ratio 
of the Euclidean length on the top of a hemisphere to the 
Euclidean radius. The paths on tops of all but the outer hemispheres 
over $V_{1}'$ 
are bounded, because of the conditions we imposed on the 
intersections 
of the base hemispheres. As for the outer hemispheres (if $V_{1}'$ is 
the orbit of more than one hemisphere) a path across the top of an 
outer hemisphere $H_{1}$ also has 
bounded hyperbolic length. But the distance to the nearest hemisphere 
over   $(S_{2})_{\geq \varepsilon _{1}}$ might be much greater, if 
$H_{1}$ has a
much larger Euclidean diameter than the nearby 
hemispheres over $(S_{2})_{\geq \varepsilon _{1}}$.
In any case, the hyperbolic length of $f_{3}(\gamma )$ is boundedly
proportional to $d_{2}$, and to 
$\vert f_{2}(\gamma 
)\vert $. Also by taking paths on the tops of hemispheres, we can 
choose
$\beta (\gamma )$ 
transverse to $\gamma $ such that  $f_{3}(\beta (\gamma ))$ is, up to 
homotopy in $S_{3}$, a union of boundedly finitely many paths of 
bounded length 
and, at most two long geodesic segments in the hemispheres over 
$V_{1}'$ 
for each crossing of $\gamma $.
 We can ensure that the 
long segments do not cancel by adjusting by a Dehn twist round 
$\gamma $, if necessary, but  so that $n_{\gamma ,\beta (\gamma 
)}([f_{2}])$ is bounded by $O(\vert f_{2}(\gamma )\vert ^{-1})$.  
So then $f_{3}(\gamma )$ 
is, up to homotopy, a  bounded distance from $(\beta (\gamma ))_{*}$. 
So then, as before, we can form $\Gamma '$ with $\# (\Gamma \cap 
\Gamma ')\leq L_{2}$, from a sequence $\Gamma _{i}$ with $\Gamma 
_{0}=\Gamma $, $\Gamma _{r}=\Gamma '$, and such that $\Gamma _{i+1}$ 
is obtained from $\Gamma _{i}$ by replacing some loop $\gamma $ with 
$\vert f_{2}(\gamma )\vert $ small by a loop $\beta (\gamma )$.
 Then there is a bounded track homotopy 
between $f_{3}$ and $f_{4}$. The homotopy is defined restricted to 
the 
set with image in $N_{\geq \varepsilon _{0}}$ exactly as before, and 
the 
sets $V_{2}$ have been 
constructed so as to ensure that there is a bounded track homotopy on 
the set $S_{3}\cap N_{<\varepsilon _{0}}$. It follows that the moduli 
of $T(f_{4}(\gamma ),\varepsilon _{0})\subset S(f_{4})$ and 
$T(f_{2}(\gamma ),\varepsilon _{0})\subset S_{2})$ are boundedly 
proportional, and that $\vert f_{2}(\beta (\gamma ))\vert '$ and 
$\vert f_{4}(\beta (\gamma ))\vert '$ are boundedly proportional. The 
bound on $d([f_{2}],[f_{4}])$ follows.

\Box 

\ssubsection{Generalised pleated surfaces.}\label{4.5}

A pleated surface $f:S\to N$ derives its metric from the hyperbolic 
metric on $N$. It will sometimes be useful to use the metric on 
$\partial N$, where $\partial N$ is the boundary if $N$ obtained by 
projecting the Poincar\'e metric from the domain of discontinuity 
$\Omega \subset 
\partial H^{3}$. This is, in fact, the case for $f_{2}:S\to S_{2}$ in 
\ref{3.9}. Write $\overline{N}=N\cup \partial N$. We consider maps 
$f:S\to N$, where $f$ maps the loops of a multicurve $\Gamma $ to 
cusps in $N$, and, for each component $\alpha $ of $S\setminus \cup 
\Gamma $, either $f\vert \alpha $ is a pleated surface, or $f\vert 
\alpha $ is a homeomorphism onto a component of $\partial N$. Then 
$[f\vert \alpha ]$ is an element of ${\cal{T}}(S(\alpha ))$ for each 
$\alpha $. If $f(\alpha )\subset \partial N$, we use the Poincar\'e 
metric on $S(\alpha )$ to define the element of 
${\cal{T}}(S(\alpha ))$. We shall refer to such a map $f$ as a 
{\em{generalised pleated surface}}. It defines an element $[f]\in 
{\cal{T}}(S(\omega ))$, where $\omega =S\setminus (\cup \Gamma )$.
Applying \ref{3.9} to each 
$S(\alpha )$ with $f(\alpha )\subset \partial N$, there is a genuine 
pleated surface a bounded $d_{\omega }$ distance away.

\section{Teichm\"uller geodesics: long thick and dominant 
definitions.}\label{5}

In this section we explain and expand some of the ideas of long thick 
and dominant (ltd) segments of geodesics in Teichm\"uller space
$\cal{T}(S)$ which were used in \cite{R1}. The theory of \cite{R1} 
was 
explicitly for marked spheres only, because of the application in 
mind, but in fact the theory works without adjustment for any finite 
type surface, given that projections $\pi _{\alpha }$ to smaller 
Teichm\"uller spaces $\cal{T}(S(\alpha ))$ for subsurfaces $\alpha $ 
of $S$ have been 
defined in \ref{2.5}. For proofs, for the most part, we refer to 
\cite{R1}.
  The basic idea is to 
get into a position to apply arguments which work along geodesics 
which never enter the thin part of Teichm\"uller space, by projecting 
to suitable subsurfaces $\alpha $ using the projections $\pi 
_{\alpha }$ of \ref{2.5}. A reader who wishes to get to the proof of 
the Ending Laminations Theorem in the case of combinatorial bounded 
geometry is advised to read to the end of the basic definition 
\ref{5.3}, and then proceed to a recommended menu from Section 
\ref{6}. 
In the theory of Teichm\"uller geodesics which is 
developed here (and earlier, in \cite{R1})
it does not seem to make sense to consider geodesics 
in the thick part of Teichm\"uller space --- which is what 
combinatorially bounded geometry means --- in strict isolation.
We use  the basic notation and theory of Teichm\"uller space 
$\cal{T}(S)$ from Section \ref{2}.

\ssubsection{Good position.}\label{5.1} 

Let $[\varphi ]\in 
\cal{T}(S)$. Let $q(z)dz^{2}$ be a quadratic 
differential on $\varphi (S)$. All quadratic differentials, as in 
\ref{2.2}, will be of total mass $1$. Let $\gamma $ be a  
nontrivial nonperipheral simple closed loop on $S$. 
Then there is a limit of isotopies of 
$\varphi (\gamma )$ to 
{\em{good position}} with respect to $q(z)dz^{2}$, with the limit 
possibly passing through some punctures. If $\gamma $ is the 
isotopy limit, then either $\gamma $ is at constant angle to 
the stable and unstable foliations of $q(z)dz^{2}$, or is a union of 
segments between singularities of $q(z)dz^{2}$ which are at constant 
angle to the stable and unstable foliations, with angle $\geq \pi $ 
between any two consecutive segments at a singularity, unless it is a 
puncture.  An 
equivalent statement is that $\gamma $ is a geodesic with respect to 
the singular Euclidean metric $\vert q(z)\vert d\vert z\vert ^{2}$. 
If two good positions do not coincide, then they bound an open 
annulus 
in $\varphi (S)$ which contains no singularities of $q(z)dz^{2}$.
See also 14.5 of \cite{R1}. 

 The 
{\em{q-d length }} $\vert \varphi (\gamma )\vert _{q}$ is length with 
respect to the quadratic differential metric for any homotopy 
representative in good position. (See 14.5 of \cite{R1}.) 
We continue, as in Section \ref{2}, to use $\vert \varphi (\gamma 
)\vert $ 
to denote the hyperbolic, or Poincar\'e, length on $\varphi (S)$ of 
the geodesic on $\varphi (S)$ homotopic to $\varphi (\gamma )$. 
If $[\varphi ]\in \cal{T}_{\geq 
\varepsilon }$ then there is a constant $C(\varepsilon )>0$ such that 
for all nontrivial nonperipheral closed loops $\gamma $,
$${1\over C(\varepsilon )}\leq {\vert \varphi (\gamma )\vert 
_{q}\over \vert \varphi (\gamma )\vert }\leq C(\varepsilon ).$$
We also define $\vert \varphi (\gamma )\vert _{q,+}$ to be the 
integral of the norm of the projection of the derivative of $\varphi 
(\gamma )$ to the tangent space of the unstable foliation of 
$q(z)dz^{2}$, and similarly for $\vert \varphi (\gamma )\vert 
_{q,-}$. So these are both majorised by $\vert \varphi (\gamma 
)\vert 
_{q}$, which is, in turn, majorised by their sum.

\ssubsection{Area.}\label{5.2}

The following definitions come from 9.4 of \cite{R1}. 
For any essential nonannulus 
subsurface  $\alpha \subset S$, $a(\alpha 
,q)$ is the area with respect to $q(z)dz^{2}$ of $\varphi (\alpha )$ 
where $\varphi (\partial \alpha )$ is in good position and bounds the 
smallest area possible subject to this restriction. If $\alpha $ is a 
loop at $x$ then $a(\alpha ,q)$ is the smallest possible area of an 
annulus of modulus $1$ and homotopic to $\varphi (\alpha )$. We are 
only interested in this quantity up to a bounded multiplicative 
constant and it is also boundedly proportional to $\vert\varphi 
(\alpha 
)\vert _{q}^{2}$ whenever $\varphi (\alpha )$ is in good position, 
and $\vert \varphi (\alpha )\vert $ is bounded. 
We sometimes write $a(\alpha ,x)$ or even 
$a(\alpha )$ for $a(\alpha ,q)$, if it is clear from the context what 
is meant.

Generalising from \ref{5.1}, there is a constant $C(\varepsilon )$ 
such that, if 
$\varphi (\alpha )$ is homotopic to a component of $(\varphi 
(S))_{\geq \varepsilon }$, then for all nontrivial nonperipheral 
non-boundary-homotopic closed loops $\gamma \in \alpha $,
$${1\over C(\varepsilon )}\vert \varphi (\gamma )\vert \leq {\vert 
\varphi (\gamma )\vert _{q}\over \sqrt{a(\alpha ,q)+a(\partial \alpha 
,q)}}\leq C(\varepsilon )\vert \varphi (\gamma )\vert .$$

Now suppose that $\ell $ is a directed geodesic segment in 
$\cal{T}(S)$
containing 
$[\varphi ]$, and that $q(z)dz^{2}$ is the quadratic differential 
at $[\varphi ]$ for $d([\varphi ],[\psi ])$ for any $[\psi ]$ in the 
positive direction along $\ell $ from $[\varphi ]$ (see \ref{2.1}.) 
Let $p(z)dz^{2}$ be the stretch of $q(z)dz^{2}$ at $[\psi ]$, and 
let $\chi $ be the minimum distortion map with $[\chi \circ \varphi 
]=[\psi ]$. Then $\chi $ maps the $q$-area element to the $p$-area 
element.
Then $a(\alpha ,q)=a(\alpha ,p)$ if $\alpha $ is a gap, but  if 
$\alpha $ is a loop, $a(\alpha ,y)$ varies for $y\in \ell $.

If $\alpha $ is a loop we also make an extra definition. We define 
$a'(\alpha ,[\varphi ], q)$ (or simply $a'(\alpha )$ if the context 
is 
clear) to be the $q$-area of the largest modulus annulus (possibly 
degenerate) homotopic to $\varphi (\alpha )$ and with boundary 
components in good position for $q(z)dz^{2}$. Then $a'(\alpha )$ is 
constant along the geodesic determined by $q(z)dz^{2}$.

\ssubsection{The long thick and dominant definition}\label{5.3}

Now we fix parameter functions $\Delta $, $r$, $s:(0,\infty )\to 
(0,\infty )$ and a constant $K_{0}$. 

Let $\alpha $ be a gap. Let $\ell $ be a geodesic segment.
We say that $\alpha $ is {\em{long, $\nu $-thick and dominant at }} 
$x$ 
(for $\ell $, and with respect to $(\Delta ,r,s)$) if $x$ is the 
centre of a segment $\ell _{1}$ in the geodesic extending $\ell $ of  
length $2\Delta (\nu )$ 
such that $\vert \psi (\gamma )\vert \geq \nu $ for all $[\psi ]\in 
\ell _{1}$ and nontrivial nonperipheral $\gamma \subset \alpha $ not 
homotopic to boundary components, but $\ell _{1}\subset {\cal 
T}(\partial 
\alpha ,r(\nu ))$ and $a(\partial \alpha ,y)\leq s(\nu )a(\alpha, y)$ 
for all $y\in \ell _{1}$. We shall then also say that $\alpha $ is 
long $\nu $-thick and dominant along $\ell _{1}$. See 15.3 of 
\cite{R1}. 

A loop $\alpha $ at 
$x$ is $K_{0}$-{\em{flat at}}  $x=[\varphi ]$ (for $\ell $) if 
$a'(\alpha )\geq K_{0}a(\alpha)$.  This was not quite the definition 
made in 
\cite{R1}  
where the context was restricted to $S$ being a punctured sphere, 
but the results actually worked for any finite type surface. The term 
arises because if $\alpha $ is $K_{0}$-flat then the metric $\vert 
q(z)\vert dz^{2}$ is equivalent to a Euclidean (flat) metric on an 
annulus homotopic to $\varphi (\alpha )$ of modulus $K_{0}-O(1)$.  
For 
fixed $K_{0}$ we may simply 
say 
{\em{flat}} 
rather than $K_{0}$-flat. In future, we shall often refer to 
prarmeter 
functions as quadruples of the form $(\Delta ,r,s,K_{0})$.

If $\alpha $ is long $\nu $-thick and dominant along a segment $\ell 
$, that 
is long thick and at all points of $\ell $, then $d_{\alpha }(x,y)$ 
is very close to $d(x,y)$ at all points of $\ell $. This is a 
consequence of the results of Section 11 of \cite{R1}. All we care 
about for the moment is that they differ by some additive constant. 
It 
is also probably worth noting (again by the results of Chapter 11 of 
\cite{R1}) that if $[\varphi ]\in \ell $ and
$\pi _{\alpha }([\varphi ])=[\varphi _{\alpha }]$, then $\varphi 
_{\alpha }(S(\alpha ))$ and the component $S(\alpha ,r(\nu )[\varphi 
])$ of 
$(\varphi (S))_{\geq 
r(\nu )}$ homotopic to $\varphi (\alpha )$ are isometrically very 
close, 
except in small neighbourhoods of some punctures, and the quadratic 
differentials $q(z)dz^{2}$ at $[\varphi ]$ for $d([\varphi 
],[\psi ])$ ($[\psi ]\in \ell$) and the quadratic differential 
$q_{\alpha }(z)dz^{2}$ at $[\varphi _{\alpha }$ for $d_{\alpha 
}([\varphi ],[\psi ])$, are very close.

\ssubsection{}\label{5.4}
Before starting to describe the usefulness of long thick and 
dominants, we need to show they exist, in some abundance.
This was the content of the first basic result about long thick and 
dominants  in 15.4 of \cite{R1}, which was stated only for $S$ being 
a puntured sphere, but the 
proof worked for a general finite type surface. 

\begin{ulemma} For some $\nu _{0}$ and $\Delta _{0}$ 
depending 
only on $(\Delta ,r,s,K_{0})$ (and the topological type of $S$), 
  the following holds. Any 
geodesic segment $\ell $ of length $\geq \Delta _{0}$ contains a 
segment $\ell '$ for which there is $\alpha $ such that:
\begin{description}
\item[.] either 
$\alpha $ is a gap  which is long $\nu $-thick and dominant along 
$\ell '$ for some $\nu \geq \nu _{0}$ and $a(\alpha 
)\geq 1/(-2\chi (S)+1)=c(S)$ (where $\chi $ denotes Euler 
characteristic,
\item[.] or  $\alpha $  is a $K_{0}$-flat loop along $\ell '$.
\end{description}
More generally  there 
is $s_{0}$ depending 
only on $(\Delta ,r,s,K_{0})$ (and the topological type of $S$) 
such that, whenever $\omega \subset S$ is such that 
$a(\partial \omega )\leq s_{0}a(\omega )$, then we can find $\alpha $ 
as above with $\alpha \subset \omega $ and 
$a(\alpha )\geq 1/(-2\chi (\omega )+1)a(\omega)$ if $\alpha $ is a 
gap.

\end{ulemma}

\noindent {\em{Proof.}} (See also 15.4 of \cite{R1}.) We consider the 
case $\omega =S$.
Write $r_{1}(\nu )=e^{-\Delta (\nu )}r(\nu )$. 
Let $g=-2\chi (S)$  and let
$r_{1}^{g}$ denote the $g$-fold iterate. We then take 
$$\nu _{0}=r_{1}^{g}(\varepsilon _{0})$$ 
for a fixed Margulis constant 
$\varepsilon _{0}$ and we define 
$$\Delta _{0}=2\sum _{j=1}^{g}\Delta (r_{1}^{j}(\varepsilon 
_{0})).$$ 
Then for some $j\leq g$, we can find $\nu =r_{1}^{j}(\varepsilon 
_{0})$ and $[\varphi ]\in \ell $ such that the segment $\ell '$ of 
length 
$2\Delta (\nu )$ centred on $y=[\varphi ]$ is contained in $\ell $ 
and such that 
for   any nontrivial nonperipheral loop 
$\gamma $, either $\vert \varphi '(\gamma )\vert \geq \nu $ for all 
$[\varphi ']\in \ell '$, or $\vert 
\varphi (\gamma )\vert \leq r_{1}(\nu )$ --- in which case 
$\vert 
\varphi '(\gamma )\vert \leq r(\nu )$ for all $[\varphi ']\in \ell 
'$. 
Suppose there are no 
$K_{0}$-flat loops at $[\varphi ]$, otherwise we are done.
For any loop $\gamma $ with $\vert \varphi (\gamma )\vert 
<\varepsilon _{0}$, 
if $\beta $ is a gap such that $\gamma \subset \partial \beta $ and 
there 
is a  component of $(\varphi (S))_{\geq 
\varepsilon _{0} }$ homotopic to $\varphi (\beta )$ and separated 
from the flat 
annulus homotopic to $\varphi (\gamma )$ by an annulus of modulus 
$\Delta _{1}$, we have, since every zero of $q(z)dz^{2}$ has order at 
most $2g$, for a constant $C_{1}$ depending only on the topological 
type of $S$,
\begin{equation}\label{5.4.1}C_{1}^{-1}a(\gamma ,[\varphi ])e^{\Delta 
_{1}}
    \leq a(\beta )\leq 
C_{1}a(\gamma ,\varphi ])e^{(2g+2)\Delta _{1}}.\end{equation}

Now let $\alpha $ be a subsurface such that $\varphi (\alpha )$ is 
homotopic to a component $S(\alpha ,\nu )$ of 
$(\varphi (S))_{\geq \nu }$ of maximal area. Then by (\ref{5.4.1}), 
we have a bound of $O(e^{(2g^{2}+2g)/\nu })$ 
on the ratio of areas of any two components of $(\varphi (S))_{\geq 
\varepsilon _{0} }$ in $S(\alpha ,\nu )$ and assuming $r(\nu )$ is 
sufficiently small given $\nu $,
$$a(\partial \alpha ,[\varphi ])\leq e^{-1/(3gr(\nu))}a(\alpha ),$$
and for all $y'\in \ell '$,
$$a(\partial \alpha ,y')\leq e^{\Delta 
(\nu)}e^{-1/(3gr(\nu))}a(\alpha 
).$$ 
Assuming $r(\nu )$ is sufficiently small given $s(\nu )$ and 
$\Delta (\nu )$, $\alpha $ is long $\nu $-thick and dominant along 
$\ell '$ for $(\Delta, r,s)$, and $a(\alpha )\geq 1/(g+1)$

The case $\omega =S$ is similar. We only need $s_{0}$ small enough 
for $a(\partial \omega )/a(\omega )$ to remain small along a 
sufficiently long segment of $\ell $.
\Box

Because of this result, we can simplify our notation. So let $\nu 
_{0}$ 
be as above, given $(\Delta ,r,s,K_{0})$. We shall simply say 
$\alpha $ is {\em{ltd}} (at $x$, or along $\ell _{1}$,  for $\ell $) 
if either $\alpha $ is a gap and long $\nu $-thick and dominant for 
some $\nu \geq \nu _{0}$, or $\alpha $ is a loop and $K_{0}$-flat. 
We shall also say that $(\alpha ,\ell _{1})$ is ltd.

\ssubsection{}\label{5.5}
We refer to Chapters 14 and 15 of \cite{R1} for a summary of all the 
results 
concerning ltd's, where, as 
already stated, the context is restricted to $S$ being a punctured 
sphere, but the results work for any finite type surface. The main 
points about ltd's are, firstly, that they are good coordinates, in 
which 
arguments which work in the thick part of Teichm\" uller space can be 
applied, and secondly that there is only bounded movement in the 
complement of ltd's. This second fact, together with \ref{4.4},
is worth scrutiny. It is, at first 
sight, surprising. It is proved in 15.14 of \cite{R1}, which we now 
state, actually slightly corrected since short interior loops in 
$\alpha $ were forgotten in the statement there (although the proof 
given there does consider short interior loops) and slightly expanded 
in the case of $\alpha $ being a loop. 

\begin{ulemma}    Fix 
long thick and dominant parameter functions 
$\Delta ,r,s,K_{0})>0$, and 
let $\nu _{0}>0$ also be given and sufficiently small. 
Then there exists  $L=L(\Delta ,r,s,K_{0},\nu _{0})$ such that the 
following 
holds. Let  $\ell $ be a geodesic segment and let $\ell _{1}\subset 
\ell $ and, given $\ell _{1}$, let  
$\alpha \subset S$ be a maximal subsurface up 
to homotopy  with the property that $\alpha \times \ell _{1}$ is 
disjoint from all ltd's  $\beta \times \ell '$ such that $\beta $ is 
either $K_{0}$-flat along $\ell '$ or $\nu $-thick long and dominant 
for some $\nu \geq \nu _{0}$, for 
$[\varphi ]\in \ell _{1}$. Suppose 
also that all components of 
$\partial \alpha $ are 
nontrivial nonperipheral. Then $\alpha $  is a disjoint union of gaps 
and loops $\beta $ such that the following hold.
\begin{equation}\label{5.5.1}\vert \varphi (\partial \beta )\vert 
\leq 
L{\rm{\ for\ all\ }}[\varphi ]\in \ell 
.\end{equation} 
If $\beta $ is a gap, then for all $[\varphi ]$, 
$[\psi ]\in \ell $ and nontrivial nonperipheral 
non-boundary-parallel closed loops $\gamma $ in $\beta $,
\begin{equation}\label{5.5.2} L^{-1}\leq {\vert \varphi 
(\gamma 
)\vert '\over 
\vert \psi (\gamma )\vert '}\leq L,\end{equation}
\begin{equation}\label{5.5.5}
    \vert \varphi (\gamma )\vert \geq L^{-1}.
    \end{equation}
If $\beta $ is a loop, then for all $[\varphi ]$, 
$[\psi ]\in \ell $,
\begin{equation}\label{5.5.3} \vert {\rm{Re}}(\pi _{\alpha 
}([\varphi ])-\pi _{\alpha }([\psi ])\vert \leq L.\end {equation}

 Also if $\gamma $ is 
 in the interior of $\alpha $, and 
$\ell _{1}=[[\varphi _{1}],[\varphi _{2}]]$, then given $\varepsilon 
_{1}>0$ there exists $\varepsilon _{2}>0$ depending only on 
$\varepsilon _{1}$ and the ltd parameter functions and flat constant 
such that
\begin{equation}\label{5.5.4}\begin{array}{l}
{\rm{If\ }}\vert \varphi (\gamma )\vert 
< \varepsilon _{2},{\rm{\ then}}\cr
{\rm{Min}}(\vert \varphi _{1}(\gamma )\vert ,\vert \varphi 
_{2}(\gamma )\vert \leq \varepsilon _{1}{\rm{,\ and}}\cr 
{\rm{Max}}(\vert \varphi _{1}(\gamma )\vert ,\vert \varphi 
_{2}(\gamma )\vert \leq  L.\cr
\end{array}\end{equation}

\end{ulemma} 

If (\ref{5.5.1}), and either (\ref{5.5.2}) and (\ref{5.5.5}), or 
(\ref{5.5.3}) 
hold for $(\beta ,\ell _{1})$, depending on whether $\beta $ is a gap 
or a loop, we say that 
$(\beta ,\ell 
_{1})$ is {\em{bounded}} (by $L$). Note that $L$ depends on the ltd 
parameter 
functions, and therefore is probably extremely large compared with 
$\Delta (\nu )$ for many values of $\nu $, perhaps even compared 
with $\Delta (\nu _{0})$. 

Here are some notes on the proof. For fuller details, see 5.14 of 
\cite{R1}. First of all, under the assumption that $\partial \alpha 
$ satisfies the condition (\ref{5.5.1}), it is shown that $\alpha $ 
is a union of $\beta $ satisfying (\ref{5.5.1}) to (\ref{5.5.5}). 
First, we show that (\ref{5.5.2}) holds for all $\gamma \subset 
\alpha $ such that $\vert \varphi _{i}(\gamma )\vert $ is bounded 
from $0$ for $i=1$, $2$, and that (\ref{5.5.4}) holds for $\alpha $.
This is done by breaking $\ell $ into three segments, with $a'(\alpha 
)$ dominated by $\vert \varphi (\partial \alpha )\vert _{q}^{2}$ on 
the two outer segments $\ell _{-}$, $\ell _{+}$, where $q(z)dz^{2}$ 
is 
the quadratic differential at $[\varphi]$ for $\ell $. The middle 
segment $[[\varphi _{-}],[\varphi _{+}]]$  has to be of bounded 
length by 
the last part of \ref{5.4}, since there are no ltd's in $\alpha $ 
along $\ell $. Then  $\vert \varphi (\partial \alpha )\vert _{q}$ is 
boundedly proportional to $\vert \varphi (\partial \alpha )\vert 
_{q,-}$ along $\ell _{-}$, and to  $\vert \varphi (\partial \alpha 
)\vert _{q,+}$ along $\ell _{+}$.  We can obtain (\ref{5.5.2}) along 
$\ell _{+}$, at least for a nontrivial $\alpha _{1} \subset \alpha $ 
for 
which  we can ``lock ''  loops $\varphi (\gamma )$, for which $\vert 
\varphi (\gamma )\vert $ is bounded,  along stable segments 
to $\varphi _{+} (\partial \alpha )$. If $\alpha _{1} \neq \alpha $ 
and 
$\gamma '\subset \partial \alpha _{1} $ is in the interior of $\alpha 
$, 
then either  
$\vert \varphi _{+}(\gamma ')\vert $ is small, or  $\vert \varphi 
_{+}(\gamma ')\vert _{q_{+}}$ is dominated by $\vert \varphi 
_{+}(\gamma 
')\vert _{q_{+},-}$, where $q_{+}(z)dz^{2}$ is the stretch of 
$q(z)dz^{2}$ at $[\varphi _{+}]$. In the case when $\vert \varphi 
_{+}(\gamma ')\vert $ is small, there is some first point $[\varphi 
_{++}]\in \ell _{+}$ for which  $\vert \varphi 
_{++}(\gamma ')\vert _{q_{++}}$ is dominated by $\vert \varphi 
_{++}(\gamma 
')\vert _{q_{++},-}$, where $q_{++}(z)dz^{2}$ is the stretch of 
$q(z)dz^{2}$ at $[\varphi _{++}]$. For this point,  $\vert \varphi 
_{++}(\gamma ')\vert $ is still small, and can be locked to a small 
segment of $\varphi _{++}(\partial \alpha )$. This means that we can 
deduce that $\vert \varphi _{2}(\gamma ')\vert $ is small, giving 
(\ref{5.5.4}). So one proceeds by induction on the topological 
type of $\alpha $,   obtaining (\ref{5.5.2}) and (\ref{5.5.4})
for $\alpha $ from that for 
$\alpha \setminus \alpha _{1}$. Then (\ref{5.5.4}) and (\ref{5.5.2}) 
imply that the set of loops with $\vert \varphi _{1}(\gamma )\vert 
<\varepsilon _{1}$ or $\vert \varphi _{2}(\gamma )\vert <\varepsilon 
_{1}$, for a sufficiently small $\varepsilon _{1}$, do not intersect 
transversally. This allows for a decomposition into sets $\beta $ 
satisfying (\ref{5.5.1}), (\ref{5.5.2}) and (\ref{5.5.5}).
 One then has to remove the hypothesis 
(\ref{5.5.1}) for $\partial \alpha $. This is done by another 
induction, considering successive gaps and loops $\alpha '$ disjoint 
from all ltd's along segments $\ell '$ of $\ell $, with $\vert 
\varphi (\partial \alpha ')\vert \leq \varepsilon _{0}$ for $[\varphi 
]\in \ell '$, possibly with $\partial \alpha '=\emptyset $. One then 
combines the segments and reduces the corresponding $\alpha '$, 
either 
combining two at a time, or a whole succession together, if the 
$\alpha 
'$ are the same along a succession of segments. In finitely many 
steps, one reaches $(\alpha ,\ell )$ finding in the process that 
$\partial \alpha $ does satisfy (\ref{5.5.1}).

As for showing that $\alpha $ satisfies (\ref{5.5.3}), that follows 
from the 
following lemma --- which is proved in 15.13 of \cite{R1}, but not 
formally stated. Note that if $\alpha $ is a loop, $a'(\alpha 
,[\varphi ])$ is constant for $[\varphi ]$ in a geodesic segment 
$\ell $, but $a(\alpha  ,q)$ is proportional to $\vert \varphi 
(\alpha )\vert _{q}^{2}$ (if $q(z)dz^{2}$ is the quadratic 
differential at $[\varphi ]$ for $\ell $), which has at most one 
minimum on the geodesic segment and  otherwise increases or decreases 
exponetially with distance along the segment, depending on whether 
$\vert \varphi (\alpha )\vert _{q}$ is boundedly proportional to 
$\vert \varphi (\alpha )\vert _{q,+}$ or  $\vert \varphi (\alpha 
)\vert _{q,-}$. 
So for any $K_{0}$, the set of $[\varphi ]\in \ell $ for which 
$a'(\alpha )\geq K_{0}a(\alpha ,[\varphi ])$ is a single segment, up 
to bounded distance. This motivates the following.  

\begin{lemma}\label{5.6} Given $K_{0}$, there is $C(K_{0})$ such that 
the following holds. Let $\ell $ be any geodesic segment. Suppose 
that $a'(\alpha )\leq K_{0}a(\alpha ,[\varphi ])$ for all $[\varphi 
]\in 
\ell $. Then for all $[\varphi ]$, $[\psi ]\in \ell $,
$$\vert {\rm{Re}}(\pi _{\alpha }([\varphi ])-\pi _{\alpha }([\psi 
]))\vert \leq C(K_{0}).$$
\end{lemma}

\noindent{\em{Proof.}}
The argument is basically given in 15.13 of 
\cite{R1}. Removing a  segment of length bounded in terms of $K_{0}$, 
$\varepsilon _{0}$ at one end, we obtain a reduced segment $\ell '$ 
such that  that $a'(\alpha )\leq \varepsilon _{0}a(\alpha ,[\varphi 
])$ for all $[\varphi ]\in \ell '$.
We use the 
quantity $n_{\alpha }([\varphi ])$ of \ref{2.5}, which is 
${\rm{Re}}(\pi _{\alpha  }([\varphi ])+O(1)$ and is given to within 
length $O(1)$ 
by $m$ minimising $\vert \varphi (\tau _{\alpha }^{m}(\zeta ))\vert 
$ for a fixed $\zeta $ crossing $\alpha $ at most twice (or a bounded 
number of times). This is the same to within $O(1)$ as the $m$ 
minimising
$\vert \varphi (\tau _{\alpha }^{m}(\zeta ))\vert 
_{q}$ for any quadratic differential $q(z)dz^{2}$. (To see this 
note that the shortest paths, in the Poincar\'e metric, across a 
Euclidean annulus $\{ z:r<\vert 
z<1\} $ are the restrictions of straight lines through the origin.)  
Assume without loss of generality that $\vert \varphi (\alpha )\vert 
_{q}$ 
is boundedly proportional to $\vert \varphi (\alpha )\vert _{q,+}$ 
for 
$[\varphi]\in \ell $, and $q(z)dz^{2}$ the quadratic differential at 
$[\varphi ]$ for $\ell $. 
The good positions of  $\varphi (\tau _{\alpha }^{m}(\zeta ))$ for 
all 
$m$ are locked together along stable segments whose qd-length is 
short in comparision with $\vert \varphi (\alpha )\vert _{q}$, if 
$\varepsilon _{0}$ is sufficiently small. So $n_{\alpha }([\varphi 
])$ 
varies by $<1$ on $\ell $, and is thus constant on $\ell '$, if 
$\varepsilon _{0}$ is sufficiently small, and hence varies by at most 
$C(K_{0})$ on $\ell $.\Box

\ssubsection{Decomposing $S\times \ell $.} \label{5.7}

Let $\alpha _{i}\subset S$ be  a gap or loop for $i=1$, $2$, isotoped 
so 
that $\partial \alpha _{1}$ and $\partial \alpha _{2}$ have only 
essential intersections, or with $\alpha _{1}\subset \alpha _{2}$ 
if $\alpha _{1}$ is a loop which can be isotoped into $\alpha _{2}$. 
Then the {\em{convex hull}} $C(\alpha _{1},\alpha _{2})$ 
of $\alpha _{1}$ and $\alpha _{2}$ is the union of $\alpha _{1}\cup 
\alpha _{2}$ and any components of $S\setminus (\alpha _{1}\cup 
\alpha 
_{2})$ which are topological discs with at most one puncture. Then 
$C(\alpha _{1},\alpha _{2})$ is again a gap or a loop. The latter 
only 
occurs if $\alpha _{1}=\alpha _{2}$ is a loop. We are only interested 
in the convex hull up to isotopy, and it only 
depends on  $\alpha _{1}$ and $\alpha _{2}$ up 
to isotopy. It is so called because, if $\alpha _{i}$ is chosen to 
have geodesic boundary, and $\tilde {\alpha _{i}}$ denotes 
the preimage of $\alpha _{i}$ in the hyperbolic plane covering $S$, 
then up to isotopy $C(\alpha _{1},\alpha _{2})$ is the projection of 
the convex hull of any component of $\tilde{\alpha _{1}}\cup 
\tilde{\alpha _{2}}$.

The following version of \ref{5.5} will be important in constructing 
the geometric model. It follows directly from the statement of 
\ref{5.5}.

\begin{ulemma} Fix ltd parameter functions $(\Delta ,r,s,K_{0})$, and 
an associated constant $\nu _{0}$ as in \ref{5.4}, and $L=L(\Delta 
,r,s,K_{0},\nu _{0})$ as in \ref{5.5}.
Let $\ell $ be any geodesic segment in ${\cal T}(S)$. Then  we can 
write $S\times \ell $ as
$$S \times \ell=\cup _{j=1}^{R} \alpha 
_{j}\times\ell _{j}$$ 
 where each $\alpha _{i}\times \ell _{j}$ is either bounded by $L$, 
 or long $\nu $-thick and dominant along $\ell _{j}$ for $(\Delta 
,r,s)$ 
 and some $\nu 
 \geq \nu _{0}$, or $K_{0}$ flat along $\ell _{j}$, depending on 
 whether $\alpha _{j}$ is a gap or a loop.

In addition the decomposition is {\em{vertically efficient}} in the 
following sense.
\begin{description} 
    \item[1.] If $\gamma \subset \partial \alpha _{j}$ or $\gamma 
    =\alpha _{j}$ for some $j$, then $\ell _{j}$ is contained in a 
    connected union $\ell '=[x,y]$ of segments $\ell _{k}$ such that
    $\gamma \subset \partial \alpha _{k}$, and 
    $\gamma $ is in the convex hull of those 
    $\alpha _{m}$ for which $(\alpha _{m},\ell _{m})$ is ltd and 
$\ell _{m}\subset \ell '$.
    
\item[2.] If $\ell _{j}$ and $\ell _{k}$ intersect precisely in 
    an endpoint, and $\alpha _{j}$ and $\alpha _{k}$ have essential 
    intersections, then there is no gap or nontrivial nonperipheral 
loop 
    $\beta \subset \alpha 
_{j}\cap \alpha _{k}$.
    
    \item[3.]   For any $\gamma $ and $\ell '=[x,y]$ as in 1, either 
$x$ is an 
    endpoint of $\ell $, or $x\in \ell _{p}$ for some ltd $(\alpha 
    _{p},\ell _{p})$ such that $\gamma $ intersects $\alpha _{p}$ 
    essentially, and similarly for $y$.
    \end{description}

\end{ulemma}

\noindent {\em{Proof.}} Choose any disjoint set of $ \alpha 
_{j}\times \ell 
_{j}$ ($1\leq j\leq R_{0}$) such that the 
complement of the 
union contains no ltd, and such that for every $(\alpha _{j},\ell 
_{j})$ and $x\in \ell _{j}$, $\gamma \subset \partial \alpha _{j}$, 
there is 
an ltd $(\alpha _{k},\ell _{k})$ with $x\in \ell _{k}$ and $\gamma 
\subset \partial \alpha _{k}$. Then 
condition 1 is satisfied. By \ref{5.5}, for every $(\alpha _{j},\ell 
_{j})$ 
which is not ltd, $\alpha _{j}$ is a disjoint union of $\beta $ such 
that $(\beta ,\ell _{j})$ is bounded. If we can refine this partition 
to satisfy conditions 2 and 3, then every $(\beta ,\ell ')$ in the 
complement of the ltds will automatically be bounded by \ref{5.5} 
(especially (\ref{5.5.4})), 
because conditions 2 and 3 will ensure that there is no $\gamma $ in 
the 
interior of $\beta $ with $\vert \varphi (\gamma )\vert <\varepsilon 
_{2}$ for 
$[\varphi ]$ an endpoint of $\ell $, at least if we take $\varepsilon 
_{1}$ small enough given the ltd parameter functions. 

  Then we modify the partition in finitely many steps, 
always keeping condition 1, until condition 2 is satisfied. We do 
this as follows. Suppose we have a partition $\cal{P}$ into ltd and 
bounded sets, satisfying 1 of the vertically efficient conditions,
  and there are $\ell _{k}$ and $\ell _{m}$ intersecting in 
precisely one point and such that ($\alpha _{k},\ell _{k})$ and 
$(\alpha _{m},\ell _{m})$ are bounded and $\alpha _{k}\cap \alpha 
_{m}$ contains a $\beta $ as is disallowed in  2 of vertically 
efficient. 
Then we can 
take $\beta $ to be a maximal union of components of $C(\partial 
\alpha _{k},\partial \alpha _{m})$. 

Then rewrite
$$\alpha _{k}\times \ell _{k}\cup  \alpha 
_{m} \times\ell _{m}= ((\alpha _{k}\setminus \beta )\times \ell 
_{k})\cup ((\alpha _{m}\setminus \beta )\times  \ell 
_{m})\cup ( \beta \times (\ell _{k}\cup \ell _{m})).$$
By \ref{5.5}, the $(\alpha ,\ell ')$ arising in this rewriting are 
still bounded {\em{for the same $L$.}} 
Since rewriting reduces the topological type of the surfaces 
involved, and no new endpoints of segments $\ell _{j}$  are 
introduced, 
finitely many rewritings gives a partition satisfying 1 and 2 of 
vertically efficient. Finally, to get 3 of vertically efficient, if 
$\gamma \times \ell '$ as in 1 of vertically efficient does not have 
endpoints as required by 3, we extend $\ell '$ through adjacent $\ell 
_{p}$ 
with $(\alpha _{p},\ell _{p})$ bounded, possibly joining up such 
segments, until endpoints are in ltds intersecting $\gamma $ 
essentially, as required. 
\Box

The pairs $(\alpha _{j},\ell _{j})$ in the above are not unique. For 
example, as already noted, it is possible for $(\alpha _{j},\ell 
_{j})$ to be both ltd 
and bounded, because the constant $L$ of \ref{5.5} is typically much 
bigger than $\Delta (\nu _{0})$, for $\nu _{0}$ 
as 
in \ref{5.4}.

\section{Long, thick and dominant ideas.}\label{6}

This is a rather long section, which is pure theory of Teichm\"uller 
geodesics, with no input from three-dimensional hyperbolic geometry. 
It does seem necessary to go through some of these results in some 
detail, where they have not previously appeared in \cite{R1}, or not 
in the same forms as given here. For understanding the proof of the 
Ending Laminations Theorem in the case of combinatorially bounded 
geometry, the  parts most obviously needed are: the first theorem in 
\ref{7.4}, subsection \ref{7.8}, and the first lemma in \ref{7.12}.
We also make use of Lemma \ref{7.2} at one point, in the case of the 
long $\nu $-thick and dominant $\alpha $ being the whole surface $S$.
However, we also, at one point, make explicit use of the theorem in 
\ref{7.6} --- which is a deduction from the 
main theorem in \ref{7.4}. These two results are about general 
Teichm\"uller geodesics, not confined to the thick part of 
Teichm\"uller space. Thus, even 
the proof in the case of combinatorial bounded geometry explicitly 
relies, at one point, on the theory of general Teichm\"uller 
geodesics. This should not be a surpise, because, historically, the 
case of combinatorial bounded geometry  is highly nontrivial. 
I regard 
the most difficult result in the whole paper as \ref{7.13}, which 
like the rest of this section, is purely about Teichm\"uller 
geodesics. This may be 
in some contrast to the experience of others who have worked on 
the Ending Laminations Theorem. I shall comment on this later. 

\ssubsection{Fundamental dynamical lemma.}\label{7.1}

The whole of the theory of ltd gaps and loops is based on a simple 
dynamical lemma which quantifies density of leaves of 
the stable and unstable foliations of a quadratic differential. This 
is basically 15.11 of \cite{R1}, where three alternative conclusions 
are given. Here is a statement assuming  the gap $\alpha $ is ltd 
at $[\varphi ]$.

\begin{ulemma} Given $\delta >0$, the following holds for suitable 
ltd 
parameter 
functions $(\Delta ,r,s,K_{0})$ and for a suitable function $L(\delta 
,\nu 
)$. Let $\alpha $ be a gap which is long $\nu $-thick 
and dominant along a segment $\ell =[[\varphi _{1}],[\varphi _{2}]]$ 
and let $[\varphi ]\in \ell $  
with $d([\varphi ],[\varphi _{1}])\geq \Delta (\nu )$. Let 
$q(z)dz^{2}$ be the 
quadratic differential at $[\varphi ]$ for $d([\varphi ],[\varphi 
_{2}])$ with stable and unstable foliations $\cal{G}_{\pm }$. Let 
$a=a(\alpha,q)$. Then 
there is no segment of the unstable foliation $\cal {G}_{+}$ of 
qd-length 
$\leq 2L(\nu ,\delta )\sqrt{a}$ with both ends on $\varphi (\partial 
\alpha )$, and 
every segment of the unstable foliation $\cal {G}_{+}$ of qd-length 
$\geq L(\nu ,\delta )\sqrt{a}$ in $\varphi (\alpha )$ intersects 
every segment of $\cal {G}_{-}$ of length $\geq \delta \sqrt{a}$. 
Similar statements hold with the role of stable and unstable reversed.
\end{ulemma}

\ssubsection{Loops cut the surface into cells.}\label{7.2}
Now we give some of the key results about long thick and dominants 
which we shall need. We start with two fairly simple results, both of 
which follow directly from \ref{7.1}. These properties 
are used several times in 
\cite{R1}, 
but may never be explicitly 
stated. The first may be reminiscent of the concept of tight 
geodesics 
in the curve complex developed by Masur and Minsky \cite{M-M2}, and 
the point may be that these occur ``naturally'' in Teichm\" uller 
space

\begin{ulemma} Given $L>0$, there is a function $\Delta _{1}(\nu )$ 
depending only 
on the topological type of $S$, such that the following holds for 
suitable parameter functions $(\Delta ,r,s,K_{0})$ 
Let $\alpha $ be a gap which is long $\nu $-thick and dominant along 
$\ell$ for $(\Delta , r, s,K_{0})$, with $\Delta (\nu 
)\geq \Delta _{1}(\nu )$. Let $y_{1}=[\varphi _{1}]$, $y_{2}=[\varphi 
_{2}]\in 
\ell $ with $d(y_{1},y_{2})\geq \Delta _{1}(\nu )$. 
Let $\gamma _{i}\subset \alpha $ with $\vert \varphi _{i}(\gamma 
_{i})\vert \leq L$, $i=1$, $2$. Then $\alpha \setminus (\gamma 
_{1}\cup \gamma 
_{2})$ is a union of topological discs with at most one puncture and 
topological annuli parallel to the boundary. Furthermore, for a 
constant $C_{1}=C_{1}(L,\nu )$, 
$$\# (\gamma 
_{1}\cap \gamma _{2})\geq C_{1}\exp 
d(y_{1},y_{2}).$$\end{ulemma}

\noindent {\em{Proof.}} Let $[\varphi ]$ be the midpoint of 
$[[\varphi _{1}],[\varphi _{2}]$ and let 
$q(z)dz^{2}$ be the quadratic differential 
for $d([\varphi ],[\varphi _{2}])$ at $[\varphi ]$. 
Because $\vert \psi  (\gamma 
_{i})\vert \geq \nu $ for all $[\psi ]\in [[\varphi _{1}],[\varphi 
_{2}])$, by \ref{5.2}, the good position  of $\varphi _{1}(\gamma 
_{1})$ satisfies 
$$\vert \varphi _{1}(\gamma _{1})\vert _{q,+}\geq C(L,\nu 
)\sqrt{a(\alpha 
,q)},$$
and similarly for $\vert \varphi _{2}(\gamma _{2})\vert _{q,-}$. So 
$$\vert \varphi (\gamma _{1})\vert _{ q,+}\geq C(L,\nu )e^{\Delta 
_{1}(\nu )/2}\sqrt{a(\alpha 
,q)},$$
$$\vert \varphi (\gamma _{2})\vert _{ q,2}\geq C(L,\nu )e^{\Delta 
_{1}(\nu )/2}\sqrt{a(\alpha 
,q)},$$
Then \ref{7.1} implies that, given $\varepsilon $, if 
$\Delta _{1}(\nu )$ is large enough given $\varepsilon $, 
$\varphi (\gamma _{1})$ cuts every segment of stable 
foliation of $q(z)dz^{2}$ of qd-length $\geq \varepsilon 
\sqrt{a(\alpha )}$ and 
$\varphi (\gamma _{2})$ cuts every segment of unstable 
foliation of $q(z)dz^{2}$ of qd- length $\geq \varepsilon 
\sqrt{a(\alpha )}$. So 
components of $\varphi (\alpha )\setminus (\varphi (\gamma _{1})\cup 
\varphi (\gamma _{2}))$ have Poincar\'e diameter
$<\nu $ if $\Delta _{1}(\nu )$ is sufficiently large, and must be
topological discs with at most one puncture or boundary-parallel 
annuli.

The last statement also follows from \ref{7.1}. 
If $d(y_{1},y_{2})<\Delta 
_{1}(\nu )$, there is nothing to prove, so now assume that 
$d(y_{1},y_{2})\geq \Delta 
_{1}(\nu )$. It suffices to bound below the number of intersections 
of $\varphi (\gamma _{1})$  and $\varphi (\gamma _{1})$. Let $L(\nu 
,1)$ be as in \ref{7.1}, and assume without loss of generality that 
$L(\nu ,1)\geq 1$. Supppose that $\Delta _{1}(\nu )$ is large enough 
that  each of $\varphi (\gamma _{1})$ and $\varphi (\gamma _{2})$ 
contains 
at least one segment which is a qd distance $\leq \sqrt{a}/L(\nu ,1)$ 
from 
unstable and stable segments, respectively, of  $qd$-length 
$\geq \sqrt{a}L(\nu ,1)$. Note that the number of singularities of 
the 
quadratic differential is bounded in terms of the topological type 
of $S$. So 
apart from length which is a  bounded multiple of $L(\nu 
,1)\sqrt{a}$, 
each of $\varphi (\gamma _{1})$ and $\varphi (\gamma _{2})$ is a 
union of such segments. Then applying \ref{7.1}, each such segment 
of $\varphi (\gamma _{1})$ intersects each such segment on 
$\varphi (\gamma _{2})$. So we obtain the result for 
$C_{1}=c_{0}L(\nu 
,1)^{-2}$, for $c_{0}$ depending only on the topological type of $S$.
\Box

\ssubsection{A partial order on ltd $(\beta ,\ell )$.}\label{7.3}

\begin{ulemma} For $i=1$, $3$, let $y_{i}=[\psi _{i}]\in \ell _{i}$, 
and 
let $\beta _{i}$ be a subsurface of $S$ with $\vert \psi 
_{i}(\partial 
\beta _{i})\vert \leq L$. Let  ltd parameter 
functions be suitably chosen given $L$. 
Let  $\ell _{2}\subset [y_{1},y_{3}]$ and  let $\beta 
_{2}\cap \beta _{i}\not = \emptyset $ for both $i=1$, $3$, and let 
$\beta _{2}$ 
be ltd along $\ell _{2}$. Then $\beta _{1}\cap \beta _{3}\not = 
\emptyset $, and $\beta _{2}$ is in the convex hull of 
$\beta _{1}$ and $\beta _{3}$. \end{ulemma}

\noindent{\em{Proof.}} This is obvious unless both $\partial 
\beta _{1}$ and $\partial \beta _{3}$ intersect the interior of 
$\beta 
_{2}$. So now suppose that they both do this.
First suppose that $\beta _{2}$ is a gap and long, $\nu $-thick and 
dominant. Let $y_{2,1}=[\psi _{2,1}]$, $y_{2,3}=[\psi _{2,3}]\in \ell 
_{2}$ 
with $y_{2,i}$ separating $\ell 
_{i}$ from $y_{2}$, with $y_{2,i}$ distance $\geq {1\over 3}\Delta 
(\nu)$ from 
the ends of $\ell _{2}$ and from $y_{2}$. If $\beta _{2}$ is a loop, 
then we can take these distances to be $\geq {1\over 6}\log K_{0}$. 
For $[\psi ]\in[y_{-}.y_{+}]$, let $\psi (\beta )$ denote the 
region bounded by $\psi (\partial 
\beta )$ and homotopic to $\psi (\beta )$, assuming $\psi (\partial 
\beta )$ is in good position with respect to the quadratic 
diferential at $[\psi ]$ for $[y_{-},y_{+}]$
Then if $\beta _{2}$ is a gap, 
$\psi _{2,1}(\partial \beta _{1}\cap \beta _{2})$ 
includes a union of segments of in approximately unstable direction, 
of Poincar\' e length bounded from $0$, and similarly for 
$\psi _{2,3}(\partial \beta _{3}\cap \beta _{2})$, with 
unstable replaced by stable. Then as in \ref{7.2}, 
 $\psi _{2}(\partial \beta _{3}\cap \beta _{2})$ 
and $\psi _{2}(\partial \beta _{1}\cap \beta _{2})$
cut $\psi _{2} (\beta _{2})$ into topological discs with at 
most one puncture and annuli parallel to the boundary. It follows 
that $\beta _{2}$ is contained in the convex hull of $\beta _{1}$ and 
$\beta _{3}$. If $\beta _{2}$ is a loop it is simpler. We replace 
$\psi (\beta _{2})$ by the maximal flat annulus $S([\psi ])$ 
homotopic 
to $\psi (\beta _{2})$, for $[\psi ]\in \ell _{2}$. Then $\psi 
_{2}(\partial \beta _{1})\cap S([\psi _{2}])$ is in approximately 
the unstable direction and $\psi 
_{2}(\partial \beta _{3})\cap S([\psi _{2}])$ in approximately the 
stable direction. They both cross $S([\psi _{2}]$, so must intersect 
in a loop homotopic to $\psi _{2}(\beta _{2})$. \Box

We  define $(\beta _{1},\ell _{1})<(\beta _{2},\ell _{2})$ 
if 
$\ell _{1}$ is to the left of $\ell _{2}$ (in some common geodesic 
segment) and $\beta _{1}\cap \beta _{2}\not = \emptyset $. We can 
make 
this definition for any segments in a larger common geodesic segment, 
and even for single points in a common geodesic segment. So in the 
same way we can define $(\beta _{1},y_{1})< (\beta _{2},\ell _{2})$ 
if $y_{1}$ is to the left of $\ell _{2}$, still with  
$\beta _{1}\cap \beta _{2}\not = \emptyset $, and so on. This 
ordering is transitive restricted to ltd's $(\beta _{i},\ell _{i})$
by the lemma.

\ssubsection{Triangles of geodesics.}\label{7.4}

The concept of long thick and dominant was mainly developed in order 
to formulate results about triangles of geodesics in ${\cal T}(S)$. 
The following theorem was proved in 15.8 of \cite{R1} in the case of 
$S$ 
being a punctured (or marked) sphere. The proof is in fact 
completely general, once the approximate product structure of the 
thin part of ${\cal T}(S)$ has been formalised, as we did in Section 
\ref{2}. Before we state the general theorem, we state it in the 
special 
case of a geodesic segment $[y_{0},y_{1}]\subset \cal{T}_{\geq 
\nu }$.

\newenvironment{theorem1}{\par\noindent{\textbf{Triangle 
Theorem}\rm{\ (special case).}}\,\,\em}{\rm}
\newenvironment{theorem2}{\par\noindent{\textbf{Triangle 
Theorem}\rm{\ (general case).}}\,\,\em}{\rm}

\begin{theorem1} There exists a function $C:(0,\infty )\to 
(0,\infty)$ such that the 
following holds. Let $[y_{0},y_{1}]\subset \cal{T(S)}_{\geq \nu }$. 
Let 
$y_{2}\in \cal{T}(S)$. Then for all $y\in [y_{0},y_{1}]$, there 
exists 
$y'\in [y_{0},y_{2}]\cup [y_{1},y_{2}]$ such that $d(y,y')\leq 
C(\nu )$. If $y'\in [y_{0},y_{2}]$ and $w\in [y_{0},y]$ then the 
corresponding $w'$ is in $[y_{0},y']$, and similarly if $y'\in 
[y_{1},y_{2}]$. \end{theorem1}

\begin{theorem2}  
There are functions $C:(0,\infty )\to (0,\infty)$, 
$\Delta _{1}:(0,\infty )\to (0,\infty)$, $\nu _{1}:(0,\infty )\to 
(0,\infty)$ 
and constants $L_{0}$, $L_{1}$
such that the following  holds for suitable parameter functions 
$(\Delta ',r',s',K_{0}')$, and for ltd parameter 
functions $(\Delta ,r,s,K_{0})$ given $(\Delta 
',r',s',K_{0}')$.

 Let $y_{0}$, $y_{1}$, $y_{2}\in {\cal T}(S)$ 
with $y_{j}=[\varphi _{j}]$. Take any $y=[\varphi 
]\in [y_{0},y_{1}]$.
\begin{description}
\item[ 1.]  Let 
$\alpha $ be a loop which is $K_{0}$-flat on 
$\ell \subset [y_{0},y_{1}]\subset {\cal T}$. 
Then $\ell $ is a union of two disjoint segments $\ell 
_{0}$ and $\ell _{1}$, and there are segments $\ell _{j}'\subset 
[y_{j},y_{2}]$ such that for all $y\in 
\ell _{j}$ there is $y'=[\varphi ']$, $\ell _{j}'$  such that 
$$\vert \varphi '(\alpha )\vert \leq L_{0},$$
$$\vert {\rm{Re}}(\pi _{\alpha 
}(y))-{\rm{Re}}(\pi _{\alpha }(y'))\vert \leq 
L_{0},$$
and either 
$\alpha $ is $K_{0}'$-flat along $\ell 
_{j}'$, or $\ell _{j}$ and $\ell _{j}'$ have 
length $\leq L_{1}$.
\item[ 2.]  
 Let $\alpha $ be a long $\nu $- thick and 
dominant gap along $\ell \subset 
[y_{0},y_{1}]\subset {\cal T}$ for  $(\Delta ,r,s)$.  Then $\ell $ is 
a 
union of two disjoint segments $\ell _{0}$ and $\ell 
_{1}$, and there are  $\ell _{j}'\subset [y_{j},y_{2}]$ such that for 
all 
$y\in \ell _{j}$ there is $y'=[\varphi ']\in \ell _{j}$,  
with $y'\in \ell _{j}'\subset [y_{j},y_{2}]$ such 
that 
$$\vert \varphi '(\partial \alpha )\vert \leq  L_{0},$$
$$d_{\alpha }(y,y')\leq C(\nu ).$$ 
and either $\alpha $ is long, $\nu '
$-thick and dominant along $\ell _{j}'$ for $(\Delta ',r',s')$ and 
some $\nu '\geq \nu _{1}(\nu )$  
or $\ell _{j}$ and $\ell _{j}'$ have length $\leq 
\Delta _{1}(\nu )$. 

\item[ 3.] If $y\in \ell _{j}$ and $y'$ 
are as in either 1 or 2 above and $\lambda \subset 
[y,y_{j}]$, with $\beta $, $\lambda $,  satisfying 
the conditions of $\alpha $, $\ell $, in 1 or 
2 above, and $\alpha \cap \beta \neq \emptyset  $, 
then $\lambda =\lambda _{j}$,   
where $\lambda  _{j}$ are defined relative to $\lambda $ 
as the $\ell _{j}$ to $\ell $.
\end{description}
\end{theorem2}

This result extends to larger cycles of geodesic segments in 
$\cal{T}(S)$. If $y_{i}\in \cal{T}(S)$ for $0\leq i\leq n$, then we 
can consider the geodesic segments $[y_{i},y_{i+1}]$ for $0\leq i<n$ 
and $[y_{0},y_{n}]$. Then for any ltd $(\alpha ,\ell )$ along 
$[y_{0},y_{1}]$, we get a corresponding result to the above
relative to a decomposition of $\ell $ into sets $\ell 
_{j}$ for $j=0$ or $1\leq j\leq n$ with $\ell _{j}\subset 
[y_{j},y_{j+1}]$ for $j>0$ and $\ell _{0}\subset [y_{0},y_{n}]$. 
Typically, one expects all but one of the $\ell _{j}$ to be empty, 
but 
they could all be nonempty. The result is generalised by considering 
a 
decomposition into triangles, for example, triangles with vertices at 
$y_{0}$, $y_{j}$ and $y_{j+1}$ for each $1\leq j<n$.

\ssubsection{Only coordinates matter.}\label{7.6} 

The ltd's which occur along a geodesic segment, up to bounded 
distance, are often determined by only some coordinates of the ends 
of the geodesic segment. This is what the following theorem says. 

\begin{utheorem} Given $L_{1}$ and suitable ltd parameter functions 
$(\Delta ,r,s,K_{0})$, there is $L_{2}$ such that the following 
holds.
    Let $y_{j}=[\varphi _{j}]$, $y_{j}'=[\varphi 
_{j}']\in \cal{T}(S)$, $j=0$, $1$.
 Suppose that there are gaps or loops 
  $\alpha _{j}$ such that
$$\vert \varphi _{j}(\alpha _{j})\vert \leq L_{1},\ \ 
\vert \varphi _{j}'(\alpha _{j})\vert \leq L_{1},$$
and
$$d_{\alpha _{j}}(y_{j},y_{j}')\leq L_{1}{\rm{\ or\ }}
\vert {\rm{Re}}(\pi _{\alpha _{j}}(y_{j})-\pi _{\alpha 
_{j}}(y_{j}'))\vert 
\leq L_{1}.$$
Then we have the following.

Let $\alpha  $ be a loop or gap which is $K_{0}$-flat or long $\nu 
$-thick and dominant along $\ell =[z_{0},z_{1}]\subset 
[y_{0},y_{1}]$, 
with essential intersections with both $\alpha _{0}$ and $\alpha 
_{1}$ and such 
that 
$$d_{\alpha _{j},\alpha }'(y_{j},z_{j})\geq L_{2}.$$
Then in the conclusion of \ref{7.4}, but considering 
$[y_{0},y_{0}']\cup 
[y_{0}',y_{1}']\cup [y_{1}',y_{1}]$ instead of $[y_{0},y_{2}]\cup 
[y_{2},y_{1}]$,
and with decompositions $\ell =\ell _{0,0}\cup \ell _{0,1}\cup \ell 
_{1,1}$ instead of $\ell _{0}\cup \ell _{1}$, we 
can take $\ell _{0,0}=\ell _{1,1}=\emptyset $, so that there is a 
corresponding segment $\ell '$ on $[y_{0}',y_{1}']$ to all of $\ell $.
    \end{utheorem}
    
\noindent{\em{Proof.}} Suppose that we cannot take $\ell 
_{0,0}=\emptyset $. Then let $\ell _{0,0}'$ be the corresponding 
segment on $[y_{0},y_{0}']$. Then by the following lemma,
for a constant $L_{3}$ depending 
only on the ltd parameter functions,
$$d_{\alpha _{0}}(y_{0},y_{0}')\geq d_{\alpha _{0},\alpha }'(y_{0},y)
+ d_{\alpha _{0},\alpha }'(y_{0}',y)-L_{3},$$
giving a contradiction. Similarly $\ell _{1,1}=\emptyset $.
\Box

\begin{lemma}\label{7.5} The following holds for suitable ltd 
parameter 
functions $(\Delta ,r,s,K_{0})$, a function   $C:(0,\infty )\to 
(0,\infty )$ 
and constants $C_{0}>0$, $C_{0}'>0$. Let $y_{i}=[\varphi 
_{i}]$, $0\leq i\leq n$ and $z_{i}$, $1\leq i<n$ be 
points on a geodesic segment in $\cal{T}(S)$, in the order $y_{0}$, 
$y_{1}$, $z_{1}$,\ldots $y_{n-1}$, $z_{n-1}$, $y_{n}$. Write $\ell 
_{i}=[y_{i},z_{i}]$ for $1\leq i<n$. Let $\alpha _{i}$ 
be gaps or loops such that $\vert \varphi _{i}(\partial \alpha 
_{i})\vert 
\leq L_{0}$ for all $i$. For $i=0$ or $n$, if $\alpha _{i}$ is a 
gap, let $\vert \varphi _{i}(\gamma )\vert \geq \varepsilon _{0}$ for 
all 
nontrivial nonperipheral non-boundary-parallel $\gamma \subset \alpha 
_{i}$. For $1\leq i<n$, let $\alpha _{i}$ be either a gap which is 
long $\nu _{i}$-thick and dominant along $\ell _{i}$, or let 
$\alpha _{i}$ be  a loop which is $K_{0}$-flat along $\ell _{i}$. 
Let $\alpha _{i}\cap \alpha _{i+1}\neq \emptyset $ 
for $0\leq i<n$.  Write 
$C_{i}=C(\nu _{i})$ for $0<i<n$ if $\alpha _{i}$ is a gap and 
$C_{i}=2C_{0}$ if $\alpha _{i}$ is a loop. Then
$$d_{\alpha _{0},\alpha _{n}}'(z_{0},z_{n})\geq 
\sum _{i=0}^{n-1}(d_{\alpha _{i}}(y_{i},z_{i})+d_{\alpha _{i},\alpha 
_{i+1}}'(z_{i},y_{i+1})-C_{i}).$$
\end{lemma}

\noindent{\em{Proof.}}
This is done by a {\em{locking technique}} which is 
used frequently in \cite{R1}. By locking, we mean the following. 
Suppose that $\beta $ and $\beta '$ are two paths in $\psi (S)$ and 
some quadratic differential $q(z)dz^{2}$ is fixed. Then we say that 
$x' \in \beta '$ is {\em{locked}} to $x\in \beta $ (for some fixed 
$\delta >0$) 
if there is  an arc of stable foliation between $x$ and $x'$ of 
Poincar\'e length
$\leq \delta $ times the injectivity radius at $\psi (S)$, also 
measured in Poincar\'e length. If $\delta $ is sufficiently small, 
depending only on the topological type of $S$, then we can equally 
well use the qd-length to measure this ratio. This means there is a 
constant $C(S)>0$ such that if an arc has length $\leq \delta \leq 
C(S)^{-1}$ times the 
injectivity radius, with both measured in the Poincar\'e metric then 
the qd-length is $\leq C(S)\delta $ times the injectivity radius, 
also 
measured in the qd-metric, and similarly with Poincar\'e metric and 
qd-metric interchanged.  Let $\chi _{t} $ denote the family of 
homeomorphisms 
obtained by scaling unstable and stable length for $q(z)dz^{2}$ by 
$e^{\pm t}$. If $\beta '$ is locked to $\beta $, then $\chi 
_{t}(\beta ')$ is locked to $\chi _{t}(\beta )$ for all $t>0$, 
because qd-length of locking segments get multiplied by $e^{-t}$, 
while the injectivity radius, measured in qd-length, cannot decrease 
by more than $e^{-t}$. A technique developed in \cite{R1} for showing 
one loop 
was much longer than another was to show that one loop had many 
points locked to each point on the other.
 
If $\alpha _{i}$ is a loop, put $\gamma _{i}=\gamma _{i}'=\alpha 
_{i}$. Now suppose that $\alpha _{i}$ is a gap. 
For some $L_{1}(\nu )$ depending only on $\nu $ (and the topological 
type of $S$),  for each 
$[\varphi ]\in \ell _{i}$ 
there is a loop $\gamma \subset \alpha _{i}$ such that
$$\nu _{i}\leq \vert \varphi (\gamma )\vert \leq L_{1}(\nu _{i}).$$
By \ref{7.1}, given $\delta >0$ there are then $\Delta 
_{1}(\nu )$ and 
$L_{2}(\nu _{i})$ depending only on $\nu $ and $L_{1}(\nu )$ such 
that 
if $d([\varphi ],z_{i})\geq \Delta _{1}(\nu _{i})$, and $[\varphi 
]\in \ell _{i}$, then for 
each $x\in \varphi (\gamma )$ and each point $x'\in \varphi (\alpha 
_{i})$ where the injectivity radius is $\geq \nu _{i}$, there 
is a stable segment in $\varphi (\alpha _{i})$ starting from $x$
and coming within $\delta \nu _{i}$ of $x'$. From now on we assume 
that $\Delta (\nu )>\Delta _{1}(\nu )$, and sufficiently  large in a 
sense to be determined.

If $\alpha _{i}$ is a loop, then define $\gamma _{i}=\gamma 
_{i}'=\alpha _{i}$. In this case, $\alpha _{i-1}\neq \alpha _{i}$, 
except possibly if $i=1$, and $\alpha _{i+1}\neq \alpha _{i}$, except 
possibly if $i=n-1$. If $\alpha _{i}$ is a gap, choose $\gamma _{i}$, 
$\gamma _{i}'\subset 
\alpha _{i}$ such that, 
$$\nu _{i}\leq \vert \varphi _{i}(\gamma _{i})\vert \leq L_{1}(\nu 
_{i}),$$
$$\nu _{i}\leq \vert \psi _{i} (\gamma _{i}')\vert \leq L_{1}(\nu 
_{i}),$$
$$C'\vert \varphi _{i+1}(\gamma _{i}'\cap \alpha 
_{i+1})\vert \geq \exp d_{\alpha _{i},\alpha 
_{i+1}}'(z_{i},y_{i+1}).$$
This last is possible by \ref{2.8}, and should be done for $0\leq 
i<n$, for a suitable constant $C'$.
In addition, choose $\gamma _{0}\subset \alpha _{0}$ so that
$$C'\vert \varphi _{1}(\gamma _{0}\cap \alpha 
_{1})\vert \geq \exp d_{\alpha _{0},\alpha _{1}}'(y_{0},y_{1}).$$

If $\alpha _{i}$ is a gap, let $z_{i}'=[\psi 
_{i}']\in \ell _{i}$ with $d(z_{i}',z_{i})=
\Delta _{1}(\nu _{i})$, for $\Delta {1}(\nu )$ to be chosen as 
follows, and if $\alpha _{i}$ is a loop, let $z_{i}'=[\psi 
_{i}']\in \ell _{i}$ with $d(z_{i}',z_{i})=
\Delta _{1}$, for $\Delta _{1}$ to be chosen as 
follows.
By \ref{6.2}, if $\alpha _{i}$ is a gap every point on $\psi 
_{i}'(\gamma _{i}')$ is locked to $\geq 
C_{2}(\nu _{i})^{-1}\exp d(y_{i},z_{i}')$ points on $\psi 
_{i}(\gamma _{i})$ along stable segments of Poincar\'e length $\leq 
L_{2}(\nu _{i})$, and every point on $\varphi _{i}(\gamma _{i})$ 
is locked to $\geq 
C'^{-1}\exp d_{\alpha _{i-1},\alpha _{i}}'
(z_{i-1},y_{i})$ points on $\varphi 
_{i}(\gamma _{i-1}')$ along stable segments of Poincar\'e length 
$\leq 
L_{2}(\nu _{i})$. Then assuming $\Delta _{1}(\nu )$ is sufficiently 
large given $L_{2}(\nu )$, and $\Delta (\nu )$ sufficiently large 
given $\Delta _{1}(\nu )$, every point on $\psi _{i}(\gamma _{i}')$ 
is 
locked to $\geq 
C_{2}(\nu _{i})^{-1}\exp d(y_{i},z_{i}')$ points on $\psi 
_{i}(\gamma _{i})$ along segments of Poincar\'e length $\leq 
C'.L_{2}(\nu _{i}).\exp -d_{\alpha _{i}}(y_{i},z_{i})<\delta 
\nu _{i}$  and 
$\geq (C_{2}(\nu _{i})C')^{-1}\exp 
d_{\alpha _{i-1},\alpha _{i}}'(z_{i-1},y_{i}).\exp d_{\alpha 
_{i}}(y_{i},z_{i}')$ points on 
$\psi _{i}(\gamma _{i-1}')$ along stable segments of Poincar\'e 
length $\leq C'.L_{2}(\nu _{i}).\exp -d_{\alpha 
_{i}}(y_{i},z_{i})\leq 
\delta \nu _{i}$. 

If $\alpha _{i}$ is a loop, $\varphi _{i}(\gamma _{i-1}')$ intersects 
the flat 
annulus homotopic to $\varphi _{i}(\alpha _{i})$ in 
$\geq C'^{-1}.C\exp 
d_{\alpha _{i-1},\alpha _{i}}'(z_{i-1},y_{i})$ segments in an 
approximately stable direction. The corresponding segments of $\psi 
_{i}'(\gamma _{i-1}')$ each contain $\geq C'^{-1}\exp d_{\alpha 
_{i}}(y_{i},z_{i}')$ disjoint segments which can be locked to $\psi 
_{i}'(\alpha _{i})$ along stable segments of  Poincar\'e length $\leq 
C'.\vert \psi _{i}'(\alpha _{i})\vert $, for suitable $C'$. Then 
$\psi 
_{i}(\gamma _{i-1}')$ contains \\ $\geq C'^{-2}\exp d_{\alpha 
_{i-1},\alpha _{i}}'(z_{i-1},y_{i}).\exp d_{\alpha 
_{i}}(y_{i},z_{i})$ segments which can be locked to $\psi (\alpha 
_{i})$ by stable segments of length $\leq \delta .\vert \psi (\alpha 
_{i})\vert $, assuming that $\Delta _{1}$ is large enough given $C'$ 
and $\delta $.

Define $C_{i}=\log C'$ if $i=0$ or $n$. If $0<i<n$, define 
$C_{i}=\Delta _{1}(\nu _{i})+\log C_{2}(\nu _{i})+3\log C'$ if 
$\alpha _{i}$ is a gap,  and $C_{i}=C_{0}'=4\log C'+\log \Delta _{1}$ 
if 
$\alpha _{i}$ is a 
loop. Define 
$$\Delta _{j}=\sum _{i=1}^{j}(d_{\alpha 
_{i-1}\alpha _{i}}'(z_{i-1},z_{i})-C_{i}),$$
$$\Delta _{j}'=\sum _{i=1}^{j}(d_{\alpha 
_{i-1}\alpha _{i}}'(z_{i-1},z_{i})-C_{i})-d_{\alpha 
_{i}}(y_{j},z_{j}).$$
We claim inductively that, assuming the ltd parameter functions are 
strong enough, for $j\geq 1$, each point  of $\psi 
_{j}(\gamma _{j}')$ is locked to $\geq \exp \Delta _{j}$ points on 
$\psi _{j}(\gamma _{0}')$ along  stable segments of Poincar\'e  $\leq 
2\delta $ times the injectivity radius in the Poincar\'e metric, and 
similarly $\varphi _{j}(\gamma _{j})$  is locked to 
$\geq \exp \Delta _{j}'$ points on 
$\varphi _{j}(\gamma _{0}')$ along  stable segments of Poincar\'e  
$\leq 
2\delta $ times the injectivity radius in the Poincar\'e metric.  As 
before, for $\delta $ sufficiently small, it suffices to prove this 
using the qd-metric. But in the qd-metric,  $\chi _{t}$ 
multiplies $qd$-length along the stable direction by $e^{-t}$, while 
the injectivity radius in the qd-metric decreases by at most a 
factor $e^{-t}$. Relative ratios of lengths of locking segments are 
preserved up to a bounded proportion. 
So locking segments between $\psi _{k}(\gamma _{i}')$ and $\psi 
_{k}(\gamma _{j}')$ have Poincar\'e length $\leq C''exp (\Delta 
_{j}-\Delta _{k})$, for $i<j<k$, for a suitable constant $C''$. 

This gives the required estimate, apart from minor adjustments if 
$\alpha _{0}$ or $\alpha _{n}$ is a loop. In those cases, given the 
definition of $d_{\alpha _{0},\alpha _{n}}'$, there is nothing to 
prove if $\alpha _{0}=\alpha _{n}$, or if $\alpha _{0}$ is not 
flat
along any segment of  $[y_{0},y_{n}]$. If $\alpha _{0}$ is  along a 
segment of $[y_{0},y_{n}]$, we can introduce another segment $\ell 
_{1}$ if necessary, renumbering, so that $\alpha _{0}$ is flat 
precisely along $\ell _{1}$.  We make similar adjustments near 
$y_{n}$ 
but otherwise the proof is exactly as above.  

\Box

\ssubsection{The graph of the qd-length function.}\label{7.20}

One of the basic technical considerations in the study of 
Teichm\"uller geodesics, as is probably already apparent, is the 
difference between the qd- and Poincar\'e metrics. The two metrics 
are 
not globally Lipschitz equivalent. But they are Lipschitz-equivalent, 
up  to scalar, on any thick part of a surface. The Lipschitz constant 
is bounded in terms of the topological type of the surface, but the 
scalar 
is completely uncontrollable. This should not be regarded as a 
problem. One simply has to look at ratios of lengths rather than at 
absolute lengths. Also, the qd-length function has a rather 
remarkable 
property. Fix a Teichmuller geodesic $\{ [\chi _{t}\circ \varphi 
_{0}]:t\in \mathbb R\} $, where $\chi _{t}$ minimises distortion 
and $d([\chi _{t}\circ \varphi _{0}],[\chi _{s}\circ \varphi 
_{0}])=\vert t-s\vert $. Let $q_{0}(z)dz^{2}$ be the quadratic 
differential at $[\varphi _{0}]$ for $d([\varphi _{0}],[\chi 
_{t}\circ \varphi _{0}]$, for $t>0$, and $q_{t}(z)dz^{2}$ the 
stretch at $[\chi _{t}\circ \varphi _{0}]$ (\ref{2.2}). Let $\vert 
.\vert _{t}=\vert .\vert _{q_{t}}$, the qd-length (\ref{5.1}). For 
any 
finite loop set $\gamma $, define
$$F(t,\gamma )=\log \vert \chi _{t}\circ \varphi _{0}(\gamma )\vert 
_{t}.$$
 By 14.7 of 
\cite{R1} (and I am sure this is well-known), 
there is a constant $C_{0}$ depending only on  the 
topological type of $S$, and a bound on the number of loops in 
$\gamma $,
and there are $c(\gamma )$, $t(\gamma )\in 
\mathbb R$ such that
\begin{equation}\label{7.20.1}
    \vert F(t,\gamma )-\vert t-t(\gamma )\vert -c(\gamma )\vert \leq 
    C_{0}.\end{equation}
The graph of the function $t\mapsto F(t,\gamma )$, for any $\gamma $, 
therefore lies within $C_{0}$ of a $V$, will the slopes of the arms 
of  the $V$ 
being $-1$ on the left and $1$ on the right, and minimum at 
$t(\gamma )$.. 

Comparision between Poincar\'e and qd-length can then be made as 
follows. Given $L_{1}>0$ there is $L_{2}\in \mathbb 
R$ such that
\begin{equation}\label{7.20.2}
    \vert \chi _{t}\circ \varphi _{0}(\gamma )\vert \leq 
L_{1}\end{equation}
    whenever 
    \begin{equation}\label{7.20.3} F(t,\gamma )-F(t,\gamma ')\leq 
L_{2}\end{equation}
for all nontrivial 
nonperipheral $\gamma '$ intersecting $\gamma $ transversely. 
Conversely, given $L_{2}\in \mathbb R$, there is $L_{1}$ such that 
(\ref{7.20.2}) holds whenever (\ref{7.20.2}) holds for all $\gamma '$ 
intersecting $\gamma $ transversely. There is a similar 
characterisation of short loops. Given $L_{2}<0$, there is $L_{1}>0$ 
(small if $L_{2}$ is negatively large)
such that, whenever (\ref{7.20.2}) holds, then (\ref{7.20.3}) holds  
for all $\gamma '$ intersecting 
$\gamma $ tranversely. Conversely, given $L_{1}>0$, there is $L_{2}$ 
(negative if $L_{1}$ is small) such that (\ref{7.20.2}) holds for 
$\gamma $, whenever (\ref{7.20.3}) holds for $\gamma $ and all 
$\gamma '$ transverse to $\gamma $.

If $\gamma $ satisfies 
\ref{7.20.3} for all transverse $\gamma '$ and $\gamma ''$ is another 
loop, disjoint from $\gamma $,  with $\vert \chi _{t}\circ 
\varphi _{0}(\gamma '')\vert $ bounded, and $F(t,\gamma )-F(t,\gamma 
'')\leq 
L_{2}$, then it is possible that $\chi _{t}\circ \varphi _{0}(\gamma 
'')$ is 
short, while $\chi _{t}\circ \varphi _{0}(\gamma )$ is not. However, 
if 
$\Gamma $ is a set of loops $\gamma $, satisfying \ref{7.20.3}, then 
any 
component $\alpha $ of the convex hull of $\Gamma $ is such that 
$\chi 
_{t}\circ \varphi _{0}(\alpha )$ is contained in a single component 
of $(\chi 
_{t}\circ \varphi _{0}(S))_{\geq \varepsilon (L_{2})}$, for a 
suitable 
$\varepsilon >0$ depending only on $L_{2}$. Conversely, if $\chi 
_{t}\circ \varphi _{0}(\alpha )$ 
is a component of $(\chi 
_{t}\circ \varphi _{0}(S))_{\geq \varepsilon }$ then we can find a 
set of 
loops $\Gamma $ with convex hull $\alpha $ such that $\gamma $ 
satisfies (\ref{7.20.3}) for all $\gamma \in \Gamma $, for a suitable 
$L_{2}=L_{2}(\varepsilon )$. 

\ssubsection{Ltd's in the projection are the same.}\label{7.7} 
We shall need the following.

\begin{utheorem}Given $L_{1}$, there is $L_{2}$, and given ltd 
parameter functions $(\Delta ',r',s',K_{0}')$, 
there are $(\Delta ,r,s,K_{0})$, $C_{0}$, $\nu _{0}$, $\nu 
_{0}'$ 
and $C:(0,\infty )\to (0,\infty )$  such that the following 
holds. Let $y_{j}=[\varphi _{j}]$, $j=0$, $1$, and let $\alpha $ be a 
loop or  gap 
with $\vert \varphi _{j}( \alpha )\vert \leq L_{1}$ or $\vert \varphi 
_{j}(\partial \alpha )\vert \leq L_{1}$ for $j=0$, $1$. Let $\ell 
=[z_{0},z_{1}]\subset [\pi _{\alpha }(y_{0}),\pi _{\alpha 
}(y_{1}]\subset \cal{T}(S(\alpha ))$, let $\beta $ be long $\nu 
$-thick and dominant along $\ell $ with respect to $(\Delta ,r,s)$ 
for some $\nu \geq \nu _{0}$, or $K_{0}$-flat along $\ell $ and let 
$\beta $, $z_{j}$, $\pi _{\alpha }(y_{j})$, $L_{1}'$ satisfy the 
conditions of \ref{7.6} with $\pi _{\alpha }(y_{j})$ replacing 
$y_{j}$. Then there is $\ell '\subset [y_{0},y_{1}]$ such that $\beta 
$ is long, $\nu '$-thick and dominant along $\ell '$ with respect to 
$(\Delta ',r',s')$ and some $\nu '\geq \nu _{0}'$ or $K_{0}'$ flat 
along $\ell '$, and for each $y\in \ell $ there is $y'\in \ell '$ 
such 
that
$$d_{\beta }(y,y')\leq C(\nu ){\rm{\ \ or\ \ }}\vert {\rm{Re}}(\pi 
_{\beta }(y)-\pi _{\beta }(y'))\vert \leq C_{0}.$$
\end{utheorem}

\noindent {\em{Proof.}} This is proved by similar techniques to 
\ref{7.4}, but since there is no precise statement like this in 
\cite{R1}, we had better give some details. We can assume $L_{1}$  is 
bounded and choose the ltd parameter functions relative to it, 
because we can then get the result for a general $L_{1}$ using 
\ref{7.6} 

 Given $\gamma \subset S$ and $[\varphi ]\in 
{\cal T}(S)$, we say that $\varphi (\gamma )$ is {\em{almost bounded 
(by $L$)}} if $\vert \varphi (\gamma 
)\vert ''\leq L$, where $\vert \varphi (\gamma 
)\vert ''$  is as in \ref{2.3}. 
A sufficient criterion 
for a loop to be almost bounded at some point on a geodesic segment
$[[\varphi _{0}],[\varphi _{1}]]$ is given by the negation of a 
necessary condition for a loop to be not bounded at any point 
of the geodesic segment, as follows. Given $L$, there is $L'$ such 
that 
if, for 
all choices of disjoint simple loops $\gamma _{0}'$, $\gamma _{1}'$ 
which both intersect $\gamma $, either 
\begin{equation}\label{7.7.1}\vert \varphi _{0}(\gamma 
)\vert \leq L\vert \varphi _{0}(\gamma _{j}')\vert {\rm{\ for\ 
}}j=0,1,\end{equation} or 
\begin{equation}\label{7.7.2}\vert \varphi _{1}(\gamma 
)\vert \leq L\vert \varphi _{1}(\gamma _{j}')\vert {\rm{\ for\ 
}}j=0,1,\end{equation}
then there is a 
point $[\varphi ]\in [[\varphi _{0}],[\varphi _{1}]]$ such that 
$\vert \varphi (\gamma )\vert ''\leq L'$. It is not clear if 
(\ref{7.7.1}) or (\ref{7.7.2}) is a 
necessary 
condition for $[\varphi ]$ to be bounded for some $[\varphi ]\in 
[[\varphi _{1}],[\varphi _{2}]]$, in general. But if 
$[\varphi ]\in [[\varphi _{0}],[\varphi _{1}]]$ 
and $\gamma \subset \alpha $, where $\beta  $ is ltd or flat  along a 
segment 
of 
$[[\varphi _{0}],[\varphi _{1}]]$ containing $[\varphi ]$, and $\vert 
\varphi (\gamma )\vert \leq L''$, then for any choice of $(\gamma 
_{0}',\gamma _{1}')$ as above for a suitable $L$ given 
$L''$, for any choice of $(\gamma _{0}',\gamma _{1}')$, one of 
(\ref{7.7.1}) or (\ref{7.7.2}) holds. This is essentially the 
content of 15.8 of \cite{R1}. 

Now let $\beta $ be a gap which is long, $\nu $-thick and dominant 
along \\ $\ell \subset [\pi _{\alpha }(y_{0}),\pi _{\alpha 
}(y_{1})]$. So 
for  
any loop $\gamma \subset \beta $ and $[\varphi ]\in \ell $ such that 
$\vert \varphi (\gamma )\vert \leq L$, we can find $[\varphi ']\in 
[y_{0},y_{1}]$ such that $\vert \varphi '(\gamma )\vert ''\leq L'$. 
But we 
actually want a bound on $\vert \varphi '(\gamma )\vert $, and that 
$\beta $ should be ltd along a segment containing $[\varphi ']$ for 
suitable parameter functions. This is done as follows. Take 
any  $w_{i}=[\psi _{i}]\in \ell $ and loop sets  $\Gamma _{i}\subset 
\beta $ 
which are cell-cutting in $\beta $ 
$i=0$, $1$, $2$ with $\vert \psi _{i}(\Gamma _{i})\vert \leq L$ such 
that $\vert \psi _{i}(\Gamma _{i})\vert \leq L$, with $w_{1}\in 
[w_{0},w_{2}]$ and 
$$L_{3}\leq d(w_{i},w_{i+1})\leq L_{4},$$
where $L_{3}$ is large enough for $\zeta _{i}\cup \zeta _{i+1}$ to be 
cell-cutting in $\beta $ for any $\zeta _{j}\in \Gamma _{j}$ for 
$i=0$, $1$, $2$, using \ref{7.2}. We shall also need $L_{3}$ large 
enough 
 for there to be no loop $\zeta \subset \beta $ with $\vert \psi 
_{i}(\zeta 
 )\vert \leq L''$ for  $i=0$, $1$, or both $i=1$, $i=2$, where $L''$ 
depends only on 
 $L$ and $L'$. This is again possible using \ref{7.2}, for $L_{3}$ 
 depending only on $L$ and $L'$. So now we fix this choice of 
$L_{3}$, 
 and $L_{4}$. Using (\ref{7.7.1}) and (\ref{7.7.2}) as above, we have 
 $w_{i}=[\psi _{i}']\in [y_{0},y_{1}]$ with $\vert \psi _{i}'(\Gamma 
 _{i})\vert ''\leq L'$, $i=0$, $1$, $2$. Now again using \ref{7.2}, 
for suitable parameter 
 functions $(\Delta '',r'',s'')$ and 
flat constant $K_{0}''$ strong enough given $L_{4}$, which bounds 
$\# (\Gamma _{0}\cap \Gamma _{1})$ and $\# (\Gamma _{1}\cap \Gamma 
_{2})$, 
 there cannot be any segment $\lambda 
\subset [w_{i}',w_{i+1}']$ and $\omega $ intersecting both $\zeta 
_{i}$ 
and $\zeta _{i+1}$ (any $\zeta _{j}\in \Gamma _{j}$) along which 
$\omega  $ is long thick and 
dominant for $(\Delta '',r'',s'')$ or $K_{0}''$-flat. Any $\omega $ 
which intersects one of $\zeta _{i}$, $\zeta _{i+1}$ intersects the 
other, since $\zeta _{i}\cup \zeta _{i+1}$ is cell-cutting in $\beta 
$. So then let $L_{5}$ be the constant $L$ given by \ref{5.5} 
relative 
to $(\Delta 
'',r'',s'')$ and $K_{0}''$. Then there is $\varepsilon >0$ bounded 
below in terms of $L_{5}$ so that if $\zeta \subset {\rm{Int}}(\beta 
)$ and 
$\vert \psi _{1}'(\zeta )\vert <\varepsilon $ then $\vert \psi 
_{i}'(\zeta )\vert <\varepsilon _{0}$ for either $i=0$ or $i=2$, by 
(\ref{5.5.4}). 
Now suppose there is such a loop $\zeta$ in the interior of $\beta $. 
Suppose without loss of generality that $\vert \psi _{0}'(\zeta 
)\vert \leq \varepsilon _{0}$.
 Then $\zeta $ has $\leq L'$ intersections with each of $\Gamma 
_{0}$, $\Gamma 
_{1}$, by the definition of $\vert .\vert ''$, since the loops of 
$\psi _{i}'(\Gamma _{i})$ have $\vert.\vert ''$-length $\leq L'$ for 
$i=0$, $1$. Since these loop sets are both cell-cutting in $\beta $, 
we 
deduce that $\vert \psi _{i}(\zeta )\vert \leq L''$ for $L''$ 
depending only on $L$ and $L'$, and $i=0$, $1$. By the choice of 
$L_{3}$ this is 
impossible. So this means we have a bound on $\vert \psi _{1}'(\Gamma 
_{1})\vert $ and we can take $[\varphi ']=w_{1}'$, if we take 
$w_{1}=y$. 

So we have $\ell '\subset [y_{0},y_{1}]$, and, for each $y\in \ell $, 
we 
have $y'\in \ell $ with $d(y,y')\leq C(\nu )$. Now we need to 
show that given $(\Delta ',r',s')$, $\ell '$ is long $\nu' $-thick 
and 
dominant along $\ell '$ for some $\nu '\geq \nu /C(\nu)$, if $(\Delta 
,r,s)$ are suitably chosen. First, we note that because $\partial 
\beta $ can be homotoped into $\Gamma _{0}$ as above, 
$\vert \psi (\partial \beta )\vert \leq C_{1}(\nu ')$ for all $[\psi 
]\in 
\ell '$, and indeed of an extension $\ell _{1}'$ of $\ell '$ at both 
endpoints, if the ltd parameter functions $(\Delta ,r,s)$ are 
sufficiently strong. The $d_{\beta }$-lengths of $\ell $ and $\ell '$ 
differ by at 
most $2C(\nu )$, and there are similar properties for $\ell _{1}'$.
The Poincar\'e length  of $\psi (\partial \beta )$ is  bounded  along 
$\ell _{1}'$. Now we need to show that the ratio $a(\partial \beta 
,[\psi ])/a(\beta )$ decreases exponentially in the middle of $\ell 
_{1}'$.  
The easiest way to see this is to make use of the 
functions $F(t,\gamma )$ of \ref{7.20}, for the geodesic 
$[y_{0},y_{1}]$. Write 
$\ell _{1}'=[[\chi _{a}\circ \varphi _{0}],[\chi _{b}\circ \varphi 
_{0}]]$.
Because the Poincar\'e length of $\psi (\partial \beta )$ is 
bounded all along $\ell _{1}'$, and we have a lower bound of $\nu '$ 
(with $\nu '=\nu '(\nu )$) on the length of 
loops $\psi (\zeta )$ for $\zeta $ in the interior of $\beta $, there 
is a function $C_{2}(\nu )$ such that
\begin{equation}\label{7.7.3}
    F(t,\partial \beta )\leq C_{2}(\nu ')+F(t,\gamma )\end{equation}
for all $\gamma $ in the interior of $\beta $ and $t\in [a,b]$.
The function $F(.,\gamma )$, for any $\gamma $ in the interior of 
$\beta $
has minimum at most $C_{2}(\nu ')$ below the minimum of 
$F(t,\partial \beta )$, if the minimum is in $[a,b]$. 
But, by comparing with $\ell _{1}$, for a 
$\gamma $ in the interior of $\beta $ for which $\chi _{t}\circ 
\varphi _{0}(\gamma )$ is bounded, and for $T={\rm{Min}}(\vert 
t-a\vert 
,\vert t-b\vert ,$, 
$$\vert \chi _{a}\circ \varphi _{0}(\gamma )\vert 
\geq C_{3}(\nu ').e^{T},$$
and similarly for $a$ replaced by $b$. It follows that, for such 
$\gamma $,
$$F(a,\partial \beta )\leq F(a,\gamma )+C_{4}(\nu ')-T,$$
and similarly with $a$ replaced by $b$. It follows that
$$F(t,\partial \beta )\leq F(t,\gamma )-T+C_{4}(\nu ')+C_{0},$$
where $C_{0}$ is the constant of (\ref{7.20.1}). Then the good 
position of $\chi _{t}\circ \varphi _{0}(\gamma )$ is bounded from 
the stable and unstable positions. 
Suppose that this $t$ is such that the good position of $\chi 
_{t}\circ \varphi _{0}(\partial \beta )$ is close to the stable 
foliation, or to the unstable foliation. These happen except on a 
bounded interval of $t$. Then the good position of
$\chi _{t}\circ \varphi _{0}(\beta )$ contains a ball of definite 
Poincar\'e radius  centred on a 
point of $\chi _{t}\circ \varphi _{0}(\gamma )$. So $\vert \chi 
_{t}\circ 
\varphi _{0}(\gamma )\vert _{t}$ is boundedly proportional to 
$\sqrt{a(\beta )}$ and, for this $t$,
$$\vert \chi _{t}\circ \varphi 
_{0}(\partial \beta )\vert _{t}\leq C_{5}(\nu ')e^{-T} \sqrt 
{a(\beta )}.$$
Then, since this is true except on a bounded interval,
it must be true for all $t\in [a,b]$, if we adjust the constant. If 
$\ell _{1}$ and 
$\ell _{1}'$  are sufficiently long, that is, if the parameter 
functions $(\Delta ,r,s)$ are 
sufficiently strong, we can ensure that the qd-length of $\psi 
(\partial \beta )$ along $\ell '$ is $\leq s(\nu )\sqrt{a(\beta )}$. 
Then  we also have  the Poincar\'e length of $\psi 
(\partial \beta )$ is $<r(\nu )$, assuming without loss of 
generality that $s$ is sufficiently strong given $r$.  The bound on 
$d_{\beta }(y,y')$ for all $y$ means that $\beta $ is $\nu '$-thick 
along $\ell '$ for some $\nu '$ depending only on $\nu $. So 
altogether, by suitable choice of $(\Delta ,r,s)$ given $(\Delta 
',r's')$ we can ensure that $\beta $ is long, $\nu '$-thick 
and dominant along $\ell '$ for $(\Delta ',r',s')$.

Finally, let $\beta $ be a loop. By \ref{5.6},  the quantity 
$n_{\beta  
}([\varphi ])$ only changes for $[\varphi ]\in \ell $.
We have 
$$n_{\beta }(z_{j})=n_{\beta }(\pi _{\alpha }(y_{j}))+O(1)=n_{\beta 
}(y_{j})+O(1)$$
So
$$n_{\beta }(y_{1})-n_{\beta }(y_{0})=n_{\beta }(z_{1})-n_{\beta 
}(z_{0})+O(1)=\vert \ell \vert +O(1).$$
The only way to achieve this is if there is $\ell 
'=[w_{0},w_{1}]\subset [y_{0},y_{1}]$ along which $\beta $ is 
$K_{0}$-flat and with $n_{\beta }(w_{j}')=n_{\beta }(y_{j})+O(1)$. 
Then assuming $L_{0}$ is suitable chosen we do indeed have, for each 
$y\in \ell $, a corresponding $y'\in \ell '$ with 
$$\vert {\rm{Re}}(\pi _{\beta }(y')-\pi _{\beta }(y)\vert 
\leq C_{0}.$$
\Box  

\ssubsection{``Orthogonal projection'' for geodesics in 
${\cal T}_{\geq \nu }$.}\label{7.8}

Now we  describe an analogue of orthogonal projection for a 
geodesic segment in $\cal{T}(S)$. In order to describe the idea, we 
first consider the definition for a geodesic segment  
$$[y_{-},y_{+}]\subset {\cal{T}}_{\geq \nu }$$
for a fixed $\nu >0$. In this case we define
$$x=x(.,[y_{0},y_{1}]):{\cal{T}}(S)\to [y_{0},y_{1}]$$
as follows. Take any $z\in {\cal{T}}(S)$. Then by the special case 
Triangle Theorem of \ref{7.4}, 
there is $y\in [y_{-},y_{+}]$, unique up to moving it a bounded 
distance in $\cal{T}(S)$, such that there exist $y'\in [y_{-},z]$ 
and $y''\in [y_{+},z]$ such that, for $C(\nu )$ as in the special 
case Triangle 
Theorem
$$d(y,y')\leq C(\nu ),\ \ d(y,y'')\leq C(\nu ).$$
We then choose such a $y$ for each $z$ and define 
$$x(z)=x(z,[y_{-},y_{+}])=y.$$ 
The function $x$ is not continuous (unless we are more careful with 
the definition, at least), but it is coarse Lipschitz, and hence, by 
a 
coarse Intermediate Value Theorem (since $x(y_{-})=y_{-}$ and 
$x(y_{+})=y_{+}$) coarsely surjective onto $[y_{-},y_{+}]$ along any 
path in 
${\cal{T}}(S)$ 
joining $y_{-}$ and $y_{+}$.
\begin{ulemma} If $[y_{-},y_{+}]\subset ({\cal{T}}(S))_{\geq \nu }$ 
then 
there 
are  $L(\nu )$ and $C_{1}(\nu )$ such that for any $y_{1}$, $y_{2}\in 
\cal{T}(S)$, if $d(x(y_{1}),x(y_{2})\geq L(\nu )$, with $x(y_{1})$ 
nearer $y_{-}$ than $x(y_{2})$, then there are 
$y_{1}'$, $y_{2}'\in [y_{1},y_{2}]$ with
$$d(y_{i}',x(y_{i}))\leq C_{1}(\nu ),$$
giving
\begin{equation}\label{7.8.1}d(y_{1},y_{2})\geq 
d(y_{1},x(y_{1}))+d(y_{2},x(y{2}))+d(x(y_{1}),x(y_{2}))-4C_{1}(\nu 
).\end{equation}\end{ulemma}

\noindent {\em{Proof.}} From the definition of $x(y_{i})$, there are 
points $y_{i,1}\in [y_{-},y_{i}]$, $y_{i,2}\in [y_{i},y_{+}]$ and 
 such that, for $j=1$, $2$, and $i=1$, $2$,
$$d(y_{i,j},x(y_{i}))\leq C(\nu ).$$
Hence, for $j=1$, $2$
$$d(y_{2,j},y_{1,j})\geq L(\nu )-2C(\nu)$$
and so, assuming $L(\nu )$ sufficiently large given $C(\nu )$, 
$y_{2,1}$ cannot be within $C(\nu )$ of any point 
on $[y_{-},y_{1}]$. So considering the triangle 
with 
vertices at $y_{1}$, $y_{2}$, $y_{-}$, and again applying the special 
case of the Triangle Theorem of \ref{7.4}, there must be 
$y_{2}'\in [y_{1},y_{2}]$ with
$$d(y_{2,1},y_{2}')\leq C(\nu ).$$
This gives the existence of $y_{2}'$, for $C_{1}(\nu )=2C(\nu )$, for 
$C(\nu )$ as in \ref{7.4}. The existence of $y_{1}'$ is similar, and 
(\ref{7.8.1}) follows.\Box

\ssubsection {``Orthogonal projection'':  two set-valued 
functions.}\label{7.9} 

Now let $[y_{-},y_{+}]$ be any geodesic segment in ${\cal{T}}(S)$. 
Before defining a function $x$ with values in ${\cal{T}}(S)$, we 
shall define two set-valued functions $T(z,+)$ and $T(z,-)$ with 
essentially complementary values in $S\times [y_{-},y_{+}]$. We shall 
sometimes use these rather than the ``orthogonal projection'' itself.

Fix ltd parameter functions  $(\Delta ,r,s,K_{0})$ for which the 
results general case of the Triangle Theorem 
of \ref{7.4} holds. By \ref{5.7}, we can 
choose a vertically efficient partition  of $S\times [y_{-},y_{+}]$ 
into sets $\alpha \times \ell $, where each $(\alpha ,\ell )$ is 
either ltd or bounded. Let $\cal{P}$ denote the set of the 
$(\alpha ,\ell )$ from the partition. 
We have the ordering of \ref{7.3}, which is transitive  on $\cal{P}$.
Fix any $z\in \cal{T}(S)$. By \ref{7.4}, for each $(\alpha ,\ell )\in 
{\cal{P}}$, $\ell $ is the disjoint union $\ell _{-}\cup 
\ell _{+}$, where one of these two segments could be empty, such that 
$\ell _{-}$ is a bounded $d_{\alpha }$ -distance from a 
corresponding segment on $[y_{-},z]$, and similarly for $(\alpha 
,\ell _{+})$, with the usual modifications if $\alpha $ is a loop.
Write $T_{\rm{ltd}}(z,-)=T_{\rm{ltd}}(z,-,[y_{-},y_{+}])$ for 
the resulting set of $(\alpha 
,\ell _{-})$ and $T_{\rm{ltd}}(z,+)$ for the set of $(\alpha 
,\ell _{+})$. Then let $T(z,-)=T(z,-,[y_{-},y_{+}])$ be the union of 
$T_{\rm{ltd}}(z,-)$ and of all bounded $(\beta ,\ell ')\in {\cal{P}}$ 
with $(\beta ,\ell ')\leq (\alpha ,\ell _{-})$ for some $(\alpha 
,\ell _{-})\in T_{\rm{ltd}}(z,-)$ and  of all $(\alpha ',\ell ')\in 
\cal{P}$ for which there is 
no $(\alpha ,\ell )\in T_{\rm{ltd}}(z,+)$ with $(\alpha ',\ell 
')>(\alpha ,\ell )$. We define $T(z,+)=T(z,+,[y_{-},y_{+}])$, and 
similarly with 
$T_{\rm{ltd}}(z,+)\subset T(z,+)$. Note that some bounded 
$(\alpha ,\ell )$ 
are likely to be in both $T(z,-)$ and $T(z,+)$, maximal elements in 
$T(z,-)$ and minimal elements in $T(z,+)$. The definitions are such 
that $S$ is the disjoint union of those  $\alpha $ such that $(\alpha 
,\ell )$ is maximal in $T(z,-)$, and similarly for $T(z,+)$ and 
the $(\alpha ,\ell )$ minimal in $T_{z,+}$. We write 
$T_{\rm{max}}(z,-)$ for the set of maximal elements in $T(z,-)$ and 
$T_{\rm{min}}(z,+)$ for the set of minimal elements in $T(z,+)$. The 
$\alpha $ with $(\alpha ,\ell )\in T_{\rm{max}}(z,-)$ are disjoint 
and 
their union is $S$, and similarly for $T_{\rm{min}}(z,+)$.

By construction, the sets $T(z,+)$ and $T(z,-)$ are coarse Lipschitz 
in $z$, in a natural sense.

\ssubsection{``Orthogonal projection'': Upper and Lower 
Boundary and $x(.,.)$.}\label{7.10}

The sets $T(z,-)$ and $T(z,+)$  fit into a more general 
framework of taking a vertically efficient partition $\cal{P}$ of 
$S\times [y_{-},y_{+}]$ into ltd and bounded sets $(\alpha ,\ell )$,
taking a maximal unordered set of 
ltd's 
in this partition, splitting the elements of this into two sets 
$E_{\rm{ltd}}(-)$ and $E_{\rm{ltd }}(+)$ and forming 
resulting sets 
$E(-)$ and $E(+)$ with maximal and minimal sets $E_{\rm{max}}(-)$ and 
$E_{\rm{min}}(+)$ respectively. We shall say that such $E(\pm )$ are 
obtained 
from 
an {\em{order splitting}} $E$ of $\cal{P}$. As in \ref{7.9}, the 
sets $E(-)$ and $E(+)$ are probably not disjoint because they can 
have some bounded $(\alpha ,\ell )$ in common.

The {\em{upper boundary}} 
of 
$E(-)$ is the set of all $(\alpha ,y)$ such that $(\alpha ,\ell )\in 
E_{\rm{max}}(-)$ and $y$ is the right endpoint of $\ell $. The lower 
boundary of $E(+)$ is defined similarly.

The upper and lower boundary can be used to define a single  element 
$ x(E(+))=x(E(-))$ of 
${\cal T}(S)$ up to bounded Teichm\"uller distance. If $E(\pm 
)=T(z,\pm )$ then we can regard this element of ${\cal{T}}(S)$ as 
the orthogonal projection of $z$ --- which we shall do sometimes, but 
not always, because it makes for worse constants, which seems an 
unnecessary complication. The $\alpha $ with $(\alpha ,y)$ in the 
upper boundary of $E_{\rm{max}}(-)$ are disjoint, with union $S$, and 
similarly for the lower boundary. Then we can define $x(E(-))$ 
up to 
bounded 
distance by defining the image under projections $\pi _{\alpha }$. So 
let $(\alpha ,y)$ be in the upper boundary of $E(-)$. If $\alpha $ is 
a gap, we define
$$\pi _{\alpha }(x(E(-)))=\pi _{\alpha }(y).$$
If $\alpha $ is a loop, we define 
$${\rm{Re}}(\pi _{\alpha }(x(E(-))))={1\over \varepsilon _{0}},\ \ 
{\rm{Im}}(\pi _{\alpha }(x(E(-))))={\rm{Im}}(\pi _{\alpha }(y)).$$
 So we are stipulating that  loops  
 $\alpha $ are not short at $x(E(-))$, for $\alpha $ in the 
decomposition. We define $x(E(+))$ similarly using the lower boundary 
of $E(+)$. Then $d(x(E(+)),x(E(-)))$ is bounded in terms of the ltd 
parameter functions and flat constant. This follows because if 
$(\alpha _{+},y_{+})$ and $(\alpha _{-},y_{-})$ are in the lower and 
upper boundaries of $E(+)$, $E(-)$ and $\alpha _{+}\cap \alpha 
_{-}\neq \emptyset $, then $y_{-}=y_{+}$ if $\alpha _{-}\neq \alpha 
_{+}$, and $d_{\alpha }(y_{-},y_{+})$ is bounded if $\alpha 
_{-}=\alpha _{+}=\alpha $. So $x(E)=x(E(-))=x(E(+))$ is well-defined 
up to bounded distance.

If $E(\pm )=T(z,\pm ,[y_{-},y_{+}])$, then we shall denote $x(E)$ by 
$x(z,[y_{-},y_{+}])$, or sometimes by $x(z)$, if the context is 
clear. So we no longer have $x(z)\in [y_{-},y_{+}]$, even up to 
bounded distance, as in the case 
of $[y_{-},y_{+}]\subset {\cal{T}}_{\geq \nu }$. But we do have this 
for suitable coordinates.

The following lemma is a general analogue of \ref{7.8}.  

\begin{lemma}\label{7.11} Choose ltd parameter functions $(\Delta 
,r,s,K_{0})$  such 
that the results on this section hold, and sufficiently strong 
given  a first set of ltd functions $(\Delta ',r',s',K_{0}')$.
There is a function $C_{1}(\nu )$ and constant $C_{1}$ depending only 
on $(\Delta 
,r,s,K_{0})$ such that the following hold.
Let $[y_{-},y_{+}]$ be any geodesic segment, and take sets 
$T(y,\pm )=T(y,\pm ,[y_{-},y_{+}])$ with respect to $(\Delta 
,r,s,K_{0})$. 

Let $y_{1}$ ,$y_{2}\in \cal{T}(S)$. Let $\ell =\ell _{1}\cap \ell 
_{2}$, 
where $(\beta ,\ell _{1})\in T(y_{1},+)$,
 $(\beta ,\ell _{2})\in T(y_{2},-)$,  with $\beta $ long $\nu$-thick 
and 
dominant or $K_{0}$-flat along $\ell $ 
 Then (as in \ref{7.8}), there is $\ell '\subset [y_{1},y_{2}]$ such 
that $\beta $ is ltd for $(\Delta ',r',s',K_{0}')$, and such that,
for all $w\in \ell $, 
there is $w'\in \ell '$ with 
$$d_{\beta }(w,w')\leq C_{1}(\nu ){\rm{\ or\ }}\vert {\rm{Re}}(\pi 
_{\beta }(w)-\pi _{\beta }(w')\vert \leq C_{1}.$$
and for all $w_{1}$, $w_{2}\in \ell '$
\begin{equation}\label{7.11.1}d(y_{1},y_{2})\geq d_{\beta 
}'(y_{1},w_{1})+d_{\beta 
}'(y_{2},w_{2})+d_{\beta }(w_{1},w_{2})-4C_{1}\end{equation}
where $C_{1}=C_{1}(\nu )$ if $\beta $ is a gap and is a fixed 
constant 
if $\beta $ is a loop.
\end{lemma}

\noindent {\em{Proof.}} 
The argument is similar to \ref{7.8}. There 
are  segments corresponding  to $(\beta ,\ell )$ on $[y_{1},y_{+}]$ 
and $[y_{2},y_{-}]$. Since there is a corresponding segment on 
$[y_{1},y_{+}]$, 
there cannot be one on $[y_{-},y_{1}]$. So considering the triangle 
with vertices $y_{1}$, $y_{2}$, $y_{-}$, and using the general 
Triangle 
Theorem of \ref{7.4}, there must be a corresponding 
segment on $[y_{1},y_{2}]$ to the one on $[y_{2},y_{-}]$, for the 
chosen parameter functions 
$(\Delta ',r',s')$ and $K_{0}'$ if $(\Delta ,r,s)$ and  $K_{0}$ are 
suitably chosen given these. (\ref{7.11.1}) then follows as  
(\ref{7.8.1}).
\Box

\ssubsection {}\label{7.12}
In hyperbolic geometry a $\kappa $-quasi-geodesic path is  
distance $O(\kappa )$ from a geodesic. No such precise result is 
available 
in Teichm\" uller geometry. But there is a simple result for 
geodesic 
segments in $\cal{T}_{\geq \nu }$ for any fixed $\nu >0$, which are 
proved analogously to corresponding results in hyperbolic space. 
Variants on 
the following lemma are possible, and some of those statements might 
be slightly simpler than the following, but this is precisely the 
form 
we shall need. An earlier proof that a quasigeodesic in 
${\cal{T}}(S)$ with endpoints joined by a geodesic in 
${\cal{T}}_{\geq 
\nu }$ is a bounded distance from that geodesic appears in 4.2 of 
\cite{Min5}.

\begin{ulemma} The following holds for any sufficiently large 
$L_{1}$. Let $z_{\pm }\in \cal{T}(S)$ with 
$[z_{-},z_{+}]\subset \cal{T}_{\geq \nu}$ and $d(z_{-},z_{+})\geq 
4L_{1}$. Let 
$x(.)=x(.,[z_{-},z_{+}])$ as in \ref{7.8}.
 Let $\{ y_{i}:0\leq i\leq n\} $ be a sequence in $\cal{T}(S)$ with 
 $d(x(y_{0}),z_{-})\leq L_{1}$, $d(x(y_{n}),z_{+})\leq L_{1}$, and 
 $d(y_{i},y_{i+1})\leq L_{1}$ 
for all $i$. Let $\{ z_{i}:0\leq i\leq n\} $ be a sequence of 
successive points on $[z_{-},z_{+}]$ with $z_{0}=z_{-}$, 
$z_{n}=z_{+}$ 
and $L_{1}^{-1}\leq d(z_{i},z_{i+1})\leq L_{1}$ for all $0\leq i<n$.  
Suppose 
 also that there is a function $K(L)$ such that, whenever
 $d(z_{i},y_{j})\leq L$ for some $i$ and $j$, then 
 $d(z_{j},y_{j})\leq K(L)$.
 Then there is $L_{2}$ depending only on $L_{1}$, $\nu $ and the 
 function $K(L)$ such that for all $i$,
 $$d(z_{i},y_{i})\leq L_{2}.$$
 \end{ulemma}
 
 \noindent {\em {Proof.}} We have $x(y,[z_{-},z_{+}])\in 
 [z_{-},z_{+}]$ for  all $y\in \cal{T}(S)$ by the definition of $x$ 
 in \ref{7.8}.  Fix $L_{2}$, to be taken sufficiently large given 
 $L_{1}$ and $\nu $. Take any $i_{1}<i_{2}$ such that, for 
$i_{1}<i<i_{2}$, 
 $$d(y_{j},z_{i})\geq L_{2}{\rm{\ for\ all\ }}j$$
 and either $i_{1}=0$ or $d(y_{j},z_{i_{1}})\leq L_{2}$ for some $j$, 
and 
 either $i_{2}=n$ or $d(y_{j},z_{i_{2}})\leq L_{2}$ for some $j$. To 
prove 
 the lemma, it suffices to bound $i_{2}-i_{1}$ in 
 terms of $L_{1}$, if $L_{2}$ is suitably defined in terms 
 of $L_{1}$,
  and to obtain a contradiction for $L_{2}$ sufficiently large 
 in terms of $L_{1}$, if $i_{2}-i_{1}=n$ and $d(z_{0},y_{0})\geq 
L_{2}$. 
 
Since $i_{2}-i_{1}\leq L_{1}d(z_{i_{1}},z_{i_{2}})$, we may suppose 
that 
 $$d(z_{i_{1}},z_{i_{2}})\geq 4K(L_{2}).$$ 
 Then since $x(.)$ is coarse Lipschitz, we can
 choose a sequence $i(j)$, $0\leq j\leq r$, such that the points 
 $x(y_{i(j)})$ occur in strictly increasing 
 order  in $[x(y_{i_{1}}),x(y_{i_{2}})]\subset [z_{-},z_{+}]$, 
 with $i(0)=i_{1}$, $i(r)=i_{2}$, and 
 \begin{equation}\label{7.12.1}
 2L_{1}\leq d(x(y_{i(j)}),x(y_{i(j+1)}))\leq 4L_{1},\end{equation} 
 $$d(z_{i_{1}},x(y_{i_{1}}))\leq K(L_{2}),
 \ \ d(z_{i_{2}},x(y_{i_{2}})\leq K(L_{2}).$$ 
 Then
 $$d(x(y_{i_{1}}),x(y_{i_{2}}))\geq {1\over 
2}d(z_{i_{1}},z_{i_{2}}).$$
 It follows from (\ref{7.12.1}) that 
 $$r\geq {i_{2}-i_{1}\over 8L_{1}^{2}}.$$

 Assuming $L_{1}$ is sufficiently large 
 given $C_{1}=C_{1}(\nu )$ of \ref{7.8}, by (\ref{7.8.1}),
 $$d(y_{i(j)},y_{i(j+1)})\geq $$
 
$$d(y_{j},x(y_{i(j)})+d(x(y_{i(j)}),x(y_{i(j+1)}))+d(x(y_{i(j+1)}),y_{i(j+1)})
 -4C_{1}$$
 $$\geq d(y_{j},x(y_{i(j)}))+d(x(y_{i(j+1)}),y_{i(j+1)}).$$
 For $1\leq j<r$, $d(y_{j},x(y_{i(j)})\geq L_{2}$. So
 $$L_{1}(i_{2}-i_{1})\geq 
 \sum _{j=0}^{r-1}d(y_{i(j)},y_{i(j+1)})$$
 $$\geq d(y_{i_{1}},x(y_{i_{1}}))+{1\over 2}(r-1)L_{2}\geq 
 {L_{2}\over 16L_{1}^{2}}(i_{2}-i_{1})-{1\over 
 2}L_{2}$$
 Now put $L_{2}=32L_{1}^{3}$. Then we have
 $$(i_{2}-i_{1})\leq 16L_{1}^{2}.$$
 So in all cases, we have a  bound on $i_{2}-i_{1}$ in terms of 
$L_{1}$,
 if $L_{2}=32L_{1}^{3}$. Also,
 $$d(y_{i_{1}},x(y_{i_{1}}))\leq L_{1}(i_{2}-i_{1}).$$
 So then, 
 $$d(y_{i_{1}},z_{i_{1}})\leq 16L_{1}^{3}+K(L_{2})$$
 and for $i_{1}\leq i\leq i_{2}$,
 $$d(y_{i},z_{i})\leq 48L_{1}^{3}+K(L_{2}).$$
 Since the righthand side is $>L_{2}$, we have this for all $i$.
 
 \Box
 
 The generalisation to ${\cal{T}}(S)$ that we shall use
is somewhat weaker than this --- necessarily so. 

\begin{ulemma}
The following holds for a constant $L_{0}$, $\Delta _{0}$ given fixed 
ltd parameter 
functions $(\Delta ,r,s,K_{0})$, related constant $\nu _{0}$ as in 
\ref{5.5},   
and given constants  $\kappa _{0}$, $L_{1}$, 
$\lambda _{0}>0$. 
 Let $d(y_{i},y_{i+1})\leq L_{1}$ for 
    all $i$ and 
$$\sum _{i=0}^{n-1}d(y_{i},y_{i+1})\leq \kappa _{0}d(y_{0},y_{n}).$$
Suppose also that $\{ (\alpha _{i},\ell _{i}):1\leq i\leq m$ is a 
totally ordered set 
of ltds along $S\times [y_{0},y_{n}]$, such that if $\alpha _{i}$ is 
a 
gap, then $\alpha _{i}$ is $\nu $-thick, 
long and dominant along $\ell _{i}$, for some $\nu \geq \nu _{0}$, 
and, if $\alpha _{i}$ is a loop, $\alpha _{i}$ $K_{0}$-flat along 
$\ell _{i}$.
Then if $[x_{1},x_{2}]\subset [y_{0},y_{n}]$ is such that 
$$\sum _{j=1}^{n}\vert \ell _{j}\cap [x_{1},x_{2}]\vert \geq \lambda 
_{0}d(y_{0},y_{n}),$$
there are  at least one $i$, $j$, $y=y_{i}$  and $w\in \ell 
_{j}\cap 
[x_{1},x_{2}]$ 
 such that
$$d_{\alpha _{j}}'(y,w)\leq L_{0}.$$

\end{ulemma}

\noindent {\em{Proof.}} Suppose the lemma is not true for 
$[x_{1},x_{2}]$. Let $C_{1}$ be as in \ref{7.11}. Removing some of 
the 
$y_{i}$, and assuming $L_{1}$ is large enough that $L_{1}/2\geq 
\Delta 
(\nu )$ for all $\nu \geq \nu _{0}$, we can 
assume that for each $0\leq i<n$, the sets $(\alpha _{j},\ell 
_{j,i})$ in $T(y_{i},+)\cap T(y_{i+1},-)$ for $\ell _{j,i}\subset 
\ell 
_{j}$ are for $j\in I(i)$ where
$$L_{1}\leq \sum _{j\subset I(i)}\vert \ell _{j,i}\vert 
\leq 
2L_{1}.$$
Of course, we no longer have the upper bound on $d(y_{i},y_{i+1})$, 
but we still have the same upper bound on $\sum 
_{i}d(y_{i},y_{i+1})$. 
By our assumption that the lemma is not true, for all $w\in \ell 
_{j}$, all $j\in I(i)$,
$$d_{\alpha _{j}}'(w,y_{i})\geq L_{0},\ \ d_{\alpha 
_{j}}'(w,y_{i+1})\geq L_{0}.$$
Let $A$ be the set of  $i$ 
such that $\ell _{j,i}\subset [x_{1},x_{2}]$ for all $j\in I(i)$.
Then by \ref{7.11} we have, for all $w\in \ell _{j,i}$, and $i\in 
A$, $j\in I(i)$,
$$d(y_{i},y_{i+1})\geq d_{\alpha _{j}}'(w,y_{i})+d_{\alpha 
_{j}}'(w,y_{i})-2C_{1}\geq 2L_{0}-2C_{1}.$$
Now we also have
$$\# (A)\geq {\lambda _{0}\over 2L_{1}}d(y_{0},y_{n}).$$
So then we get, assuming $L_{0}\geq 2C_{1}$, as we may do,
$$\sum _{i\in A}d(y_{i},y_{i+1})\geq 
{\lambda _{0}L_{0}\over 2L_{1}}d(y_{0},y_{n}).$$
This gives a contradiction if $L_{0}$ is large enough that
$${\lambda _{0}L_{0}\over 2L_{1}}>\kappa _{1}.$$
\Box

Note that we do not get such a strong result as in the case when 
$[y_{0},y_{n}]\subset \cal{T}_{\geq \nu }$ for some $\nu $ --- and 
actually, our assumption in \ref{7.11} was a little weaker than this. 
The reason 
is that in the case $[y_{0},y_{n}]\subset \cal{T}_{\geq \nu }$ we can 
deduce a bound on $d(y_{i},y_{j})-d(x(y_{i}),x(y_{j}))$ 
from bounds on $d(y_{i},x(y_{i}))$ and $d(y_{j},x(y_{j}))$. But what 
we have above is (effectively) a bound on $d_{\alpha 
}(y_{i},x(y_{i}))$ for some $i$ and for some subsurface $\alpha $. 
But 
we used the bound on $d(y_{0},y_{n})$ to get this, and would need to 
bound $d(y_{i},y_{j})$ to proceed further So there is no real 
analogue of 
the quasi-geodesic-implies-geodesic result for Teichm\"uller 
geodesics 
in general. Instead, it seems to be possible to obtain results for 
families of paths following certain rules. The paths  
through  pleated surfaces that will be used do follow such a set of 
rules, as we shall see.

\ssubsection{A Chain of ltd's.}\label{7.13}

We shall use the following generalization of \ref{5.4}. The notation 
$(\alpha _{j},\ell _{j})$ is used, the same notation as in \ref{5.7}. 
In subsequent sections the $\{ (\alpha _{j},\ell _{j}):1\leq j\leq 
R_{1}\} $ produced below will be  a subset of the set in \ref{5.7}.

The proof of the result is different in character from that of 
\ref{5.4}, being, essentially, a construction of a zero measure 
Cantor 
set, while \ref{5.4} obtained a set with a lower bound on area. This 
result is in any case more sophisticated, because it uses 
\ref{5.5} --- and hence also \ref{5.4} --- in the course of the 
proof. This result can be regarded as a parallel to the existsence of 
a tight geodesic in the curve complex used by Minsky et al.. For 
reasons which are not entirely clear to me, but which may be 
significant, this result appears to be much harder to prove.

\begin{utheorem} Fix long thick and dominant parameter 
functions and flat constant $(\Delta ,r,s,K_{0})$. 
Then there exist $\Delta _{0}$ and $\nu 
_{0}$  depending 
only on $(\Delta ,r,s,K_{0})$
and the topological type of $S$ such that the following holds. 
Let $[y_{0},y_{T}]=[[\varphi _{0}],[\varphi _{T}]]$ be any geodesic 
segment 
in 
${\cal T}(S)$ 
of length $T\geq \Delta _{0}$, parametrised by length. 
Then 
there exists a sequence $(\alpha _{i},\ell _{i})$ ($1\leq i\leq 
R_{1}$) of ltd's  with $(\alpha _{i},\ell _{i})<(\alpha _{i+1},\ell 
_{i+1})$ 
for $i<R_{0}$, where the ordering $<$ is as in \ref{7.3} and such 
that  each segment of $[y_{0},y_{T}]$ of 
length 
$\Delta _{0}$ intersects some $\ell _{i}$. 

\end{utheorem}

\noindent {\em{Start of Proof.}}

We shall prove the theorem by showing that, if $\Delta _{0}$ is 
sufficiently large, for some loop $\gamma \subset S$
and some $\xi \in \gamma $, for every segment $\ell $ of 
length $\Delta _{0}$, along $[y_{0},y_{T}]$, there is an ltd $(\alpha 
,\ell ')$
with $\ell '\subset \ell $ and $\varphi _{0}(\xi )\in \varphi 
_{0}(\alpha )$. It 
then follows that the intersection   of all such $\varphi _{0}(\alpha 
)$ is 
nonempty and therefore any two such $\alpha $ intersect essentially.
It is not quite the case that any sequence $(\alpha _{i},\ell _{i})$ 
for $i\leq j$ is extendable, but, essentially,  if 
$(\alpha _{j},\ell _{j})$ is chosen so that $\alpha _{j}$ does not 
intersect any bounded gaps (in the sense of \ref{5.5}) for distance 
$\Delta _{0}^{1/n}$ beyond $\ell _{j}$, for a suitable $n$ depending 
only on the topological type of $S$,
then further extension becomes possible.  

 Let $y_{t}=[\chi _{t}\circ \varphi _{0}]$ be such that 
$d(y_{0},y_{t})=t$ with stretch $q_{t}(z)dz^{2}$ at $y_{t}$ and 
$d(y_{-},y_{+})=t$ and with $\chi _{t}$ minimising distortion. We use 
this parameterisation for the whole geodesic contaniing 
$[y_{0},y_{T}]$. We 
shall 
write $\vert .\vert _{t}$ for $\vert .\vert _{q_{t}}$ and $\vert 
.\vert _{t,+}$ for the unstable length $\vert .\vert _{q_{t},+}$.
We fix a vertically efficient decomposition of $S\times 
[y_{0},y_{T}]$ as 
in \ref{5.6}, in which every $(\beta ,\ell )$ is either ltd or 
bounded by $L$ in the sense of \ref{5.5}, with $L$ as in \ref{5.5}.
In particular, for $[\varphi ]\in \ell $:

\begin{equation}\label{7.13.2}\vert \varphi (\partial \beta )\vert 
\leq 
L.\end{equation}

In fact, (\ref{7.13.2}) holds for all $[\varphi ]\in \ell '(\partial 
\beta )$, 
where $\ell '(\partial \beta )=
\{ [\chi _{t}\circ \varphi _{0}]:t_{0}(\partial \beta )\leq t\leq 
t_{0}'(\partial \beta )\} $ is the possibly larger interval 
containing $\ell $ 
such that $\partial \beta $ does not intersect $\alpha $ 
transversally, for any $\alpha $ which is ltd along a segment of 
$\ell '$. Here, $t_{0}(\partial \beta )<0$ and $t_{0}'(\partial 
\beta )>T$ are allowed. Also, for $y_{t}=[\chi _{t}\circ \varphi 
_{0}]\in \ell $, 
by the last 
part of \ref{5.4} we have, if $(\beta ,\ell )$ is bounded by $L$,

\begin{equation}\label{7.13.1}\vert \chi _{t}\circ \varphi 
_{0}(\partial \beta 
)\vert _{t}^{2}\geq  s_{0}a'(\beta )\end{equation} 
for a constant 
$s_{0}$ depending only on the ltd parameter functions
 and the topological 
type of $S$. Here, $a'(\beta )$ is as in \ref{5.2}.

Now for any loop or finite set of (not necessarily disjoint) loops 
$\gamma $, 
 let $F(t,\gamma )$ be the qd-length function as in \ref{7.20}, with 
minimum $t(\gamma )$.
Now suppose that $(\beta ,\ell )$ is bounded by $L$ and that $\ell $ 
is of 
length $\geq \Delta _{0}$. Then from (\ref{7.13.2}), (\ref{7.13.1}) 
and (\ref{7.20.1}) 
we obtain, for a constant $C_{1}$ depending only on the topological 
type of $S$,
\begin{equation}\label{7.13.4}a'(\beta )\leq C_{1}L^{2}e^{-\Delta 
_{0}}.\end{equation}
Also, by definition,
$$t_{0}(\partial \beta )\leq t(\partial \beta )\leq t_{0}'(\partial 
\beta ).$$

\begin{lemma}\label{7.14} The following holds for a suitable constant 
$L_{1}$ depending only on the topological type of $S$, 
if $\Delta _{0}$ is sufficiently large given a constant $L$. Let 
$\beta $ be 
bounded by $L$ along $\ell \subset [y_{0},y_{T}]$, where $\ell $ has 
length $\geq \Delta _{0}$.    
Then there is 
$(\beta _{1},\ell 
_{1})$ bounded by 
$L_{1}$, such that $\beta \subset \beta _{1}$, 
$\beta $ is  properly contained in $\beta _{1}$ if 
$t(\partial \beta)\leq t_{0}(\partial \beta )+\sqrt{\Delta _{0}}$
and $\ell _{1}\subset 
\ell $ has length at least $\Delta _{0}^{1/3}$.
\end{lemma}

By 
(\ref{7.13.4}) with $\Delta _{0}$ replaced by $\Delta _{0}^{1/3}$, 
the union of all such $\beta _{1}$ must also be properly 
contained in $S$, assuming $\Delta _{0}$ is sufficiently large given 
the ltd parameter functions. Then applying the lemma a number of 
times 
which is bounded in terms of the topological type of $S$, with 
$\Delta 
_{0}$ replaced by $\Delta _{0}^{1/n}$ for different $n$, we obtain 
the following.

\begin{corollary}\label{7.15} The following holds for a suitable 
constant 
$L_{1}$ as above, 
if $\Delta _{0}$ is sufficiently large given $L$, and for $k$ which 
is 
bounded in terms of the topological type of $S$. Let $\ell \subset 
[y_{0},y_{T}]$ be any interval of length $\geq \Delta _{0}$. Then 
there is $\ell '\subset \ell $ of length $\Delta _{1}\geq \Delta 
_{0}^{1/k}$ 
such that, for any $\ell ''\subset \ell '$ of length $\geq \Delta 
_{1}^{1/3}$, the union  of $\beta $ which are bounded by $L$ along 
$\ell '$ is the same as the union of $\beta $ which are bounded by 
$L_{1}$ along $\ell ''$, and for all such $\beta $,
$t(\partial \beta )\geq t_{0}(\partial \beta )+\Delta 
_{1}^{1/2}$.
\end{corollary}

\ssubsection{Proof of \ref{7.14}.}\label{7.16}
The lemma is proved by using the characterisation of bounded 
Poincar\'e  
Poincar\'e length of \ref{7.20}, in particular, the discussion on 
comparing length of non-transverse loops. 
Fix $L_{1}$, depending only on the 
topological type of $S$, so that, for any $[\varphi ]\in 
{\cal{T}}(S)$,
a maximal multicurve  $\Gamma $ with $\varphi (\Gamma )\vert \leq 
L_{1}$ exists. Also fix a Margulis constant $\varepsilon _{0}$.
We consider the function
$$G(t,\partial \beta )={\rm{Max}}(F(t,\partial \beta  
)-F(t,\Gamma ))=F(t,\partial \beta )-{\rm{Min}}F(t,\Gamma ),$$
where, the minimum is taken over loop sets  $\Gamma $ such that
 $\gamma $ satisfies 
\ref{7.20.3} for all $\gamma \in \Gamma $, and with specified
 convex hull for $\Gamma $. The convex hull is a  finite union of 
gaps and loops
 $\alpha $ such that, for each $\alpha $, 
 $\chi _{t}\circ \varphi _{0}(\alpha )$ is homotopic to a component 
of 
$(\chi _{t}\circ \varphi _{0}(S))_{\geq 
\varepsilon _{0}}$. The closure of the 
union of the $\alpha $ contains $\beta $, and 
every such $\alpha $ either intersects $\beta $, or shares a boundary 
component with $\beta $. The possibility that $\alpha $ is disjoint 
from $\beta $ only happens if $t(\partial \beta )<t_{0}(\partial 
\beta )+\sqrt{\Delta _{0}}$. Any boundary component $\gamma _{1}$ 
which is 
shared by $\beta $ and $\alpha $ is such that $F(\gamma _{1},t)$ is 
maximal among $\gamma 
_{1}\subset \partial \beta $.
 Assuming $L_{2}$ is large enough given the topologial type 
of $S$, for each $t$, there is at least one $\Gamma $ satisfying 
these 
conditions,  but only boundedly many can give the minimum value of 
$G(t,\partial \beta )$, with bound depending on $L_{2}$. We write 
$\Gamma _{t}$ for  a choice of $\Gamma $ such that
$$\vert G(\partial \beta ,t)-(F(\partial \beta ,t)-F(\Gamma 
_{t},t))\vert \leq C(\varepsilon _{0}),$$
for a suitable constant $C(\varepsilon _{0})$ depending only on 
$\varepsilon _{0}$ and the topological type of $S$.

For $t$ such that $\vert \chi _{t}\circ f_{0}(\partial \beta )\vert 
\leq L$, that is, $t_{0}(\partial \beta )\leq t\leq t_{0}'(\partial 
\beta )$, $G(t,\partial \beta )$ is bounded above in terms of $L$. 
For $t$ such that $\vert \chi _{t}\circ f_{0}(\partial \beta )\vert 
\geq \varepsilon $, $G(t,\partial \beta )$ is bounded below in terms 
of $\varepsilon $. 
So by our choice of $t_{0}(\beta )$, 
$t_{0}'(\beta )$, $G(t,\partial \beta )$ is bounded above and below 
in terms of $L$ and $\varepsilon _{0}$ (our fixed Margulis constant).
 By the restrictions we have put on the choice of 
$\Gamma $, $G(t,\partial \beta )$ is also bounded below in terms of 
$L$ for all $t\in [t_{0}(\partial \beta ),t_{0}'(\partial \beta )]$, 
if $t(\partial \beta )\geq t_{0}(\partial \beta )+\sqrt{\Delta 
_{0}}$. 

So in all cases, by the choice of $\Gamma _{t}$
for $t\in [t_{0}(\partial \beta ),t_{0}'(\partial 
\beta )]$ there is a bound on the number of different $\Gamma _{t}$ 
in 
terms of $L$ and $\varepsilon _{0}$ if 
$t(\partial \beta )\geq t_{0}(\partial \beta )+\sqrt{\Delta _{0}}$, 
and if  $t(\partial \beta 
)\leq t_{0}(\partial \beta )+\sqrt{\Delta _{0}}$ 
then the number of different sets of 
loops
$\Gamma _{t}$ for $t\in [t_{0}(\partial \beta ),t_{0}'(\partial 
\beta )]$ is $\leq C_{1}(L,\varepsilon _{0})\sqrt{\Delta _{0}})$. 
We see this as follows. The 
graph of $F(t,\partial \beta )$ is an approximate V shape, with 
slope $+1$ and $-1$ to the right and left of the minimum. The graph 
of 
each $F(t,\Gamma )$  is also 
such a $V$ shape, and the graph for $\Gamma =\Gamma _{s}$, 
$s\in [t_{0}(\partial \beta 
),t_{0}'(\partial \beta )]$, must be above that of $t\mapsto 
F(t,\partial 
\beta )-C_{2}(L,\varepsilon _{0})$ for $t\in [t_{0}(\partial \beta 
),t_{0}'(\partial \beta )]$. At 
$t=t_{0}(\partial \beta )$, $t=t_{0}'(\beta )$,
$G(t,\partial \beta )$ is bounded above and below.  Then minima of 
the 
functions $F(t,\Gamma )$ move progressively to the right, and the 
right branches of the $V$'s move down, but all of them within a 
rectangle with sides slopes $-1$ and $+1$, where the $-1$ sides have 
width $\leq C_{3}(L,\varepsilon 
_{0})$ if 
$t(\partial \beta )\geq t_{0}(\partial \beta )+\sqrt{\Delta _{0}}$, 
and width  $\sqrt{\Delta _{0}}+C_{3}(L,\varepsilon 
_{0})$ if $t(\partial \beta 
)\leq t_{0}(\partial \beta )+\sqrt{\Delta _{0}}$.
 For each 
$t$, there are only a bounded number of choices for $\Gamma _{t}$. 
The 
minima cannot get too close, because then there would be too many 
choices for $\Gamma _{t}$, for a slightly larger $L_{2}$. So the 
number of different loop sets $\Gamma _{t}$ which arise must be 
bounded by $C_{1}(L,\varepsilon 
_{0})$ or $C_{1}(L,\varepsilon 
_{0})\sqrt{\Delta _{0}}$ in the respective cases. So there
 must be $\ell _{1}\subset \ell $ of length  $>\Delta 
_{0}^{1/3}$ along which the choice of $\Gamma _{t}$ can be 
chosen to be constant.  Let $\beta _{1}$ be the union of the gaps in 
the  convex hull of $\Gamma _{t}$. Then $\beta _{1}$ is bounded by 
$L_{1}$ along $\ell _{1}$, where $L_{1}$ is independent of $L$, and 
$\beta _{1}$ strictly contains $\beta $ if $t(\partial \beta 
)\leq t_{0}(\partial \beta )+\sqrt{\Delta _{0}}$, as required. \Box

\ssubsection{More on transfer between Poincar\'e and 
qd-length.}\label{7.17}

So now, we consider 
$(\beta ,\ell )$ bounded by $L_{1}$, with $\ell $ of length $\geq 
\Delta _{0}^{1/k}$ and 
$$t(\beta )>t_{0}(\beta )+\Delta _{0}^{1/2k}.$$

  For 
$t<t(\beta )$, the $q_{t}$-length of unstable segments across $\chi 
_{t}\circ f_{0}(\beta )$ is $\leq C_{1}L_{1}e^{t-t(\beta )}$, 
for $C_{1}$ depending only on the topological type of $S$. We now 
need to interpret and strengthen this statement in terms of 
Poincar\'e length, under 
the assumption that $\ell $ has length $\geq \Delta _{0}^{1/k}$ and 
that $\Delta _{0}$ is sufficiently large given the ltd parameter 
functions and flat constant, and hence also sufficiently large 
given $s_{0}$ and $L_{1}$. If $-{1\over 2}\Delta 
_{0}^{1/2k}<t-t(\beta 
)<-{1\over 4}\Delta 
_{0}^{1/2k}$, then we can be sure that the Poincar\'e length of 
unstable 
segments
is $\leq C_{2}e^{t-t(\beta )}$, because the 
$q_{t,+}$-length  is exponentially smaller than the $q_{t,-}$-length, 
and the Poincar\'e length of $\chi _{t}\circ f_{0}(\partial \beta )$ 
is bounded. We also have the following.

\begin{ulemma}
 For $p$ depending only on $L_{1}$, whenever $t-t(\beta 
)<-{1\over 4}\Delta _{0}^{1/2k}$, the following holds. Let $\gamma 
\subset \chi _{t}\circ \varphi _{0}(S)$ be a segment of unstable 
foliation 
whose diameter is boundedly proportional to the injectivity radius at 
any point of $\gamma $. Then all components of 
$\gamma \cap \chi _{t}\circ 
f_{0}(\beta )$ have Poincar\'e length $\leq -\exp (\Delta 
_{0}^{1/3k})$ times the injectivity radius at $\gamma $ and at least 
one in any $p$ complementary components in $\gamma $   has Poincar\'e 
length 
 at least $\exp (\Delta 
_{0}^{1/3k})$ times the 
length of any unstable segment in $\gamma \cap \chi _{t}\circ 
f_{0}(\beta )$.
\end{ulemma}

\noindent{\em{Proof.}}
Because we are restricting to a set whose diameter is boundedly 
proportional to the injectvity radius
we only need to prove this for $-{1\over 2}\Delta 
_{0}^{1/2k}<t-t(\beta 
)<-{1\over 4}\Delta 
_{0}^{1/2k}$, if we replace $\exp (\Delta 
_{0}^{1/3k})$ by $\exp (\Delta 
_{0}^{3/7k})$. For if we can do this, it is true for this range of 
$t$ with respect to $q_{t}$-length, with $\exp (\Delta 
_{0}^{3/7k})$ replaced by $\exp (\Delta 
_{0}^{2/5k})$. Then  since unstable length is multiplied by $e^{t}$ 
under application of $\chi _{t}$ for $t<0$, and we are restricting to 
a set whose diameter is boundedly proportional to the injectivity 
radius, we have the result for all $t<-{1\over 4}\Delta 
_{0}^{1/2k}$, with respect to $q_{t}$-length. But then we also have 
the result for Poincar\'e length, if we replace $\exp (\Delta 
_{0}^{2/5k})$ by $\exp (\Delta 
_{0}^{1/3k})$, because the image under $\chi _{-t}$ of an unstable 
segment of short Poincar\'e length remains inside the injectivity 
radius of the surface, and so Poincar\'e length and $q_{t}$-length 
are 
still comparable. So now we suppose 
that  $-{1\over 2}\Delta _{0}^{1/2k}<t-t(\beta )<-{1\over 4}\Delta 
_{0}^{1/2k}$. Suppose that  $p>e^{3L'}$. Then if we take a string of 
$p$ 
successive 
segments, two of them must intersect $\chi _{t}\circ \varphi 
_{0}(\partial 
\beta )$ at Poincar\'e distance $<e^{-2L_{1}}$ apart. 
The segments of intersection 
with $\chi _{t}\circ \varphi _{0}(\partial 
\beta )$ have Poincar\'e length $<\exp (-{1\over 5}\Delta 
_{0}^{1/2k})$. 
Then if all the components 
outside $\chi _{t}\circ \varphi _{0}(\beta )$ have length 
$<e^{-4L_{1}}$ we 
can make a loop from such segments of Poincar\'e length $<e^{-L_{1}}$ 
intersecting $\chi _{t}\circ \varphi _{0}(\partial \beta )$ 
transversally, 
assuming that $\Delta _{0}$ is sufficiently large. 
This is impossible.\Box

\ssubsection{Proof of \ref{7.13}: construction of the 
sequences.}\label{7.18}
For $1\leq i\leq R_{0}$, some $R_{0}$, we shall 
find   sequences  $t_{i}$,  $\alpha _{i}$, $\zeta _{i}$ such that the 
following hold, for a constant $C_{1}$.
\begin{description}
\item[1.] $\alpha _{i}$ is ltd at $y_{t_{i}}$ for $[y_{0},y_{T}]$.
\item[2.] $\zeta _{i+1}\subset  \zeta _{i}\subset f_{0}(\alpha _{i}$.
\item[3.] For all $i\geq 1$, 
$\chi _{t_{i}}(\zeta _{i})$ is an unstable 
segment whose Poincar\'e length is boundedly proportional, with 
bounds given by $C_{1}$, to the injectivity radius at that point of 
$\chi _{t_{i}}\circ \varphi _{0}(S)$. If $\alpha _{i}$ is a loop, 
$\chi _{t_{i}}(\zeta _{i})$ is an unstable 
segment 
in  a component of $(\chi _{t_{i}}\circ 
\varphi _{0}(S))_{< \varepsilon _{0}}$ homotopic to 
$\chi _{t_{i}}\circ \varphi _{0}(\alpha _{i})$, in the  part which is 
flat with respect to the $q_{t_{i}}$-metric.
\item[4.] $t_{1}\leq \Delta _{0}$ and $t_{R_{0}}\geq T-\Delta _{0}$.
\item [5.] For all $i<R_{0}$, $t_{i}+\Delta _{0}^{1/3k}<t_{i+1}\leq 
t_{i}+\Delta _{0}$.
\item[6.] For all $i$, $\chi _{t_{i}}(\zeta _{i})\cap \chi 
_{t_{i}}\circ \varphi _{0}(\beta )=\emptyset $ for any $(\beta ,\ell 
)$ such that $(\beta ,\ell 
)$ is  bounded by $L_{1}$,  and $\ell $ is of length $\geq \Delta 
_{0}^{1/3k}$  with $\ell $ starting ending within 
${1\over 2}\Delta _{0}^{1/k}$ of the right of $y_{t_{i}}$.
\end{description}
This completes the proof of \ref{7.13},  
since we then have 
$$\zeta _{j}\cap \varphi _{0}(\alpha _{i})\neq \emptyset $$ 
for all $i\leq j$.

To choose $(\alpha _{1},\ell _{1})$, let $\ell $ be the initial 
setment 
of $[y_{0},y_{T}]$ of length $\Delta _{0}$. Apply \ref{7.15} to 
find $\ell '$, $\Delta _{1}=\Delta _{0}^{1/n}$ for some $n$
and take any $\ell ''\subset \ell '$ of length $\Delta 
_{1}^{1/3}$ in the first half of $\ell '$. Assuming that the 
ltd parameter functions are sufficiently strong 
given $L_{1}$, 
that a component of a $\partial \beta $ cannot be contained in the 
interior of an ltd gap $\alpha $, and cannot be an ltd loop. 
Then by \ref{5.5} we can 
find  an ltd $(\alpha _{1},\ell _{1})$ with $\ell _{1}\subset \ell '$ 
and with $\alpha $ of sufficiently large qd-area not to be contained 
in any $\beta $, and hence disjoint from all the $\beta $ which are 
bounded 
along $\ell ''$.  We can then choose $\zeta _{1}\subset 
f_{0}(\alpha _{1})$ so that 
$\chi _{t_{1}}(\zeta _{1})$ has Poincar\'e length boundedly 
proportional to the injectivity radius at any point of 
$\chi _{t_{1}}(\zeta _{1})$, in the flat-metric part of $\chi 
_{t_{1}}\circ 
f_{0}(\alpha _{1})$ if $\alpha _{1}$ is a loop.

So now suppose that we have found $\zeta _{i}$, $t_{i}$ and $\alpha 
_{i}$ 
and we look for $\zeta _{i+1}$ by looking for an appropriate subset 
of $\chi _{t_{i}}(\zeta _{i})$.  Then we claim that the set of 
$I(\beta )$ with $(\beta ,\ell )$ bounded by $L_{1}$, $\ell \subset 
[y_{t_{i}},y_{T}]$ of 
length $\geq \Delta _{0}^{1/k}$  is not all of  $\chi _{t_{i}}(\zeta 
_{i})$, 
 if 
$\Delta _{0}$ is sufficiently large in terms of the ltd parameter 
functions. We only need to consider $\beta $ with $t_{0}(\partial 
\beta )\geq t_{i}+{1\over 2}\Delta _{0}^{1/k}$, by condition 6. Then, 
as in the case $i=1$, by \ref{7.15}, for $\Delta _{1}=\Delta 
_{0}^{1/2k^{2}}$, we only 
need to prove this for $\ell $ of length $\geq \Delta _{1}$ and 
\begin{equation}\label{7.13.7}t_{i}+\Delta _{1}\leq t_{0}(\partial 
\beta )\leq \Delta _{0},\end{equation}
\begin{equation}\label{7.13.8}
    t(\partial \beta )\geq t_{0}(\partial \beta )+\sqrt{\Delta 
    _{1}}.\end{equation}
 First, 
 the number of such  $\beta $ is $\leq C\Delta _{0}$ for $C$ 
depending 
 only on the topological type of $S$. Define 
 $$I(\beta )=\chi _{t_{i}}(\zeta _{i})\cap \chi _{t_{i}}\circ \varphi 
_{0}(\beta ).$$
Then, restricting to this set of $\beta $ satisfying (\ref{7.13.7}) 
and 
(\ref{7.13.7}), by \ref{7.17}, the sum of the Poincar\'e lengths of 
the intervals of the set $I(\beta )$ is $\leq C_{1}\Delta 
_{0}pe^{-\Delta _{1}}$ times the injectivity radius at any point of 
$\chi _{t_{i}}(\zeta _{i})$. Since $\chi _{t_{i}}(\zeta _{i})$ has 
Poincar\'e length $\geq C_{1}^{-1}$ times the injectivity radius, 
the complement of the union of 
the $I(\beta )$ is nonempty. Note that this argument would not work 
without the upper bound in (\ref{7.13.7}). The calculation just done 
is the 
first step in constructing a zero measure Cantor set, albeit of 
Hausdorff dimension close to $1$ if $\Delta _{0}$ is large.

Now  
choose an interval in $[y_{t_{i}},y_{t_{i}+\Delta _{0}}]$ as $\ell 
''$ 
in \ref{7.15}, with $\ell ''$ in the first half of $\ell '$, again 
for 
$\ell '$ as in \ref{7.15}. The complement in $\zeta _{i}$ of the 
$I(\beta )$ 
for the 
finitely many $\beta $ for
which $\beta $ is bounded by $L_{1}$ along $\ell '$ is nonempty,by 
\ref{7.16}, and  the complement of this set of $\beta $ has qd-area 
bounded from $0$ and  has 
nonempty intersection with
convex hull of 
the $\alpha $ which are ltd along segments of $\ell ''$, by 
\ref{5.5}, 
assuming that $\Delta _{0}$ is sufficiently large. 
Assuming that the ltd parameter functions are sufficiently strong 
given $L_{1}$, 
that a component of a $\partial \beta $ cannot be contained in the 
interior of an ltd gap $\alpha $, and cannot be an ltd loop. So the 
complement of the union of the $\beta $ is a union of components of 
the 
convex hull of the ltd's. Each 
component $J$ of the complement in $\zeta _{i}$
of these $I(\beta )$  must intersect $\chi _{t_{i}}\circ 
\varphi _{0}(\alpha )$ for some 
ltd  $(\alpha ,\ell )$, because otherwise  we can make a nontrivial 
loop out 
of the boundaries of the $\chi _{t_{i}}\circ \varphi _{0}(\partial 
\beta _{i})$ and $J$ 
which is disjoint from all the $\chi _{t_{i}}\circ \varphi 
_{0}(\alpha )$, 
and homotoping that loop to good position, it remains disjoint from 
$J$, and separates the convex hull from $J$. Take any such $(\alpha 
,\ell )$ 
to be $(\alpha _{i+1},\ell _{i+1})$, and $\zeta _{i+1}$ so that $\chi 
_{t_{i}}(\zeta _{i+1})$ is a component of $J\cap \chi _{t_{i}}\circ 
\varphi _{0}(\alpha _{i+1})$. Then $\chi _{t_{i+1}}(\zeta _{i+1})$ 
has 
Poincar\'e length bounded from $0$, even if we restrict to the 
intersection with $(\chi _{t_{i+1}}\circ \varphi _{0}(S))_{\geq 
\varepsilon 
_{0}}$. Then by the choice of $\ell ''$, $(\alpha 
_{i+1},\ell _{i+1})$ and $\zeta _{i+1}$ have all the required 
properties.

\Box  

\section{Geometric model manifolds.}\label{7}

In this section, we construct the geometric model $M=M(\mu 
(e_{1}),\cdots 
\mu (e_{n}))$ for any homeomorphism type of relative Scott core with 
ends $e_{i}$. We start by 
constructing geometric models in the geometrically finite cases, when 
$\mu (e_{i})\in {\cal{T}}(S(e_{i}))$ 
for all $i$. The simplest case of combinatorial bounded geometry is 
dealt with first, in \ref{6.1}. The other geometrically finite cases 
follow. All models are geometrically finite until \ref{6.12}.
The geometrically infinite case, when
 $\mu (e_{i})\in 
{\cal{O}}_{a}(S(e_{i}))$, for some, or all, $i$,
is constructed by taking geometric limits, 
for which, of course, we need to know that if $\mu ^{k}(e_{i})\to \mu 
(e_{i})$ then the corresponding geometric models $M(\mu 
^{k}(e_{1}),\cdots 
\mu ^{k} (e_{n}))$ converge geometrically, up to Lipschitz 
equivalence, 
to a single limit. This is dealt with in \ref{6.12} in the case of 
combinatorial bounded geometry, and in the general case in 
\ref{6.13}. The geometrical model in the geometrically infinite case 
is a manifold which is also a metric space, defined up to a bounded 
Lipschitz equivalence, and boundedly Lipschitz equivalent to a 
Riemannian manifold whose ends are topological products.

\ssubsection{The combinatorially bounded geometry geometrically 
finite  
Kleinian surface case.}\label{6.1}

We assume that 
the relative Scott core is homeomorpic to $S_{d}\times 
[0,1]$, where $S_{d}$ is the horodisc deletion of a finite type 
surface $S$, and that the end invariants $y_{0}$, $y_{u}\in 
{\cal{T}}(S)$ are such 
that 
$$[y_{0},y_{u}]=\{ y_{t}:t\in [0,u]\} \subset \cal{T}_{\geq \nu }$$
for some fixed $\nu 
>0$. We are parametrising by length, so that $d(y_{s},y_{t})=t-s$ for 
any $s<t$. We can, of course, extend the parametrisation, by $\mathbb 
R$, to the geodesic in ${\cal{T}}(S)$ containing $[y_{0},y_{u}]$. We 
let $S_{t}$ be the hyperbolic surface for $y_{t}$, with 
hyperbolic metric $\sigma _{t}$ on $S_{t}$. We fix a smooth manifold 
structure on $S$. We take any continuous family of 
$C^{1}$ diffeomorphisms $\varphi _{t}:S\to S_{t}$ such that $[\varphi 
_{t}]=y_{t}$ and such that if $\vert t-s\vert <1$, 
$$\varphi _{t}\circ \varphi _{s}^{-1}:S_{s}\to S_{t}$$
has norm of derivative bounded by $K$ for some fixed $K$. Here, the 
norm of the derivative is taken with respect to the metrics $\sigma 
_{t}$ and $\sigma _{s}$. We can also choose the family $\varphi _{t}$ 
so that
$$(z,t)\mapsto \varphi _{s}^{-1}\circ \varphi _{t}(z):S\times \mathbb 
R\to S$$
is $C^{1}$ (or even $C^{\infty }$). Then the model metric  $\sigma 
_{M}$ on $M=S\times [0,u]$ 
is
$$\sigma _{M}=\varphi _{t}^{*}\sigma _{t}+dt^{2}.$$
Fix $s$ and consider the map
$$(z,t)\mapsto (\varphi _{s}(z),t):S\times [0.u]\to S_{s}\times 
[0,u].$$
Then the model metric is transformed to the metric 
$$(\varphi _{s}\circ \varphi _{t}^{-1})_{*}\sigma _{t}+dt ^{2},$$
which, on $S_{s}\times [s-1,s+1]$, is boundedly equivalent to $\sigma 
_{s}+dt^{2}$, because of the bound on the derivative of $\varphi 
_{s}\circ \varphi _{t}^{-1}$. It follows that the exact choice of the 
family $\varphi _{t}$ is unimportant. Subject to the above 
constraints, any choice gives the same geometric model up to 
Lipschitz 
equivalence.

In this special case, it may be preferred to use the singular 
Euclidean  metric on $S_{t}$ coming from the quadratic differential 
$q_{t}(z)dz^{2}$ on $S_{t}$ for the geodesic and to use $\varphi 
_{t}=\chi _{t}\circ \varphi _{0}$ for some $\varphi _{0}$ and for 
$\chi _{t}:S_{0}:S_{t}$ minimising distortion. In this case , we have 
a singular Euclidean structure on $S$ such that the metric on 
$S\times [0,u]$ is
$$e^{2t}dx^{2}+e^{-2t}dy ^{2}+dt^{2}.$$
Although singular, the metric induced by this Riemannian metric is 
boundedly Lipschitz equivalent to the previous one. See also 
\cite{Min1}.

\ssubsection{Geometrically finite Kleinian surface case: $S_{j,t}$ 
and properties of $\varphi 
_{t}$.}\label{6.2}

We fix a Margulis constant $\varepsilon _{0}$ for both dimensions $2$ 
and $3$. 

We  again consider a geodesic segment $[y_{0},y_{u}]=\{ y_{t}:t\in 
[0,u]\} $ and let $S_{t}$, $\sigma _{t}$ be as in \ref{6.1}. We fix 
a vertically efficient deomposition as in \ref{5.7} into sets $\alpha 
_{j}\times \ell _{j}$, $1\leq j\leq R$. Write
$$I_{j}=\{ t:y_{t}\in \ell _{j}\} .$$ 
Then the model manifold $M=M(y_{0},y_{u})$ is given topologically by
$$M=S\times [0,u]=\cup _{j=0}^{R}\alpha _{j}\times I_{j},$$
where each $\alpha _{j}$ is either a gap --- not including the 
boundary --- or a loop. It remains to define the model Riemannian 
metric $\sigma _{M}$. This is done by decomposing $M$ into pieces 
corresponding to each $(\alpha _{j},\ell _{j})$

Now for each $j$ and $t\in I_{j}$ we define  a subsurface 
$S_{j,t}\subset 
S_{t}$ which is 
homotopic to $\varphi _{t}(\alpha _{j})$, for any $[\varphi 
_{t}]=y_{t}$. 
If $\alpha _{j}$ is a loop, then we know from \ref{5.5} that $\vert 
\varphi _{t}(\alpha _{j})\vert \leq L$. Then there is $\varepsilon 
(L)$, 
continuous in $L$, such that 
closed $\varepsilon (\vert \varphi _{t}(\gamma _{i})\vert 
)$-neighbourhoods 
$T(\varphi _{t}(\gamma _{i}))$ of 
nonintersecting geodesics $\varphi _{t}(\gamma _{i})$ are disjoint, 
for any 
hyperbolic 
surface (\cite{Bus} Chapter 4). In fact, 
we can take $T(\gamma )$ to be the 
$\varepsilon _{0}$-Margulis tube $T(\varphi _{t}(\gamma ),\varepsilon 
_{0})$, 
if $\vert \varphi (\gamma )\vert \leq {1\over 2}\varepsilon 
_{0}$.  

We define 
$$S_{j,t}=T(\varphi _{t}(\alpha _{j}))$$
if $\alpha _{j}$ is a loop. If $\alpha _{j}$ is a gap, then $S_{j,t}$ 
is the 
closure of the component of the 
complement of 
$$\cup \{ T(\varphi _{t}(\gamma )):\gamma \subset \partial \alpha 
_{j}\} $$
which is homotopic to $\varphi _{t}(\alpha _{j})$. Note that the 
vertically efficient conditions ensure that, if $\alpha _{j}$ is a 
loop, $S_{j,t}$ has modulus 
bounded below in terms of the ltd parameter functions for all $t\in 
I_{j}$. If $\alpha  $ is any loop  such that $\alpha =\alpha _{j}$ 
or $\alpha \subset \partial \alpha _{j}$  for at least one 
$j$, 
then let $I(\alpha )=[s(\alpha ),t(\alpha )]$ be the connected  union 
of the   $I_{j}$ for all such $j$. For such a $j$ for which $\alpha 
_{j}$ is a loop, define $S_{j,t}=S_{\alpha ,t}$. 
Then, except when $s(\alpha )=0$, or similarly $t(\alpha )=u$,  
$S_{\alpha ,s(\alpha )}$ and $S_{\alpha ,t(\alpha )}$ also have 
modulus bounded above. This is because of
 the properties of vertically efficient in \ref{5.7}. These 
properties 
 imply that we 
also have $s(\alpha )\in 
I_{k}$ (or similarly $t(\alpha )\in I_{k}$) for some ltd gap or loop 
$\alpha 
_{k}$. Here, $\alpha _{k}$ either  contains $\alpha $ in its 
interior, or intersects $\alpha $, transversely. So  there is 
a loop $\gamma \in \alpha _{k}$ which 
intersects $\alpha $ transversely, and with $\vert \varphi 
_{t}(\gamma )\vert $ bounded in terms of the ltd parameter functions. 

Now we choose  a continuous family of $C^{1}$ diffeomorphisms, 
$\varphi _{t}$, with $[\varphi 
_{t}]=y_{t}$, and 
satisfying the following properties. 
In all the following let $y_{s}$, $y_{t}\in \ell _{j}$.
\begin{description}
    \item[1.] $\varphi _{t}\circ \varphi _{s}^{-1}(S_{j,s})=S_{j,t}$ 
    and  $\varphi _{t}\circ \varphi _{s}^{-1}(\partial 
    S_{j,s})=\partial S_{j,t}$.
    \item[2.] $\varphi _{t}\circ \varphi _{s}^{-1}:S_{j,s}\to 
    S_{j,t}$ has bounded derivative with respect to the norm induced 
    by the metrics $\sigma _{s}$, $\sigma _{t}$  if $(\alpha 
_{j},\ell 
    _{j})$ is bounded,
   or if  $\vert s-t\vert \leq 1$ 
    and $(\alpha _{j},\ell _{j})$ is ltd.
    \item[3.] For any loop $\alpha $ which is $\alpha _{j}$ for at 
    least one $j$, $\partial S_{\alpha ,s}$, 
    $\varphi _{t}\circ \varphi _{s}^{-1}:\partial S_{\alpha ,s}\to 
\partial 
    S_{\alpha ,t}$ has constant derivative with respect to the length 
induced 
    by $\sigma _{s}$, $\sigma _{t}$.
    \item[4.] For each loop $\alpha $ which is $\alpha _{j}$ for at 
    least one $j$, for a fixed loop $\beta 
    (\alpha )$ intersecting $\alpha $ at least  once and at 
most twice, 
    (as in \ref{2.5}) and chosen intersection point $x$, let $\beta 
    _{t}$ be the closed geodesic on $S_{t}$ which is isotopic to  
    $\varphi _{t}(\beta (\alpha ))$, and let $x_{t}$ be the 
    image of $x$ under such an isotopy which also isotopes $\varphi 
    _{t}(\alpha )$ to a closed geodesic. Let $x_{1,t}$ and 
$x_{2,t}$ be 
    the endpoints in $\partial S_{\alpha ,t}$ of the  geodesic 
    arc of $\beta _{t}$ containing $x_{t}$. Then for all $s$, 
    $t\in I_{j}$, and $i=1$, $2$,
    $$\varphi _{t}\circ \varphi _{s}^{-1}(x_{i,s})=x_{i,t}.$$
    \end{description}
    
There are many such choices of $\varphi _{t}$, because the conditions 
on $\varphi _{t}$ are conditions on $\varphi _{t}\circ \varphi 
_{s}^{-1}$ restricted to $S_{j,s}$, for $s<t$, $s$, $t\in I_{j}$ 
for some $j$. But we cannot, as in the bounded geometry case, take 
$\varphi _{t}\circ \varphi _{s}^{-1}$ to minimise distortion for all 
$s<t$.

If $\alpha _{j}$ is a gap, let
$$W_{j}=\cup _{t\in I
_{j}}\varphi _{t}^{-1}(S_{j,t})\times \{ t\} $$
Thus, $W_{j}$ is a submanifold of $S\times [0,u]$ with piecewise 
smooth boundary which is 
homeomorphic to $\alpha _{j}\times I_{j}$. We also define
$$W_{j,t}=\varphi _{t}^{-1}(S_{j,t})\times \{ t\} .$$
If $\alpha $ is a loop 
which is 
$\alpha _{j}$ for at least one $j$, let
$$T(\alpha )=\cup _{t\in I(\alpha )} \varphi _{t}^{-1}S_{\alpha 
,t}\times \{ t\}
 $$
Then $T(\alpha )$ is solid torus,
called a {\em{model Margulis tube}}. No metric has been specified 
yet. 

 \ssubsection{Geometrically finite Kleinian surface case with two 
ends: metric on the 
complement 
of model Margulis tubes.}\label{6.3}

Let $\alpha _{j}$ be a gap. We define the model Riemannian metric 
$\sigma _{M}$ on 
$W_{j,t}$ by
$$\sigma _{M}=\varphi_{t}^{*}(\sigma _{t})+(c_{j}(z,t))^{2}dt^{2},$$
where, if we do not mind discontinuities in the metric on the union 
of all such $W_{j,t}$, we can take $c_{j}(z,t)=c_{j}$ to be constant, 
and
$$c_{j}=\begin{array}{ll} 1&{\rm{\ if\ }}(\alpha _{j},\ell 
_{j}){\rm{\ is\ ltd,}}\cr {1\over  \vert \ell _{j}\vert +1}& 
{\rm{\ if\ }}(\alpha _{j},\ell 
_{j}){\rm{\ is\ bounded.}}\end{array}$$
If we want the metric to be continuous or smooth, we can adjust 
$c_{j}(z,t)$ in a neighbourhood of $\partial W_{j}\cap \partial 
W_{k}$ 
whenever $\alpha _{j}$ and $\alpha _{k}$ are both gaps.
 As in the case of 
combinatorial bounded geometry, the precise definition of $\varphi 
_{t}$ 
does not change the metric up to Lipschitz equivalence, because the 
map 
$$(z,t)\mapsto (\varphi _{s}(z),t):(W_{j},\sigma _{M}) \to 
(S_{j,s}\times 
I_{j},\sigma _{s}+dt^{2})$$
is boundedly Lipschitz
if $(\alpha _{j},\ell _{j})$ is bounded, and boundedly Lipschitz 
restricted to the set of 
$(z,t)\in W_{j}$ with $\vert s-t\vert \leq 1$ if $(\alpha _{j},\ell 
_{j})$ is ltd.

\ssubsection{Margulis tubes.}\label{6.4}

 For any
three-dimensional
$\varepsilon _{0}$-Margulis 
tube, the metric induced by 
the hyperbolic metric on the tube boundary is Euclidean, and 
 Margulis tubes are also parametrised by $H^{2}$, at least if the 
core 
 loop is sufficiently short compared to $\varepsilon _{0}$.
 Indeed, if we choose coordinates so that in 
the upper half-space model of hyperblic space $H^{3}$, the Margulis 
tube has lift
$$\{ (re^{t+i\theta},e^{t}):0\leq r\leq R,\ \theta , t\in \mathbb 
R\} ,
$$
and the hyperbolic isometry corresponding to the core loop is the 
map $(z,u)\mapsto (e^{\lambda }z,e^{{\rm{Re}}(\lambda }u)$
with $\lambda =o(\varepsilon 
_{0})$. Then the metric on the boundary of the tube is 
$$R^{2}d\theta ^{2}+(R^{2}+1)dt^{2}.$$
and the corresponding element of the Teichm\"uller space of the 
torus, $H^{2}$, 
is 
$${2\pi i\over \lambda }(1+R^{-2})^{1/2},$$
and 
$$R={\varepsilon _{0}\over \vert \lambda \vert }(1+f(\lambda )),$$
for a $C^{1}$ function $f$ with $f(\lambda )=0$.
We write
$${2\pi i\over \lambda }=B+iA.$$
The length of the core loop is then ${\rm{Re}}(\lambda )$, and 
$${\rm{Re}}(\lambda )={2\pi A\over \vert B+iA\vert ^{2}},$$
and
$$R={\varepsilon _{0}\vert B+iA\vert \over 2\pi }(1+f(\lambda )).$$
See also \cite{Min0} 
Section 6. 

\ssubsection{Geometrically finite Kleinian surface case with two 
ends: 
the metric on model Margulis tubes.}\label{6.5}

It remains to define the metric $\sigma _{M}$ on the model Margulis 
tube $T(\alpha )$ if $\alpha =\alpha _{j}$ for at least one $j$.
We have 
$$\partial T(\alpha )=\cup _{k=1}^{4}\partial _{k}T(\alpha )$$
where $\partial _{1}T(\alpha )$ and $\partial _{2}T(\alpha )$ are the 
components of $\partial _{v}T(\alpha )$ (vertical boundary), where 
$$\partial _{v}T(\alpha )=\cup _{t\in I(\alpha )}\varphi 
_{t}^{-1}(\partial S_{\alpha ,t}),$$
and 
$$\partial _{3}T(\alpha )=\varphi _{s(\alpha )}^{-1}S_{\alpha 
,s(\alpha )},\ \ \partial _{4}T(\alpha )=\varphi 
_{t(\alpha )}^{-1}S_{\alpha ,t(\alpha )}.$$
We define the metric restricted to $\partial _{3}T(\alpha )$ to be 
$\varphi _{s(\alpha )}^{*}(\sigma _{s(\alpha )})$ and similarly for 
$\partial _{4}T(\alpha )$. On $\partial 
_{v}T(\alpha )$, the 
metric has already been defined from the 
definitions of the  
metric on the sets $W_{k}$ for gaps $\alpha _{k}$. We choose the 
numbering of $\partial _{1}T(\alpha )$ and $\partial _{2}T(\alpha )$ 
so that $x_{k,t}\in \partial _{k}T(\alpha )$ for $t\in I(\alpha )$, 
for $x_{k,t}$ as in 4 of \ref{6.3}. 
Both components $\partial _{k}S_{\alpha ,t}$ ($i=1$, $2$) of 
$\partial S_{\alpha ,t}$ (using the index $k$ in the same way) have 
the 
same length $\lambda _{\alpha }(t)$ in the $\sigma _{t}$ metric, by 
the definition of $S_{\alpha ,t}$, since $\lambda _{\alpha }(t)$ is a 
function of $\vert \varphi _{t}(\alpha )\vert $. Choose an 
orientation on $S_{\alpha ,t}$. For 
$x\in \partial _{k}S_{\alpha ,s}$ let $\lambda _{\alpha }(x,s)$ 
denote the 
length 
of the shortest 
positively oriented segment from $x_{k,s}$ to $x$, so that
$$g _{\alpha }(z)=\exp (2\pi i \lambda _{\alpha }(x,s)/\lambda 
_{\alpha }(s))$$
is smooth on $\partial _{1}S_{\alpha ,s}$, $\partial _{2}S_{\alpha 
,s}$. 
Then, by the properties of \ref{6.3}, in particular 3 and 4 of 
\ref{6.3}, the 
map 
$$(x,t)\mapsto (g _{\alpha }\circ\varphi _{s}(x),t):\partial 
_{v}T(\alpha )\to \{ e^{2\pi i\theta }:\theta \in \mathbb 
R \} \times I(\alpha )$$
pushes forward $\sigma _{M}$ on $\partial 
_{k}T(\alpha )$ ($k=1$, $2$) to
$$d\theta ^{2}+(c_{\alpha ,k}(t)^{2})dt ^{2}$$
on $S^{1}\times I(\alpha )$, for a suitable function $c_{\alpha 
,k}(t)$.
Even if we chose the $c_{k}(x,t)$ with some 
dependence on $x$ in order to smooth the metric, the choice can be 
made so that the  $c_{\alpha ,k}(t)$ depend only on $t$, because for 
each $t$,
$\partial S_{\alpha ,t}\subset\partial S_{m,t}$ for at least one gap 
$\alpha _{m}$, and intersects no other $\partial S_{p,t}$
 for a gap $\alpha _{p}$, unless $t$ is an endpoint of $I_{m}$, in 
which case 
 there is exactly one other gap $\alpha _{p}$ with 
 $\partial S_{\alpha ,t}\subset\partial S_{p,t}$.
We can then make a further change of variable
$$C_{\alpha ,k}(t)=\int _{s(\alpha )}^{t}c_{\alpha ,k}(u)du.$$
Then
$$(e^{i\theta },t)\mapsto (e^{i\theta },C_{\alpha ,k}(t)):S^{1}\times 
I(\alpha )\to S^{1}$$
pushes forward $d\theta ^{2}+(c_{\alpha ,k}(t)^{2})dt ^{2}$ to 
$d\theta 
^{2}+dt ^{2}$ on $S^{1}\times [0,A_{k}(\alpha )]$ where 
$A_{k}(\alpha )=C_{\alpha ,k}(t(\alpha ))$. So $\sigma _{M}$ on 
$\partial _{k}T(\alpha )$ 
is  
a Euclidean metric, up to isometry. 

The annuli $\partial _{3}T(\alpha )=S_{\alpha ,s(\alpha )}$ and 
$\partial _{4}T(\alpha )=S_{\alpha ,t(\alpha )}$ are 
conformally equivalent to Euclidean annuli $S^{1}\times 
[0,A_{3}(\alpha )]$ and $S^{1}\times [0,A_{4}(\alpha )]$. The 
uniformising maps restricted to boundaries have constant derivative 
with respect to length. If the modulus of $S_{\alpha ,s(\alpha )}$ 
is bounded above and below then the uniformising map transforms 
$\sigma 
_{s(\alpha )}$ to a metric which is boundedly Lipschitz equivalent to 
the Euclidean metric $d\theta 
^{2}+dt^{2}$ on $S^{1}\times [0,A_{3}(\alpha )]$. The modulus is 
bounded below in terms of $L$. It 
is also bounded above if $\vert \varphi _{s(\alpha )}(\alpha )\vert $ 
is bounded above, which is true if $s(\alpha )>0$ by the definition 
of vertically efficient, as already noted. Similar statements hold 
for $\partial _{4}T(\alpha )$.

So now if we match up boundaries, we have a conformal map 
$$\psi _{\alpha }:\partial T(\alpha )\to (S^{1}\times [ 0,A])/\sim ,$$
where 
$$A=A(\alpha )={1\over 2\pi }\sum _{k=1}^{4}A_{k}(\alpha ),$$
and where $\sim $ is defined by 
$$(0,e^{2\pi i \theta })\sim (e^{2\pi i(B+\theta )},2\pi A)$$
for some $B=B(\alpha )\in \mathbb R$.  We now determine $B(\alpha )$. 
The meridian is
$$\beta _{s(\alpha )}\cup \{ x_{1,t}:t\in I(\alpha )\} 
\cup \{ x_{2,t}:t\in I(\alpha )\} \cup \beta _{t(\alpha )},$$
for $x_{i,t}$ and $\beta _{t}$ as in 4 of \ref{6.2}. Thus, $\beta 
_{t}$ 
is an arc of the geodesic homotopic to $\varphi _{t}(\beta 
(\alpha ))$, for a fixed loop $\beta (\alpha )$ intersecting $\alpha 
$ 
at most twice (as in \ref{6.2}). Then $n_{\alpha }(y_{s(\alpha )})$ 
is, 
by definition, the value of $n$ which minimises 
$\vert \varphi _{s(\alpha )}\tau _{\alpha }^{n}(\beta 
(\alpha ))\vert $, where here (as earlier) $\tau _{\alpha }$ denotes 
oriented Dehn twist around $\alpha $. A similar statement holds for 
$t(\alpha )$. Then because $S_{\alpha ,t}$ has 
modulus bounded from $0$ in terms of the ltd parameter functions, the 
number of twists of $\beta _{s(\alpha )}$ round the geodesic 
homotopic to 
$\varphi _{s(\alpha )}(\alpha )$ is within $C$ of $-n_{\alpha 
}(y_{s(\alpha )})$, where $C$ 
depends only on the ltd parameter functions. A similar statement 
holds 
for $S_{\alpha ,t(\alpha )}$, $\beta _{t(\alpha )}$, $n_{\alpha 
}(y_{t(\alpha )})$. It follows that 

\begin{equation}\label{6.5.1}B={\rm{Re}}(-\pi _{\alpha }
    (y_{t(\alpha )})+\pi _{\alpha 
}(y_{s(\alpha )}))=-n_{\alpha }(y_{t(\alpha )})+n_{\alpha 
}(y_{s(\alpha 
)})
+O(1),\end{equation}
where the $O(1)$ term depends on the ltd parameter functions.

It is natural to take generators of the fundamental group of 
$\partial 
T(\alpha )$ as follows. Take  the homotopy class of $\partial 
_{k}T(\alpha )$ for any $k$, equivalently of the 
components of $\varphi _{t}^{-1}(S_{\alpha ,t})$ (any $t\in I(\alpha 
)$) as 
first generator and the meridian as second generator. As usual, 
identify the Teichmuller space of the torus with the upper 
half-plane $H^{2}$, using these generators. The corresponding point 
in $H^{2}$ is then $B+iA$.
Note that $A$ is bounded from $0$ in terms of the ltd parameter 
fuctions, since $A_{3}$ and $A_{4}$ are. 

Now suppose that
$${\rm{\ either\ }}s(\alpha )>0{\rm{\  or\ }} \vert \varphi 
_{0}(\alpha 
)\vert \geq \varepsilon _{0},$$
and 
$${\rm{\ either\ }}t(\alpha )<u{\rm{\  or\ }} \vert \varphi 
_{u}(\alpha 
_{j})\vert \geq \varepsilon _{0}.$$
The map $\psi _{\alpha }$ is 
a bounded Lipschitz equivalence from $\sigma _{M}\vert \partial 
T(\alpha )$ to the Euclidean metric $d\theta ^{2}+dt ^{2}$. 
Now we define the metric $\sigma _{M}$ on $T(\alpha)$ using 
a model Margulis tube $T(B+iA)$ such that the corresponding point in 
$H^{2}$ for $\partial T(B+iA)$ is exactly $B+iA$ whenever $B+iA$ is 
sufficiently large 
given $\varepsilon _{0}$, and within a bounded distance of it 
otherwise. 
This is possible, becase the map $f$ of \ref{6.4} is boundedly 
$C^{1}$.
Let $\sigma $ be the hyperbolic 
metric on $T(B+iA)$. We can assume that map $\psi _{\alpha }$ already 
defined on $\partial T(\alpha )$ has image $\partial T(B+iA)$, and 
under the current assumption it is boundedly biLipschitz on $\partial 
T(\alpha )$. Then extend diffeomorphically on the interior to
$$\psi _{\alpha }: T(\alpha )\to T(B+iA)$$
Then we define
\begin{equation}\label{6.5.2}
    \sigma _{M}=\psi _{\alpha }^{*}\sigma {\rm{\ on\ 
    Int}}(T(\alpha )).\end{equation}
If we wish to keep the metric continuous and smooth we can smooth it 
in the preimage of a small neighbourhood of $\partial T(\alpha )$ in 
$T(\alpha )$. 

Now suppose that 
$$t(\alpha )=u,\ \ \ \vert \varphi _{u}(\alpha )\vert <\varepsilon 
_{0},$$ and 
$${\rm{\ either\ }}s(\alpha )>0{\rm{\  or\ }} \vert \varphi 
_{0}(\alpha 
)\vert \geq \varepsilon _{0}.$$
Write 
$$\varepsilon _{2}=\vert \varphi _{u}(\alpha )\vert .$$
Then we identify $\partial 
_{1}T(\alpha )\cup \partial _{2}T(\alpha )\cup \partial _{3}T(\alpha 
)$ with an 
annulus in the tube  boundary of a Margulis tube $T(B+iA')$ 
corresponding to a 
point $B+iA'\in H^{2}$ for some 
$A'\geq A''=A_{1}+A_{2}+A_{3}$ and for $B$ as above. The length of 
the core loop will be 
$\varepsilon _{1}$. Then, as in \ref{6.4}, the relation between 
$\varepsilon 
_{1}$ and $B+iA'$ is
$$\varepsilon _{1}={2\pi A'\over \vert B+iA'\vert ^{2}}.$$
We want to determine $\varepsilon _{1}$  
so that the Poincar\'e metric on $S_{\alpha ,t(\alpha )}$
is the same, up to bounded distortion, as that of an annulus in the 
Margulis tube with boundary components on the boundary of the 
Margulis tube, and separated in the boundary of the Margulis tube an 
an annulus of area $A''$ on one side, and an annulus of area $A'-A''$ 
on the other. 
If
$$\varepsilon _{2}\leq {2\pi .2A''\over \vert B+2iA''\vert ^{2}},$$
then we define $\varepsilon _{1}=\varepsilon _{2}$. If
$${2\pi .2A''\over \vert B+2iA''\vert ^{2}}\leq \varepsilon _{2}\leq 
{2\pi \over \vert B+2iA''\vert }$$
then we define
$$A'=2A''.$$
The reason for this is that the radius of injectivity within diameter 
$1$ of a core loop of complex length $\lambda $ varies between $\vert 
\lambda \vert $ and ${\rm{Re}}(\lambda )$. Finally, if 
$${2\pi \over \vert B+2iA''\vert }\leq \varepsilon _{2}$$ 
then 
$A''<A'<2A''$, so that $(A'-A'')/A'<{1\over 2}$, and we define $A'$ 
inplicitly by the equation
$$\log ((A'-A'')/A')+\log \vert A'+iB\vert +\log \varepsilon _{2}=0.$$

Now we take an annulus $Y(B,A',A'')$ in $\partial T(B+iA')$ which is 
proportion 
$A''/A'$ of the total area and bounded by geodesics in the 
Eucidean metric in the homotopy class in $T(\alpha )$ of the core 
loop and 
of Euclidean length $\varepsilon _{0}(1+O(\varepsilon _{0}))$. The 
homotopy class of these geodesics in $\partial T(B+iA')$ is 
determined 
uniquely, at least if $\varepsilon _{1}=o(\varepsilon _{0})$, and in 
the other case it does not matter which homotopy class we choose with 
these properties. Now we define an annulus $C(B,A',A'')\subset 
T(B+iA')$ with 
the same boundary as $Y(B,A',A'')$. If $\varepsilon _{2}.\vert 
B+2iA''\vert >2\pi $, we take $C(B,A',A'')$ to be the union of the 
shortest 
geodesics in 
$T(B+iA')$ joining points in different components of $\partial 
Y(\alpha )$. 
These geodesics all have length $-2\log (\varepsilon _{2})+O(1)$ and 
their union is indeed an annulus. The annulus comes within a bounded 
distance of the core loop if $\varepsilon _{2}.\vert B+2iA''\vert  $ 
is 
bounded.  If $\varepsilon _{2}.\vert B+2iA''\vert  \leq 2\pi $, then 
we take 
$C(B,A',A'')$ to be the 
union of two annuli, formed by taking the union of 
geodesics in $T(B+iA')$ from each component of $\partial Y(B,A',A'')$ 
to the 
core loop of $T(B+iA')$, meeting both $\partial Y(B,A',A'')$ and the 
core loop perpendicularly. In this case, $C(B,A',A'')$ is not a 
smooth manifold if $A'>2A''$, 
but this is  what we want. In this case, the convex hull boundary  of 
the
corresponding hyperbolic manifold has a sharp angle at the 
geodesic $(\alpha )_{*}$, where $\alpha _{*}$ is the closed geodesic 
in the free homotopy class represented by $\alpha $.

 We also have a piecewise smoooth diffeomorphism
$$\psi _{\alpha }:\partial T(\alpha )\to C(B,A',A'')\cup 
Y(B,A',A''),$$
smooth on the interior of $\partial _{k}T(\alpha )$ for all $k$
which is boundedly biLipschitz, which  
maps  $\partial _{4}T(\alpha )$ to 
$C(B,A',A'')$ and $\partial _{1}T(\alpha )\cup \partial _{2}T(\alpha 
)\cup \partial _{3}T(\alpha )$ to $Y(B,A',A'')$. Let 
$T(B,A',A'')$ be the part of $T(B+iA')$ bounded by $C(B,A',A'')\cup 
Y(B,A',A'')$. when we extend to $\psi _{\alpha }$ to a diffeomorphism 
on the 
interior
$$\psi _{\alpha }:T(\alpha )\to T(B,A',A'')$$
and then we again use (\ref{6.5.2}) to define $\sigma _{M}$ on 
$T(\alpha )$.

The case when 
$$s(\alpha )=0,\ \ \ \vert \varphi _{0}(\alpha )\vert <\varepsilon 
_{0},$$ 
and 
$${\rm{\ either\ }}t(\alpha )<u{\rm{\  or\ }} \vert \varphi 
_{u}(\alpha 
)\vert \geq \varepsilon _{0}$$
is treated exactly similarly, with $s(\alpha )$ and $t(\alpha )$ 
interchanged, 
and also $0$ and $u$, $A_{3}$ and $A_{4}$, and 
$\partial _{3}T(\alpha )$ and $\partial _{4}T(\alpha )$. 
Finally we consider the case when
$$s(\alpha )=0,\ \ \ \vert \varphi _{0}(\alpha )\vert <\varepsilon 
_{0},$$
$$t(\alpha )=u,\ \ \ \vert \varphi _{u}(\alpha )\vert <\varepsilon 
_{0}.$$
In this case we define
$$\varepsilon _{2}=\vert \varphi _{T}(\alpha )\vert ,$$
$$\varepsilon _{3}=\vert \varphi _{0}(\alpha )\vert ,$$
and suppose without loss of generality that $\varepsilon _{2}\leq 
\varepsilon _{3}$. Then we define
$$A''=A_{1}+A_{2}+{1\over \varepsilon _{3}}.$$
 Then we define $\varepsilon _{1}$ and $A'$ exactly as before.

The set $Y(B,A',A'')$ now has two components whose areas are $A_{1}$ 
and 
$A_{2}$ chosen Euclidean distance $1/\varepsilon _{3}$ apart in 
$\partial T(B+iA')$ with respect to the Euclidean metric on $\partial 
T(B+iA')$. The 
set $C(B,A',A'')$ also has two components, one of which is a union of 
geodesic segments in $T(B+iA')$ joining the components of $\partial 
Y(B,A',A'')$ 
which are 
Euclidean distance $1/\varepsilon _{3}$ apart, in which case these 
geodesics have length $2\log (\varepsilon _{0}/\varepsilon 
_{3})+O(1)$. 
The other component is defined as before, depending on whether 
$\varepsilon _{2}.\vert 
B+2iA''\vert >2\pi $ or $\varepsilon _{2}.\vert 
B+2iA''\vert \leq 2\pi $ Then, as before, we take $T(B,A',A'')$ to be 
the subset of 
$T(B+iA')$ 
bounded by $Y(B,A',A'')$ and $C(B,A',A'')$ and pull back the metric 
using the formula of (\ref{6.5.2}).

The metric $\sigma _{M}$ has now been completely determined, up to 
bounded distortion, in the case of two geometrically finite ends.
If the model Margulis tube is homotopic to $\{ t\} 
\times \gamma $ then we call the core loop $\gamma _{**}$ and shall 
sometimes call the 
solid torus $T(\gamma _{**})$ rather than $T(\gamma )$ as above.

\ssubsection{Geometrically finite Kleinian surface case: more than 
two ends.}\label{6.9}

We now consider the case of a model relative Scott core $M_{c}$  
being 
homeomorphic 
to $S_{d}\times [0,1]$, where $S_{d}$ is a compact surface with 
boundary,
but where the boundary component homeomorphic to $S_{d}\times \{k\} $ 
might be a union of several components of $\partial_{d} M_{c}$ and of 
annuli $\partial 
M_{c}\cap \partial M_{d}$, for $k=0$, $1$. This is an example of 
something which happens regularly. It is sometimes convenient to 
group ends of $M_{d}$ together, by taking a connected 
submanifold $W$ of 
$M_{c}$, and considering the group of ends in each component of 
$M_{d}\setminus W$. We only do this for $W$ with the following 
properties: each component of $\partial W\cap \partial M_{d}$ is an 
essential annulus in $\partial M_{d}$, a component of $\partial 
W\setminus \partial M_{d}$ cannot be homotoped into $\partial M_{d}$, 
and 
if $S_{d}$ is a component of $\partial W\setminus \partial M_{d}$ 
bounding the closure 
$U$ of a nonempty component of $M_{c}\setminus W$, then $(U,S_{d})$ 
is 
homeomorphic to $(S_{d}\times [0,1],S_{d})$.  In such cases, as in 
this case, 
we define 
$M_{d,W}$ to be the union of $M_{d}$ and components of $M\setminus 
M_{d}$ which are disjoint from $W$. Then if $U'$ is the closure of 
the 
component of 
$M_{d,W}\setminus W$ bounded by $S_{d}$, $(U',S_{d},\partial 
M_{d,W}\cap U')$ is homeomorphic to $(S_{d}\times [0,\infty 
),S_{d},\partial S_{d}\times [0,\infty ))$. So $U'$ conntains a 
unique 
end of $M_{d,W}$.

In the current case, $W$ is homeomorphic to $S_{d}\times [0,1]$, and 
$M_{d,W}$ has two ends, which we call $e_{\pm }$.
Then $S_{d}$ is the horodisc deletion of a finite type surface 
$S$, and we have two isotopically disjoint  multicurves $\Gamma 
(e_{\pm })$, one of which could be empty,  and the end invariants 
of $M_{d}$ 
split into two sets 
$(y_{1},\cdots y_{m})$ and $(y_{m+1},\cdots y_{p})$, which are the 
end 
invariants $\mu (e_{\pm })$ of $e_{-}$, $e_{+}$ respectively. Here, 
$y_{i}\in {\cal{T}}(S(\alpha _{i}))$ and $\alpha _{i}$ ($1\leq i\leq 
m$) are the gaps of $\Gamma (e_{-})$, while $\alpha _{i}$ ($m+1\leq 
i\leq 
p$) are the gaps of $\Gamma (e_{+})$. In this case we make the model 
manifold $M(\mu (e_{-}),\mu (e_{+}))$
by slightly modifying the model $M(y(e_{-}),y(e_{+}))$ of \ref{6.2} 
to 
\ref{6.5}, (replacing $y_{0}$ and $y_{u}$ by $y(e_{-})$ and 
$y(e_{+})$)
where 
$y(e_{-})\in {\cal{T}}(S)$ is defined using the $y_{i}$ for $1\leq 
i\leq 
m$ and $y(e_{+})$ is defined using the $y_{i}$ for $m+1\leq i\leq p$.

Fix a Margulis constant $\varepsilon _{0}$. We 
define $y(e_{-})$ up to bounded distance by defining 
projections $\pi _{\alpha _{i}}(y(e_{-}))$ for $1\leq i\leq m$ and 
$\pi 
_{\gamma }(y(e_{-}))$ for $\gamma \in \Gamma (e_{-})$. So we take
$$\pi _{\alpha _{i}}(y(e_{-}))=y_{i},\ 1\leq i\leq m,$$
$$ \pi _{\gamma }(y(e_{-}))={1\over \varepsilon _{0}},\ \gamma \in 
\Gamma _{0}.$$
The definition of the imaginary part of $\pi _{\gamma }$ needs a 
choice of loop transverse 
to $\gamma $ intersecting it at most twice (see \ref{2.5}). But in 
fact this choice is irrelevant, because the imaginary part only 
influences the geometry of the Margulis tube, which we are about to 
remove. If 
$y(e_{-})=[\varphi _{0}]$ then we have chosen to make $\vert \varphi 
_{0}(\gamma )\vert =\varepsilon _{0}(1+o(1))$. We define $y(e_{+})$ 
similarly, using $\Gamma (e_{+})$ and $y_{i}$ for $m<i\leq p$.

The model manifold $M(\mu (e_{-}),\mu (e_{+}))$ is then 
obtained from the model \\ $M(y(e_{-}),y(e_{+}))$ constructed in 
\ref{6.2} to 
\ref{6.5} by simply leaving out the models for parts of Margulis 
tubes 
$T(\gamma _{**})$ for $\gamma \in \Gamma (e_{-})\cup \Gamma (e_{+})$, 
and 
gluing in, instead, pieces of horoballs with the right 
boundary piece. We only do this when $\gamma \in \Gamma (e)$ for just 
one of $e=e_{+}$, $e_{-}$. So we have a single annulus $A$, with a 
Euclidean metric. We simply glue in $A\times [0,\infty )$ with the 
product Euclidean metric.

\ssubsection{Independence of ltd parameter functions.}\label{6.6}

The model manifold 
is independent, up to Lipschitz equivalence with bound depending on 
the ltd parameter functions, of the decomposition of 
\ref{5.7}, and of the long thick and dominant parameter functions 
used. To see this, we only need to compare the models given by 
vertically efficient partitions for two sets 
of ltd parameter functions  $(\Delta _{1},r_{1},s_{1},K_{1})$, 
$(\Delta _{2},r_{2},s_{2},K_{2})$ where $(\Delta 
_{2},r_{2},s_{2},K_{2})$
is sufficiently stronger than $(\Delta _{1},r_{1},s_{1},K_{1})$ for 
the following to hold. For all $\nu $,  
$r_{2}(\nu )\leq \nu _{1})$, where every ltd gap  in the 
decomposition 
for 
$(\Delta _{1},r_{1},s_{1},K_{1})$ is long $\nu $-thick and dominant 
for some $\nu \geq \nu _{1}$. 
We can choose this $\nu _{1}$ so that , if $(\alpha ,\ell )$ is 
bounded for 
$(\Delta _{1},r_{1},s_{1},K_{1})$, then $\vert \varphi (\gamma )\vert 
\geq 
\nu _{1} $ for all $\gamma $ in the interior of 
$\alpha $. Such a $\nu _{1}$ exists by \ref{5.4}, 
and the definition 
of 
bounded in \ref{5.5} and the results discussed there.  So if $(\alpha 
_{i},\ell _{i})$ is in the partition for $(\Delta 
_{i},r_{i},s_{i},K_{i})$ and $\ell _{1}\cap \ell _{2}=\emptyset $, we 
have $\alpha _{1}\subset\alpha _{2}$ or $\alpha _{1}\cap \alpha 
_{2}=\emptyset $. If $\alpha _{1}\subset \alpha _{2}$, then the 
vertically efficient condition for the second partition  implies that 
$\ell _{1}\subset \ell _{2}$. So then the 
partition for $(\Delta _{1},r_{1},s_{1},K_{1})$ is a refinement of 
the partition for $(\Delta _{2},r_{2},s_{2},K_{2})$. We can assume 
that the same family $\varphi _{t}$ is used to define both models,  
because, as already noted in \ref{6.3}, the choice of family, subject 
to the 
conditions of \ref{6.2}, does not affect the model. 

Let $(\alpha 
_{1},\ell _{1})$ be  ltd for  $(\Delta _{1},r_{1},s_{1},K_{1})$ with 
$(\alpha _{1}\times \ell _{1})$ in the first partition. 
Let $\alpha _{2}\times \ell 
_{2}$ be 
 from the second partition, with $\alpha _{1}\times \ell _{1}\subset 
 \alpha _{2}\times \ell _{2}$. If 
  $(\alpha _{2},\ell _{2})$ is bounded for $(\Delta 
_{2},r_{2},s_{2},K_{2})$ then we have a bound on $\vert \ell 
_{1}\vert $ 
in terms of $(\Delta _{2},r_{2},s_{2},K_{2})$ and also a bound on the 
number of such $(\alpha _{1},\ell _{1})$ in terms of 
$(\Delta _{2},r_{2},s_{2},K_{2})$. Write $W_{1}$ $W_{2}$ for the 
corresponding subsets of $S\times [0,u]$, as defined in \ref{6.2},
so that $W_{1}\subset W_{2}$. Then the scaling factor $c_{1}(t)$ is 
bounded, since $\vert \ell _{1}\vert $ is bounded. If $\alpha _{1}$ 
is 
a loop, the corresponding Margulis tube has bounded geometry, because 
the quantity $B_{1}$ of \ref{6.4} is bounded in terms of 
${\rm{Re}}(\pi _{\alpha 
_{j}}(y_{t_{j}})-\pi _{\alpha _{j}}(y_{s_{j}}))$, if $\ell 
_{j}=[y_{s_{j}},y_{t_{j}}]$, by (\ref{6.5.1}). So, in both models, 
the 
model metric on $W_{2}$ is of bounded geometry, with bound depending 
on the ltd parameter functions. Now suppose  the same situation,  
but with $(\alpha _{2},\ell _{2})$ ltd with respect to $(\Delta 
_{2},r_{2},s_{2},K_{2})$. If $\alpha _{1}$ is a gap then the metric 
on $W_{1}$ is defined in exactly the same way in both models. If 
$\alpha _{1}$ is a loop properly contained in $\alpha _{2}$, then 
$\vert \varphi (\alpha ) \vert \geq \nu _{2}$ for all $[\varphi ]\in 
\ell _{1}$, for $\nu _{2}$ such that for every $(\alpha ,\ell )$ for 
in the second partition with $\alpha $ a gap, $\alpha $ is long $\nu 
$-thick and dominant for some $\nu \geq \nu _{2}$. Then $W_{1}$ has 
bounded geometry in both models. So the geometry of all pieces is 
boundedly equivalent in both models. 

\ssubsection{Consistency with previous models and 
examples.}\label{6.7}

The model metric on Margulis tubes in \ref{6.5} is effectively 
exactly the same as the model metric on Margulis tubes in Minsky's 
punctured torus paper \cite{Min0}. The situation here is slightly 
more general, because the quantity $A=A(\alpha )$ of \ref{6.5} can be 
arbitrarily large, although it is always bounded from $0$. In the 
punctured torus case it is bounded above also, because in the 
punctured torus case, if $(\alpha ,\ell )$ is ltd, then either 
$\alpha $ is the punctured torus itself, or $\alpha $ is a loop.

This model is also consistent with Rafi's examples \cite{Raf}, 
\cite{Raf2}.
 In one of 
Rafi's examples, adapting to the present notation, $S$ is a closed 
surface of genus two, and $[y_{0},y_{u}]$ is the union of two 
segments $\ell _{1}$, $\ell _{2}$ along which loops $\gamma _{1}$ and 
$\gamma _{2}$ 
respectively are $K_{0}$-flat, apart from  a bounded segment in the 
middle. The convex hull of $\gamma _{1}$ and $\gamma _{2}$ could be 
(for 
example) $S\setminus \gamma _{3}$ where $\gamma _{3}$ is not 
$K_{0}$-flat along any segment and not adjacent to $\alpha $ for any 
ltd $(\alpha ,\ell )$. In the 
corresponding hyperbolic  manifold with  ending data $(y_{-},y_{+})$ 
(or even close to this) the geodesic loop $(\gamma _{3})_{*}$ is not 
short. But for $[\psi ]\in \ell _{1}\cup \ell _{2}$ bounded from 
$y_{\pm }$, $\vert \psi (\gamma _{3})\vert $ is small, and $\to 0$ as 
the length of both $\ell _{1}$ and $\ell _{2}$ $\to \infty $. Rafi 
proves this, and there is also a proof of something similar in 
15.21 of \cite{R1}. 
However, in the 
model manifold $M(y_{-},y_{+})$, the geodesic $(\gamma _{3})_{**}$ is 
not short. A model Margulis tube is inserted, but it has bounded 
geometry.

\ssubsection{Model manifolds 
for compression bodies.}\label{6.10}

We now have a geometric model for any hyperbolic $3$-manifold $N$
for which $(N_{d},\partial N_{d})$ is homeomorphic to $(S_{d}\times 
{\mathbb R},\partial S_{d}\times \mathbb R)$ for the horodisc 
deletion 
of a finite type surface $S$. The  model in the general case is 
obtained by gluing together interval bundle models. The data for a 
model is given by topological type and ending lamination data. The 
ending lamination data  for an end $e$ in ${\cal{T}}(S(e))\cup 
{\cal{O}}_{a}(S(e))$ is given relative to an identification of the 
corresponding boundary component of the relative Scott core with 
$S(e)$. 
The definition of  ${\cal{O}}_{a}(S(e))$ is given relative to a base 
metric on $S(e)$.
We are now going to discuss this identification and base metric 
briefly, 
in the case 
of relative compression bodies.

There are a number of definitions of compression body in the 
literature, all of them equivalent. See for example \cite{Ot1}, 
\cite{Can}, \cite{Bon2}. The relative version which is used here is 
adapted from the definition of \cite{Sou}. A relative compression 
body is, in 
fact, a compression body in the usual sense, but the decomposition is 
slightly different.

A {\em{relative compression body}} $(W,\partial W\cap A)$ is a 
compact manifold with boundary $\partial W$, with $\partial 
W\setminus A$ a finite disjoint union of surfaces 
wih boundary such that all components of  $A$ are closed 
annuli, and such that there is a subsurface $S_{d,0}$ of $\partial 
W\setminus A$ with the 
following properties. All  components of  $\partial W\setminus (A\cup 
S_{d,0})$ are 
incompressible, 
and the inclusion $S_{d,0}\to W$ is surjective on $\pi _{1}$.  
The surface $S_{d,0}$ is called the {\em{exterior surface}}.  We 
only consider manifolds which have a chance of being homotopy 
equivalent to hyperbolic manifolds, so we assume also that $W$ is 
$K(\pi ,1)$. 

Now we assume that $S_{d,0}$ is compressible in 
$W$, since 
otherwise we are in the interval bundle case. Let $S_{0}$ be the 
surface of which $S_{d,0}$ is the horodisc deletion.
We can then find a multicurve $\Gamma _{0}$ on $S_{0}$ 
such that the loops of $\Gamma _{0}$ bound disjoint embedded discs in 
$W$, and such that the closure of each component of the 
complement of the union of the discs and $\partial M_{c}$ is either a 
three-ball or an interval bundle homeomorphic to $S_{d,i}\times 
[0,1]$ 
for some component $S_{d,i}$ of $\partial W\setminus A$, with 
the intersection of $A$ with the interval bundle identifying with 
$\partial 
S_{d,i}\times [0,1]$. It is natural to choose our base 
metric on $S_{0}$ so that, for a fixed choice of $\Gamma _{0}$, the 
loops of $\Gamma _{0}$ are of bounded length, and of length 
bounded from $0$. We shall always do this. We therefore fix 
$z_{0,0}=[\varphi _{0}]\in {\cal {T}}(S_{0})$ such that $\vert 
\varphi 
_{0}(\gamma )\vert $ is bounded and bounded from $0$ for $\gamma \in 
\Gamma _{0}$.
The gaps of $\Gamma _{0}$ are 
$\alpha _{i}$, $1\leq i\leq p$, where $p\geq m$. In $M$ we can attach 
a disjoint disc to each loop of $\Gamma _{0}$. This embeds $\alpha 
_{i}$ in a surface $S_{d,i}$, which is a union of $\alpha _{i}$ and a 
number of discs. If $1\leq i\leq m$ then $S_{d,i}$ is homotopic in 
$W$ to an incompressible component of $\partial W\setminus 
A$, and $S_{d,i}$ is the horodisc deletion of a surface 
$S_{i}$.  
For $m<i\leq p$, $\alpha _{i}$ and the attached discs bound a ball in 
$W$. For $1\leq i\leq m$ 
we define $z_{0,i}\in {\cal{T}}(S_{i})$ to have the same conformal 
structure as  $\pi _{\alpha _{i}}(z_{0,0})$, but forgetting 
the 
puncture in each disc. 

Now suppose that $M$ is to be the model for a hyperbolic 
$3$-manifold, with 
$M_{d}$ as the model for the horoball deletion, and 
that $W\subset M$ is a submanifold of the model for the relative 
Scott core, 
as suggested at the start of \ref{6.9}. Let $M_{d,W}$ be as defined 
in \ref{6.9}. The ends of $M_{d,W}$ are then sets of ends of $M_{d}$. 
We fix geometrically finite ending invariants for the ends of 
$M_{d}$. This gives corresponding ending invariants $\mu (e_{i})$  
for 
$M_{d,W}$. The model manifold for this end, with these invariants, is 
obtained by modifying a model $M(z_{0,i},y_{i})$ to a model 
$M(z_{0,i},\mu (e_{i}))$ for a certain point 
$y_{i}\in {\cal{T}}(S_{i})$, as described in \ref{6.2} to \ref{6.9}. 
Then the 
geometric model 
$M(\mu (e_{0}),\cdots \mu (e_{m}))$ will be made from the models
$$M_{k}=M(z_{0,k},\mu (e_{k})),$$
$0\leq k\leq m$.  It seems 
simpler to construct 
the geometric manifold for the unquotiented 
$(z_{0,k},y_{k})$, and to note 
that geometric models for different 
choices  in a single ${\rm{Mod}}_{0}(\partial M_{k},M_{k})$ orbit are 
homeomorphic. 

It only remains to show how to glue together the different models 
$M_{k}$ together. 
We can assume that these submanifolds are all disjoint in $M$, by 
moving them apart slightly, and we 
have to extend the metric to the complement. Let $D(\gamma )$ be the 
disc in $M_{c}$ attached to  the loop in $\partial M_{0}$ 
corresponding to  
$\varphi _{0}(\gamma )$, for $\gamma \in \Gamma _{0}$. We also 
assume 
that these discs are all disjoint, and have interiors disjoint from 
all $M_{i}$. We define the metric on 
each of these discs to be the pullback of the Euclidean metric on the 
unit disc, where the map on the boundary is length-preserving with 
respect to the metric on  $ M_{0}$. Then we label the components of 
the 
complement in $M$ of the $M_{k}$ ($0\leq k\leq m$) by $M_{k}'$ for 
$1\leq k\leq p$, where $M_{k}'$ is between $M_{k}$ and $M_{0}$ if 
$k\leq m$, and $M_{k}'$ is a ball if $k>m$.  So now we have a metric 
on the boundary of each 
complementary component $M_{k}'$, $1\leq k\leq p$, and the closure of 
each $M_{k}'$ is homeomorphic to $S_{k}\times [0,1]$ 
if $1\leq j\leq m$.  In the case of 
$k\leq m$, the metric on the boundary component which is also in 
$\partial M_{k}$ is a hyperbolic 
metric by construction, corresponding to $y_{k}\in 
{\cal{T}}(S_{k})$. 
For the other boundary 
component, the metric determines a conformal structure, and hence an 
element of ${\cal{T}}(S_{k})$ which we claim is a bounded distance 
from $z_{0,k}$. To see this, we consider the definition of $\pi 
_{\alpha 
_{k}}$. The 
corresponding hyperbolic surface is a punctured surface homeomorphic 
to $S(\alpha _{k})$. But the definition in \ref{2.5} initially 
describes 
this hyperbolic surface in terms of its conformal structure, 
attaching 
a punctured disc to the subset of $\varphi _{0}(S_{0})$ homotopic 
to 
$\varphi _{0}(\alpha _{k})$ and bounded by the geodesics homotopic 
to 
$\varphi _{0}(\partial \alpha _{k})$. So we might just as well 
regard this as a 
conformal structure on the surface, with the punctures corresponding 
to $\partial \alpha _{k}$ removed. We also have a metric, taking the 
metric on the component of  $\varphi _{0}(S_{0}\setminus \Gamma 
_{0})$ 
homotopic to $\varphi _{0}(\alpha _{k})$ and bounded by the 
geodesics homotopic to $\varphi _{0}(\partial \alpha _{k})$
and pushing forward 
the usual Euclidean  unit disc metric on the added discs. This 
metric is boundedly Lipschitz equivalent, under a piecewise smooth 
diffeomorphism $\varphi 
_{k}$, to the unique hyperbolic metric 
which it is conformally equivalent to, which is the metric 
corresponding to $z_{0,k}$. So we can take a metric on $M_{k}'$ which 
is boundedly equivalent to the pullback of product metric $\sigma 
_{k}+dt^{2}$ on $S_{k}\times [0,1]$, extending the metric already 
defined on the boundary. For $k>m$, choose the pullback under a 
homeomorphism $\varphi 
_{k}'$ of a metric on the round unit ball which is boundedly 
equivalent 
to the Euclidean metric, extending the metric already define on the 
boundary. The metric has some discontinuities since the sets 
$\partial M_{k}'$ 
are only piecewise smooth, but we can make the metric smooth and 
continuous by small local perturbations, without changing the metric 
up to coarse Lipschitz equivalence. 

\ssubsection{Model for the Scott core.}\label{6.11}

 Choose $W\subset 
M_{c}$  as in \ref{6.9}. We can choose $W=M_{c}$, but other choices 
might be more useful. An model for $W$ has already 
been given in the case when $W$ is a relative compression body 
\ref{6.10}. As we shall see, the general case is similar, with the 
model being made by gluing together models for interval bundles. 
First we need to review standard methods for decomposing $W$ into 
simpler pieces. We are assuming that $M_{c}$ has nonempty 
nonspherical boundary, 
and so the same is true of $W$. Therefore it is Haken, and can be 
decomposed into simpler pieces. A decomposition is given in \cite{Hem}
in Chapter 13, with most of the basic work done in Chapter 6. Here is 
an adaptation which gives some extra properties relating to the 
annuli 
and tori in $\partial W\cap \partial M_{d}$. For the moment, we write 
$A$ for the disjoint union of annuli and tori in $\partial W\cap 
\partial M_{d}$. 

\begin{ulemma}
 Let $W$ be a compact connected oriented aspherical 
$3$-manifold 
 in which any torus is parallel to the boundary, and with nonempty 
 boundary. 
    Let $A$ be 
a disjoint union of essential annuli in $\partial W $, which 
are all homotopically distinct.
Then $W$ is obtained 
from a finite sequence $W_{i}$, $0\leq i\leq n$ such that $W_{n}=W$, 
$W_{0}$ is a disjoint union of balls and interval bundles. Any 
interval bundle component $V$ is homeomorphic to $S\times [0,1]$ in 
such a way  that $\partial S\times [0,1]$  is 
homeomorphic to a subset of $A\cap \partial V$ and $S\times \{ 0\} $ 
identifies with an incompressible 
subsurface of $W$. There is a 
surjective map $j_{i}:W_{i}\to W_{i+1}$ which is injective restricted 
to the interior of $W_{i}$, at most two-to one, and pairs up disjoint 
subsurfaces of $\partial W_{i}$ to form $W_{i+1}$. Each component of 
$\partial W_{i}$ contains at least one maximal connected subsurface 
which is mapped by $j_{i}$ to an incompressible surface $S_{1}$ in 
a component $V$ of $W_{i+1}$.  

The following properties will hold for $S_{1}$. From now on, by abuse of notation, we identify 
components of $W_{i}$ with their images in $W_{j}$, for any $j>i$.
\begin{description}
    \item[1.] $\partial S_{1}\neq \emptyset $.

\item[2.] $S_{1}$ is incompressible in $V$: a simple 
loop which is homotopically nontrivial in $S_{1}$ does not bound a 
disc in 
$V$.
    \item[3.] $S_{1}$ is boundary incompressible in the following 
sense.
  There is no union of a homotopically nontrivial arc in $S_{1}$ 
between points of $\partial S_{1}$, and an arc in $\partial 
V\setminus A$, which bounds a disc in $V$.

\item[4.] If $S_{1}$ is not a disc, and $S_{2}$ is any connected 
    union of components of $\partial V\setminus A$ and $A$ 
    which are all intersected by $\partial S_{1}$, then $S_{2}$ is 
    incompressible in $V$.

    \item[5.] If $S_{1}$ is a disc, and $i+1<n$, then $V$ is obtained 
from 
    some  component $V'$ of $W_{i+2}$ by attaching incompressible 
surfaces to $\partial V'$, 
    none of which is a disc.

 \item[6.] If $S_{1}$ is a disc, and $i+1<n$, then $S_{1}$ is 
indecomposable in the 
following sense. Let $S_{2}$ be any component of $\partial (\partial 
V\cap 
\partial V')$, where $V'$ is the component of $W_{i+2}$ containing 
$V$. Isotope $\gamma =\partial S_{1}$ to have only essential 
intersections with $\partial S_{2}$. Then up to homotopy preserving 
intersections with $\partial S_{2}$,
$\gamma \neq \gamma _{1}*\zeta *\gamma _{2}*\overline{\zeta }$, 
where $\gamma _{1}$ and $\gamma _{2}$ are both nontrivial in 
$\partial 
V$ but trivial in $V$
\end{description}

\end{ulemma}

\noindent {\em{Proof.}} 
 The distinguished set of annuli in 
$\partial W$ plays a role.  Rules 1 to 
3 are straightforward, and probably standard. The other 
rules are formulated in such 
a way that we can bound geometry in section \ref{8}, and rule 3 will 
also be used in section \ref{8}. Rule 5 says 
we should use a compressing disc whenever it is possible to do so, in 
keeping with Rule 4. Rule 5 says we should put in all such discs at 
once. Rule 6 is a decomposability condition: as indeed is the more 
standard Rule 3.

Suppose that $W_{j}$ has been constructed with the required 
properties for $j\geq i+1$ (by abuse of notation we are assuming $n$ 
has already been determined) and let $V$ be a component of $W_{i+1}$.

We start by taking a maximal set of compressing discs attached to 
each 
component of $\partial V$, subject to Rule 4. We can do this 
satisfying 
Rule 6 if 
$i+1<n$, since each disc is a sum of indecomposable ones. So now we 
assume that
there are no compressing discs which can be attached to $\partial V$

  Then we apply 6.8 and 6.5 
 of \cite{Hem}. Rules 2,4, and 5 is automatically satisfied. We can decompose 
 to satisfy Rule 3 if necessary.

\ssubsection{Interval bundles in the Scott core.}\label{6.14}
In the case we are interested in, when $W$ is homotopy equivalent to 
a 
hyperbolic $3$-manifold $N$, $W$ contains no essential tori except 
possibly parallel to the boundary, if $A$ has toroidal 
components. The 
decomposition then gives more information, as the following lemma 
shows.

\begin{ulemma} Continue with the notation of \ref{6.11}. 
 Suppose that $W$ has no essential tori, except possibly parallel to 
 components of $\partial W$ which are in $\partial M_{d}$. Let 
 $\Gamma _{0}$ be a union of multicurves from the boundaries of the 
 compressing discs of sets $W_{i}\setminus W$ on {\em{compressible }}
 components of $\partial W\setminus 
 A$. Let 
 $\Gamma _{0}'$ be the union of a maximal multicurves on these 
 boundary components  which are noncollapsing, 
 on these compressible boundary components, such that  
 $\Gamma _{0}\cup \Gamma _{0}'$ cuts each such component of $\partial 
 W\setminus A$ into discs with at most one puncture and 
 annuli parallel to the boundary.  Let 
$$\Sigma =\Gamma _{0}'\cup A\cup \cup _{1\leq j\leq n}\partial 
(\partial W_{j}\cap 
\partial 
W_{j-1}).$$ 
Define $\Sigma _{0}=\Sigma $. Inductively, define 
$\Sigma _{i+1}$ to be the union, over all  interval bundle components 
$V$ of 
$W_{0}$, and homeomorphism to $S\times [0,1]$, of annuli and 
rectangles with 
alternate sides  in $A$ which are 
  images under the homeomorphism of arcs
 and loops $\alpha \times [0,1]$ where the image of $\alpha \times 
 \{ j\} $ is already such an arc or loop in $\Sigma _{i}$, for some 
 $j=0$, $1$, with $\partial \Sigma _{i}$ replaced by $\Sigma $ 
 if $i=0$. Then for some $r$, each component $V$ of 
 $W_{0}\setminus 
\Sigma _{r}$ is either a ball, or homeomorphic to an interval bundle
$S_{d}\times [0,1]$, where the sets homeomorphic  $S_{d}\times 
\{ 0,1\} $ are in incompressible components of  $\partial W\setminus 
A$. In the case when $V$ is a ball
 $\partial 
V\setminus \Sigma _{r}$ is a union of topological discs and annuli 
parallel to the boundary.\end{ulemma}

In future, we shall write $\Sigma '$ for the set $\Sigma _{r}$.

\noindent{\em{Proof.}} 
Inductively, the component  of $W_{0}\setminus \Sigma _{i}$ 
have the same properties as those listed for components of  $W_{0}$. 
Thus, each component $V$ is either a ball or homemorphic to $S\times 
[0,1]$ under a 
homeomorphism sending $\partial V\cap A$ to a subset of $\partial 
S\times [0,1]$. There is a bounded $r$ for which nontrivial interval 
bundle components of $W\setminus \Sigma _{r}$ do not decompose. 
Write $\Sigma _{r}=\Sigma '$. The 
only way for a component  $V$ of $W\setminus \Sigma '$ to be such 
an interval bundle is if it is part of a stack of such interval 
bundles,  such that each such interval bundle $(V,\partial V\cap A)$ 
is homeomorphic 
to $(S\times [0,1],\partial S\times [0,1])$ for some surface $S$. The 
surface $S$ must be the same throughout the stack, and the boundary 
component homeomorphic to $S\times \{ 1\} $ on one matches up with 
the boundary component homeomorphic to $S\times \{ 0\} $ on the next. 
Since there are only finitely many components, either the stack 
closes 
up to give a $S^{1}$-bundle over $S$, or there are two end bundles 
with 
boundary in $\partial W\setminus A$. The first possibility can only 
occur if $S$ has boundary which is not in $A$, 
since $W$ is connected and has boundary which is not in $A$. But then 
we get a contradiction to $W$ having no essential tori, except 
parallel to the boundary. So only the second possibility occurs. In 
this case there is just one manifold on the stack and the components 
homeomorphic to $S\times [0,1]$ are in $\partial W\setminus A$. 

The statement about $\partial \Sigma '$ cutting $\partial V$ into 
cells follows, because if not, there is a boundary of a compressing 
disc in the boundary of some $W_{i}$ containing $V$, which is 
disjoint from $\Sigma _{r}$. This is impossible. For all $i$, the 
construction is such that  all compressing discs automatically 
intersect at least one loop in  $\Sigma $ transversally.   

\Box

\ssubsection{The model for the non-interval-bundle part of the 
$W$.}\label{6.15}

Lemmas \ref{6.11} and \ref{6.14} give a decomposition of the core $W$ 
into a union $W''$ of interval bundles, and the complement $W'$. The 
common boundary is a union of embedded annuli between components of 
$\partial W\setminus A$. These common annuli are taken 
to 
be metrically $S^{1}\times [0,1]$, up to bounded distortion. 
Otherwise, the metrics are constructed completely separately. Here, 
we 
treat the noninterval bundle piece $W'$. The construction is 
inductive. 
We take a metric such that $\Sigma '$ of \ref{6.14} has bounded 
length. We 
transfer this metric to each component of $\partial V$, for each 
component $V$ of $W_{i}$, for each $i$. Because of the cell cutting, 
this metric is unique up to bounded distortion, and does exist. 
This means that we 
have assigned an element  $y=y(S)=[\varphi ]\in {\cal{T}}(S)$ to each 
component $S$ of 
$\partial V$, for each component $V$ of $W_{i}$, up to bounded 
distance. 
We have a model manifold $M(y)$ homeomorphic to $S\times [0,1]$, 
simply by taking $\varphi _{*}\sigma +dt^{2}$ where $\sigma $ is the 
Poincar\'e metric on $\varphi (S)$. For the moment, we call this 
metric $\sigma _{y(S)}$. 
So 
now we need to obtain a metric on $W'$ from the metrics on each 
$\partial V$, by gluing together the spaces $M(y)$. The idea is 
exactly 
the same as in \ref{6.10} --- which 
was one reason for doing that first. Perturb the  $W_{i}$ so that 
the components of $W_{i-1}$ inside $W_{i}$ have boundaries in the 
interior of $W_{i}$. Attach incompressible surfaces to $\partial 
W_{i}$ in $W_{i}$, cutting $W_{i}$ up into the component pieces and 
disjoint from the perturbed components of $W_{i-1}$ inside. The 
region 
between $V$, and the 
components of $W_{i-1}$ inside, minus the cutting surfaces, is a 
union 
of  open 
interval bundles. Each corresponding closed interval bundle is 
homeomorphic to $S\times [0,1]$, where $S$ is the corresponding 
component of $\partial V'$, for a component $V'$ of $ 
W_{i-1}$ inside $V$. We can then map $S\times [0,1]$ to this region 
by a map 
$\Phi $ which is a diffeomorphism between the interiors, and take the 
metric 
$\Phi _{*}(\sigma _{y(S)})$ on the region. This metric has 
discontinuities as the 
boundary, but is the right metric up to bounded distortion. We can 
make it continuous and smooth by perturbing at the boundary. If some 
of the compressing surfaces are discs, we can make the metric on the 
surfaces inside by using the conformal structure on the surface 
outside, as explained in \ref{6.10}.

This model is made up of finitely many bounded interval bundles. It 
might be of interest to vary the topological type, and, instead of 
gluing together models $M(y)$, use models $M(y_{1},y_{2})$ for 
$[y_{1},y_{2}]$ possibly a long geodesic segment in ${\cal{T}}(S)$. 
It seems likely that the whole theory developed here could carry over 
to produce geometrical models for hyperbolic manifolds, and locally 
uniform biLipschitz constants, in this setting, under suitable 
conditions. 

\ssubsection{The model manifolds for compressible ends and ends 
without 
incompressible interval bundle bridges. }\label{6.16}

Curiously, a model manifold is easier to construct for a 
{\em{compressible}} 
end $e$ of $M_{d,W}$ than for an incompressible one. 
In fact, we have touched on the 
construction in \ref{6.10}. Let $S_{d}(e)$ be the corresponding 
component of $W\setminus \partial M_{d}$, the horodisc deletion of 
$S(e)$. Let $\Gamma _{0}(e)$ be the multicurve of 
boundaries of compressing discs in components of $\partial 
(S_{d}(e)\cap 
\partial W_{n-1})$. As in \ref{6.10}, we simply choose 
$z_{e,0}=[\varphi 
_{e,0}]$ so that $\vert \varphi _{e,0}(\Gamma _{0})\vert $ is 
bounded, but 
this is also the same, up to bounded distortion, as requiring that 
the arcs of $\Sigma `$ on $S_{d}(e)$ (as in \ref{6.14}) have 
bounded length.  If 
the geometrically finite ending invariant for this component of 
$M_{d,W}\setminus W$ is 
$y(e)$ then the model for this end is $M(z_{e,0},y(e))$.

Now suppose that $e$ is an  incompressible end, and is 
{\em{not}} the end of an interval bundle bridge. Then again, the point $z_{e,0}\in 
{\cal{T}}(S)$ is determined by making all arcs in $\Sigma '$ on 
$S$ of bounded length, and we again take the model to be  
$M(z_{e,0},y(e)),$ if the assocated end invariant is 
$(y(e))$. 

\ssubsection{The model for the interval bundle part of $W$ and 
associated ends}
\label{6.17}

The interval bundle part of the Scott core is a disjoint union of 
interval bundles each of which forms a bridge between two  components 
of 
$\partial W\setminus \partial M_{d}$. The model manifold on these 
interval bundles is dependent on the ending invariants of the 
associated ends, at least in the case of {\em{incompressible}} ends. 
Also, the model manifolds for the associated ends are interdependent.

We first deal with interval bundles between $S_{1}$ and $S_{2}$ 
where at least one of $S_{d,1}$, $S_{d,2}$ is {\em{compressible}}. 
Suppose that $S_{d,1}$ is compressible.
Suppose that the surface that is identified is $\alpha $. Then 
$\alpha $ itself must be incompressible. The corresponding model 
manifold is just the portion $\alpha \times [0,1]$ of $M(z_{0})$, 
with 
the model metric, where $z_{0}$ is the chosen basepoint on $S_{1}$, 
as in \ref{6.16}. If we choose to use the chosen basepoint on 
$S_{2}$, 
it does not matter. The metric is the same up to bounded distortion, 
depending only on the topological type. The choice of basepoints 
depends only on the topological type, not on the end invariants. 

Now we deal with models for  other bridging interval bundles, and the 
end model manifolds that they bridge between. As in \ref{6.16}, we 
only need to choose a  basepoint in ${\cal{T}}(S(e))$ for each end. 
We 
can then define the model manifolds for the ends as in \ref{6.16}, 
and 
the model for the bridging manifold as above. For each  
incompressible end $e'$ let $\omega (e,e')$ be the maximal subsurface 
of $S(e)$ which is homotopic to a subsurface of $S(e')$. By this, we  homotopic in $N$ and not in $S(e)$, in the case $e=e'$. We write 
$\omega (e',e)$ for this surface on $S(e')$.
Let 
 $$\beta  =\beta (e)=S(e)\setminus \cup _{e'}\omega  (e,e').$$

We start off with an 
initial choice $z_{e,0}=[\varphi _{e,0}]\in {\cal{T}}(S(e))$. We 
choose $\pi _{\beta }(z_{e,0})$ so that $\vert \varphi _{e,0}(\Sigma 
')\vert $, for $\Sigma '$ as in \ref{6.14}, is bounded, as before. 
Outside of $\beta $, the choice 
of 
$z_{e,0}$ is arbitrary for the moment. We also need to fix a 
geometrically finite ending invariant $(y(e))$. We 
let $y_{e,+}$ be the element of ${\cal{T}}(S(e))$ defined using the 
$(y(e))$, in the same way as $y_{-}$ in \ref{6.9}. 
We recall that the model for the end $e$ is defined by removing some 
pieces of Margulis tube and replacing them by pieces of horoball.

Now we need  to define $z_{e,0}'\in {\cal{T}}(S(e))$. Fix 
sufficiently strong ltd parameter functions $(\Delta 
_{1},r_{1},s_{1},K_{1})$
and a vertically efficient ltd-bounded decomposition of $S(e)\times
[z_{e,0},y_{e,+}]$.  For each incompressible end $e'$ for which
$\omega (e,e')\neq \emptyset $, we define $E_{\rm{ltd}}(e,e',-)$ to be
set of all ltd $(\alpha ,\ell )$ in the decomposition of $S(e)\times
[z_{e,0},y_{e,+}]$ with $\alpha \subset \omega (e,e')$.  Then since
$\vert \varphi _{e,0}(\partial \omega (e,e'))\vert $ is bounded in
terms of topological type, by \ref{7.2}, for sufficiently strong ltd
parameter functions , if $(\alpha ',\ell ')<(\alpha ,\ell )$, $\alpha
\subset \omega (e,e')$, and $(\alpha ',\ell ')$ is ltd, we also have
$\alpha '\subset \omega (e,e')$.  It follows that if $E(e,e',-)$ is
the set of all $(\alpha ' ,\ell ')$ with $(\alpha ',\ell ')<(\alpha
,\ell )$ for some $(\alpha ,\ell )\in E_{\rm{ltd}}(e,e',-)$ and
$E(e,e',+)$ is all the other $(\alpha ,\ell)$, then
$E(e,e')=(E(e,e',-),E(e,e',+))$ is an order splitting of the ltd
decomposition of $S(e)\times [z_{e,0},y_{e,+}]$ in the sense of
\ref{7.9}.  Then we define
 $$z(e,e')=x(E(e,e')).$$
 We remark that these conditions are consistent in the case $e=e'$.
 
 We define
 $$\pi _{\beta }(z_{e,0}')=\pi _{\beta }(z_{e,0}).$$
  Then, for $\omega  =\omega (e,e')$, define
 $$\pi _{\omega  }(z_{e,0}')=x(\pi _{\omega 
 }(z_{e,0}),[\pi _{\omega }(z(e,e')),\pi _{\omega }(z(e',e))]).$$
 Here, $x(.,.)$ is as in \ref{7.10}. This definition is such that 
 $d_{\omega (e,e')}(z_{e,0}',z_{e',0}')$ 
 is bounded in terms of \\ $(\Delta _{1},r_{1},s_{1},K_{1})$. Also, a 
ltd  decomposition for $S(e)\times [z_{e,0}',y_{e,+}]$ is given by a 
subset of 
the decomposition for $S(e)\times [z_{e,0},y_{e,+}]$.

Then our model for the end $e$ is obtained by altering 
$M(z_{e,0}',y_{e,+})$ just as in \ref{6.9}. The model for the bridge 
interval bundle between the ends $e$ and $e'$ is obtained from $M(\pi 
_{\omega 
}(z_{e,0}),\pi _{\omega 
}(z_{e,0}))$ by removing a model  horoball along the boundary 
components 
of $\partial \omega $.  That is, if $[z_{e,0}',y_{e,+}]=\{ y_{t}:t\in 
\{0,u]\} $, we remove $\cup _{t}\varphi _{t}^{-1}(H(\partial \omega 
,\varepsilon _{0}))$, where $H(\partial \omega 
,\varepsilon _{0})$ is the union of components of 
$(S_{t})_{<\varepsilon 
_{0}}$ homotopic to $\varphi _{t}(\partial \omega )$. In general, we 
 need to glue in a model Margulis tube in the place of each one --- 
 which might or might not be bounded. 
A bounded strip on the boundary comes from the non-interval bundle 
part of the model. We determine which Margulis tube to use from the 
geometry on the boundary, as usual.

From now on we redefine $z_{e,0}$ so that $[\pi _{\omega 
}(z_{e',0}'),\pi _{\omega }(z_{e',0})]\cup [\pi _{\omega 
}(z_{e,0}),\pi 
_{\omega }(z_{e,0}')]$ is a bounded distance, coordinatewise, from 
$[\pi _{\omega }(z_{e',0}'),\pi _{\omega }(z_{e,0}')]$, whenever $e$ 
and 
$e'$ are incompressible ends with $\omega =\omega (e,e')\neq 
\emptyset$. 
To do this, we leave $\pi _{\beta }(z_{e,0})$ as before. But for each 
pair $(e,e')$, we choose $x(e,e')=x(e',e)\in 
[\pi _{\omega }(z_{e',0}'),\pi _{\omega }(z_{e,0}')]$, and define 
$z_{e,0}$ by
$$\pi _{\omega }(z_{e,0})=\pi _{\omega }(z_{e',0})=\pi _{\omega 
}(x(e,e')),$$
$$\pi _{\partial \omega }(z_{e,0})=\pi _{\partial \omega 
}(z_{e',0})=\pi 
_{\partial \omega 
}(x(e,e')),$$
whenever $\omega =\omega (e,e')=\omega (e',e)$.

\ssubsection{Model for the whole manifold.}\label{6.18}

The model for a hyperbolic manifold with core $N_{c}$ and 
geometrically finite 
invariants $y(e)$ for each end is obtained by gluing together 
models for $N{c}$ and models $M(z_{e,0},y(e))$ for each end. The 
boundaries are the same up to bounded distortion and we can remove 
discontinuities by perturbation, without changing the metric up to 
coarse Lipschitz equivalence. The model for the core, as we have 
seen, is formed by gluing together interval bundles. Rather than 
using a model for $N_{c}$ and each of the ends, it is sometimes more 
convenient to take a model for some $W\subset N_{c}$ and glue this 
together with models for each end of $N_{d,W}$. This is the case, for 
example, when $N$ is homeomorphic to $S\times \mathbb R$, but has 
more than two ends.

\ssubsection{The combinatorially bounded geometry geometrically 
infinite  
Kleinian surface case.}\label{6.12}

Now we consider the combinatorially bounded geometry  case of $(\mu 
_{-},\mu _{+})\in 
{\cal{T}}(S)\cup {\cal{GL}}_{a}(S)$ and convergence of the Kleinian 
surface models $M(y_{n,-},y_{n,+})$ to $M(\mu _{-},\mu _{+})$ as 
$y_{n,-}\to \mu _{-}$, $y_{n,+}\to \mu _{+}$. From now  on in this 
subsection, we assume 
that $\mu _{+}\in {\cal{GL}}_{a}(S)$ and $\mu _{-}\in 
{\cal{GL}}_{a}(S)\cup {\cal{T}}(S)$. 

Let $i$ be the intersection number 
of 
\ref{3.5}. We consider the condition that, for at least one 
normalised transverse invariant measure on $\mu _{+}$, for some 
$c=c(\mu 
_{+})>0$, 
\begin{equation}\label{6.12.1}i(\mu _{+},\gamma )\geq c(\mu 
_{+})\vert \gamma \vert 
^{-1}\end{equation}
 for all simple closed 
loops $\gamma $, and similarly for $\mu _{-}$. Here, $\vert \gamma 
\vert $ is measured with respect to a fixed hyperbolic metric on $S$, 
as is used in \ref{3.5}. It can be shown that 
this condition implies that $\mu _{+}$ has only one transverse 
invariant measure up to scalar, and similarly for $\mu _{-}$. It  
probably helps to regard $\mu _{+}$ as a measured foliation with 
transverse invariant measure equivalent to Lebesgue measure (which it 
is, up to measure isomorphism). The 
proof of the fundamental dynamical lemma \ref{7.1} shows that the 
combinatorially bounded geometry condition implies that every leaf of 
length $L$ comes within  distance  $c_{1}/L$ of every point 
in the measured foliation, 
for every $L>0$, for $c_{1}$ depending only on $c(\mu _{+})$. It 
follows 
that given any invariant set $E$ under the foliation, every point is 
a positive density point for $E$, and therefore $E$ has full measure 
and the transverse invariant measure is ergodic. But if every 
transverse invariant measure is ergodic, there is only one, up to 
scalar. 

Fix a basepoint $y_{0}=[\varphi _{0}]\in \cal{T}(S)$, We claim that 
if 
$y_{n,+}=[\varphi _{n,+}]\to \mu _{+}$
and, for all $n$ 
\begin{equation}\label{6.12.2}[y_{0},y_{n,+}]\subset 
({\cal{T}}(S))_{\geq \nu },\end{equation}
then (\ref{6.12.1}) holds for a $c(\mu _{+})$ depending only on 
$\nu $ and $y_{0}$. For let $\zeta _{n,+}$ be a loop such that 
\begin{equation}\label{6.12.3}C_{1}^{-1}\vert \varphi 
_{0}(\zeta _{n,+})\vert \leq \exp 
d(y_{0},y_{n,+})\leq C_{1}\vert \varphi 
_{0}(\zeta _{n,+})\vert \end{equation} and 
\begin{equation}\label{6.12.4}\vert \varphi _{n,+}(\zeta _{n,+})\vert 
\leq C_{1}.\end{equation}
This is possible for a suitable $C_{1}$ by \ref{2.3}, and the fact 
that 
$y_{n,+}\in (\cal{T}(S))_{\geq \nu }$. Enlarging $C_{1}$ we then also 
have, for all $y=[\varphi ]\in [y_{0},y_{n,+}]$, 
$$C_{1}\vert \varphi 
(\zeta _{n,+})\vert \geq \exp 
d(y_{n,+},y). $$
Then any limit of $\zeta 
_{n,+}/\vert \zeta _{n,+}\vert $ (taking $\vert .\vert $ with respect 
to a fixed hyperbolic metric on $S$) has zero intersection with 
$\mu _{+}$, and must be $\mu _{+}$ by arationality. For any simple 
closed 
loop $\gamma $, choose $y=[\varphi ]\in [y_{0},y_{n,+}]$ such that 
$\vert 
\varphi (\gamma )\vert $ is minimal over all such $y$. Then $\vert 
\varphi (\gamma 
)\vert $ is bounded from $0$, and the good position of $\varphi 
(\gamma )$ is such that along most of its length it is bounded from 
the stable and unstable foliations of the quadratic differential for 
$d(y,y_{n,+})$. So for suitable $C_{1}$,
$$C_{1}\vert \varphi _{0}(\gamma )\vert \geq \exp d(y_{0},y).$$
Using this and   
\ref{7.2}, since $S$ is ltd along $[y_{0},y_{n,+}]$,
$$i(\gamma ,\zeta _{n,+})=\# (\varphi (\gamma )\cap \varphi (\zeta 
_{n,+}))\geq C_{2}^{-1}\vert \varphi (\zeta _{n,+})\vert \geq 
(C_{1}C_{2})^{-1}\exp d(y,y_{n,+})$$
$$\geq C_{3}\exp (d(y_{0},y_{n,+})-d(y_{0},y))\geq C_{4}{\vert \zeta 
_{n,+}\vert \over \vert \gamma \vert }.$$
So (\ref{6.12.2}) implies (\ref{6.12.1}). 

We claim that the converse is also essentially true. So suppose that 
(\ref{6.12.1}) holds and that $y_{n,+}\to \mu _{+}$. Choose $\zeta 
_{n,+}$ so that (\ref{6.12.3}) is replaced by 
\begin{equation}\label{6.12.5}C_{1}^{-1}{\vert \varphi 
_{0}(\zeta _{n,+})\vert \over \vert \varphi _{n,+}(\zeta _{n,+})\vert 
'}\leq \exp 
d(y_{0},y_{n,+})\leq C_{1}{\vert \varphi 
_{0}(\zeta _{n,+})\vert \over \vert \varphi _{n,+}(\zeta 
_{n,+})\vert },\end{equation}
and that (\ref{6.12.4}) is 
replaced by
\begin{equation}\label{6.12.6}\vert \varphi _{n,+}(\zeta _{n,+})\vert 
''
\leq C_{1}.\end{equation}
Here, $\vert .\vert '$ and $\vert .\vert ''$ are as in \ref{2.3}. 
This time we do not have a lower bound on $\vert \varphi _{n,+}(\zeta 
_{n,+})\vert $. But choose $y_{n,+}'=[\varphi _{n,+}']\in 
[y_{0},y_{n,+}]$ such that (\ref{6.12.6}) holds for $\varphi _{n,+}'$ 
replacing $\varphi _{n,+}$, and also
\begin{equation}\label{6.12.7}\vert \varphi _{n,+}'(\zeta 
_{n,+})\vert 
\geq C_{1}^{-1}.\end{equation}
Then we claim that
\begin{equation}\label{6.12.8}[y_{0},y_{n,+}']\subset 
({\cal{T}}(S))_{\geq \nu }.\end{equation}
Suppose this is 
not true and that there is a loop $\gamma $ and $y=[\varphi ]\in 
[y_{0},y_{n,+}']$ with $\vert \varphi (\gamma )\vert \leq \nu $. We 
take $[w_{1},w_{2}]=[[\psi _{1}],[\psi _{2}]]\subset 
[y_{0},y_{n,+}']$ 
such that 
$$\vert \psi _{j}(\gamma )\vert \leq C_{1},$$  
and 
$$d(w_{1},w_{2})\geq L(\nu ),$$
where $L(\nu )\to \infty $ as $\nu 
\to 0$. In fact, we can take $L(\nu )\geq C_{1}^{-1}\log (1/\nu )$.  
Then, using \ref{2.3.2}, 
$$i(\gamma ,\zeta _{n,+})\leq C_{2}\exp d(w_{2},y_{n,+}')
\leq C_{2}\exp -L(\nu )\exp 
(d(y_{0},y_{n,+}')-d(w_{1},y_{0}))$$
$$\leq C_{3}\exp -L(\nu ){\vert \varphi 
_{0}(\zeta _{n,+})\vert \over \vert \varphi _{0}(\gamma 
)\vert }.$$
 This contradicts 
(\ref{6.12.1})

Similar arguments work for $\mu _{-}$ and sequences $y_{n,-}$, 
$y_{n,-}'$, if either $\mu _{-}\in {\cal{GL}}_{a}(S)$ satisfies 
(\ref{6.12.1}) or (\ref{6.12.2}) holds for the $y_{n,-}$. So now if 
we assume that either $\mu _{-}\in {\cal{GL}}_{a}(S)$ and satisfies 
(\ref{6.12.1}) and $y_{n,-}$, $\zeta _{n,-}$, $y_{n,-}$ are defined 
similarly to $y_{n,+}$, $\zeta _{n,+}$, $y_{n,+}'$, or $\mu _{-}\in 
{\cal{T}}(S)$ and $y_{n,-}=y_{n,-}'=\mu _{-}$ for all $n$. Then for a 
suitable $\nu >0$,
$$[y_{0},y_{n,-}']\in (\cal{T}(S))_{\geq \nu },$$
and hence applying \ref{7.4} to the triangle with vertices at 
$y_{0}$, $y_{n,+}'$, $y_{n,-}'$, for suitable $\nu $, for all $n$,
$$[y_{n,-}',y_{n,+}']\subset (\cal{T}(S))_{\geq \nu }.$$

Now we claim that if $\mu _{-}$, $\mu _{+}\in {\cal{GL}}_{a}(S)$ 
satisfy (\ref{6.12.1}) and $\mu _{+}\neq \mu _{-}$ then for a 
constant $C=C(\mu _{+},\mu _{-})$, for all sufficiently large $n$,
\begin{equation}\label{6.12.9}d(y_{n,-}',y_{0})+d(y_{0},y_{n,+}')\leq 
d(y_{n,-}',y_{n,+}')+C.\end{equation}
The reason is simply that $\mu 
_{+}\neq \mu _{-}$ means $i(\mu _{+},\mu _{-})>0$ by arationality,
 and hence for all 
sufficiently large $n$ and constants $C_{j}$ depending on 
$\mu _{+}$, $\mu _{-}$,  
$$\exp d(y_{n,-}',y_{n,+}')\geq C_{1}^{-1}i(\zeta _{n,+},\zeta 
_{n,-})\geq C_{2}^{-1}\vert \zeta _{n,+}\vert .\vert \zeta 
_{n,-}\vert $$
$$\geq C_{3}^{-1}\vert \varphi _{0}(\zeta _{n,+})\vert .
\vert \varphi _{0}(\zeta _{n,-})\vert \geq C_{4}^{-1}\exp 
(d(y_{n,-}',y_{0})+d(y_{0},y_{n,+}')).$$
It follows from \ref{7.4} and the fact that the geodesic segments 
joining $y_{0}$, $y_{n,+}'$, $y_{n,-}'$ are in $({\cal{T}}(S))_{\geq 
\nu }$ that, for a constant $C'$ depending only on $\mu _{\pm }$, for 
all sufficiently large $n$,
\begin{equation}\label{6.12.10}d(y_{0},x(y_{0}))\leq C',\end{equation}
where $x(y_{0})$ denotes the orthgonal projection (\ref{7.8}) to 
$[y_{n,-}',y_{n,+}']$ for any $n$, that is, $y_{0}$ is distance $\leq 
C'$ from some point on $[y_{n,-}',y_{n,+}']$. It then follows from 
\ref{7.4} that each point on $[y_{n,-}',y_{0}]\cup [y_{0},y_{n,+}']$ 
is distance $\leq C''$ from some point on $[y_{n,-}',y_{n,+}']$, for 
$C''$ independent of $n$. So then $y_{0}$ is a distance $\leq C'$ 
from some point on $[y_{n,-},y_{n,+}]$, for all sufficiently large 
$n$, and any $y'\in [y_{n,-}',y_{n,+}']$ is distance $\leq C''$ from 
some point in $[y_{n,-},y_{n,+}]$. 

This finally puts us in a position to prove geometric 
convergence, up to bounded coarse Lipschitz equivalence, of suitably 
based
models $(M(y_{n,-},y_{n,+}),x_{n})$, $(M(y_{n,-}',y_{n,+}'),x_{n}')$,
if $y_{n,+}\to \mu _{+}$ where $\mu _{+}$ satisfies (\ref{6.12.1}) and
either similar properties hold for $\mu _{-}$, or $y_{n,-}=\mu _{-}\in
{\cal{T}}(S)$ for all $n$.  In all cases we have $y_{n,0}\in
[y_{n,-},y_{n,+}]$ and $y_{n,0}'\in [y_{n,-}',y_{n,+}']$ with
$$d(y_{0},y_{n,0})\leq C',\ \ d(y_{0},y_{n,0}')\leq C'.$$
We translate the vertical coordinate so that 
$$M(y_{n,-},y_{n,+})=S\times [-u_{n,-},u_{n,+}],\ \ 
M(y_{n,-}',y_{n,+}')=S\times [-u_{n,-}',u_{n,+}'],$$
$$d(y_{n,-},y_{n,0})=u_{n,-},\ \ d(y_{n,+},y_{n,0})=u_{n,+},$$
$$d(y_{n,-}',y_{n,0}')=u_{n,-}',\ \ d(y_{n,+}',y_{n,0}')=u_{n,+}'.$$
We take $x_{0}$ to be any fixed point in $S\times \{ 0\} $,  
$x_{0}=x_{n}=x_{n}'$ for all $n$. 

Then to prove geometric convergence 
of the models up to bounded coarse Lipschitz equivalence, it suffices 
to
prove that, for a constant $C_{0}$, for any $\Delta >0$, for all
sufficiently large $k$ and $n$, if $y\in [y_{k,-},y_{k,+}]$ and
$d(y_{0},y)\leq \Delta $, there is $y'\in [y_{n,-},y_{n,+}]$ with
\begin{equation}\label{6.12.11}d(y,y')\leq C_{0},\end{equation}
and similarly for $y_{k,\pm }$, $y_{n,\pm }$ replaced by 
$y_{k,\pm }'$, $y_{n,\pm }'$.  For suppose we have this.  Then
replacing $y_{k,\pm}$ by suitable points in the original geodesic
segment $[y_{k,-},y_{k,+}]$ if necessary, and replacing $n$ by a
subsequence if necessary, we can assume that (\ref{6.12.1}) holds for
all $y\in [y_{k,-},y_{k,+}]$, and for all $n\geq k$.  We can also
assume that $d(y_{k,+},y_{n,+})$ and $d(y_{k,-},y_{n,-})$ is bounded
from $0$ for all $n>k$.  Then we can construct a boundedly coarse
biLipschitz map $\varphi _{k,n}$, from $M(y_{k,-},y_{k,+})$ to a
subset of $M(y_{n,-},y_{n,+})$, fixing $x_{0}$, for any $k>n$, 
and such that the distance, in the
Riemannian metrics, between $\varphi _{k,n}(\partial
M(y_{k,-},y_{k,+}))$ and $V_{k,n}=\partial M(y_{n,-},y_{n,+})\setminus \varphi
_{k,n}(\partial M(y_{k,-},y_{k,+}))$, is bounded from $0$.  Then we
can make a Riemannian manifold $V_{p}$ by gluing together
$M(y_{1,-},y_{1,+})$ and $V_{n,n+1}$ for $1\leq n<p$, taking a metric which is
the model metric for $M(y_{1,-},y_{1,+})$ on $M(y_{1,-},y_{1,+})$ and
the model metric for $M(y_{n+1,-},y_{n+1,+})$ on
$V_{n,n+1}$, except near $\partial \varphi
_{n,n+1}(M(y_{n,-},y_{n,+}))$.  on all of $V_{n,n+1}$, the metric is taken
boundedly equivalent to the metric on $M(y_{n+1,-},y_{n+1,+})$. 
The manifold $(V_{p},x_{0})$ has a based submanifold which is
naturally diffeomorphic to $(M(y_{n,-},y_{n,+}),x_{0})$, for each
$n\leq p$, with bounds on the derivative and derivative inverse with
respect to the Riemannian metrics, and the geometric limit $\lim
_{p\to \infty }(V_{p},x_{0})$ exists as a based Riemannian manifold. 
Each set $V_{p+1}\setminus V_{p}=V_{p,p+1}$ is homeomorphic to
$S\times I$ where $I$ is the union of one or two intervals, depending
on whether $\mu _{-}\in {\cal{T}}(S)$ or $\mu _{-}\in
{\cal{GL}}_{a}(S)$.  It follows that $\lim _{p\to \infty }V_{p}$ is
homeomorphic to $S\times [0,\infty )$ or $S\times \mathbb R$.  The
geometric limit $\lim _{p\to \infty }V_{p}$ depends only on $\mu
_{\pm}$ up to coarse biLipschitz equivalence, not on the precise
sequences $y_{n,\pm }$, by (\ref{6.12.11}), and not on the precise
choice of metric near $\partial V_{p}$ for any $p$.

It suffices to
prove (\ref{6.12.11}) for $y_{k,\pm }'$ and $y_{n,\pm }'$, since every
point in $[y_{k,-}',y_{k,+}']$ is a bounded distance from a point in
$[y_{k,-},y_{k,+}]$, and $d(y_{0},y_{k,+}')\to \infty $ as $k\to
\infty $ (because $\vert \zeta _{k,+}\vert \to \infty $), and
similarly for $y_{k,-}'$ if $\mu _{-}\in {\cal{GL}}_{a}(S)$.  So now
let $y\in [y_{k,-}',y_{k,+}']$ with $d(y_{0},y)\leq \Delta $.  Suppose
that $y$ is not a bounded distance $\leq C_{0}$ from a point in
$[y_{n,-}',y_{n,+}']$.  Then by \ref{7.4} for quadrilaterals, and
$(\cal{T}(S))_{\geq \nu }$, taking vertices $y_{n,\pm }'$, $y_{k,\pm
}'$, for suitable $C_{0}$ (given by \ref{7.4}), $y$ must be a distance
$\leq C_{0}$ from $y'\in [y_{n,+}',y_{k,+}']\cup 
[y_{k,-}',y_{n,-}']$. 
We can assume without loss of generality that $y'=[\varphi ']\in
[y_{n,+}',y_{k,+}']$.  If $n$ and $k$ are both large given $\Delta $
then $d(y_{k,+}',y')$ and $d(y_{n,+}',y')$ are both large.  Then by
\ref{7.2}, $$i(\zeta _{k,+},\zeta _{n,+})=\# (\{ \varphi '(\zeta
_{k,+})\cap \varphi '(\zeta _{n,+})\} )\geq C_{1}^{-1}\vert \varphi
'(\zeta _{k,+})\vert .\vert \varphi '(\zeta
_{n,+})\vert $$
$$\geq C_{2}^{-1}e^{-2\Delta }\vert \varphi _{0}(\zeta 
_{k,+})\vert .\vert \varphi _{0}(\zeta 
_{n,+})\vert \geq C_{3}^{-1}e^{-2\Delta }\vert \zeta 
_{k,+}\vert .\vert \zeta 
_{n,+}\vert .$$
This cannot be true for arbitrarily large $k$ and $n$, because we 
would 
then deduce that $i(\mu _{+},\mu _{+})>0$. So geometric convergence 
of the models is complete in this case.

\ssubsection{Geometric convergence of models in the general 
geometrically infinite case}\label{6.13}

All geometric models in the geometrically finite case are obtained by 
gluing together finitely many models $M(y_{-},y_{+})$, with 
insignificant modifications to some model Margulis tubes in the 
ending models with extra cusps in the ends. For the 
geometrically infinite models, we simply want to take geometric 
limits 
of models of the form $M(y_{0},y_{n,+})$ for some exterior models, or 
a geometric limit of a single sequence of models 
$M(y_{n,-},y_{n,+})$, for a suitable choice of 
basepoint and with $y_{n,+}$ convergent to some point in $\mu _{+}\in 
\partial{\cal{T}}(S)$, where $\partial {\cal{T}}(S)$ is the 
modification of the Thurston boundary described in \ref{3.6}, and 
either the same is true for $\mu _{-}$, or $\mu _{-}\in 
{\cal{T}}(S)$.  

So  let $y_{n,+ }=[\varphi _{n,+}]\to \mu _{+}$ and $y_{n,-}\in \mu 
_{-}$. If $\mu _{\pm }$ are both arational geodesic laminations on 
$S$, we 
assume that $i(\mu _{-},\mu _{+})>0$. If $\mu _{+}'$ and $\mu _{-}'$ 
denote the lamination parts of 
$\mu _{\pm}$, we also assume that no 
closed loop $\gamma \subset S $ in the closure of the support of $\mu 
_{+}'$ 
or $\mu _{-}'$ satisfies $i(\gamma ,\mu _{+}')=i(\gamma ,\mu 
_{-}')=0$. The latter condition is to ensure that any geometric limit 
is connected (but is not actually a necessary condition for this).
We need to show 
that, for a suitable choice of base-point $x_{n}\in 
M(y_{n,-},y_{n,+})$,
$(M(y_{n,-},y_{n,+}),x_{n})$ has a single geometric limit up to 
bounded coarse  biLipschitz equivalence. 

 Fix a basepoint 
$y_{0}=[\varphi _{0}]$ as in \ref{6.12}. As in \ref{6.12}, we can 
find loops $\zeta _{n,\pm }$ and $y_{n,\pm }'=[\varphi _{n,\pm }']$ 
converging to $\mu _{\pm }$  
such that (\ref{6.12.6}) and (\ref{6.12.7}) hold. If $\mu _{+}$ has 
at least one lamination component, we can also assume that $\zeta 
_{n,+}$ has nonempty intersection with the support of at least one 
minimal component of $\mu 
_{+}$, and similarly for $\zeta _{n,-}$, $\mu _{-}$. Let 
$x(y_{0},[y_{n,-},y_{n,+}])$ be the orthogonal projection of 
$y_{0}$ 
relative to 
$[y_{n,-},y_{n,+}]$, as in \ref{7.10}.  We claim that, in order to 
show
geometric convergence, up to bounded coarse Lipschitz equivalence, it
suffices to show that for a suitable constant $L_{1}$, given $\Delta
>0$, for all sufficiently large $k$,

\begin{equation}\label{6.13.1}d(y_{0},x(y_{0},[y_{k,-},y_{k,+}]))\leq
L_{1},\end{equation} 

and for all sufficiently large $k$, and $n$, and
all ltd $(\alpha ,\ell )$ for $[y_{k,-},y_{k,+}]$ with $d_{\alpha
}'(y,y_{0})\leq \Delta $ for $y\in \ell $, there is $(\alpha ,\ell ')$
which is ltd along $\ell '\subset [y_{n,-},y_{n,+}]$, possibly with
respect to different ltd parameter functions (which is enough, by
\ref{6.6}), and for all $y\in \ell $, there is $y'\in \ell '$ with

\begin{equation}\label{6.13.2}d_{\alpha }(y,y')\leq L_{1}{\rm{\ or\
}}\vert {\rm{Re}}(\pi _{\alpha }(y)-\pi _{\alpha }(y'))\vert \leq
L_{1},\end{equation}
depending on whether $\alpha $ is a loop or a
gap.  We see that these suffice as follows. So suppose that both 
(\ref{6.13.1}) and (\ref{6.13.2}) hold. 

We first look for suitable basepoints $x_{0}=(z_{0},0)$ in the models. 
This means looking for a gap or loop $\alpha $ such that $z_{0}\in
\varphi _{0}^{-1}(S_{\alpha ,0})$, that is, $S_{\alpha ,0}$ is
homotopic to $\varphi (\alpha )$ and bounded by $S_{\partial \alpha
,0}$ in the terminlogy of \ref{6.2}.  We can drop the first few terms
of the sequence if necessary, and assume that there is $(\alpha ,\ell
)$ which is ltd along $[y_{1,-},y_{1,+}]$, $d_{\alpha }'(y,y_{0})\leq
\Delta _{1}$ for all $y\in \ell $, for $\Delta _{1}$ depending on the
ltd parameter functions and $y_{1,\pm }$, and (\ref{6.13.2}) holds for
all $n$.  Similarly, (\ref{6.13.2}) also holds for $y_{n,0}'=[\varphi
_{n,0}']\in [y_{n,-}',y_{n,+}']$ replacing $y_{n,0}$, again because
the ltd's along $[y_{0},y_{n,\pm }']$ are a subset of those along
$[y_{0},y_{n,\pm }]$ up to bounded distance and hence the ltd's along
$[y_{n,-}',y_{n,+}']$ are a subset along $[y_{n,-},y_{n,+}]$, up to
bounded distance.  Conversely if we have (\ref{6.13.1}) and
(\ref{6.13.2}) with $y_{n,\pm }'$ replacing $y_{n,\pm }$ and
$y_{n,0}'\in [y_{n,-}',y_{n,+}']$ replacing $y_{n,0}$, then we have
(\ref{6.13.1}) and (\ref{6.13.2}) for $y_{n,\pm }$.  So suppose that
we have all of these.  Let $\varphi _{n,0}$ and $\varphi _{n,0}'$ be
homeomorphisms which are part of families $\varphi _{n,t}$, $\varphi
_{n,t}'$, relative to $[y_{n,-},y_{n,+}]$ and $[y_{n,-}',y_{n,+}']$,
satisfying the properties of \ref{6.2}.  Translate the vertical
coordinate of the model as in \ref{6.12}, defining $u_{n,\pm }$,
$u_{n,\pm }'$ as there.  If $\alpha $ is a gap and long $\nu
_{n}$-thick and dominant at $y_{n,0}$ then we take $x_{n}=(z_{n},0)$
so that $\varphi _{n,0}(x_{n})$ is in the $\nu _{n}$-thick part of
$\varphi _{n,0}(S)$ which is homotopic to $\varphi _{n}'(\alpha )$. 
Let $x_{n}'=(z_{n}',0)$ be similarly defined relative to $\varphi
_{n,0}'$.  We can, and do, choose $\varphi _{n}$ and $\varphi _{n}'$
so that $z_{n}=z_{n}'=z_{0}$.  If $\alpha $ is a loop, we can choose
it so that it is transverse to a lamination component of each of $\mu
_{\pm }$.  Then we claim that we have a lower bound on $\vert \varphi
(\alpha )\vert $ for $[\varphi ]\in [y_{n,-}',y_{n,+}']$ and for all
$n$, and hence similarly for $[y_{n,-},y_{n,+}]$.  For from
(\ref{6.13.1}) we have, for $x(y_{0})=x(y_{0},[y_{n,-}',y_{n,+}'])$,
for a constant $C_{0}$ depending only on the ltd parameter functions,
\begin{equation}\label{6.13.3}d(y_{n,-}',x(y_{0}))+d(x(y_{0}),y_{n,+}')\leq
d(y_{n,-}',y_{n,+}')+C_{0},\end{equation} which can be seen by using
the bound \ref{2.8} for $d(y_{n,-}',y_{n,+}')$ in terms of
$${\rm{Max}}_{\beta }(d_{\beta }'(x(y_{0}),y_{n,-}')+d_{\beta
'}(x(y_{0}),y_{n,+}')).$$ 
So then for  $\zeta _{n,\pm }$, $y_{n,\pm }'=[\varphi _{n,\pm }']$, 
$y_{n,\pm }=[\varphi _{n,\pm }]$ 
satisfying (\ref{6.12.6}) and (\ref{6.12.7}) as above, 
suppose that we do not have a lower bound on $\vert \varphi (\alpha 
)\vert $.
Let  $w_{j}=[\psi _{j}]$, $j=1$, $2$, $[w_{1},w_{2}]\subset 
[y_{n,-}',y_{n,+}']$, $\vert \psi _{j}(\alpha )\vert \leq C_{1}$, 
$d(w_{1},w_{2})\geq \Delta _{2}$, where $\Delta _{2}$ can be taken 
large if some $\vert \varphi (\alpha )\vert $ is small. We have  
$$i(\alpha ,\zeta 
_{n,-}).i(\alpha ,\zeta _{n,+})\leq C_{2}^{-1}\exp 
((d(w_{1},y_{n,-}')+d(w_{2},y_{n,+}'))$$
$$\leq C_{3}^{-1}e^{-\Delta _{2}}\exp 
((d(y_{0},y_{n,-}')+d(y_{0},y_{n,+}'))
\leq C_{4}^{-1}e^{-\Delta _{2}}\vert \zeta _{n,-}\vert \vert \zeta 
_{n,+}\vert $$
If this is true for arbitrarily large $\Delta _{2}$ and hence also 
arbitrarily large $n$, then taking limits, we obtain 
$$i(\alpha ,\mu _{+}').i(\alpha ,\mu _{-}')=0,$$
or $i(\alpha ,\mu _{+})=0$ if $y_{n,-}=y_{n,-}'=\mu _{-}$. By 
the choice of $\alpha $, this is impossible. 

Next, we can extend (\ref{6.13.2}) to the set of all $(\alpha ,\ell )$ 
in a vertically efficient decomposition for $S\times 
[y_{k,-},y_{k,+}]$ with $d_{\alpha }'(y,y_{0})\leq \Delta $ for $y\in 
\ell $.For given such an $(\alpha ,\ell )$,  we can find a 
corresponding set for $S\times
[y_{n,-},y_{n,+}]$, if $n$ is sufficiently large, by taking upper and 
lower boundaries (\ref{7.10})
of sets of ltd's below and above $\alpha \times \ell $ in $S\times
[y_{k,-},y_{k,+}]$, and taking the corresponding upper and lower 
boundaries in
$S\times [y_{n,-},y_{n,+}]$.  We can thicken slightly so that the
upper and lower boundaries are disjoint.  Then the two pieces
corresponding to $(\alpha ,\ell )$ in $M(y_{k,-},y_{k,+})$ and
$M(y_{n,-},y_{n,+})$ are again boundedly coarse Lipschitz equivalent. 
Model Margulis tubes in each are determined by the metrics on their
boundaries.  So each model Margulis tube in $M(y_{k,-},y_{k,+})$ which
is completely encased by pieces corresponding to ltd or bounded  
$(\alpha ,\ell )$, with $d_{\alpha }'(y,y_{0})\leq \Delta $ for $y\in 
\ell $, 
is boundedly Lipschitz equivalent to a model Margulis tube in
$M(y_{n,-},y_{n,+})$.

Next, in analogy to what was done in \ref{6.12}, we change the
sequences $y_{n,\pm }$ so that (\ref{6.13.1}) and (\ref{6.13.2}) hold
for all $(\alpha ,\ell )$ in the decomposition for
$[y_{k,-},y_{k,+}]$, and for all $n>k$. For a suitable $\Delta _{n}$, let
$(E_{n}(-,+),E_{n}(+,+))$ be an order splitting for $[y_{0},y_{n,+}]$
(\ref{7.10}) so that $d_{\alpha }'(y,y_{0})\leq \Delta _{n}$ for all
$y\in \ell $ and $(\alpha ,\ell)\in E_{n}(-,+)$ and $d_{\alpha
}'(y,y_{0})\geq \Delta _{n}$ for $y\in \ell $ and $(\alpha ,\ell )\in
E_{n}(+,+)$.  Restricting to a subsequence if necessary, we can assume that
(\ref{6.13.2}) holds for $\Delta =\Delta _{k}$ and all $n>k$ and
similarly for the sequence $y_{n,-}$.  Replace $y_{n,+}$ by
$x(E_{n}(-,+))=x(E_{n}(+,+))$ in the notation of \ref{7.10}, and
similarly for $y_{n,-}$.  Then, as in \ref{6.12}, we have a sequence of
maps $\varphi _{k,n}$ from $(M(y_{k,-},y_{k,+}),x_{0})$ into
$(M(y_{n,-},y_{n,+}),x_{0})$ for $k\leq n$.  We can choose $\varphi
_{k,n}$ to have uniformly bounded derivative and inverse derivative
with respect to the Riemannian metrics on $M(y_{k,-},y_{k,+})$ and
$M(y_{n,-},y_{n,+})$.  Then, as in \ref{6.12}, we form the sequence
$V_{n}$ by gluing together pieces $M(y_{1,-},y_{1,+})$ and
$$V_{n,n+1}=M(y_{n+1,-},y_{n+1,+})\setminus 
\varphi _{n+1,n}(M(y_{n,-},y_{n,+})),$$
so that
the Riemannian metric on $V_{n,n+1}$, except in a neighbourhood of
\\ $\varphi _{n,n+1}(M(y_{n,-},y_{n,+}))$, is the Riemannian metric on
$M(y_{n+1,-},y_{n+1,+})$, and is boundedly equivalent to this metric
everywhere on $V_{n,n+1}$.  

If $\mu _{+}$ is reducible, that is,
$i(\mu _{+},\gamma )=0$ for at least one nontrivial nonperipheral
closed loop, then there are some model Margulis tube boundaries which
intersect $\partial V_{n}$ for all $n$.  We denote by $\partial
_{h}V_{n}$ the complement in $\partial V_{n}$ of any such Margulis tube
boundaries.  In order for the limit to be a Riemannian manifold and a
topological product, we need the distance between $\partial _{h}V_{n}$ and
$\partial _{h}V_{n+1}\setminus \partial _{h}V_{n}$ to be bounded from $0$,
replacing the original sequence for a sufficiently fast increasing
subsequence if necessary.  We concentrate on the boundary
corresponding to $\mu _{+}$.  (There is boundary corresponding to $\mu
_{-}$ only if $\mu _{-}\notin {\cal{T}}(S)$.)  Write $\partial
_{h,+}V_{n}$ for the union of boundary components corresponding to
$y_{n,+}$.  Fix $k$ so that $\partial _{h,+}V_{k}$ contains components
corresponding to all Teichm\"uller space components of $\mu _{+}$.  It
suffices to show that the distance in the model Riemannian metric of
$\partial _{h,+}V_{n}\setminus \partial _{h,+}V_{k}$ from $x_{0}$ $\to
\infty $ as $n\to \infty $, and the same for any model Margulis tubes
intersecting $\partial _{h,+}V_{n}\setminus \partial _{h,+}V_{k}$.  We
see this as follows.  If not, then there is $\Delta >0$, and, for each
$n$, a path from $x_{0}$ to $\partial _{h,+}V_{n}\setminus \partial
_{h,+}V_{k}$ of length $\leq \Delta $.  Restricting to a subsequence
if necessary, we can assume that the path always passes through the
same model Margulis tubes and the same sets $W_{j}$, in the notation
of \ref{6.2}, corresponding to sets $\alpha \times \ell $ in the
vertically efficient decomposition for $S\times [y_{n,-},y_{n,+}]$,
for some sufficiently large $n$.  So then the path for $n$ can be
assumed to end in $(\gamma ,u_{n,+})$ for a nontrivial nonperipheral
closed loop $\gamma \subset S$ such that $\vert \varphi _{n,+}(\gamma
)\vert \leq \Delta '$ for all $n$.  There is no relation between
$\Delta $ and $\Delta '$, but $\Delta '<+\infty $, because $\Delta '$
is determined by the $W_{j}$ passed through, correspinding to sets
$\alpha \times \ell $.  Then $\gamma $ lies in a geodesic lamination
component of $\mu _{+}$, but $i(\mu _{+},\gamma )=0$, which is a
contradiction.

Now we prove (\ref{6.13.1}).  We 
consider the sets of ltds 
$T(y_{0},+)$ and $T(y_{0},-)$ of \ref{7.9}. If (\ref{6.13.1}) does 
not 
hold for a sufficiently large $L_{1}$ given $L_{2}$, then, by 
\ref{5.5}, there must 
be some 
maximal totally ordered set  of ltds $(\alpha _{i},\ell _{i})$, 
$1\leq i\leq 
m$ along 
$[y_{0},x(y_{0})]$
with $\ell _{m}$ nearest to $x(y_{0})$, and 
$$\sum _{i}\vert \ell _{i}\vert \geq L_{2}.$$
Let $w_{m}=[\xi ]$ be the  end of $\ell _{m}$ nearer to $x(y_{0})$, 
so 
that 
$d_{\alpha _{n}}'(w_{m},x(y_{0}))\leq L_{0}$, for $L_{0}$ depending 
only on the ltd parameter functions.  Let $\gamma \subset \alpha 
_{m}$ 
be a loop such that $\vert \xi (\gamma )\vert $ is bounded. 
Then by \ref{7.5}, for  a constant $C_{5}$ bounded in terms of $\vert 
\xi (\gamma )\vert $ and  the ltd 
parameter functions,  
$$i(\gamma ,\zeta _{k,+})\leq C_{5}\vert \xi (\zeta _{k,+}\cap \alpha 
_{m}\vert \leq C_{5}^{2}e^{-L_{2}/2}\vert \varphi _{0}(\zeta  
_{n,+})\vert .$$
So if this is true for arbitrarily large $k$,
$$i(\gamma ,\mu _{+}')\leq e^{-L_{2}/4}.$$
But similarly 
$$i(\gamma ,\mu _{-}')\leq e^{-L_{2}/4}.$$
This contradicts our assumption on $\mu _{\pm}'$. So 
now we have (\ref{6.13.1}) for all sufficiently large $k$. 

So now we need to show (\ref{6.13.2}), using a generalisation of the 
technique used to prove (\ref{6.12.11}). So let $k$ be sufficiently 
large that (\ref{6.13.1}) holds, and let $(\alpha ,\ell )$ be ltd 
for $[y_{k,-},y_{k,+}]$ with $d_{\alpha }'(y_{0},y)\leq \Delta $ for 
$y\in \ell $. We again use 
\ref{7.4}, but this time the quadrilateral case, for the 
quadrilateral 
with vertices at $y_{k,\pm }$ and $y_{n,\pm }$. So fix $\Delta >0$, 
and let $(\alpha ,\ell )$ be long $\nu $-thick and dominant or 
$K_{0}$-flat (having fixed ltd parameter functions) along 
$[y_{k,-},y_{k,+}]$ within $d_{\alpha }'$ distance $\Delta $ of 
$y_{0}$. By assumption, $\alpha $ does not intersect any loop $\zeta 
$ with $i(\mu _{+},\zeta )=0$, or $i(\mu _{-},\zeta )=0$. Let 
$y=[\varphi ]\in \ell $. Suppose that (\ref{6.13.2}) 
does not hold for $y'\in [y_{n,-},y_{n,+}]$ for $n$ sufficiently 
large. By \ref{7.4}, there must then be 
$y'=[\varphi ']\in [y_{k,+},y_{n,+}]\cup[y_{k,-},y_{n,-}]$. We assume 
without 
loss 
of generality that
$y'\in [y_{k,+},y_{n,+}]$. Precisely, we have
$$d_{\alpha }(y,y')\leq C(\nu ){\rm{\ or\ }}\vert {\rm{Re}}(\pi 
_{\alpha }(y)-\pi _{\alpha }(y')\vert \leq K_{0}.$$
Suppose that $\alpha $ is a gap which is long $\nu $-thick and 
dominant. Fix a loop $\gamma \subset \alpha $ with $\vert \varphi 
_{0}(\gamma )\vert \leq L_{0}$, with $L_{0}$ depending only on 
the topological type of $S$. Now by \ref{7.2} since $\varphi '(\zeta 
_{n,+}\cap \alpha )$ and 
$\varphi '(\zeta _{k,+}\cap \alpha )$ are close to the stable and 
unstable foliations of the quadratic differential for $d(y',y_{n,+})$ 
respectively, for a constant $C_{2}(\nu )$,
$$i(\zeta _{n,+},\zeta _{k,+})\geq \# (\varphi (\alpha \cap \zeta 
_{k,+}\cap \zeta _{n,+}))\geq (C_{2}(\nu ))^{-1}i(\zeta _{n,+},\gamma 
).i(\gamma ,\zeta _{k,+}).$$
If this is true for arbitrarily large $k$, and $n$, taking limits, 
this 
means that
$$0=i(\mu _{+}',\mu_{+}')\geq C_{2}(\nu )^{-1}i(\mu _{+}',\gamma 
)^{2}>0,$$
which is a contradiction. So now we have (\ref{6.13.2}),
and the proof of geometric 
convergence of the models is finished in this general case.

\section{ Model-adapted families of pleated surfaces. }\label{8}

Let $N$ be a three-dimensional hyperbolic manifold with finitely 
generated fundmantal group. Let $\overline{N}$ 
be the union of $N$ and the quotients by the covering group of the 
complement of the limit set in $\partial H^{3}$. Let $W$ be a 
submanifold of the relative 
Scott core $N_{c}$ as in \ref{6.9}, and $N_{d,W}$ as 
defined there. Thus, each component of $N_{d,W}\setminus W$ is a 
neighbourhood of a unique end $e$ of $N_{d,W}$, and  the closure 
is homeomorphic to 
$S_{d}(e)\times [0,\infty )$, where $S_{d}(e)$ is the bounding 
component of $\partial W\setminus \partial N_{d}$. In this section, we
 construct a family of  pleated 
surfaces in $N$,  given a map $f_{e,+}:S(e)\to 
\overline{N}$ homotopic to inclusion of $S(e)$ in $N$, for each end 
$e$ of 
$N_{d,W}$. Here, $f_{e,+}$ is either a pleated surface, or a map to 
$\overline{N}\setminus N$, or a mixture of both.
This family of pleated surfaces is made up of a sequence of pleated 
surfaces for each end, and a family of pleated surfaces for $W$. The 
family of pleated surfaces for $W$ is a family of pleated surfaces 
for the noninterval bundle part of $W$ and a sequence for each 
interval bundle in $W$, using the decomposition of \ref{6.11}. The 
family for the noninterval bundle part of $W$ is actually independent 
of the choice of ending pleated surfaces. We shall prove 
that the geometry of this family on the noninterval bundle part  
depends only 
on the topological type of $(W,N)$, in the case when all ends are 
incompressible, and on the topological type and a constant $c_{0}$ 
if some ends are compressible. This result was proved by Thurston \cite{T3} in the case of incompressible boundary, with comments on what was needed to extend  to the case of compressible boundary. The proof of ``bounded window frames'' given here, in the case of incompressible boundary, is different from that in \cite{T3}.
The general result is mostly proved in \ref{8.12}, with a key hypothesis left to be proved in Section \ref{10}.

For the case of $N$ combinatorial bounded geometry, 
it is only necessary to read to the end of \ref{8.3}, and for the 
case of  $N$ being an interval bundle, to the end of \ref{8.4}.

\ssubsection{Sequences of multicurves and pleated 
surfaces with particular properties.}\label{8.2}

First, we consider sequences of maximal multicurves and 
pleated 
surfaces with certain properties. Suppose that $[z_{0},y_{+}]\subset 
{\cal{T}}(S)$. We suppose that $S$ is embedded in $N$
For each such $[z_{0},y_{+}]$, we shall choose an increasing 
sequence $\{ z_{i}=[\varphi _{i}]:0\leq i\leq n\} $ of points in 
$[z_{0},y_{+}]$ with $y_{+}=z_{n}$, with $d(z_{i},z_{i+1})\leq 1$ 
and a sequence of maximal  multicurves 
$\Gamma _{j}$ on $S$ such that the following holds for a suitable 
constant $\kappa _{0}$, and an integer $r_{0}$. 

  \subsubsection{}\label{8.2.1}  $\Gamma _{j}$ is a noncollapsing 
  maximal multicurve. 

 \subsubsection{}\label{8.2.2}  $\vert 
 \gamma \vert \geq \varepsilon _{0}$ for all $\gamma \in \Gamma 
_{j}$, unless $\gamma \in \cap _{k=0}^{n}\Gamma _{k}$.  
   
\subsubsection{}\label{8.2.3} Either $\# (\Gamma _{j}\cap \Gamma 
_{j+1})\leq r_{0}$, or 
there is $\zeta _{j}\in \Gamma _{j}$, and a loop $\gamma _{j}$
such that $i(\zeta ,\gamma _{j})=0$ for $\zeta \neq \zeta _{j}$, 
$0<i(\zeta _{j},\gamma _{j})\leq r_{0}$, $\Gamma _{j}\cup \{ \gamma 
_{j}\} \setminus \{ \zeta _{j}\} $ and $\Gamma _{j}\cup \{ \tau 
_{\gamma _{j}}(\zeta _{j})\} \setminus \{ \zeta _{j}\} $ are 
noncollapsing, and 
$\Gamma _{j+1}=(\Gamma _{j}\setminus  \{ 
\zeta 
    _{j}\} )\cup \{ \tau _{\gamma _{j}}^{\pm n_{j}}(\zeta _{j}\} $.
    
 \subsubsection{}\label{8.2.4}  
  
$$\sum _{p =0}^{n-1}\log (\# (\Gamma _{p}\cap \Gamma _{p+1})+1)\leq 
\kappa _{0}d(z_{0},z_{n}).$$

\vskip 1 true cm

Then for each $j$, we let $f_{j}:S\to N$ be a pleating surface whose 
pleating 
locus includes $\Gamma _{j}$ and homotopic to the embedding of $S$ in 
$N$.  The extra property which we shall 
require is the following, which is automatic in the case when $S$ 
is incompressible in $N$.

\subsubsection{}\label{8.2.5} For $f=f_{j}$ or $f_{j+1}$ (not 
necessarily both),  whenever  $\vert f(\gamma )\vert 
\leq D_{0}$ for $\gamma $ nontrivial  in $S$, then 
$f(\gamma )$ is nontrivial in $N$.
\vskip 1 true cm

\subsubsection{Consequences for the  pleated 
surfaces.}\label{8.2.6}

\noindent {\em{Bounded distance in Teichm\"uller space.}}
If we do have \ref{8.2.1} to \ref{8.2.5}, 
then by \ref{8.2.2}, \ref{4.4},\ref{4.7}, for a constant $\kappa 
_{0}'$, depending 
only on 
depending the topological type of 
$S$,and a constant $\kappa _{1}$ 
depending 
only on 
depending only on $r_{0}$ and on the topological type of 
$S$,
$$d([f_{j}],[f_{j+1}])\leq \kappa _{0}'(\log (\# (\Gamma _{j}\cap 
\Gamma _{j+1})+1),$$
$$\sum _{j=0}^{n-1}d([f_{j}],[f_{j+1}])\leq \kappa 
_{1}d(z_{0},z_{n}).$$

\noindent {\em{Bounded distance between impressions}}
Also, by \ref{4.3}, there is a homotopy in $N$ 
between the impressions of $f_{j}$ and $f_{j+1}$ with homotopy tracks 
of length $\leq L_{0}$ in the case when $i(\Gamma _{j},\Gamma 
_{j+1})\leq r_{0}$, and also in the case
  $\Gamma _{j+1}=(\Gamma _{j}\setminus  \{ 
\zeta _{j}\} )\cup \{ \sigma _{\gamma _{j}}^{n_{j}}(\zeta _{j}\} $, 
when  
$\vert (\gamma 
_{j})_{*}\vert $ is bounded from $0$. In general in this case, 
 we can interpolate 
pleated 
surfaces $f_{j,k}$ between $f_{j}$ and $f_{j+1}$, $0\leq k\leq 
n_{j}$, 
$f_{j,0}=f_{j}$, $f_{j,n_{j}}=f_{j+1}$, $f_{j,k}$ has pleating 
locus     
$(\Gamma _{j}\setminus  \{ \zeta 
    _{j}\} )\cup \{ \tau _{\gamma _{i}}^{k}(\zeta _{i})\} $,
and there is a homotopy in $N$ between $f_{j,k}$ and $f_{j,k+1}$ 
whose homotopy tracks have length $\leq L_{0}$.  The hypotheses of 
\ref{8.2.3}) 
ensure that all the pleating loci are noncollapsing

\subsection{Sequence of pleated surfaces: combinatorial 
bounded geometry Kleinian surface case}\label{8.3}

We show that we can find maximal multicurves $\Gamma _{i}$ satisfying 
\ref{8.2.1}-\ref{8.2.4} with certain properties in the  case of a 
hyperbolic manifold $N$ 
with $N_{c}$ homeomorphic 
to $S_{d}\times [0,1]$ and with ending invariants $\mu _{\pm}$ of 
combinatorial bounded geometry, that is, satisfying (\ref{6.12.1}). 
We 
also assume for the moment that both ends are geometrically infinite, 
that is, that $\mu _{\pm }\in {\cal{GL}}_{a}(S)$.
We identify $S_{d}$ with a component of $\partial N_{c}$, so that 
$S_{d}$ 
is embedded in $N_{d}$.

In this case, since $S_{d}$ is incompressible in 
$N_{d}$, 
\ref{8.2.5} and \ref{8.2.6} are 
automatically satisfied. We start by choosing homotopic  pleated 
surfaces 
$f_{\pm }:S\to N$ as in \ref{3.8},  with $[f_{\pm 
}]=y_{\pm }$ close to $\mu _{\pm}$ and loop sets $\Gamma _{\pm }$ 
such that 
\begin{equation}\label{8.3.0}
    \vert f_{+}(\Gamma _{+})\vert \leq L_{0},\ \ \vert f_{-}(\Gamma 
_{-})\vert \leq L_{0}.\end{equation}
 Then we can choose 
$\zeta _{\pm }\subset \Gamma _{\pm }$ similarly to $\zeta _{n,\pm }$ 
in \ref{6.12}, and then use these to define $y_{\pm }'=[\varphi 
_{\pm }']\in [y_{-},y_{+}]$  as in \ref{6.12}, with
\begin{equation}\label{8.3.1}[y_{-}',y_{+}']\subset 
({\cal{T}}(S))_{\geq \nu },\end{equation}
and maximal multicurves $\Gamma _{\pm }'$ such 
that $\vert \varphi _{+}'(\Gamma _{+}')\vert $ is bounded, and 
similarly for $\vert \varphi _{-}'(\Gamma _{-}')\vert $, where 
$\Gamma _{+}\cap \Gamma _{+}'$ contains at least one loop, as does 
$\Gamma _{-}\cap \Gamma _{-}'$. We then take $f_{+}':S\to N$ to be a 
pleated surface with pleating locus $\Gamma _{+}'$, and similarly for 
$f_{-}'$.
We can extend to cover the case of $\mu _{\pm }\in 
{\cal{T}}(S)\cup {\cal{GL}}_{a}(S)$ by choosing $f_{+}$ as in 
\ref{3.9} if $\mu _{+}\in {\cal{T}}(S)$, and similarly if 
$\mu _{-}\in {\cal{T}}(S)$.

Now  we shall choose a sequence which satisfies 
\ref{8.2.1}-\ref{8.2.6} for $r_{0}$ 
 and $\kappa _{0}$ depending only on $\nu $, with $z_{0}=y_{-}'$ and 
$y_{+}=y_{+}'$. We choose points 
$z_{i}=[\varphi _{i}]\in [y_{-}',y_{+}']$, and maximal multicurves 
$\Gamma _{i}'$, $0\leq i\leq 
n$, such that $z_{0}=y_{-}'$, $z_{n}=y_{+}'$, $\Gamma 
_{0}'=\Gamma _{-}'$, $\Gamma _{n}'=\Gamma _{+}'$, $d(z_{i},z_{i+1})$ 
is 
bounded above by $1$ and bounded below by ${1\over 2}$ and, for a 
constant 
$L_{0}=L_{0}(\nu )$ (enlarging the previous $L_{0}$ if necessary)
$$\vert \varphi _{i}(\Gamma 
_{i}')\vert \leq L_{0}.$$
Then 
the bound on $d(z_{i},z_{i+1})$ means that for a constant $C_{1}$ 
depending only on $\nu $ and the topological type of $S$, 
$$\# (\Gamma 
_{i}'\cap \Gamma _{j}')\leq C_{1}e^{C_{1}\vert i-j\vert }.$$
Now by \ref{4.2} we can replace each $\Gamma _{i}'$ by a loop set 
$\Gamma _{i}$ such that, enlarging $L_{0}(\nu )$ if necessary,
\begin{equation}\label{8.3.2}\vert \varphi _{i}(\Gamma _{i})\vert 
\leq 
L_{0},\end{equation}
and for all $\gamma \in \Gamma _{i}$,
\begin{equation}\label{8.3.3}\vert \gamma _{*}\vert \geq \varepsilon 
_{0}.\end{equation}
(\ref{8.3.2}) implies a bound on $\# (\Gamma _{i}\cap \Gamma _{i}')$, 
and on 
$\# (\Gamma _{i}\cap \Gamma _{j})$, when $\vert i-j\vert $ is bounded.
We then take $f_{j}$ to be a pleated surface with pleating locus 
including $\Gamma _{j}$.

Then conditions \ref{8.2.1} to \ref{8.2.6} are  satisfied. In 
condition 
\ref{8.2.3}, only the 
first alternative, $\# (\Gamma _{j}\cap\Gamma _{j+1})\leq r_{0}$, 
holds.
Also, by \ref{4.3}, \ref{4.4} applied to $f_{n_{+}}$ and $f_{+}$, 
$\vert f_{n_{+}}(\Gamma _{+}')\vert 
\leq L_{0}(\nu )$, assuming that $L_{0}(\nu )$ is large enough given 
$r_{0}(\nu )$, and similarly for $f_{n_{-}}$ and $\Gamma _{-}'$. 
Since $\Gamma _{+}\cap \Gamma _{+}'\neq \emptyset $, and assuming 
$L_{0}(\nu )$ is large enough given $\nu $, it follows that 
 $f_{n_{+}}(S)$ either comes within a distance 
$L_{0}(\nu )$ of $f_{+}(S)$, or of a Margulis tube 
intersected by $f_{+}(S)$. 

\ssubsection{Sequence of pleated surfaces: interval 
bundle case.}\label{8.4}
First we consider the case when $\partial N_{c}$ is homeomorphic to 
an interval 
bundle $S_{d}\times [0,1]$ and $N_{d}$ has just two ends $e_{\pm }$. 
If $e_{+}$ is geometrically infinite, we choose $\Gamma _{+}$, 
$[f_{+}]$ as in \ref{3.8}, with 
\begin{equation}\label{8.4.0}\vert f_{+}(\Gamma _{+})\vert \leq 
L_{0}.\end{equation}
If $e_{+}$ is geometrically finite, we 
choose $\Gamma _{+}$ as $\Gamma $ in \ref{3.9}, and $f_{+}$ as 
$f_{2}$ of \ref{3.9}. So for suitable $L_{0}$, once again, 
(\ref{8.4.0}) holds. We choose $f_{-}$ similarly.  
So $f_{\pm }$ are injective on $\pi _{1}$. 
We  shall now find a sequence of maximal multicurves $\Gamma _{i}$,
$0\leq i\leq n$, satisfying conditions 
\ref{8.2.1} to 
\ref{8.2.4} {\em{corresponding to}} 
$[y_{-},y_{+}]$. This means that, enlarging $L_{0}$ if necessary,
 depending only 
on topological type, we 
shall have
\begin{equation}\label{8.4.1}
    d([f_{0}],[f_{-}])\leq L_{0},\ \ d([f_{n}],[f_{+}])\leq 
    L_{0}.\end{equation}

In later sections we shall also use this method to 
define a sequence of maximal multicurves associated to a geodesic 
segment $[z_{0},y_{+}]\subset {\cal{T}}(S)$, where $y_{+}$ will be a 
pleated surface but $z_{0}$ will in general not be. There will then 
be a bound on $d([f_{n}],[f_{+}])=d([f_{n}],y_{+})$, but in general 
there will be no a priori bound on $d(z_{0},[f_{0}])$ (although such 
a bound will be obtained eventually). We can then 
choose 
a sequence of pleated surfaces $f_{i}:S\to N$ obtained 
from $\Gamma _{i}$ so 
that \ref{8.2.5} and \ref{8.2.6} are satisfied. 

We proceed as follows. We assume that ltd parameter 
functions 
have been fixed so that all relevant results in Sections \ref{5} and 
\ref{6} work. Write $y_{-}=z_{0}$. We fix a vertically efficient 
deomposition of $S\times [z_{0},y_{+}]$ into sets $\alpha _{j}\times 
\ell 
_{j}$. 
For each $j$ such that $\alpha _{j}$ is a gap, we choose points 
$z_{i,j}$ and loop sets $\Gamma 
_{i,j}$ on $S(\alpha _{j})$, $0\leq i\leq n(j)$, exactly as in 
\ref{8.3}, with $S(\alpha _{j})$ replacing $S$. In condition 
\ref{8.2.3} the first of the two alternatives becomes 
$$\# (\Gamma 
_{i,j}\cap \Gamma _{i,j+1})\leq r_{0}(\nu _{j}),$$
amd we can ensure that, if $z_{i,j}=[\varphi _{i,j}]$, then for a 
function $L_{0}(\nu )$ independent of the ltd parameter functions,
$$\vert \varphi _{i,j}(\Gamma _{i,j})\vert \leq L_{0}(\nu _{j}).$$
The vertically efficient 
conditions ensure that, for each $j$, there is 
$\nu _{j}>0$ bounded below in terms of the ltd parameter functions 
such that $\vert \varphi (\gamma )\vert \geq \nu _{j}$ for all loops 
$\gamma \subset \alpha _{j}$ which are nontrivial nonperipheral and 
not 
boundary-homotopic, and $[\varphi ]\in \ell $. If $(\alpha _{j},\ell 
_{j})$ is ltd, then we have in addition that $\alpha _{j}$ is long 
$\nu 
_{j}$-thick and dominant along $\ell _{j}$. In any case, we have
$$\ell _{j}\subset ({\cal{T}}(S(\alpha _{j})))_{\geq \nu _{j}},$$
so that the same argument as in \ref{8.3} works. 
So then we need to use these pieces to produce the points $z_{i}$ and 
multicurves $\Gamma _{i}$. First, we choose a sequence 
$z_{i}^{1}=[\varphi 
_{i}^{1}]$ of successive points on $[y_{-},y_{+}]$, such 
that $d(z_{i}^{1},z_{i+1}^{1})\leq L_{0}'$, for $L_{0}'$ depending 
only 
on the topological type, and where $d_{\alpha 
_{j}}(z_{i}^{1},z_{p,j})\leq L_{0}(\nu _{j})$
 for some $p$, if $z_{i}^{1}\in \ell _{j}$ such that $\alpha _{j}$ is 
a gap, and enlarging $L_{0}(\nu )$ if necessary, but 
still independent of the parameter functions.
 We then take $\Gamma _{i}^{1}$ to be the union of all the $\Gamma 
_{p,j}$ and $\partial \alpha _{j}$ such that $\alpha _{j}$ is a gap, 
$z_{i}^{1}\in \ell _{j}$ and $\pi 
 _{\alpha _{j}}(z_{i}^{1})=z_{p,j}$. We can also assume that 
$$\vert \varphi _{i}^{1}(\Gamma _{i}^{1})\vert \leq L_{0}'',$$
where $L_{0}''$ depends only on topological type of $S$ and on 
the ltd parameter functions. 

The 
properties \ref{8.2.1} to \ref{8.2.4} probably do not hold for this 
choice of $\{ 
z_{i}^{1}\} $ and $\{ \Gamma _{i}^{1}\} $. In particular, \ref{8.2.2} 
probably does not hold, because it is likely that $\vert (\partial 
\alpha 
_{j})_{*}\vert <\varepsilon _{0}$, and that there there is a badly 
bent annulus with core mapping to $(\partial \alpha 
_{j})_{*}$ for many of the pleated surfaces in the sequence. These 
are, however, the only loops in the pleating loci which can be too 
short, because all others have been removed. If we 
rectify this, then we need to be careful about \ref{8.2.3} and 
\ref{8.2.4}. At any rate, any modifications necessary concern the 
loops $\partial \alpha _{j}$. Our new sequence will be obtained 
from $\Gamma _{i}^{1}$ by replacing the loops $\partial \alpha _{j}$ 
only. So fix a
 loop $\gamma $ which is a component of  $\partial \alpha _{j}$ for 
at 
least one $\alpha _{j}\times \ell _{j}$. 
Then the $i$ for which $\Gamma _{p,j}\subset \Gamma _{i}^{1}$ for 
some $j$, $p$ with $\gamma \subset \partial \alpha _{j}$ form an 
interval $\{ i:m_{1}(\gamma )\leq i\leq m_{2}(\gamma )\} $. Let 
$n_{\gamma }(z)=n_{\gamma ,\beta (\gamma )}(z)$ (as in \ref{2.6}) be 
defined relative to some fixed 
loop $\beta (\gamma )$ with one or two transverse intersections 
with $\gamma $. Then we are going to choose $z_{i}^{2}=[\varphi 
^{2}]=[\varphi ^{1}\circ \tau ]$. Here, $\tau $ is a composition of 
Dehn twists $\tau _{\gamma }^{p_{i,\gamma }}$ for a set of disjoint 
loops $\gamma $ and integers 
$p_{i,\gamma }$. We choose $p_{i}=p_{i,\gamma }$ so that 
$n_{\gamma }(z_{i}^{2})$ is monotone in 
$i$ for $m_{1}(\gamma )\leq i\leq m_{2}(\gamma )$, 
$z_{i}^{2}=z_{i}^{1}$ for $i=m_{1}(\gamma )$, $m_{2}(\gamma )$, and 
$\vert n_{\gamma }(z_{i}^{2})-n_{\gamma }(z_{i+1}^{2})\vert \leq 
r_{0}$ 
except for at most one $i$. This then gives
\begin{equation}\label{8.4.2}\sum _{i=m_{1}(\gamma )}^{m_{2}(\gamma 
)-1}
\vert \log (\vert n_{\gamma }(z_{i}^{2})-n_{\gamma 
}(z_{i+1}^{2})\vert 
+1)\vert \leq \kappa _{2}d(z_{m_{1}(\gamma )}^{2},z_{m_{2 
}(\gamma )}^{2}).\end{equation}
So to choose $\Gamma _{i}^{2}$, choose any loop $\zeta _{i}$ such 
that $\zeta _{i}$ has one or  two transverse intersections with 
$\gamma $, and 
$$\vert \varphi _{i}^{1}(\zeta _{i})\vert ''\leq L_{0}'$$
where $\vert .\vert ''$ is as in \ref{2.3}, and $L_{0}'$, as before, 
depends only on the topological type of $S$. It is possible to make 
such a choice.  Then by choice of $p_{i}$, we can ensure the 
conditions on $n_{\gamma }(z_{i}^{2})$. Since $n_{\gamma }([\varphi 
\circ \tau _{\gamma }^{m}],\beta )=n_{\gamma }([\varphi ],\beta )+m$ 
for any choice of $\beta $, the condition on the $n_{\gamma 
}(z_{i}^{2})$ translates to: $\vert p_{i}- 
p_{i+1}\vert \leq r_{0}$ for all but at most one $i$, and 
(\ref{8.4.2}) holds 
with $p_{i}$ replacing $n_{\gamma }(z_{i}^{2})$. This then gives 
\ref{8.2.3} and \ref{8.2.4}. To get \ref{8.2.2}, note that by 
\ref{4.1}, the bound $\vert p_{i}-p_{i+1}\vert \leq r_{0}$ is 
compatible 
with $\vert (\tau _{\gamma }^{p_{i}}(\zeta _{i}))_{*}\vert \geq 
\varepsilon _{0}$ 

We now consider (\ref{8.4.1}), which involves twists round loops of 
$\Gamma _{\pm }$. 
We might not take $f_{0}=f_{-}$, $f_{n}=f_{+}$, because we want 
\ref{8.2.2} satisfied for $\Gamma _{0}$, $\Gamma _{n}$. 
By \ref{4.4} we shall have (\ref{8.4.1}) for $f_{0}$ if we choose 
$\Gamma _{0}$ 
so that there is no badly bent annulus for $(f_{0},f_{-})$. 
Fix any 
loop $\gamma \in \Gamma _{-}\setminus \Gamma _{+}$ with 
$\vert \gamma _{*}\vert =
\vert f_{-}(\gamma )\vert <\varepsilon 
_{0}$ and fix $\beta =\beta (\gamma )$ with one or two intersections 
with 
$\gamma $, 
and no other intersections with $\Gamma $. By the method used in 
\ref{4.1}, the set of $k$ for which $(\tau _{\gamma 
}^{k}(\beta ))_{*}$ does not intersect a bounded neighbourhood of 
$\gamma _{*}$ lie in an interval $I$ of integers of length $\leq 
1/\vert 
\gamma _{*}\vert $. So we choose $\Gamma _{0}$  to include $\tau 
_{\gamma 
}^{k}(\beta )$ for some $k$ outside this range, and if possible so 
that $\vert f_{0}(\tau _{\gamma 
}^{k}(\beta ))\vert ''$ is bounded. Then we proceed as for any 
$\gamma $ above. Then (\ref{8.4.2}) will hold 
for $\gamma $,  because, for this $\gamma $, 
$d(z_{m_{1}(\gamma )}^{2},z_{m_{2 
}(\gamma )}^{2})$ is bounded below by a multiple of $-\log 
\vert \gamma _{*}\vert $--- unless $m_{2}(\gamma )=n$. 

If both $e_{\pm }$ are geometrically finite, and $[f_{\pm 
}]=\mu (e_{\pm })$ 
are the corresponding invariants, and 
there is a loop $\gamma $ such that $\vert f_{+}(\gamma )\vert \leq 
\varepsilon _{0}$, $\vert f_{-}(\gamma )\vert \leq 
\varepsilon _{0}$, then   we can 
choose the  loop sets $\Gamma _{0}$, $\Gamma _{n}$ such that, if 
$\beta _{\pm }$ are elements of $\Gamma _{0}$, $\Gamma _{n}$ 
intersecting 
$\gamma $, then 
$$n_{\gamma ,\beta _{+}}([f_{+}])-n_{\gamma ,\beta _{-}}([f_{-}])\leq 
L_{2}\vert \log \vert f_{+}(\gamma )\vert -\vert f_{-}(\gamma )\vert 
\vert .$$
We see this from the proof of \ref{3.9}. The loops $\beta _{\pm 
}(\gamma )$ 
simply had to be adjusted  by Dehn twists which excluded intervals of 
integers of lengths $O(1/\vert f_{\pm }(\gamma )\vert )$. So if 
adjustment is necessary, we can choose the two adjustments on the 
same 
side of the smaller excluded interval.
Then \ref{8.4.2} still holds.

Now we consider modifications for general interval bundles, the case 
considered in \ref{6.9}. In this case, $f_{\pm }:S\to \overline{N}$ 
are generalised pleated surfaces in the sense of \ref{4.5}.  We again 
choose maximal multicurves $\Gamma _{\pm }$ and $f_{\pm }$ so that 
all closed loops in the pleating locus of $f_{+}$ are in 
$\Gamma _{+}$, and similarly for $f_{-}$, $\Gamma _{-}$, and so that 
(\ref{8.4.0}) is satisfied. Now $[f_{+}]$ is an element of 
${\cal{T}}(S(\alpha _{+}))$, where $\alpha _{+}=S\setminus (\cup 
\Gamma _{+})$, and similarly for $f_{-}$, $\Gamma _{-}$.
In order to completely determine 
$y_{\pm }=[f_{\pm }]\in {\cal{T}}(S)$, we also need to define
$\pi _{\gamma }(y_{+})$ for any $\gamma \in \Gamma _{+}$ 
and  
$\pi _{\gamma }(y_{-})$ for any $\gamma \in \Gamma _{-}$. It is 
convenient to make a choice which minimises $d(y_{-},y_{+})$ up to an 
additive constant. Having defined $y_{\pm }$, the rest of the 
construction is exactly as before.

\ssubsection{General case: sequence of pleated surfaces corresponding 
to ends, and to bridges between incompressible ends.}\label{8.5}
 
let $N$, $W$, $N_{d,W}$ $\overline{N}$ be as in the introduction 
to Section \ref{8}.
Let $e$ be an end of $N_{d,W}$.
The sequence for the end $e$ is determined by $[f_{e,+}]=y_{e,+}$, 
$\Gamma _{+}(e)$, and by 
another element $z_{e,0}=[\varphi _{e,0}]\in {\cal{T}}(S(e))$, which 
has to be 
determined.  It is determined as in \ref{6.16} and \ref{6.17}, 
depending on whether $e$ is compressible or incompressible.  The map
$f_{e,+}:S\to \overline{N}$, is a generalised pleated surface as in 
\ref{8.4}, and there is a maximal multicurve
$\Gamma _{+}(e)$ which includes all the closed loops in the pleating 
locus of $f_{e,+}$, and such that
$$\vert f_{e,+}(\Gamma _{+}(e))\vert \leq L_{0}.$$
If $f_{e,+}$ maps some nonperipheral loops on $S$ to 
cusps in $N$, then, similarly to \ref{8.4}, we define 
$y_{e,+}=[f_{e,+}]$ 
as an element of ${\cal{T}}(S(e))$ for each $\gamma \in \Gamma 
_{+}(e)$.
Given a choice of $z_{e,0}$,
we define $\pi _{\gamma }(y_{e,+})$  for $(\gamma \in \Gamma 
_{+}(e)$ so 
as to minimise $d(z_{e,0},y_{e,+})$, up to a bounded additive 
constant.

Once we have fixed $y_{e,+}$ and $z_{e,0}$, as in \ref{8.4}, if $e$ 
is 
incompressible, we can 
choose a sequence of maximal multicurves 
$\Gamma _{e,j}$, $0\leq j\leq n_{+}(e)$ and  a sequence of pleated 
surfaces $f_{e,j}$ with pleating locus containing $\Gamma _{e,+}$ so 
that \ref{8.2.1} to \ref{8.2.6} are satisfied. In analogy to 
\ref{8.4.1} 
we shall require
$$d([f_{e,n}],y_{e,+})\leq L_{0}.$$
But it is important to note that we do {\em{not}} attempt to bound 
$d([f_{e,0}],z_{e,0})$, nor $\vert (\Gamma _{e,0})_{*}\vert $. 
Instead, we choose $\Gamma _{e,0}$ so that
$$\vert \varphi _{e,0}(\Gamma _{e,0})\vert ''\leq L_{0},$$
where $L_{0}$ depends only on the topological type. We can usually 
choose $z_{e,0}=[\varphi _{e,0}]\in ({\cal{T}}(S(e))_{\geq 
\varepsilon 
_{0}}$ and hence use $\vert .\vert $ rather than $\vert .\vert ''$. 
In fact, we can always do this if there are no interval bundle 
bridges between incompressible ends. Even in the case where there are 
such bridges, We can choose $z_{e,0}$, consistent with 
the choice of $z_{e,0}$ in \ref{6.17} so that
 $z_{e,0}=[\varphi _{e,0}]\in ({\cal{T}}(S(e))_{\geq \varepsilon 
_{0}}$, unless some loop $\gamma $ in the common subsurface $\omega 
(e,e')=\omega (e',e)$ is such that both $\vert \varphi _{e,0}(\gamma 
)\vert <\varepsilon _{0}$ and $\vert \varphi _{e',0}(\gamma 
)\vert <\varepsilon _{0}$.

The choice of $z_{e,0}$, $z_{e',0}$ was made in \ref{6.17} so that 
$\pi _{\omega }(z_{e,0})=\pi _{\omega }(z_{e',0})$ whenever $e$, $e'$ 
are incompressible ends with $\omega =\omega (e,e')\neq \emptyset $. 
However, we also defined $z_{e,0}'$ and $z_{e',0}'$, essentially 
with $[\pi _{\omega }(z_{e,0}'),\pi _{\omega }(z_{e',0}')]$ as large 
as possible. We shall later make use of a sequence of maximal 
multicurves and pleated surfaces 
for $[\pi _{\omega }(z_{e,0}'),\pi _{\omega }(z_{e',0}')]$. that is, 
for the model manifold $M(\pi _{\omega }(z_{e,0}'),\pi _{\omega 
}(z_{e',0}'))$.We call these sequences $\Gamma _{i}(e,e')$,
$f_{e,e',i}$  ($0\leq i\leq n(e,e')$), consisting  of  maximal 
multicurves 
$\Gamma _{i}(e,e')$ on $\omega 
(e,e')$, $\partial \omega (e,e')\subset \Gamma _{i}(e,e')$ 
and corresponding pleated surfaces $f_{e,e',i}$
with domain $\omega (e,e')$. 
constructed from the interval  We can construct such a sequence as 
before, 
satisying \ref{8.2.1} to \ref{8.2.4}. We also choose $\Gamma 
_{e,0}$, $\Gamma _{e',0}$ so that 
$$\Gamma _{e,0}\cap \omega (e,e')=\Gamma _{e',0}\cap \omega (e,e'),$$
$$f_{e,0}=f_{e',0}{\rm{\ on\ }}\omega (e,e').$$
At the moment, we have no 
upper or lower bound 
on $\vert (\partial \omega (e,e'))_{*}\vert $. So it is not yet 
possible to verify \ref{8.2.6}, even though \ref{8.2.5} is automatic.
We shall see later in this section that 
$\vert (\partial \omega (e,e'))_{*}\vert $ is bounded above. 
But it may not be bounded below. So since the sequence $\Gamma 
_{e,i}$ 
is chosen so that $\vert \gamma _{*}\vert \geq 
\varepsilon _{0}$ for all $\gamma \in \Gamma _{i}(e)$, it may not be 
possible to make $\Gamma _{n(e,e')}(e,e')\subset \Gamma _{0}(e)$. But 
we 
can ensure that for a suitable $r_{0}$,
$$\# (\Gamma _{n(e,e')}(e,e')\cap \Gamma _{0}(e))\leq r_{0}.$$

If $e$ is a compressible end, then we can try to carry out the same 
construction of sequences of multicurves and pleated sequences as for 
an incompressible end. But it is not clear that the multicurves are 
noncollapsing (\ref{8.2.1}). Nor is it clear that \ref{8.2.5} to 
\ref{8.2.6} are satisfied. For this, Lemmas \ref{8.7} and \ref{8.8} 
below are relevant. Lemmas \ref{8.11} and \ref{8.9} will be used to 
construct the 
sequence of pleated surfaces for the core.

\begin{lemma}\label{8.7} Let $N$ be a hyperbolic $3$-manifold with 
finitely generated fundamental group. Let $S$ be a finite type 
surface, possibly with boundary, with $(S_{d},\partial S_{d})$ 
embedded in $(N_{d},\partial 
N_{d})$, and bounding an essential submanifold $W$ of $N$. Fix a Margulis constant $\varepsilon _{0}$ and integer 
$r$ and a constant $L_{0}$. There is 
$L=L(r,(S,N), \varepsilon _{0},L_{0})$ such 
that the following 
holds.  Let
$\Gamma _{0}\subset S$ be a maximal multicurve. Let 
$f_{0}:S\to N$ be a pleated surface homotopic to inclusion, with 
pleating locus $\Gamma _{0}$,  and such that any nontrivial 
nonperipheral component of $\partial S$ is mapped to a geodesic of 
length $\leq L_{0}$. Let $A$ be the union of badly bent annuli for 
$f_{0}$.
Suppose that $\gamma '\subset S$ is a simple closed  nontrivial  
loop  such that $\gamma '$ bounds a disc 
 $D_{1}$ in $W$, with interior 
disjoint from $S$ and $\# (\gamma '\cap \Gamma 
_{0})\leq r$. 

Then there is $\gamma $ such 
that $\gamma $ is nontrivial  in 
$S$ and bounds an embedded disc in $W$, and 
$\vert f_{0}(\gamma )\setminus A\vert 
\leq L$ and $f_{0}(\gamma )$ is homotopic to the union in $N$ of 
$f_{0}(\gamma )\setminus A$ and finitely many geodesic segments in 
$N$ of length $\leq L$.

\end{lemma}

\noindent {\em{Proof.}} 
Since $f_{0}(\gamma ')$ has $\leq r$ 
intersections with $f_{0}(\Gamma _{0})$, we can find a union of $\leq 
8r$ geodesic segments in $N$ which is homotopic to $f_{0}(\gamma ')$ 
and a bounded distance from $f_{0}(\gamma ')$. To do this, an arc in 
a 
component of $S\setminus \Gamma _{0}$, with endpoints in $\Gamma 
_{0}$, can be homotoped, keeping endpoints in $\Gamma _{0}$, to 
either 
a single segment in the pleating locus of $f_{0}$, which is contained 
in a geodesic segment homotopic at both ends to loops of $\Gamma 
_{0}$, and two geodesic arcs which can be taken arbitrarily short in 
$S(f_{0})$, or a union of two asymptotice, a segment along a loop of 
$\Gamma _{0}$, and up to four short arcs joing these up and joining 
to 
endpoints in $\Gamma _{0}$. So we can form a loop homotopic 
to $f_{0}(\gamma ')$ which is contained in the union of $\leq 2r$ 
arcs 
of $f_{0}(\Gamma _{0})$, $\leq 2r$ geodesics arcs in the image of the 
pleating locus of $f_{0}$ which are asymptotic to $f_{0}(\Gamma 
_{0})$ at each end, and arbitrarily short arcs joining these. Now 
complete this union of $\leq 8r $ geodesics to a triangulation of 
$D_{1}$. 
So then we have a finite union of $\leq 8r-2$ geodesic triangles in 
$N$. This 
union bounds a disc $D_{1}$, where
 $D_{1}$ is a union of $\leq 8r-2$ 
topological discs of 
diameter $\leq L_{1}$ connected by long thin pieces with two boundary 
components, which we call ``rectangles''. Note that a sequence of 
short arcs across triangles, avoiding the ``thick'' parts of the 
triangles, must have $\leq 24r-6$ arcs between intersections with 
$\partial D_{1}$, because otherwise we can find a  closed 
loop in $D_{1}$ whose image in $N$ can be homotoped to a closed 
geodesic. This would imply the closed loop was nontrivial, which is 
impossible since $S_{1}$ is a disc.   We can choose this 
decomposition into bounded diameter pieces and thin rectangles
so that a rectangle always connects two bounded 
diameter pieces, simply by adding a bounded diameter piece at an end.
We choose one of the bounded diameter topological discs $D_{2}$ which 
has 
only one boundary component. This is possible: in any partition of a 
disc by finitely 
many disjoint arcs, there is at least one complementary component 
with at most one arc in its boundary. If $D_{2}=D_{1}$ we are done. 
If not, there is exactly one rectangle adjoined to $D_{2}$.   If 
there is a 
nontrivial loop in $S(f_{0})$ with  length $\leq L_{1}$ in 
$S(f_{0})\setminus A$
which is trivial in $N$, then we are done. If there is no such loop,
then we apply the Short Bridge 
Arc Lemma  of \ref{3.3.4} to the arc in 
$\partial D_{2}$ which is adjacent to the rectangle. The hypothesis 
of \ref{3.3.4} is satisfied, because the length of $\partial 
D_{2}\cap \partial D_{1}$ is bounded.
Then we can replace the arc in common with the rectangle with an arc 
in 
$S(f_{0})$ of bounded length.  Then either $D_{2}$ is the disc we 
require, 
or if its boundary is trivial in $S(f_{0})$, we can remove it and 
repeat the process with the remainder of the disc. After repeating 
the process a bounded number of times, we have a disc with boundary 
in $S(f_{0})$ of length $\leq L$ for suitable 
$L$. In both cases, identifying $S(f_{0})$ with $S$ embedded in $N$, we can use the Loop Theorem \cite{Hem} by another one, with boundary contained in the boundary of the first length (and hence again with image under $f_{0}$ of length $\leq L$) and embedded in $W$.

\Box

\begin{lemma}\label{8.11} Take the same hypotheses on $N$, $S$, $W$
$f_{0}$, $\Gamma _{0}$ as in \ref{8.7},
but suppose that  every component of the surface $S$, which is not 
assumed to be connected, is incompressible in $W$, and 
instead of bounding a disc, $\gamma '$ bounds a surface  
$S_{1}$ in $W$ such that
$\# (\partial S_{1}\cap\Gamma _{0})\leq r$, $S_{1}$  
is
incompressible and boundary incompressible (in the sense of 
\ref{6.11} Rules 2 and 3). 

Then one of the following holds.
\begin{description}
\item[1.] $\vert f_{0}(\partial S_{1})\setminus A\vert \leq L$, and 
$f_{0}(\partial S_{1})$ is homotopic in $N$ to the union of 
$f_{0}(\partial S_{1})\setminus A$ and finitely many geodesic 
segments 
in $N$ of length $\leq L$. In particular, $(\partial S{1})_{*}\leq L$.
\item[2.] There is an essential embedded annulus $S_{2}$ in $W$ with boundary in $S$, such that  $\vert f_{0}(\partial S_{2})\setminus A\vert \leq L$. In particular, $(\partial S{2})_{*}\leq L$.
\end{description}

It follows that the lengths of the boundary components of maximal interval bundles in $W$ are bounded. 
\end{lemma}

{\em{Remark}}The final statement leads another proof of Thurston's ``bounded window frame'' theorem in \cite{T3} -- but using the full force of the Annulus Theorem.

\noindent{\em{Proof.}}  We start as in \ref{8.7}: taking a bounded 
union of $\leq 8r$ geodesic segments and extending to a triangulation 
of $S_{1}$. The number of triangles needed is $\leq 8r+2k$, where 
$-k$ is the Euler 
characteristic  of $S_{1}$. The corresponding geodesic triangles in 
$N$ determine 
a pleated surface $f_{1}:S_{1}\to N$, but not the boundary components 
may consist of finitely many geodesic segments, and not be completely 
geodesic. This pleated surface gives $S_{1}$ the structure of a 
complete hyperbolic surface $S(f_{1})$, but the boundary only 
consists of 
finitely many geodesic segments, and is only piecewise geodesic. 
The universal cover of $S(f_{1})$ then identifies a subset  with a 
closed convex subset of the hyperbolic plane, and the covering group 
is a discrete group of hyperbolic isometries of the hyperbolic plane. 
So Margulis' Lemma holds \cite{T}. For a suitable Margulis 
constant $\varepsilon _{0}>0$, nontrivial geodesics of length 
$<\varepsilon _{0}$ are disjoint and separated by distance $-\log 
\varepsilon _{0} -O(1)$ from all closed geodesics of length $\leq 1$. 
But 
in fact, Margulis' Lemma extends in this case of a hyperbolic surface 
with piecewise geodesic boundary. Nontrivial arcs of length $\leq 
\varepsilon _{0}$ between boundary components, and in different 
homotopy classes, are also disjoint. A nontrivial arc of length 
$<\varepsilon _{0}$ and nontrivial loops of length $<\varepsilon 
_{0}$ are  
separated by rectangles and cylinders of length $-\log \varepsilon 
_{0}
-O(1)$. 
So we have a decomposition of $S(f_{1})$ into a finite  
union 
of $\leq 2k+8r$ bounded diameter pieces, which, this time, may be 
connected by 
cylinders around $\leq 2k$ short closed loops, as well as $\leq 
2k+8r$ 
long thin rectangles. 
 In 
order to bound $\vert (\partial S_{1})_{*}\vert $, it suffices to 
show there are no rectangles over a given length.  

Suppose that $R$ is such a rectangle, with long sides in  $f_{0}(S(f_{0}))$ If  $R$  intersects 
$N_{<\varepsilon }$ for a suitable $T(f,\gamma _{*},\varepsilon)$, and then we deduce from the Radius of Injectivity Lemma \ref{3.3.3}  that $\gamma $ is freely homotopic in $N$ to $f_{0}(\gamma _{1})$ and $f_{0}(\gamma _{2})$  for loops $\gamma _{1}$ and $\gamma _{2}$ which are not freely homotopic in $S$. Then  we can choose the annulus $f_{1}(S_{2})$ homotopic to $\gamma $, giving conclusion 2.  So now we assume that the sides of the rectangle do not 
intersect $f_{0}((S(f_{0}))_{<\varepsilon })$  for a suitable 
$\varepsilon $ depending only on the original Margulis constant $\varepsilon _{0}$ the constant $C(\varepsilon _{0})$ of the Radius of Injectivity Lemma. Parametrize nearest points on the long sides of the rectangle by $(\zeta _{1}(t),\zeta _{2}(t))$, where $t$ is a length parameter, with $t\in [0,T]$, say. Then $\{ (\zeta _{1}(t),\zeta _{2}(t)):t\in [0,T]\} $ lies in a compact subset of $(S(f_{0}))_{\geq \varepsilon }\times ((S(f_{0}))_{\geq \varepsilon }$, of diameter bounded in terms of $\varepsilon $ using the product of the hyperbolic metric on $S(f_{1})$. Since $S$ is incompressible, If the lift of the rectangle has boundary components in the same component of the lift of $f_{0}(S)$, and $T$ is sufficiently large given $\varepsilon $, we can apply the Short Bridge Arc Lemma to deduce that an arc across the rectangle can be homotoped to a short arc in $S$. So such rectangles can either be removed (if they are boundary-homotopic in $S_{1}$), or discounted, using the boundary incompressibility of $S_{1}$. 

So now we assume that the lift of the rectangle has boundary components in different lifts of $f_{0}(S)$.
 So if $T$ is sufficiently large given $\varepsilon $, for any $t_{0}\in [0,T]$, we can find $t_{1}$ $t_{2}$ with $\vert t_{2}-t_{1}\vert \geq 1$ but $\vert t_{i}-t_{0}\vert $ bounded above in terms of $\varepsilon $ for $i=1$, $2$, and
such that the hyperbolic distance in $S(f_{0})$ between $f_{0}(\zeta _{j}(t_{1}))$ and $f_{0}(\zeta _{j}(t_{2}))$ is $<\varepsilon /10$ for $j=1$, $2$, and the distance in $N$ between $f_{0}(\zeta _{1}(t_{j})$ and $f_{0}(\zeta _{2}(t_{j})$ is also $<\varepsilon /10$  for $j=1$, $2$. Then we can join $\zeta _{j}(t_{1})$ and $\zeta _{j}(t_{2})$ by a short arc $\beta _{j}$ to give a nontrivial closed loop $\alpha _{j}$ on $S(f_{0})$, which we identify with $S$, embedded in $N$. Let $\gamma _{j}$ denote the short arc across the rectangle from $\zeta _{1}(t_{j}$ to $\zeta _{2}(t_{j})$. Then $f_{0}(\gamma _{1}*\beta _{2}*\overline{\gamma _{2}}*\overline{\beta _{1}})$ has length $\leq 2\varepsilon /5$. Since it is in $N_{\geq \varepsilon }$ is must be trivial. So $\gamma _{1}*\beta _{2}*\overline{\gamma _{2}}*\overline{\beta _{1}}$ is trivial in $W$. So the loops $\alpha _{1}$ and $\alpha _{2}$ are freely homotopic in $W$ but not in $S_{1}$.  By the proof of Waldhausen's  Annulus Theorem (Theorem 3) in \cite{C-F}, we can find a boundary incompressible embedded annnulus with boundary  in $S(f_{0})$  arbitrarily close to  a subset of $\alpha _{1}\cup \alpha _{2}$.Taking the image of this boundary under $f_{0}$ This gives alternative 2 above. If $T$ is bounded for all such rectangles $R$, then we obtain alternative 1.

The bound on the boundary of maximal interval bundles is achieved by applying the above with $S_{1}$ an annulus, if there is one, and then replacing by an annulus with bounded boundary length, and then repeating the process until a maximal set of boundary incompressible annuli with bounded boundary lengths has been constructed. This uses the fact that a  sufficiently long geodesic segment in $H^{2}$ which has endpoints  a bounded distance from a segment joining $x_{0}$ and $g^{n}.x_{0}$ must project to have self-intersections in $S(f_{0})$. This needs to be applied to $g$ in the conjugacy class of bounded annuli boundaries, in order to construct a next annulus with bounded boundary, if there is another annulus homotopically disjoint from a set already constructed. 
\Box
    
In \ref{8.7},
when 
$S_{1}$ is a disc, we only know that there is some 
disc, possibly different, whose boundary has bounded length. In the 
case $W\neq W_{n}$, we can strengthen this, as follows. This will mean 
that if we have good information on the geometry of $W=W_{n}$, then 
 we 
can ensure bounds on $\partial (\partial W_{i}\cap \partial W_{i+1})$.

\begin{lemma}\label{8.9} Let $N$, $S$ be as in \ref{8.7}.  Let 
 $S=S_{1}\cup S_{2}$, and such that $\gamma _{1}$, $\gamma _{2}$ and 
 $\gamma _{1}*\gamma _{2}$ are nontrivial in $N$, whenever $\gamma 
 _{i}$ is a closed loop with endpoints on $\partial S_{1}=\partial 
 S_{2}$ which is nontrivial in $S_{i}$.  Let a 
Margulis constant $\varepsilon _{0}$ and another constant $L_{1}$ be 
given.

Then there exists $L=L((S,N),\varepsilon _{0},L_{1})$ 
 such that the 
following 
holds. Let $\gamma \subset S$ bound a disc in $N$ such that $\gamma $ 
is nontrivial in $S$ and indecomposable, 
in the sense that $\gamma $ 
is not homotopic in $S$ to  $\gamma _{3}*\zeta  *\gamma 
_{4}*\overline{\zeta }$, 
where $\gamma _{3}*\zeta  *\gamma _{4}*\overline{\zeta }$ has only 
essential 
intersections with $\partial S_{1}=\partial S_{2}$ and $\gamma _{j}$ 
is 
homotopically nontrivial in $S$ and trivial in $N$ for $j=3$, $4$.

Let $f:S\to N$ be a pleated surface homotopic to the identity, and 
with  $\vert f(\partial S_{2})\vert \leq L_{1}$. Let $A$ be the union 
of badly 
bent annuli for $f$. Then $\vert f(\gamma )\setminus 
A\vert \leq L$ 
and $f(\gamma )$ is homotopic in $N$ to the union of $f(\gamma 
)\setminus 
A$ and a finite union of geodesic segments in $N$ of length $\leq 
L$.\end{lemma}

\noindent{\em{Proof.}} We use the argument of \ref{8.7}, splitting 
$f(\gamma)$ up into bounded diameter bits and long thin 
rectangles. If 
there is a long thin 
rectangle, between two different discs, there is either one with 
boundary in 
$f(S_{i})$ for one of $i=1$, $2$, or one with one 
boundary component in each of $f(S_{1})$, $f(S_{2})$. 
In the first case, since $S_{i}$ is incompressible, we can 
apply the  Short Bridge 
Arc Lemma of \ref{3.3.4} to deduce that  the arc across the rectangle 
can be homotoped into $f(S_{i})$. So we have two discs, neither of 
which can have boundary completely in $S_{i}$, because $S_{i}$ is 
incompressible. So we get a contradiction to 
indecomposablity, if there 
are no such rectangles and we have a bound on $\vert f(\gamma 
)\setminus A\vert $. The alternating bits in $N$ are bounded because 
of 
the way $f(\gamma )$ was split up. In the second case, the same 
argument works 
as 
in \ref{8.11}. If there is a long thin rectangle then, as in 
\ref{8.11}, there must be a long subrectangle in $N_{\geq \varepsilon 
_{0}}$. Then, as in \ref{8.11} we can find  
closed loops on $S_{1}$ and $S_{2}$ which are not multiples 
of  loops in $\partial S_{1}=\partial S_{2}$, and such that the 
product is trivial  in $N$, 
contradicting our assumption.\Box

The following lemma is useful for 
obtaining information about a pleated surface corresponding to a 
compressible component of $W$. It will be applied in conjunction with 
\ref{8.7}.

\begin{lemma}\label{8.8} Let $N$ be a hyperbolic 
$3$-manifold with finitely generated fundamental group, and let 
$W\subset N_{c}$ be such that each component of $N_{d}\setminus W$ is 
homeomorphic to the interior of an interval bundle. Fix such a 
component, with corresponding component $S_{d}$ of $\partial _{d}W$, 
where $S_{d}$ is the horodisc deletion of $S$.  
Let  $z_{0}=[\varphi _{0}]\in 
{\cal{T}}(S)$. 
Let $w_{0}=[\xi _{0}]\in {\cal{T}}(S)$. Let $\Gamma _{+}\subset S$ 
be a multicurve of loops which are all nontrivial 
in $N$, with 
$\vert \xi _{0}(\Gamma _{+})\vert \leq L_{0}$.
Let $\Gamma _{+}$ satisfy  
\begin{equation}\label{8.8.1}\begin{array}{l}
    i(\zeta ,\mu )\geq c_{0}\vert 
\zeta 
\vert {\rm{\ for\ at\ least\ one\ }}\zeta \in \Gamma _{+}\cr 
{\rm{\ if\ }}i(\gamma ',\mu )\leq c_{0}\vert \gamma '\vert ,\cr\end{array}\end{equation}
whenever $\gamma '$ is a simple  closed loop
 which is nontrivial in $S$ but trivial in $N$, 
$\mu $ is a  geodesic lamination on $S$, 
with a normalised transverse invariant measure.

Let $y_{0}=[\psi _{0}]$ satisfy 
$$\vert \psi _{0}(\gamma )\vert ''\leq L_{0},$$
where $\gamma $ is a closed loop 
which is nontrivial on $S$ but 
bounds a disc in $N$. 

Let $x(.)=x(.,[y_{0},w_{0}])$, as in \ref{7.10}. 

Then for a constant $D_{2}=D_{2}(L_{0},z_{0},c_{0})$, 
which is locally bounded in $z_{0}$ and  also depends 
on the topological type of $(S_{d},W)$ and 
suitable ltd parameter functions,
$$d(z_{0},x(z_{0}))\leq D_{2}.$$

\end{lemma}

\noindent {\em{Proof.}} Suppose for contradiction that 
$d(z_{0},x(z_{0})\geq D_{2}$. 
Then we can find $\alpha $ which is ltd 
along a segment of $[x_{\alpha }(z_{0}),z_{0}]$, and for some 
$x_{\alpha 
}(z_{0})=[\chi _{0}]\in [y_{0},w_{0}]$, for $C_{1}$ depending only on 
the 
ltd parameter functions,
$$\vert \chi _{0}(\partial \alpha )\leq 
C_{1},$$
$$d_{\alpha }(x(z_{0}),x_{\alpha }(z_{0}))\leq \log C_{1},$$
$$d_{\alpha }'(x_{\alpha }(z_{0}),z_{0})\geq 
D_{2}-\log C_{1}.$$

 Now for 
suitable $L_{0}$ we can choose $\zeta _{2}$ such 
that $\zeta 
_{2}\subset \alpha $ and 
$$  
\vert \chi _{0}(\zeta _{2})\vert \leq 
L_{0}.$$
Then by \ref{7.5}, or, at least, by the locking technique employed in 
\ref{7.5}, enlarging the constant $C_{1}$ if necessary but still
only depending on the  
the ltd parameter functions,
$$\vert \varphi _{0}(\gamma)\vert \geq 
C_{1}^{-1}i(\gamma ,\zeta _{2})\exp d_{\alpha }'(x_{\alpha 
}(z_{0}),z_{0}).$$
This is trivially satisfied if $\gamma $ does not intersect 
$\zeta _{2}$.
It uses the definition of $x(z_{0})$, that is, that $x_{\alpha 
}(z_{0})$ is a bounded $d_{\alpha }$ distance, depending only on the 
ltd parameter functions, from each of $[z_{0},w_{0}]$ and 
$[z_{0},y_{0}]$. Similarly, for any $\zeta \in \Gamma _{+}$,
$$\vert \varphi _{0}(\zeta )\vert \geq C_{1}^{-1}i(\zeta 
,\zeta _{2})\exp d_{\alpha }'(x_{\alpha }(z_{0}),z_{0}).$$
Now $\vert 
\varphi _{0}(\zeta )\vert $ is boundedly proportional to $\vert \zeta 
\vert $ for any loop $\zeta $, with bound depending locally uniformly 
on $z_{0}$.  So  we have
$$i(\gamma ,\zeta _{2})\leq C_{1}^{3}e^{-D_{2}}.\vert \gamma \vert 
,$$
and similarly with $\zeta $ replacing $\gamma $ for any $\zeta \in 
\Gamma _{+}$.  Putting 
$\gamma =\gamma '$,  and  $\mu =\zeta _{2}/\vert \zeta 
_{2}$, we obtain a contradiction to (\ref{8.8.1}), 
if $D_{2}$ is large enough given $L_{0}$, $c_{0}$ and locally on 
$z_{0}$.

\Box

\ssubsection{Sequence of pleated surfaces corresponding to an 
incompressible end.}\label{8.6}

We shall always use the following hypothesis on $\Gamma _{e,+}$, 
which 
depends on a constant $c_{0}$. The 
assumption is a bit stronger than required in (\ref{8.8.1}). We can 
manage with just (\ref{8.8.1}) for much of what follows, but not all, 
so we might as well fix on the stronger assumption now. 
\begin{equation}\label{8.6.1}\begin{array}{l}
    i(\zeta ,\mu )\geq c_{0}\vert 
\zeta 
\vert {\rm{\ for\ at\ least\ one\ }}\zeta \in \Gamma _{+}\cr 
{\rm{\ if\ }}i(\mu ,\mu 
')=0{\rm{\ and\ }}i(\mu '\gamma ')\leq c_{0}\vert \gamma \vert ,\cr\end{array}\end{equation}
whenever $\gamma '$ is a simple  closed loop
 which is nontrivial in $S$ but trivial in $N$, 
 and $\mu $, $\mu '$ are  geodesic laminations on $S$, 
with  normalised transverse invariant measures.
With this assumption, fix any $z_{e,0}=[\varphi _{e,0}]\in 
{\cal{T}}(S(e))$. 
We can find sequences  of multicurves $\Gamma _{e,j}$ on $S(e)$ 
and 
$f_{e,j}:S(e)\to N$ homotopic to inclusion of $S(e)$ in $N$ 
($0\leq j\leq n_{+}(e)$) satisfying \ref{8.2.1} to \ref{8.2.4}, with 
$$\vert \varphi _{e,0}(\Gamma _{e,0})\vert \leq L_{0},$$
where $L_{0}$ is bounded in terms of $z_{e,0}$ and $c_{0}$, but does 
not depend further on $y_{e,+}$. We see this as follows. By 
replacing $z_{e,0}$ by a point on $[z_{e,0},y_{e,+}]$ sufficiently 
far 
from $z_{e,0}$ given $D_{0}$, by \ref{8.7} and \ref{8.8} we can 
ensure 
that 
\begin{equation}\label{8.6.2}\vert \varphi (\gamma )\vert ''\geq 
    D_{0}\vert \gamma \vert \end{equation}
for 
all $[\varphi ]\in [z_{e,0},y_{e,+}]$ and any $\gamma $ which is 
nontrivial on $S(e)$ but trivial in $N$. Now for some constant 
$L_{0}$ depending on an initial choice of ltd parameter functions,
we can choose the sequences $z_{i}^{2}=[\varphi _{i}^{2}]$ and 
$\Gamma 
_{i}=\Gamma _{i}^{2}$ as in \ref{8.4} so that for 
$$\vert \varphi _{i}^{2}(\Gamma _{i}^{2})\vert ''\leq L_{0}.$$
If $D_{0}$ is sufficiently large, none of the loop sets can be 
collapsing. 
For if 
some 
one is, then for some $i$ and $\gamma \in \Gamma _{i}^{2}$ if 
we take any loop $\zeta $ intersecting $\gamma $ at most twice and 
disjoint 
from other elements of $\Gamma _{i}^{2}$, not separating loops with 
the same image in $N$, $\tau _{\gamma }^{n}(\zeta )$ is trivial in 
$N$ for all $n$. Replacing $\zeta $ 
by $\tau _{\gamma }^{n}(\zeta )$ for a suitable $n$, we can assume 
that $\vert \varphi _{i}^{2}(\zeta )\vert ''$ is bounded in terms of 
$L_{0}$. 
This contradicts (\ref{8.6.2}) if $D_{0}$ is large enough given 
$L_{0}$. 
\ref{8.8} also implies that (\ref{8.2.5}) and (\ref{8.2.6}) are 
satisfied for $j\geq j_{0}$, if $j_{0}$ is the largest integer such 
that 
$$d(x([f_{j}],[z_{e,0},y_{e,+}]),z_{e,0})\leq D_{2}.$$
We shall use this in Section \ref{10}.

Then by 
\ref{8.8}, if $D_{0}$ is sufficiently large given $D_{2}$ and 
$D_{1}$, $\vert f_{j}(\gamma )\vert \geq D_{1}$ whenever 
$T([f_{j}],+,[z_{0},y_{+}])$ (in the notation of \ref{7.9}) does not 
contain $(\alpha _{1},\ell _{1})$. This is automatically true for 
$j=n_{+}$.

\ssubsection{Family of pleated surfaces for the non-interval bundle 
part of the core}\label{8.12}

So far, we have constructed a sequence of maximal multicurves and 
pleated surfaces for each end $e$ of $N_{d,W}$. Conditions 
\ref{8.2.1} to \ref{8.2.4} are satisfied for all ends. Conditions 
\ref{8.2.5} and \ref{8.2.6} are also satisfied for incompressible 
ends. 
It is not  clear if \ref{8.2.5} is satisfied for compressible ends. 
If 
it is, then \ref{8.2.6} is also satisfied. We are now going to 
construct a family of pleated surfaces for $W$, under a temporary 
assumption 
which we shall prove in Section \ref{10}. The temporary assumption we 
make is 
the following. It is, of course, unnecessary if all ends are 
incompressible.

\subsubsection{Assumption on compressible ends}\label{8.12.1}
For any compressible end,
$$\vert f_{e,0}(\Gamma _{e,0})\vert \leq L_{1}.$$

\begin{utheorem}  Assume that \ref{8.12.1} holds. Let $W_{i}$ be a 
decomposition of $W$ satisfying conditions 1 to 7 of \ref{6.11}. Then 
we can find a family of maximal multicurves and corresponding pleated 
surfaces for the non-interval-bundle part $W'$ of $W$, consisting of  
sequences of two multicurves and  
pleated sequences for each component $S$ of $\partial V$, each 
non-ball 
component $V$ of $W_{i}$, each $i\geq 1$, such that  the following 
hold, for a suitable constant $L_{0}$, and for $L_{2}$ sufficiently 
large given $L_{0}$ and $L_{1}$.
\begin{description} 
     
     \item[1.]For the two multicurves$\Gamma $, 
 $\Gamma '$ in the sequence 
 for $S$ $\# (\Gamma \cap \Gamma ')\leq L_{0}$. 

 \item[2.]If $i=n$ and $S=S(e)$ and $e$ is a compressible end then 
the second
multicurve of the 
sequence for $S(e)$ is $\Gamma _{0}(e)$, and the corresponding 
pleated 
surface is $f_{e,0}$. (These are the first elements in the sequence 
for the end $e$.)

\item[3.]If $S$ is incompressible, then the first multicurve of the 
sequence includes $\partial (S\cap S')$ for any component $S'$ of 
$\partial V'$, for $V'$ a component of $W_{i-1}$, $V'\subset V$. For 
any such $S'$,
$$\vert \partial (S\cap S'))_{*}\vert \leq L_{2}.$$

\item[4.] Whether $S$ is compressible or incompressible, with $i<n$, 
the second multicurve in the sequence for $S$ coincides 
with the 
first multicurve for $S''$ on $S\cap S''$, where $S''$ is a component 
of $\partial V''$ for a component $V''$ of $W_{i+1}$, $V\subset 
V''$, $S\cap S''\neq \emptyset $. If $V_{2}$ is another 
component of $W_{i}$ in $V''$ and $S_{2}$ is the component of 
$\partial V_{2}$ meeting $S$, then the second multicurves in the 
sequences for $S$ and $S_{2}$ coincide on $S\cap S_{2}$.

\item[5.] If $S$ is compressible, the first multicurve in the 
sequence 
for $S$ can be 
written in the form $\Gamma _{1}=\Gamma _{1,1}\cup \Gamma _{1,2}$, 
where $\Gamma 
_{1,1}\cup \Gamma _{1,2}'$ is also a maximal multicurve, and 
$\Gamma _{1,2}'$ is the set of distinct 
isotopy classes in $S$ of $\partial (S\cap S')$ for any component 
$S'$ of 
$\partial V'$, for $V'$ a component of $W_{i-1}$, $V'\subset V$. 

\item[6.] Let $\Sigma '$ be the graph in $W$ of \ref{6.14}. All the 
maximal multicurves in the sequence for $S$ are made up 
of $\leq r_{0}$ arcs of $\Sigma '$ on $S$. For any pleated surface 
$f:S\to N$ in the sequence, and any $\gamma \subset S$ made up of 
$\leq r_{0}$ arcs of $\Sigma '$ if $A$ is the set of badly bent 
annuli 
for $f$, then
$$\vert f(\gamma )\setminus A\vert \leq L_{2},$$
and $f(\gamma )$ is homotopic in $N$ to the union of $f(\gamma 
)\setminus A$ and finitely many arcs of total length $\leq L_{2}$.

\item[7.] There is a map $f_{W'}:W'\to N$ which is homotopic to 
inclusion, 
such that $f\vert \partial W'\cap S_{d}(e)=f_{e,0}$ for any end $e$,
and such that for any arc $\gamma $ in $\Sigma '$, or between arcs of 
$\Sigma '$ which are homotopic then $\vert f(\gamma )\vert \leq 
L_{2}$. 
\end{description}
 
\end{utheorem}

\ssubsection{Proof of \ref{8.12}: 1-5 and a bit.}\label{8.13}

The construction is inductive, working down. To start with we only 
aim 
to satisfy 1 to 4, and the following extension of 5:
\begin{description}
    \item[5 and a bit] If $S$ is compressible, the first multicurve 
in the sequence can be 
written in the form $\Gamma _{1,1}\cup \Gamma _{1,2}$, where $\Gamma 
_{1,1}\cup \Gamma _{1,2}'$ is also a maximal multicurve, and 
$\Gamma _{1,2}'$ is the set of distinct 
isotopy classes in $S$ of $\partial (S\cap S')$ for any component 
$S'$ of 
$\partial V'$, for $V'$ a component of $W_{i-1}$, $V'\subset V$. For 
any  $\gamma \in \Gamma _{1,2}$, $\gamma \in \Gamma _{1,1}$ adjacent 
to a loop of $\Gamma _{1,2}$ 
(equivalently of $\Gamma _{1,2}'$
$$\vert \gamma _{*}\vert \leq L_{2}.$$
\end{description}

     So suppose that $S$ is a 
component of $\partial V$, $V$ a component of $W_{i}$. If $i<n$, 
suppose that the conditions are already satisfied for all components 
of $V''$, where $V''$ is the component of $W_{i+1}$ containing $V$. 
If $i=n$ and $S=S(e)$ for compressible $e$, then 2 is satisfied by 
definition. If  $V$ is incompressible, then we obtain 3 
from \ref{8.11}, for any $i$. Then 1 is true, for $L_{0}$ depending 
only on the topological type 
of $W'$. Satisfying 4 inductively is no problem.

So now we need to show that 5 and a bit can be satisfied. 
 If $i=n$, and $S=S(e)$ for a compressible end $e$, let $\Gamma _{0}$ 
be 
the second multicurve in the sequence for $S=S(e)$. By construction, 
$\Gamma _{0}$ is noncollapsing. If $i+1<n$, let 
$S''$ be component of $\partial V''$, $V''\subset W_{i+1}$,  
which intersects $S$, and let $\Gamma _{0}=\partial (S\cap S'')$.
According to rule 6 of \ref{6.11}, $S''$ is incompressible in $N$, 
and so the loops of 
$\Gamma _{0}$ are nontrivial in both $N$ and $S$. According to Rule 
4 of \ref{8.12}, the second multicurve $\Gamma _{2}$ in the sequence 
for $S$ 
should include $\Gamma _{0}$. We can choose $\Gamma _{2}$ so that 
this is true.  Using \ref{4.1}, we can also ensure that $\vert 
\gamma _{*}\vert \geq \varepsilon _{0}$ for $\gamma \in \Gamma 
_{2}\setminus \Gamma 
_{0}$.  Now let
$\Gamma _{1,2}'$ be the set of boundaries of compressing discs 
attached 
to $S$ to form the components of $W_{i-1}$ in $V$. By Rule 6 of 
\ref{6.11}, $\Gamma _{1,2}'$ bounds a maximal set of disjoint 
compressing discs. 
Then every loop on $S''$ disjoint 
from $\Gamma _{1,2}'$ is nontrivial in $V''$, and hence in $N$. 
(It can be proved inductively that there is no 
compressing disc attached to the exterior of $V''$.) 
Let $T\subset S'$ be 
the surface with boundary which is the convex hull of $\Gamma 
_{1,2}'$ 
and $\Gamma _{0}$. Then all loops in  $\partial T$ (if $\partial 
T\neq \emptyset $) are 
nontrivial in $N$, since otherwise $\# (\Gamma _{1,2}')$ is not 
maximal. 

Let $f_{2}$ be the second pleated surface in the sequence for $S$. 
Then the 
conditions of \ref{8.9} are satisfied by rule 6 of \ref{6.11}. We can 
therefore use \ref{8.9} to  bound $\vert f_{2}(\Gamma 
_{1,2}')\setminus A\vert $, 
and the 
length of connecting geodesic segments in $N$ where $A_{2}$ is 
the union of badly bent annuli for $f_{2}$.  We already have a bound 
on 
$\vert f(\Gamma _{0})\vert =\vert (\Gamma _{0})_{*}\vert $ by 
\ref{8.9}. So let $\tau =\Gamma _{0}\cup \Gamma _{1,2}'$. We have a 
bound on $\vert f(\tau )\vert $.
Then choose $\Gamma _{1,1}$ with $\partial T\subset \cup \Gamma 
_{1,1}$ 
such that $\Gamma _{1,2}'\cup \Gamma _{1,1}$ is a maximal multicurve 
on 
$S$. 
Since $\# (\Gamma _{1,2}')$ is maximal, $\Gamma _{1,1}$ is 
noncollapsing. Then by the lemma in \ref{3.8}, we can extend 
$\Gamma _{1,1}$ to a 
maximal noncollapsing maximal multicurve $\Gamma _{1}=\Gamma 
_{1,1}\cup \Gamma 
_{1,2}$ on $S$.  Since $\partial T\subset \Gamma _{1,1}$ and $\Gamma 
_{1,2}\subset T$, the loops of $\Gamma _{1,2}$ and $\Gamma _{1,1}$ 
adjacent to $\Gamma _{1,2}'$ 
(equivalently $\Gamma _{1,2}$) are in $\tau $ up to homotopy. So, 
simply depending on the topology of $\tau $ (that is, of the 
decomposition into the $W_{i}$, and, ultimately, of $W$) we can bound 
the number of arcs of $\tau $ making up the loops of $\Gamma _{1,2}$, 
and the arcs of $\Gamma _{1,1}$ adjacent to $\Gamma _{1,2}$. So we 
have 
5 and a bit. We can also ensure, by \ref{4.1}, that $\vert \gamma 
_{*}\vert 
\geq \varepsilon _{0}$ for $\gamma \in \Gamma _{1,1}\setminus 
\partial T$.

\ssubsection{Proof of 6 and 7 of \ref{8.12}}\label{8.14}
We have bounds on $\vert f(\gamma )\vert $ whenever $f:S\to N$ is a 
pleated surface in the family, and $\gamma $ is an arc in $\Sigma 
\cap S$. The arcs in $\Sigma '$ are obtained from $\Sigma $ by 
sucessive homotopies between the boundaries of the interval bundles 
which make up $W$. The elements of $\Sigma $ are themselves obtained 
from homotopy images, this time, from a certain set of closed simple 
loops. So it suffices, to bound the length, to find homotopies 
between 
pleated surfaces in the family with the same domain, with bounded 
homotopy tracks, outside badly bent annuli round short loops in the 
pleating locus. At the join between two or three interval bundles in 
$W$, two or three surfaces have parts of the domain in common, with a 
closed geodesic in common to the pleating locus of the two or three 
surfaces. When this happens, one surface $S$ is a boundary component 
of  $\partial V$ for some component $V$ of $W_{i}$, some $i$, and the 
other one or two surfaces are boundary components of components of 
$W_{i-1}$. 

We consider first the case when $S$ is compressible. Let $f_{1}$ and 
$f_{2}$ be the first and second pleated surfaces with domain $S$, as 
in \ref{8.13}. In this case, the domain of any $S'$ which is attached 
to $S$, and a component of the boundary of some component of 
$W_{i-1}$, 
identifies with a subsurface of $S$ bounded by loops in $f_{1}$. Now 
we consider the homotopy between $f_{1}$ and $f_{2}$. 
If $S$ is 
compressible, the condition of \ref{4.3} which replaces 
injectivity-on-$\pi _{1}$ is probably not satisfied. But for $f_{2}$, 
the loop set $\Gamma _{0}$ constructed in \ref{8.13} 
decomposes $S$ into incompressible surfaces, which, for the moment we 
call $S_{j}'$, $1\leq j\leq r$. Essentially, we shall apply the Short 
Bridge Arc Lemma to these. For $f_{1}$, we have a 
decomposition of $S$ into surfaces $S_{1,1}$ and $S_{1,2}$ with 
$\Gamma 
_{1,2}\cup \Gamma _{1,2}'\subset S_{1,2}$ and $\partial 
S_{1,2}\subset \Gamma 
_{1,1}$, and $S_{1,1}$ is incompressible. In fact, we can further 
decompose $S_{1,1}$ using $\Gamma _{1,1}$, using the $S_{j}'$, 
into surfaces $S_{1,1,j}$ 
which are incompressible. 
The surface $S_{2}$ is in $T$ 
and $\Sigma '\cap T$ can be homotoped into $\tau $. 
We apply \ref{8.9} to $f_{1}$ and deduce 
bounds on $\vert f_{1}(\tau)\setminus A_{1}\vert $, and on the 
geodesic 
arcs homotopic  in $N$ to components of $f_{1}(\tau )\cap A_{1}$ 
where $A_{1}$ is the union of any badly bent annuli on $S$ for $f$. 
If there are any badly bent annuli, they are round loops of $\Gamma 
_{1,1}$ adjacent to $\Gamma _{1,2}$, each of which actually gives two 
loops in the pleating locus for $f_{1}$. Note that this then 
automatically transfers to the next surface down. Let $A_{2}$ be the 
set of badly bent annuli for $f_{2}$. We have chosen the pleating 
locus of $f_{2}$ so that these only occur round loops $\partial 
(S_{j}'\cap S)$, where $S'$ are the boundary components of $V'$ which 
$S$ meets, where $V'$ is the component of $W_{i+1}$ containing $V$.

 First, we bound $\vert 
f_{2}(\partial S_{1,j}\cap S_{k}')\setminus A_{2}\vert $ 
  and the geodesic segments 
homotopic to $f_{2}(\partial S_{1,j}\cap S_{k}')\cap A_{2}$ 
using the the Short Bridge Arc Lemma \ref{3.3.4} applied to 
the surfaces $S_{k}'$. Similarly we bound
$\vert f_{1}(\Gamma _{0}\cap S_{1,j})\setminus A_{1}\vert $ and 
 and the geodesic segments 
homotopic to $f_{1}(\Gamma _{0}\cap S_{1,j})\cap A_{1}$. Taking the 
unions over all $S_{k}'$, or over $S_{1,1}$ and $S_{1,2}$, we obtain 
bounds on 
$\vert f_{2}(\partial S_{1,j})\setminus A_{2}\vert $, 
$\vert f_{2}(\partial S_{1,j})\cap A_{2}\vert $,
$\vert f_{1}(\Gamma _{0})\setminus A_{1}\vert $, 
$\vert f_{1}(\Gamma _{0})\cap A_{1}\vert $.
Then we use 
\ref{3.3.4} to homotope 
$f_{2}\vert \tau \setminus (f_{1}^{-1}(A_{1})\cup f_{2}^{-1}(A_{2}))$ 
to 
$f_{1}\vert \tau \setminus (f_{1}^{-1}(A_{1})\cup 
f_{2}^{-1}(A_{2}))$, 
and hence bound $\vert f_{1}(\tau )\setminus A_{1}\vert $, and the 
union of geodesic segments homotopic to $f_{1}(\tau )\cap A_{1}$.
 Similarly 
if $\tau '\subset \Sigma '$ and $\vert f_{2}(\tau ')\setminus 
A_{2}\vert $ is known to be 
bounded and bounds are known on the geodesic segments homotopic to 
$f_{2}(\tau '\cap A_{2}$, we can homotope $f_{2}\vert \tau '\setminus 
( f_{2}^{-1}(A_{2})\cup f_{2}^{-1}(A_{2}))$ to $f_{1}\vert \tau 
'\setminus 
(f_{1}^{-1}(A_{1})\cup f_{2}^{-1}(A_{2}))$, and hence bound $\vert 
f_{1}(\tau' )\setminus A_{1}\vert $, and the union of geodesic 
segments homotopic to $f_{1}(\tau ')\cap A_{1}$. Assuming as we may 
do 
that $\# (\Gamma _{1}\cap \Gamma _{2})$ is bounded in terms of the 
topological type, we can then extend the homotopty, using \ref{3.3.4} 
as before, to produce a bounded track homotopy beween $f_{1}$ and 
$f_{2}$ outside $f_{1}^{-1}(A_{1})\cup f_{2}^{-1}(A_{2})$.

We can also transfer lengths to a surface $f_{0}:S_{0}\to N$, where 
$S_{0}$ is a component of $\partial V'$, $V'$ a component of 
$W_{i-1}$, $S_{0}$ attached to $S$, and incompressible. If this 
happens, pairs of loops in 
$\Gamma _{1,1}$ on $S$ identify in $N$, and also identify with a loop 
in the pleating locus of $f_{0}$. The second maximal multicurve on 
$S_{0}$ is given by $\Gamma _{1,1}$ after identification. We use all 
loops of $\Gamma _{1,1}$ if $S_{0}$ is connected. In general, the 
loops of $\Gamma _{1,1}$, after identifications, give the union of 
multcurves for the boundary components of the components of $W_{i-1}$ 
attached to $S$.  Let $A_{0}$ be the union of badly 
bent annuli for $f_{0}$. These,  again, have as cores a subset of the 
loops of $\Gamma _{1,1}$ with badly bent annuli in $A_{1}$ The 
subsurface we are calling 
$S_{1,1}$ glues up along matching loops in $S_{1,1}$ to give 
$S_{0}$. Suppose that $\tau _{1}=\tau _{1,1}\cap \tau _{1,2}$, 
$\tau _{2}=\tau _{1,2}\cup \tau _{2,2}\subset S_{1,1}$
have endpoints in $\partial 
S_{1,1}=\partial S_{1,2}$, with all four graphs meeting along a loop 
in $\partial S_{1,1}$, $\tau _{1,1}\cup \tau _{1,2}$ identifying with 
a graph on $S_{0}$. Suppose also , from looking at $f_{1}(\tau _{1})$ 
and $f_{1}(\tau _{2})$, we already have bounds on the length of 
geodesic 
segments making up $\tau _{1}$, $\tau _{2}$. Suppose also that $\tau 
_{1,2}$ and $\tau _{2,2}$ lie in parts of $S$ where the muticurve 
loops from $\Gamma _{1,2}$ are identified in $N$, and that they are 
homotopic in $N$, under homotopy preserving endpoints. Then 
$\tau_{1,1}\cup \tau _{1,2}$ identifies with a graph $\tau _{0}$ on 
$S_{0}$. 
Then applying 
 \ref{3.3.4} to $f_{0}$, and using the bounds on $\tau _{1}$, $\tau 
 _{2}$, we can bound $\vert f_{0}(\tau )\setminus A_{0}\vert $ and 
the union of geodesic segment homotopic to $f_{0}(\tau )\cap A_{0}$ 
via homotopy preserving segment endpoints. 

If $S$ is incompressible, the transfers of lengths are rather 
similar. \ref{3.3.4} can be applied directly to transfer lengths 
between $f_{1}$ and $f_{2}$, which makes things slightly easier. 
Transfer of lengths around common loops is slightly different, in 
that 
$S$ is incompressible and the one or two surfaces from boundary 
components of components of $W_{i-1}$ might be compressible or 
incompressible. Compressible surfaces, as before, have to be 
decomposed into incompressible subsurfaces in order to apply 
\ref{3.3.4}. 

Now we construct the map $f:W'\to N$ which is homotopic to inclusion. 
We now have a connected graph of geodesic segments, each of which is 
either has the endpoints of $f(\gamma )$ for some $f:S\to N$ in the 
family of pleated surfaces for $W$, and $\gamma \subset \Sigma '\cap 
S$, or is the homotopy track between $f_{1}(x)$ and $f_{2}(x)$ for 
some endpoint $x$ of a maximal arc in $\Sigma '\cap S$, where $S$ is 
the domain of $S_{1}$ and $S_{2}$. All these geodesic segments have 
length $\leq L_{1}$ for some suitable $L_{1}$ which, ultimately, 
depends only on the topological type of $W'$, and the constant 
$c_{0}$. 
This graph then gives a $1$-complex for the  cell-complex 
homeomorphic to 
$W$,  in which the components are the components of $W_{0}$. We 
refine the cell decomposition to a decomposition into tetrahedra, 
keeping the same vertices, but adding edges. Again, the corresponding 
geodesic segments in $N$ have length bounded in terms of $L_{1}$ and 
the topological type of $W'$. Then we map these topological 
tetrahedra 
in the topological decomposition of $W'$ to hyperbolic tetrahedra in 
$N$. This gives the required homotopy equivalence, except on the 
boundary. On $S_{d}(e)\cap \partial W'$ we define $f=f_{e,0}$. Since 
$\vert f_{e,0}(\tau )\vert \leq L_{2}$ for any arc $\tau \in 
S(e)\cap \Sigma '$, the arc is a bounded distance from the geodesic 
with the same endpoints, and we can homotope $f_{e,0}\vert S_{d}(e)$ 
to the map which takes the $1$-complex to a union of geodesic 
segments, by a homotopy with bounded tracks. So putting these 
together, we obtain a homotopy equivalence with the desired 
properties. 
\Box

\section{Combinatorially Bounded geometry case.}\label{9}

\ssubsection{Proof of 
\ref{1.3} in the combinatorial 
bounded geometry Kleinian surface case.}\label{9.1}

We now prove our main theorem \ref{1.3} in the case of 
topological type $S\times \mathbb R$ for a finite type surface $S$, 
with end invariants $\mu (e_{\pm })\in {\cal{GL}}_{a}(S)$ and with 
combinatorially bounded geometry, and therefore, exactly two ends. 
The other two main 
theorems 
follow from 1.3. We choose pleated surfaces $f_{\pm }:S\to N$ 
approaching the ends, 
as in \ref{3.8}, where the pleating loci include
maximal multicurves $\Gamma _{\pm }$ where 
$\vert f_{+}(\Gamma _{+})\vert \leq L_{0}$, and similarly for 
$f_{-}$, $\Gamma _{-}$. We then fix on a choice of loop $\zeta _{\pm 
}\in \Gamma _{\pm }$. As we have seen in \ref{6.12}, we can find 
$y_{\pm }'=[\varphi _{\pm }']\in [y_{-},y_{+}]$ such that $\vert 
\varphi _{+}'(\zeta _{+})\vert \leq L_{0}$, similarly for $y_{-}'$, 
$\zeta _{-}$ and for a suitable 
$\nu >0$,
$$[y'_{-},y'_{+}]\subset ({\cal{T}}(S))_{\geq \nu }.$$
For $z_{0}=y_{-}'$ and $y_{+}=y_{+}'$, we then choose a sequence 
$\Gamma _{j}$ of 
maximal multicurves, as in \ref{8.3}, and $f_{j}$ a pleated surface 
with 
pleating locus including $\Gamma _{j}$, for $-n_{-}\leq j\leq n_{+}$ 
such that $\vert f_{n_{+}}(\zeta _{+})\vert \leq L_{0}(\nu )$, and 
similarly for $f_{n,{-}}$, $\zeta _{-}$, and so that \ref{8.2.1} to 
\ref{8.2.6} are 
satisfied. We then know by \ref{7.6} that  every point on  
$[y'_{-},y'_{+}]$ is a bounded distance $\leq C(\nu )$ from a 
point on $[[f_{n_{-}}],[f_{n_{+}}]]$, as noted in \ref{8.3}.

By \ref{8.3} we can find a sequence $f_{i}:S\to N$ of pleated 
surfaces with pleating locus of $f_{i}$ 
including $\Gamma _{i}$ for $-n_{-}\leq i\leq n_{+}$, such that 
\ref{8.2.1}-\ref{8.2.6} hold --- of which \ref{8.2.5} is redundant in 
this case. 
 
   It is convenient to choose the indexing of the $z_{j}$ as 
$-n_{-}\leq i\leq n_{+}$  --- rather than $0\leq j\leq n_{+}$ as in 
\ref{8.3} ---
 so that for the corresponding sequence $\{ 
z_{i}:-n_{-}\leq i\leq n_{+}\} $ we have $z_{0}$ in a bounded subset 
of $\cal{T}(S)$. We know that there is such a point $z_{0}\in 
[y_{-}',y_{+}']$ by \ref{6.12}. By definition, $y_{-}'=z_{-n_{-}}$ 
and $y_{+}'=z{n_{+}}$.  But at the moment we have no bound on 
$d(z_{i},[f_{i}])$ or on $d(z_{i},[f_{i}])$ for a general $i$. But by 
\ref{8.2.6} we do have, for a suitable $L_{0}$, for all $p$,
$$d([f_{p}],[f_{p+1}])\leq L_{0}.$$

Now we consider $x(.,[y_{-}',y_{+}'])$.   Now we claim that the 
hypotheses of \ref{7.12} hold for $[y_{-}',y_{+}']$ replacing 
$[z_{-},z_{+}]$. 
So suppose that  for some 
$p$ and $i$,
and $d([f_{p}],z_{i})\leq L$, and $d([f_{p}],z_{j})\leq r'+1$ for 
$\vert 
i-j\vert 
\leq r'$. By the definition of the sequence in \ref{8.3}, 
$z_{i}=[\varphi _{i}]$ 
and $\Gamma _{i}$ have the properties that, for a constant $L_{0}$, 
 $\vert \varphi _{i}(\Gamma _{i})\vert \leq L_{0}$. So for a constant 
$L_{0}'=L_{0}'(L,L_{0},r')$, and $\vert j-i\vert \leq r'$,
$$\vert (\Gamma _{j})_{*}\vert \leq \vert f_{p}(\Gamma _{j})\vert 
\leq 
L_{0}'.$$
  So 
we also have a bound on $\vert 
f_{j}(\Gamma _{j})\vert $. But then the bound on 
$d([f_{p}],[f_{p+1}])$ for all $p$ gives,
for a suitable constant $L_{0}''=L_{0}''(L_{0},L,r')$, for  
$\vert j-i\vert \leq r'$,
$$\vert f_{i}(\Gamma _{j})\vert \leq L_{0}''.$$
If $r'$ is large enough given $\nu $, by \ref{7.2},
$\cup \{ 
\Gamma _{j}:\vert i-j\vert \leq r'\} $ is cell-cutting. 
 So then the 
same loops are bounded at $z_{i}$ as at $[f_{i}]$. 
So we have a bound on $d([f_{i}],z_{i})$. So the hypotheses of 
\ref{7.12} 
are satisfied. So then by \ref{7.12}, for all 
$j$, and a 
suitable constant $L_{1}$, we have
\begin{equation}\label{9.1.3}
    d([f_{j}],z_{j})\leq L_{1}.
    \end{equation}
In particular we have, still for a constant $L_{1}$ depending only 
on $\nu $,
$$\vert (\Gamma _{j})\vert _{*}\leq L_{1}.$$
\Box

\ssubsection{Combinatorially bounded geometry Kleinian surface case: 
no Margulis tubes.}\label{9.2}

Fix a basepoint $w_{0}\in N$. Let $\Delta _{1}>0$ be given. Then we 
can choose $f_{\pm }$ so that the $\Delta _{1}$-neighbourhood of 
$w_{0}$ is 
strictly 
between $f_{-}(S)$ and $f_{+}(S)$. To do 
this, we use the Bounded 
Diameter Lemma of \ref{3.3}. There is  $D_{0}$ given by the Bounded 
Diameter Lemma, and an integer $k_{0}$ bounded in terms of 
topological 
type, such that $f_{n_{+}}(S)$ is connected to $f_{+}(S)$ by a 
chain of at 
most $k_{0}$ Margulis tubes and sets of diameter $\leq D_{0}$. The 
same is true for $f_{-}(S)$, $f_{-n_{-}}(S)$. So we choose 
$f_{+}$ and 
$f_{-}$ so that $f_{\pm}(S)$ cannot be connected to the $\Delta _{1}$ 
neighbourhood of $w_{0}$ by such chains. 
Then we claim that for some $\varepsilon _{1}>0$ depending only on 
$\nu $, the $\Delta _{1}$-neighbourhood  of a fixed basepoint 
$w_{0}\in 
N$ does not intersect any $\varepsilon _{1}$-Margulis tube.
We see this as follows. For each $i$, $f_{i}(S)$ 
separates 
$N$. In particular, this is true for $i=n_{+}$ and $i=-n_{-}$.
The point $w_{0}\in N$ is in a  
component $N_{1}$ of the complement in $N$ of 
$f_{n{+}}(S)\cup f_{n_{-}}(S)$ whose boundary meets both $f_{n{+}}(S)$
and $f_{n{-}}(S)$.  The union of the homotopies between $f_{i}$ and 
$f_{i+1}$ is a
homotopy between $f_{n_{-}}$ and $f_{n_{+}}$.  The basic principle we
use is:
\subsubsection{}\label{9.2.1}
   $N_{1}$ is contained in the image of the homotopy between 
$f_{n_{-}}$ and $f_{n_{+}}$. 
    Hence each point in $N_{1}$ is in the image of a homotopy between 
    $f_{i}$ and $f_{i+1}$ for some $i$.
  
For if this is not true, we have a homotopy equivalence between 
$(S\times
[0,1],S\times \{ 0,1\} )$ and $S\times
[0,1]\setminus \{ x\} ,S\times \{ 0,1\} )$, for any internal point 
$x$. 

But we now
know, by the bound on $d([f_{i}],z_{i})$ that $[f_{i}]\in
({\cal{T}}(S))_{\geq \varepsilon _{2}(\nu )}$ for some constant
$\varepsilon _{2}(\nu )>0$ depending only on $\nu $.  So since
$f_{i}:S(f_{i}))\to N$ is decreasing, $f_{i}(S)$ has diameter $\leq
L_{1}=L_{1}(\nu )$.  Also, by \ref{8.2.5}, the diameter of the
homotopy between $f_{i}$ and $f_{i+1}$ is bounded by $L_{0}(\nu )$,
taking $L_{0}(\nu )$ sufficiently large.  So for a suitable
$\varepsilon _{1}=\varepsilon _{1}(\nu )>0$, the image of this
homotopy does not intersect any $\varepsilon _{1}$-Margulis tube.  So
$N_{1}$ is disjoint from all $\varepsilon _{1}$-Margulis tubes.  Now
if $\Delta _{1}$ is sufficiently large given $\nu $, $L_{0}'$ and
$\Delta _{2}$, for $d(z_{0},z_{i})\leq \Delta _{2}$, $\vert \gamma
_{*}\vert \geq \varepsilon _{1}$ for all $\gamma $ such that $\vert
\varphi _{i}(\gamma )\vert \leq L_{0}'$.

\ssubsection{Combinatorially bounded geometry Kleinian surface case: 
lower bounds between pleated surfaces.}\label{9.3}

So now we consider $f_{i}:S\to N$ for $i$ such that 
$d(z_{0},z_{i})\leq \Delta _{2})$. Then we have bounds on 
$d([f_{i}],[f_{i+1}])$, on the geometry of $S(f_{i})$ and also on the 
diameter of the homotopy in $N$ between $f_{i}$ and $f_{i+1}$. The 
following lemma will be applied more generally - always with the 
bounded geometry assumptions given here. It is proved in 
\cite{Min1}, essentially the same proof as given here.

\begin{ulemma} There is an integer $k_{0}=k_{0}(L_{1},\varepsilon 
,C)$ such 
that the 
following holds. Let $f_{i}:S\to N$ ($-m_{-}\leq i\leq m_{+}$) be a 
sequence of homotopic  pleated surfaces such that the pleating locus 
of each $f_{i}$ contains a different maximal multicurve $\Gamma 
_{i}$, 
such that the following hold.
\begin{description}

    \item[1.] 
$f_{i}(S)\subset N_{1}\subset N_{\geq \varepsilon }$ for all $i$. 
\item[2.] $(f_{i})_{*}:\pi 
_{1}(S)\to \pi _{1}(N_{1})$ is injective.
\item[3.]  $[f_{i}]\in  ({\cal{T}}(S))_{\geq \varepsilon }$ for all 
$i$. 
\item[4.] $d([f_{i}],[f_{i+1}])\leq L_{0}$ 
\item[5.] $\vert (\Gamma _{i})_{*}\vert \leq L_{0}$.
\item [6.] There is a homotopy between 
$f_{i}$ and $f_{i+1}$ with tracks of length $\leq L_{0}$ and with 
image contained in $N_{1}$. 

\end{description}
Then if 
$f_{i}(S)$ and $f_{j}(S)$ come within $C$, we have $\vert i-j\vert 
\geq k_{0}$.\end{ulemma}

\noindent {\em{Proof.}} First, note that, given $\Delta $, there is 
$n(\Delta )=n(\Delta ,L_{1},\varepsilon  )$ such that for any $i$ 
with 
$m_{-}+n(\Delta )\leq i\leq m_{+}-n(\Delta )$ the distance between 
$f_{i}(S)$ and $f_{j}(S)$ is $>\Delta $ for at least one $j$ with 
$\vert i-j\vert <n(\Delta )$. For if not we have a large number of 
distinct simple closed geodesics of length $\leq L_{1}$ in a set of 
duameter $\Delta $ is so large that some nontrivial loop must have 
length $<\varepsilon $, which is a contradiction. Now fix a $\Delta $ 
sufficiently large given $L_{1}$ and $C$.

Now suppose for a sequence of hyperbolic manifolds $N_{p}$ all 
satisfying the same combinatorial bounded geometry condition and with 
corresponding pleated surface families $f_{i}^{p}$, that 
$f_{i}^{p}(S)$ 
and $f_{j}^{p}(S)$ come within distance $C$ 
in $N_{p}$ for an arbitrarily large $\vert i-j\vert $ as $p\to \infty 
$. Then because the 
distance between $f_{k}^{p}(S)$ and $f_{k+1}^{p}(S)$ is bounded by 
$L_{0}$ 
for all $k$, and $p$, we can find a sequence $\Delta _{p}$ with 
$\lim _{p\to \infty }\Delta _{p}=+\infty $, and we can find $i(p)$, 
$j(p)$, $k_{q}=k_{q}(p)$ for $1\leq q\leq 4$ such that
$$i(p)-k_{q}(p)\leq n(2\Delta ),\ q=1,\ 2,$$
$$\lim _{p\to \infty }\vert i(p)-j(p)\vert =+\infty ,$$
with $f_{k_{1}}^{p}(S)$ and $f_{k_{2}}^{p}(S)$ on opposite sides of 
$f_{i}^{p}(S)$ distance $\geq \Delta $ away from $f_{i}^{p}(S)$, 
$f_{k_{3}}^{p}(S)$ and $f_{k_{4}}^{p}(S)$ on opposite sides of 
$f_{i}^{p}(S)$ distance $\geq \Delta _{p}$ away $f_{i}^{p}(S)$, but 
$f_{i}^{p}(S)$ and $f_{j}^{p}(S)$ within $C$.

 By \cite{F-H-S}, there 
is $h_{q}^{p}$ homotopic to $f_{k_{q}}^{p}$ such that $h_{q}^{p}:S\to 
N$ is an 
embedding and $h_{q}^{p}(S)$ is in  a small neighbourhood of 
$f_{k_{q}}^{p}(S)$ and distance $\geq \Delta $ from $f_{i}^{p}S$ for 
$q=1$, $2$. 
Then by 
\cite{Wal}, the  submanifold of $N_{p}$ containing $f_{i}^{p}(S)$ and 
bounded 
by $h_{1}^{p}(S)$ and $h_{2}^{p}(S)$ is homeomorphic to 
$S\times [0,1]$, and similarly for the manifold bounded by 
$h_{3}^{p}(S)$ 
and $h_{4}^{p}(S)$. (To 
get into the context of \cite{Wal}, we need to consider the horodisc 
deletion $S_{d}$, but this is a trivial matter.)
Composing with an element of the mapping class group ${\rm{Mod}}(S)$, 
we can assume the $f_{i}^{p}$ all lie in a compact set. Then by the 
bound 
on $n(\Delta )$, the same is true of the $f_{k_{q}}^{p}$ and the 
$h_{q}^{p}$ for $q=1$, $2$. 
So we can take a geometric limit of $(N_{p},f_{i}^{p}(S))$. The limit 
of $f_{i}^{p}$ is a pleated surface $f:S\to N'$, where $N'$ is the 
limit of $N_{p}$, taking a basepoint in $f_{i}^{p}(S)$. Then $f$ is 
injective on $\pi _{1}(S)$, because if $f(\gamma )$ is trivial, then 
for 
sufficiently large $p$ given $\gamma $, $f_{i}^{p}(\gamma )$ is 
freely homotopic to a short loop in 
$N_{p}$ in the region bounded by $f_{k_{3}}^{p}(S)$ and 
$f_{k_{4}}^{p}(S)$. But there are no short loops in this region. We 
can also assume that $f_{k_{1}}^{p}$ and $f_{k_{2}}^{p}$ converge, 
 to pleated surfaces which are homotopic  to $f$. These limiting 
surfaces are homotopic to embeddings and to each other, since 
$f_{k_{1}}^{p}$ and 
$f_{k_{2}}^{p}$ are, and by \cite{F-H-S}, there are embeddings 
homotopic to them with images in arbitrarily small neighbourhoods of 
$f_{k_{1}}^{p}(S)$, $f_{k_{2}}^{p}(S)$. So then, applying 
 \cite{F-H-S} and \cite{Wal} again,  the $\Delta /2$
 neighourhood of $f(S)$ in $N'$ is contained in a submanifold 
 homeomorphic to $S\times [0,1]$, with $f$  homotopic to 
$z\mapsto (z,{1\over 2})$ under the homeomorphism. Now let $\gamma 
_{p}$ be a loop in the pleating locus of $f_{j}^{p}$, which by choice 
of $\Delta $ is in the submanifold of $N'$ homeomorphic to $S\times 
[0,1]$. Then, taking a 
subsequence, we can assume that 
$f_{j}^{p}(\gamma _{p})$ converges to a geodesic $\gamma _{*}$ in 
$N'$ which must be homotopic to $f(\gamma )$ for some $\gamma \in S$, 
because $N'$ is homeomorphic to $S\times 
[0,1]$. 
Then for large $p$, we have two distinct closed geodesics 
$\gamma _{*}$ and $(\gamma _{p})_{*}$  in $N_{p}$ 
of length $\leq L_{1}$ and with distance apart tending to $0$ as 
$p\to \infty $. This is impossible. 
\Box     

\ssubsection{Combinatorially bounded geometry Kleinian surface case: 
construction of biLipschitz map.}\label{9.4}

For 
$\Delta >0$, we let $p_{\pm}$ be the largest integers such that 
$d(z_{0},z_{i})\leq \Delta $ for $-p_{-}\leq i\leq p_{+}$. Then we 
construct a map $\Phi =\Phi _{-p_{-},p_{+}}$ from 
$M(z_{-p_{-}},z_{p_{+}})$ to $N$. The map $\Phi $ maps  $S\times \{ 
t_{j}\} $ to 
$f_{j}(S)$, if 
$d(z_{0},z_{j})=t_{j}$ for $j>0$ and $-t_{j}$ for $j<0$. Since we 
have a bound on $d(z_{j},[f_{j}])$ we can also choose $\Phi $ so  
that this  is boundedly Lipschitz with respect to the Riemannian 
metric on $M(z_{-p_{-}},z_{p_{+}})$. We can then assume without loss 
of generality that $\Phi (z,t_{j})=f_{j}(z)$. The distance between 
$f_{j}(z)$ and $f_{j+1}(z)$ in $N$ is then bounded, because a bounded 
loop in the homotopy class of a bounded geodesic is a bounded 
distance from that geodesic. So we can extend $\Phi $ to $S\times 
[t_{j},t_{j+1}]$ to map $z\times  [t_{j},t_{j+1}]$ to the geodesic 
between $f_{j}(z)$ and $f_{j+1}(z)$. Then $\Phi $ is coarse 
Lipschitz, since it is so restricted to $S\times [t_{j},t_{j+1}]$ for 
all $j$. It is also coarse biLipschitz, for the following reason. 
By \ref{9.3}, if $f_{i}(s)$ and $f_{j}(S)$ are a bounded distance 
apart, then we have a bound on $\vert i-j\vert $. So we only need to 
ensure that  if $\tilde{S}$ and $\tilde{N}$ 
denote the universal covers of $S$, $N$, a sufficiently long path in 
$\tilde{S}\times [a,b]$ will map  under 
$\tilde{\Phi}:\tilde{S}\times [t_{p_{-}},t_{p_{+}}]\to \tilde{N}$ to 
a path which is homotopic, via homotopy preserving endpoints, to a 
geodesic segment of length $\geq 1$ 
 This is  true by  the same argument as in \ref{9.3}. For if this is 
 not true, we can extend length by a bounded amount, and can
 find a sequence of closed geodesics such that 
 $\vert \gamma _{p}\vert \to \infty $, in the metric on 
 $S(f_{i}^{p})$, 
 and yet $f_{i}^{p}(\gamma _{p})\subset N_{p}$ is homotopic to a 
 closed geodesic $(\gamma _{p})_{*}$ of length $\leq L_{2}$, where 
 $L_{2}$ is bounded in 
 terms of $L_{1}$. Normalise as before by applying the modular group, 
  and take a subsequence so that $f_{i}^{p}(\gamma )$ converges for 
all 
  $\gamma $. Since $(\gamma _{p})_{*}$ lies in a fixed compact set it 
  must be $\gamma _{*}$ for some fixed $\gamma $, for all 
  sufficiently large $p$. But this again 
  contradicts Injectivity-on-$\pi _{1}$.

So $\Phi _{-p_{-},p_{+}}$ is uniformly coarse biLipschitz. The 
definition restricted to $S\times [t_{j},t_{j+1}]$ depends only on 
$j$, 
not on $p_{\pm}$ Taking limits as $\Delta \to \infty $, that is as 
$p_{+}\to +\infty $ and $-p_{-}\to -\infty $, we obtain a coarse 
biLipschitz map $\Phi :M(\mu _{-},\mu _{+})\to N$, which is onto, 
because 
for $p_{\pm}$ sufficiently  large given $\Delta_{1}$, the image 
contains the $\Delta _{1}$-neighbourhood of $w_{0}\in N$.

\Box

\section{Lipschitz bounds.}\label{10}

In this section, we do the groundwork for the construction of a 
biLipschitz map from model $M$ to hyperbolic manifold $N$ in complete 
generality. We follow the same general strategy as in the Kleinian 
surface case with combinatorial bounded geometry, which was carried 
out in Section \ref{9}. The main theorem from which Lipschitz bounds 
are deduced is \ref{10.1} 
which is an analogue of the work carried out in 
\ref{9.1}. To simplify the statement, we give two versions, 
the interval bundle case and then the general case.  
Note, in particular, the resemblance between (\ref{9.1.3}) and 
(\ref{10.1.1}). Theorem \ref{10.1} does, of course, take 
more work than \ref{9.1}. 
Theorem \ref{10.2} is the start of an induction to prove it. 
The two theorems are proved in \ref{10.3} to \ref{10.6}. 
There is a also an analogy between \ref{9.2} and 
\ref{10.7}, the respective bounds on Margulis tubes. In this general 
case, there are, of course, Margulis tubes. \ref{10.8} gives a 
relationship betwen the geometry of model Margulis tubes in $M$ and 
Margulis tubes in $N$. The rest of the section is the basic work 
which 
will be needed to get a map from $M$ to $N$ which is coarse 
biLipschitz, as well as Lipschitz, moreover, a map for which the 
biLipschitz constants will be locally uniform in the ending 
lamination 
data.  \ref{10.10} is a fairly straight 
analogue of \ref{9.3}. 

Throughout this section, we fix a compact connected  $3$-dimensional 
\\ submanifold-with-boundary $W$ of $N_{c}$ as in \ref{6.9} and as in 
the 
introduction to Section \ref{8}. As before, we define $N_{d,W}$ to be 
the union of $N_{d}$ and of components of $N\setminus N_{d}$ 
deisjoint 
from $W$. The conditions are such that each end $e$ of $N_{d,W}$ has 
a 
neighbourhood in $N_{d,W}\setminus W$ bounded by a component 
$S_{d}(e)$ of $W\setminus \partial N_{d}$, and $S_{d}(e)$ is the 
horodisc deletion of a finite type surface $S(e)$. We denote by 
$\overline{N}$ the closure of $N$ in the projection of $H^{3}\cup 
\Omega $, where $\Omega $ is the domain of discontinuity of the 
covering group on $\partial H^{3}$.

\begin{theorem}\label{10.1}
 
For each end $e$ of $N_{d,W}$, let
$f_{e,+}:S(e) \to \overline{N}$ be a generalised pleated surface as 
in \ref{4.5}, and homotopic to 
inclusion of $S(e)$ in $N$. Let  $\Gamma _{+}(e)$  be a maximal 
multicurve 
which includes all closed loops in the pleating locus of $f_{e,+}$, 
and suppose that
\begin{equation}\label{10.1.3}\vert f_{e,+}(\Gamma _{+}(e))\vert \leq 
L_{0}.\end{equation}

Let $z_{e,0}=[\varphi _{e,0}]\in {\cal{T}}(S)$. 
Write $y_{e,+}=[f_{e,+}]$. If $\Gamma _{+}(e)$ includes some loops 
which are cusps in $N$, let $y_{e,+}$ be defined as an element of 
${\cal{T}}(S(e))$ so as to minimise $d(z_{e,0},y_{e,+})$ up to an 
additive constant, as in \ref{8.5}.

   Fix ltd parameter functions
$(\Delta _{1},r_{1},s_{1},K_{1})$ and a vertically efficient
ltd-bounded decomposition of 
$S(e)\times [z_{e,0},y_{e,+}]$ with respect to these.
Then the following holds for a constant $L_{1}$, and for 
sufficiently strong 
ltd parameter functions $(\Delta _{1},r_{1},s_{1},K_{1})$, where 
$L_{1}$ depends only on the parameter functions and  $L_{0}$. 

\noindent {\em{Interval Bundle Case.}}

Denote the ends by $e_{\pm }$ and $f_{\pm }=f_{e_{\pm },+}$, 
$y_{\pm }=[f_{\pm }]$, $\Gamma _{+}(e_{\pm})=\Gamma _{\pm }$. Also, 
in 
this case, $z_{e_{+},0}=y_{-}$

For any 
$(\alpha ,\ell )$ in the decomposition of $S(e)  \times 
[y_{-},y_{+}]$, there is 
a pleated surface $f:\alpha \to N$ whose pleating locus incudes 
$\partial \alpha $, such that
$$\vert f(\partial \alpha )\vert \leq L_{1},$$
\begin{equation}\label{10.1.1}d_{\alpha }([f],y)\leq L_{1}{\rm{\ or\ 
}}\vert {\rm{Re}}(\pi 
_{\alpha }([f])-\pi _{\alpha }(y))\vert \leq L_{1},\end{equation}
depending on whether $\alpha $ is a gap or a loop.

\noindent{\em{General case.}}
Suppose that there is a constant $c_{0}>0$  with the following 
property.
Whenever $e$ is a compressible end, $\gamma '$ is a simple  closed 
loop
 which is nontrivial in $S(e)$ but trivial in $N$, and 
$\mu $ and $\mu '$ are geodesic laminations on $S$, with a normalised 
transverse 
invariant measure on $\mu $,
$$i(\mu ,\mu ')=0{\rm{\ and\ }}i(\mu ',\gamma ')\leq c_{0}\vert \gamma \vert ,$$ 
then for at least one $\zeta \in \Gamma _{+}(e)$,
\begin{equation}\label{10.1.2}
    i(\zeta ,\mu )\geq c_{0}\vert \zeta \vert .
\end{equation}

\noindent {\em{Compressible ends.}}
 Let $e$ be a compressible end.
Then (\ref{10.1.1}) holds as above, but 
with $[z_{e,0},y_{e,+}]$ replacing $[y_{-},y_{+}]$, but 
for a constant $L_{1}$, which, this time, 
depends, 
in addition,  $c_{0}$ and locally on $z_{e,0}$. 

\noindent {\em{Incompressible ends.}}
Let $e$ be an {\em{incompressible}} end.  
 Then (\ref{10.1.1}) holds as in the 
Interval Bundle case but with $[z_{e,0},y_{e,+}]$
replacing $[y_{-},y_{+}]$, where $z_{e,0}$ is chosen relative to the 
set of all 
$y_{e',+}$ as described in \ref{6.17}, and also with 
$[\pi _{\omega }(z_{e',0}'),\pi _{\omega }(z_{e,0}')]$ replacing 
$[y_{-},y_{+}]$, whenever $\omega =\omega (e,e')$ is the 
maximal (possibly empty) subsurface of $S(e)$ which is homotopic to a 
subsurface 
$\omega (e',e)$ of $S(e')$, and $z_{e,0}'$, $z_{e',0}'$ are as in 
\ref{6.17}.

\end{theorem}

In the combinatorial bounded geometry Kleinian surface case, the 
interval bundle version of \ref{10.1}, is Theorem \ref{9.1}. 
We shall want to apply the 
interval bundle case of the theorem to a sequence of manifolds, all 
of 
the same interval bundle type. We shall want to apply the general 
case of the theorem 
to a sequence of manifolds of the same topological type
with the same choices of $L_{0}$, $z_{e,0}$, $z_{e',0}$ all along the 
sequence, and of $c_{0}$ in the presence of compressible ends, but 
with 
$\Gamma _{+}(e)$ varying. We shall then 
use the theorem and the follow-up results to obtain uniform 
biLipschitz bounds along the sequence. We shall deal later with the 
question of satisfying the hypotheses of this theorem.

Theorem \ref{9.1} was proved using the first lemma of \ref{7.12}.
The idea of the proof in this general case is to use the second lemma 
of \ref{7.12} to start an inductive proof. The 
induction of \ref{10.1}
is started by proving the following theorem, which is a 
subset of \ref{10.1}.

\begin{theorem}\label{10.2}

Let $f_{\pm }:S\to N$ be homotopic generalised pleated surfaces. Let 
$\Gamma _{\pm }$ be maximal multicurves containing all closed loops 
in 
the pleating loci of $\Gamma _{\pm }$ and let $\vert 
f_{+}(\Gamma _{+})\vert \leq L_{0}$, $\vert 
f_{-}(\Gamma _{-})\vert \leq L_{0}$.
Let $\omega \subset S$ be a subsurface with $\partial \omega $ in 
the pleating locus of $f_{\pm }$ with 
$$\vert (\partial \omega )_{*}\vert \leq L_{0}.$$
Let $(\alpha _{i},\ell _{i})$ ($1\leq i\leq R_{0}$ or $1\leq i\leq 
R_{e,0}$) 
be a totally chain of ltds for $\omega \times [y_{-},y_{+}]$. Thus, 
the 
$(\alpha _{i},\ell _{i})$ are as in 
\ref{7.13}, 
 with $d_ {\omega }$-distance $\leq \Delta _{0}$ between the end of 
$\ell 
_{i}$ and the start of $\ell _{i+1}$, for a fixed constant $\Delta 
_{0}$ depending only on the ltd parameter functions

\noindent {\em{Incompressible surface case}} 
Let $f_{\pm }$ be injective on $\pi _{1}$.

\noindent {\em{Compressible surface case}}
Let $S=S(e)$ for a compressible end. Let $\Gamma _{+}$ satisfy 
(\ref{10.1.2}). Suppose that 
for $z_{e,0}$ as in \ref{10.1} and for $x(.,.)$ as in \ref{7.10},
$$d(z_{e,0},x(z_{e,0},[[f_{-}],[f_{+}]]))\geq D_{2}.$$

Then (\ref{10.1.1})
holds simply for $(\alpha ,\ell )=(\alpha _{i},\ell _{i})$, $1\leq 
i\leq R_{0}$, and for a constant $L_{1,1}$ replacing 
$L_{1}$.

\end{theorem}

\ssubsection{Outline proof of \ref{10.2} in the Interval Bundle 
Case.}\label{10.3}

The idea is to use the second lemma of \ref{7.12}. The difficulty, as 
was pointed out at the time, compared with the first lemma of 
\ref{7.12}, is that, although we can get good bounds straight away at 
some points, we cannot then proceed immediately to other points, 
because a bound on $d_{\alpha _{j}}(x(y_{j}),y_{j})$ for gaps 
$\alpha 
_{j}$ does not imply a bound on $d(y_{j},x(y_{j}))$. 
So the idea, given $(\alpha _{i_{0}},\ell _{i_{0}})$ 
and $y\in \ell _{i_{0}}$, is to get good bounds progressively closer 
to $y$, successively using \ref{7.12}, modifying points to get into a 
position to apply \ref{7.12} again.
For the moment, we concentrate on the interval bundle case. We shall 
consider the general case later.

So the aim is  to find an 
integer $m_{1}$, and  constants $L_{1,1,1}$, $L_{1,1,2}$ depending 
only on the $(\Delta 
_{1},r_{1},s_{1},K_{1})$ and constant $\kappa _{1}$ of 
\ref{8.2.6}, and  
finite sequences  of geodesic segments 
$[y_{m,-},y_{m,+}]$, $[z_{m,-},z_{m,+}]$, and sequences of ltd's
 $(\alpha _{m,\pm },\ell _{m,\pm })$ ($0\leq m\leq m_{1}$)  
with $y_{0,-}=y_{-}=z_{0,-}$, $z_{0,+}=y_{0,+}=y_{+}$, 
$y_{m,\pm }=[f_{m,\pm }]$, 
where $f_{m,\pm}$ is a  pleated surface,  and  
with the following properties. We assume without loss of generality 
that $f_{+}=f_{-}$ off $\omega $.
\begin{description}

    \item[1.] Either $\alpha _{m,+}=\omega $ and $y_{m,+}=y_{+}$ or, 
for 
some $n$, $\alpha _{m,+}=\alpha _{n}$  $\ell _{m,+}$ is a segment of  
$\ell _{n}$ including 
the left endpoint, $z_{m,+}$ is the right end of $\ell _{m,+}$. The 
pleating locus of $f_{m,+}$ includes $\partial \omega \cup \partial 
\alpha _{m,+}$ 
Similar properties hold for $\alpha _{m,-}$, $y_{m,-}$, $f_{m,-}$, 
$z_{m,-}$. For all $m$, $f_{m,+}=f_{m,-}=f_{+}$ off $\omega $.
\item[2.]
$$y\in [z_{m,-},z_{m,+}].$$

\item[3.] 
\begin{equation}\label{10.3.1} \begin{array}{l}
    \vert \varphi _{m,+}(\partial \alpha _{m,+})\vert \leq L_{1,1,1}
{\rm{\ and,\ if\ }}y_{m,+}\neq y_{+},\cr d_{\alpha 
_{m,+}}(y_{m,+},z_{m,+})\leq 
L_{1,1,1}{\rm{\ or\ }}\cr 
\vert {\rm{Re}}(\pi _{\alpha _{m,+}}(y_{m,+})-\pi 
_{\alpha _{m,+}}(z_{m,+}))\vert \leq L_{1,1,1},\cr 
\end{array}\end{equation}
depending on whether $\alpha _{m,+}$ is a gap or a loop,  and 
similarly with $+$ replaced by $-$,and left and right 
interchanged, where $L_{1,1,1}$ depends only on $(\Delta 
_{1},r_{1},s_{1},K_{1})$. Also, if $\alpha _{m,+}$ is a 
gap, the pleating locus 
of 
$f_{m,+}$ in $\alpha 
_{m,+}$ is a maximal multicurve $\Gamma _{m,+}$ with $\vert \varphi 
_{m,+}(\Gamma _{m,+})\vert \leq L_{1,1}$ where $[\varphi 
_{m,+}]=z_{m,+}$, and similarly for $f_{m,-}$, $\alpha _{m,-}$, 
$z_{m,-}$, $y_{m,-}$, if $y_{m,-}\neq y_{-}$.

\item [4.] The following two cases depend on whether or not the 
hypotheses of the second lemma of \ref{7.12} are satisfied.
\begin{description}
    \item [Case 1.]
If
$$d(y_{m,-},y_{m,+})\leq 2d_{\alpha 
_{m,-}.\alpha _{m,+}}'(y_{m,-},
y_{m,+})$$
then either $\alpha _{m+1,+}=\alpha 
_{m,+}$ and $f_{m+1,+}=f_{m,+}$ on $\alpha _{m,+}$ or 
$$d_{\alpha _{m+1,+},\alpha _{m,+}}'(y_{m+1,+},y_{m,+})\geq {1\over 8}
d_{\alpha _{m,-}\alpha _{m,+}}'(y_{m,-},
y_{m,+}),$$
and similarly with $+$ replaced by $-$, and 
\begin{equation}\label{10.3.3}d_{\alpha _{m+1,-}\alpha 
_{m+1,+}}'(y_{m+1,-},
y_{m+1,+})\leq {7\over 8}d_{\alpha _{m,-}\alpha 
_{m,+}}'(y_{m,-},y_{m,+}),\end{equation}
and
\begin{equation}\label{10.3.4}d(y_{m+1,-},y_{m+1,+})\leq
2\kappa _{1}d(y_{m,-},y_{m,+}).\end{equation}

\item [Case 2.] If 
$$d(y_{m,-},y_{m,+})\geq 2d_{\alpha 
_{m,-},\alpha _{m,+}}'(y_{m,-},
y_{m,+})$$
then 
$\alpha _{m+1,+}=\alpha 
_{m,+}$, $f_{m+1,+}=f_{m,+}$ on $\alpha _{m,+}$, $\alpha 
_{m+1,-}=\alpha _{m,-}$, $f_{m+1,-}=f_{m,-}$ on $\alpha _{m,-}$ and
\begin{equation}\label{10.3.5}d(y_{m+1,-},y_{m+1,+})\leq 
{3\over 
4}d(y_{m,-},y_{m,+}).\end{equation}
\end{description}

\item[5.]
$$d(y_{m_{1},-},y_{m_{1},+})\leq L _{1,1,2}^{3},$$
but for $m<m_{1}$
$$d(y_{m,-},y_{m,+})\geq L_{1,1,2}^{2}.$$

\end{description}

If we can satisfy 1-4 then we can obviously find $m_{1}$ such that 5 
holds, because 4 implies that 
$$d_{\alpha _{m,-},\alpha _{m,+}}'(y_{m,-},y_{m,+})$$
is nonincreasing in $m$ and eventually decreasing.Then we obtain 
(\ref{10.1.1}) as 
follows. Choose cell-cutting loop sets $\Gamma _{m_{1},+}\subset 
\alpha _{m_{1},+}$ and 
$\Gamma _{m_{1},-}\subset \alpha _{m_{1},-}$ such that $\vert 
f_{m_{1},+}(\Gamma 
_{m_{1},+})\vert \leq L_{1,1,3}$ (for $L_{1,1,3}$ depending on 
$L_{1,1,1}$, $L_{1,1,2}$), 
and 
similarly for $f_{m_{1},-}$, $\alpha _{m_{1},-}$, $\Gamma 
_{m_{1},-}$. 
Because of the bound on $d([f_{m_{1},-}],[f_{m_{1},+}])
=d(y_{m_{1},-},y_{m_{1},+})$, 
we have 
$$\vert f_{m_{1},+}(\Gamma _{m_{1},-}\cup \partial \alpha 
_{m_{1},-})\vert \leq L_{1,1,4}$$
for $L_{1,1,4}$ depending  on $L_{1,1,i}$, 
$i=1$, $2$, $3$.  By \ref{7.2} and 
\ref{7.3}, $\alpha $ is in the convex hull of $\alpha _{m_{1},+}$ 
and $\alpha _{m_{1},-}$. By the bound on $d(y_{m_{1},-},y_{m_{1},+})$ 
and (\ref{10.3.1}), if we take a maximal multicurve $\Gamma '$ on 
$\alpha  $ such that 
$\vert \varphi (\Gamma ')\vert \leq 
L_{0}$, each loop of $\Gamma '$ can be made out of a finite number of 
arcs of 
$$\Gamma _{m_{1},+}\cup \partial \alpha _{m_{1},+}
\cup \Gamma _{m_{1},-}\cup \partial \alpha _{m_{1},-}.$$
So then we have a bound on $\vert f_{m_{1},+}(\Gamma ')\vert $. Then 
we have a bound on $\vert (\Gamma ')_{*}\vert $. By \ref{4.1}, we can 
assume that $\vert \gamma _{*}\vert $ is bounded from $0$, keeping 
the property that $\vert \varphi (\Gamma ')\vert \leq 
L_{0}$, for $L_{0}$ depending only on the topological type. Extend 
$\Gamma '$ to a maximal multicurve on $S$ with a bound on number of 
intersections with $\Gamma _{m_{1},+}$, depending on $L_{1,1,4}$. Let 
$f$ 
have pleating locus $\alpha $. We have bounds on $d_{\alpha }
([f_{m_{1},+}],[f])$, and hence on 
$d_{\alpha }([f],y)$ or $\vert {\rm{Re}}(\pi _{\alpha }([f])-\pi 
_{\alpha 
}(y))\vert $, (depending on whether $\alpha $ is a gap or a loop), 
and 
on $\vert f(\Gamma )\vert $.

It is crucial to the construction that Case 2 can occur for many 
successive $m$ with no adverse effect on $d_{\alpha 
_{m,+}}(y_{m,+},z_{m,+})$, because $z_{m,+}$ is constant for a string 
of such $m$, and  by the observation at the end of \ref{3.1}, 
$d_{\alpha _{m,+}}(y_{m,+},z_{m,+})$ then depends only on 
$f_{m,+}\vert 
\alpha _{m,+}$: the change in $\pi _{\alpha _{m,+}}([f_{m,+}])$ will 
be slightly affected by $f_{m,+}\vert S\setminus \alpha _{m,+}$, 
but only by a uniformly bounded amount.

\ssubsection{Construction of $y_{m+1,\pm }$ from $y _{m,\pm 
}$.}\label{10.4}

For the 
moment we are considering the Interval Bundle Case.
We need to consider Cases 1 and 2 of \ref{10.3}.  We choose 
ltd parameter functions $(\Delta _{1},r_{1},s_{1},K_{1})$ such that 
all 
the results of Section \ref{6} 
hold. 
Let $\Gamma _{i}$, $f_{i}$  $(0\leq i\leq n_{+}$) be the maximal 
multicurve and 
pleated surface sequences constructed in 
\ref{8.4}, and thus satisfying
\ref{8.2.1}-\ref{8.2.6}. 

We 
proceed to the construction of $y_{m+1,\pm }$ from $y_{m,\pm }$. 
First we consider Case 1.

Let $(\alpha ,y)$ be as in the statement of \ref{10.1} and write 
$y=[\varphi ]$.
Write $y_{m,0}$ for the point in $[y_{m,-},y_{m,+}]$ such that
$$\vert \varphi (\partial \alpha )\vert \leq C(\nu ){\rm{\ or\ 
}}\vert \varphi (\alpha )\vert \leq L_{1,1,1},$$
$$d_{\alpha }(y_{m,0},y)\leq C(\nu ),{\rm{\ or\ }}\vert 
{\rm{Re}}(\pi 
_{\alpha }(y_{m,0})-\pi _{\alpha }(y))\vert \leq L_{1,1,1},$$
depending on whether  $\alpha $ is a gap which is long $\nu $-thick 
and dominant or a loop which is $K_{1}$-flat. Such a point 
$y_{m,0}$ exists by \ref{7.6}, assuming that $L_{1,1,2}$ is 
sufficiently large given the ltd parameter functions. Here, 
$L_{1,1,1}$ plays the role of $L_{1}$ in \ref{7.4} and \ref{7.6}, and 
$L_{1,1,2}$ plays the role of $L_{2}$ in \ref{7.6}. Also by 
\ref{7.6}, there are segments 
$\ell _{i}'\subset [y_{m,0},y_{m,+}]$ a bounded $d_{\alpha 
_{i}}$-distance from $\ell _{i}$ (or similarly if $\alpha _{i}$ is a 
loop: precise statment in \ref{7.6}) for all $(\alpha _{i},\ell 
_{i})$ with $(\alpha ,y_{m,0})\leq (\alpha _{i},\ell _{i})\leq 
(\alpha 
_{m,+},y_{m,+})$, except within distance $L_{1,1,2}$ of the 
endpoints, assuming $L_{1,1,2}$ sufficiently large. Then by 
\ref{7.5} and the properties of the 
$(\alpha _{i},\ell _{i})$ we have, for $\delta _{0}$ depending only 
on $\Delta _{0}$ of \ref{7.13} and the ltd parameter functions, 
$$\delta _{0}d_{\alpha _{m,-},\alpha _{m,+}}'(y_{m,-},y_{m,+})\leq 
\sum _{k}\vert \ell _{k}'\cap [y_{m,-},y_{m,+}]\vert ,$$
$$\delta _{0}d_{\alpha _{m,0},\alpha _{m,+}}'(y_{m,0},y_{m,+})\leq 
\sum _{k}\vert \ell _{k}'\cap [y_{m,0},y_{m,+}]\vert ,$$
$$\delta _{0}d_{\alpha _{m,-},\alpha _{m,0}}'(y_{m,-},y_{m,0})\leq 
\sum _{k}\vert \ell _{k}'\cap [y_{m,-},y_{m,0}]\vert .$$
 
Now we assume without 
loss of generality that 
$$\sum _{k}\vert \ell _{k}'\cap [y_{m,0},y_{m,+}]\vert \geq {1\over 3}
\sum _{k}\vert \ell _{k}'\cap [y_{m,-},y_{m,+}]\vert ,$$
because we have either this, or the corresponding statement with $+$ 
and $-$ interchanged. In this case, we take $\alpha _{m+1,-}=\alpha 
_{m,-}$, $y_{m+1,-}=y_{m,-}$. Then we apply the second lemma of 
\ref{7.12} with $x_{1}=y_{m,0}$, and $x_{2}\in \ell _{j}$ with 
$(\alpha ,x_{1})<(\alpha _{j},x_{2})\leq (\alpha _{m,+},y_{m,+})$, and
$${1\over 4}\sum _{k}\vert \ell _{k}'\cap [y_{m,-},y_{m,+}]\vert \leq 
\sum _{k} \vert \ell _{k}'\cap [x_{1},x_{2}]\vert .$$
In order to apply this strictly, we shall need to interpolate extra 
points between $[f_{k}]$ and $[f_{k+1}]$ where $\Gamma _{k}$ and 
$\Gamma _{k+1}$ are related by a possibly unbounded Dehn twist 
rather than $\# (\Gamma _{k}\cap \Gamma _{k+1})$ being bounded. 
In this case we shall need to use 
some 
of the surfaces $f_{k,m}$ between $f_{k}$ and $f_{k+1}$ as in 
\ref{8.2.6}, to keep 
the property of successive $d(y_{j},y_{j+1})$ being bounded, which is 
the hypothesis of \ref{7.12}.

Then we can find $y_{i}=[f_{i}]$ and $w\in \ell _{j}\cap 
[x_{1},x_{2}]$ 
satisfying 
the conclusion of the lemma in \ref{7.12}. Then $\alpha 
_{m+1,+}=\alpha _{j}$ and $y_{m+1,+}=y_{i}$, $z_{m+1,+}=w$, 
$f_{m+1,+}=f_{i}$, $\Gamma _{m+1,+}=\Gamma _{i}$
satisfy 
the conditions of \ref{10.3}, in particular (\ref{10.3.1}). If 
$\alpha _{j}$ is a gap, which is long, $\nu $-thick and dominant 
aong $\ell _{j}$, then we get the 
bound on $\vert (\Gamma _{i})_{*}\vert $ in the same way as in 
Section \ref{9}. That is, we have a bound on $d_{\alpha 
_{j}}(z_{k},w)$ for some $z_{k}\in \ell _{j}$ as in \ref{8.4}, 
and hence a bound on $d_{\alpha _{j}}(y_{i},z_{k})$. This means that 
we get a bound on $\vert f_{i}(\Gamma _{k'})\vert $ in terms of 
$\nu $ and $p=p(\nu )$ 
for $k'$ in an interval $I$ of  $
\leq p(\nu )$ integers containing $k$. 
If $\Delta _{1}(\nu )$ grows sufficiently fast with $\nu $, as we can 
assume, then we can choose such an interval $I$ with $z_{k'}\in 
\ell _{j}$ for all $k'\in I$, and so that $\cup _{k'\in I}\Gamma 
_{k'}$ cuts $\alpha _{j}$ into cells and annuli parallel to the 
boundary. The bounds on $\vert f_{i}(\Gamma _{k}')\vert $ for all 
$k'\in I$ 
then give a bound, depending on $\nu $, on $d_{\alpha 
_{j}}([f_{k}],[f_{i}])$.  We can then assume that $f_{i}=f_{k}$.
We get condition \ref{10.3.4} by the properties of the $f_{i}$, in 
particular \ref{8.2.4} and \ref{8.2.6}

Now we consider case 2. 

The point $y_{m+1,\pm }$ will be 
obtained from a sequence $y_{m,i,\pm}$ with $y_{m,0,\pm }=y_{m,\pm}$. 
Also, 
$y_{m,i,\pm }=[f_{m,i,\pm }]$ for a pleated surface  
$f_{m,i,+}$ with $f_{m,0,\pm }=f_{m,\pm }$. For each $i$, we shall 
have either 
$y_{m,i,-}=y_{m,i+1,-}$ or $y_{m,i,+}=y_{m,i+1,+}$, but the other one 
will be different. There will be a finite increasing 
sequence of subsurfaces $\alpha _{m,i,\pm }$ 
($i\geq 0$) such that the pleating locus of $f_{m,i,+}$ 
includes $\alpha _{m,i,+}$, and either $\alpha _{m,i+1,+}$ is 
strictly 
larger than $\alpha _{m,i,+}$, or $\alpha _{m,i+1,-}$ is strictly 
larger than $\alpha _{m,i,-}$. Then for some $r\leq -2\chi (S)$, 
where $\chi $ denotes Euler characteristic, we 
have $\alpha _{m,r,\pm }=S$ and $f_{m+1,\pm }=f_{m,r,\pm }$. We 
shall always have  $f_{m,i+1,+}=f_{m,i,+}$ 
on $\alpha 
_{m,i,+}$, 
and similarly with $+$ replaced 
by $-$. For $i<r$ $\alpha _{m,i,+}$ will be a union of $\beta $ such 
that $(\beta ,\ell ')$ is an ltd for $[y_{m,-},y_{m,+}]$.

 For a suitable 
$L_{1,1,5}$, we shall always take 
$$d(y_{m,i+1,-},y_{m,i+1,+}) \leq 
d(y_{m,i,-},y_{m,i,+})+L_{1,1,5}$$
$$\leq d(y_{m,-},y_{m,+})+(i+1)L_{1,1,5},$$
and  
$$d_{\alpha _{m,i+1,-},\alpha _{m,i+1,+}}'(y_{m,i+1,-},y_{m,i+1,+})$$
$$\leq {\rm{Max}}\left( d_{\alpha _{m,i,-},\alpha 
_{m,i,+}}'(y_{m,i,-},y_{m,i,+}),{1\over 
2}d(y_{m,i,-},y_{m,i,+})\right) $$
$$\leq {\rm{Max}}\left( d_{\alpha _{m,-},\alpha 
_{m,+}}'(y_{m,-},y_{m,+}),{2\over 3}d(y_{m,-},y_{m,+})\right) $$
$$\leq 
{2\over 3}d(y_{m,-},y_{m,+}).$$
 We then take $y_{m+1,+}=y_{m,r,+}$ 
and $y_{m+1,-}=y_{m,r,-}$ 
for the first $r$ such that
\begin{equation}\label{10.4.1}{8\over 9}d(y_{m,r,-},y_{m,r,+}))\leq 
d_{\alpha _{m,r,-},\alpha 
_{m,r,+}}'(y_{m,r,-},y_{m,r,+})\leq {2\over 
3}d(y_{m,-},y_{m,+}).\end{equation}
This will certainly be true if $\alpha _{m,r,+}=S$, which 
happens for some $r\leq -2\chi (S)$, but may happen 
earlier. This will ensure (\ref{10.3.5}), as required, provided that 
$L_{1,1,5}$ is sufficiently large given $L_{1,1,2}$ and $\chi (S)$.

Now we need to consider how to define $y_{m,i+1,\pm}$ and $\alpha 
_{m,i+1,\pm }$ when (\ref{10.4.1}) does not hold for $i$ replacing 
$r$.
Then we have two cases to consider: when
\begin{equation}\label{10.4.2}  d_{\alpha 
_{m,i,-},\alpha 
_{m,i,+}}'(y_{m,i,-},y_{m,i,+})\leq d_{\alpha 
_{m,i,-}}'(y_{m,i,-},y_{m,i,+})-{1\over 18}d(y_{m,i,-},y_{m,i,+}), 
\end{equation}

and the case when (\ref{10.4.2})  does not hold, but
\begin{equation}\label{10.4.3}d_{\alpha 
_{m,i,-}}'(y_{m,i,-},y_{m,i,+})\leq 
d(y_{m,i,-},y_{m,i,+}) -{1\over 18}d(y_{m,i,-},y_{m,i,+}).
\end{equation}

If (\ref{10.4.2}) holds, we take 
$y_{m,i+1,-}=y_{m,i,-}$, and need to define $y_{m,i+1,+}$. 
So now consider the case of (\ref{10.4.2}) holding. 
The case of (\ref{10.4.3}) is similar, with $S$ replacing $\alpha 
_{m,i,-}$ and $\alpha _{m,i,-}$ replacing $\alpha _{m,i,+}$.

We use the geodesic segment $[y_{m,i,-},y_{m,i,+}]$. 
We consider the set $\cal{A}$ of all $(\beta ,\ell )$ such that 
$\ell \subset [y_{m,i,-},y_{m,i,+}]$ and  
$$(\alpha _{m,i,-},y_{m,i,-})<(\beta ,\ell ),\ \ \beta \cap \alpha 
_{m,i,+}=\emptyset .$$
If $(\beta _{1},\ell _{1}')<(\beta _{2},\ell _{2}')$ and $(\beta 
_{1},\ell _{1}')\in \cal{A}$ then $(\beta 
_{2},\ell _{2}')\in \cal{A}$ also, by \ref{7.3} applied to each of 
the 
triples 
$$\{ (\alpha _{m,i,-},y_{m,i,-}),
(\beta _{1},\ell _{1}'),(\beta _{2},\ell _{2}')\} ,\ \  
\{ (\beta _{1},\ell _{1}'),(\beta _{2},\ell _{2}'), 
(\alpha _{m,i,+},y_{m,i,+})\} .$$
Now for $z\in 
[y_{m,i,-},y_{m,i,+}]$, 
let $C(z)$ denote the convex hull (\ref{5.7}) of the $\beta $ with 
$(\beta ,\ell ')\in \cal{A}$ for some $\ell '$ with $\ell '\cap 
[y_{m,i,-},z]\neq \emptyset $. Then $C(z)$ increases with $z$ and 
takes only finitely many values. Let $z_{m,i,j,+}$ be the successive 
points on $[y_{m,i,-},y_{m,i,+}]$ such that $C(z)=\omega _{m,i,j,+}$
is constant for $z\in [z_{m,i,j-1,+},z_{m,i,j,+}]$ with 
$z_{m,i,0,+}=y_{m,i,-}$, $z_{m,i,t,+}=y_{m,i,+}$, $t=t(m,i,+)$, where 
$t$ is bounded by $-2\chi (S)$. Suppose that
$(\beta ,\ell ')$ is ltd for $\ell '\subset 
[z_{m,i,j-1,+},z_{m,i,j,+}]$ and 
$\beta \cap \omega _{m,i,j,+}\neq \emptyset $. Then $\beta \cap \beta 
'\neq \emptyset $ for any $(\beta ',\ell '')\in {\cal{A}}$ with 
$\ell ''\subset [z_{m,i,j-1,+},z_{m,i,j,+}]$. There is at least one 
such $(\beta ',\ell '')$ with $(\beta ',\ell '')\leq (\beta ,\ell ')$.
So $(\beta ,\ell ')\in {\cal{A}}$ 
by the order-closure property, and hence $\beta \subset \omega 
_{m,i,j,+}$.  It 
follows that for a constant $L$ depending only on the ltd parameter 
functions, $\vert \varphi (\partial \omega _{m,i,j,+})\vert \leq L$ 
for all $[\varphi ]\in [z_{m,i,j-1,+},z_{m,i,j,+}]$. 
Then, by (\ref{10.4.2}),
$$d_{\alpha _{m,i,-}}'(y_{m,i,-},y_{m,i,+})\geq 
{1\over 18}d(y_{m,i,-},y_{m,i,+}).$$
Now
$$d_{\alpha _{m,i,-}}'(y_{m,i,-},y_{m,i,+})$$
$$={\rm{Max}}(d_{\alpha 
_{mi,-},\alpha _{m,i,+}}'(y_{m,i,-},y_{m,i,+}),d_{\alpha 
_{m,i,-},S\setminus\alpha _{m,i,+}}'(y_{m,i,-},y_{m,i,+}))+O(1)$$
$$=
d_{\alpha 
_{m,i,-},S\setminus\alpha _{m,i,+}}'(y_{m,i,-},y_{m,i,+})+O(1).$$
Using (\ref{10.4.2}) for the last equality. So
 $d_{\alpha _{m,i,-}}'(y_{m,i,-},y_{m,i,+})$ is bounded above and 
 below, 
 up to an additive constant depending only on the 
ltd parameter functions, by
$$\sum _{j=1}^{t}d_{\omega _{m,i,j,+}}(z_{m,i,j-1,+},z_{m,i,j,+}).$$
So for some $j$ and $\lambda _{1}=(.01)/(-\chi (S))$, we have
$$d_{\omega _{m,i,j,+}}(z_{m,i,j-1,+},z_{m,i,j,+})\geq 
{1\over 20}\lambda _{1}d(y_{m,i,-},y_{m,i,+}),$$
$$\sum _{k=1}^{j-1}d_{\omega 
_{m,i,k,+}}(z_{m,i,k-1,+},z_{m,i,k,+})<{1\over 3}
d(y_{m,i,-},y_{m,i,+}).$$
Now by \ref{7.7}, the ltd's $(\beta ,\ell ')$ along 
$[z_{m,i,j-1,+},z_{m,i,j,+}]$ with 
$\beta \subset \omega _{m,i,j,+}$ are in natural correspondence with 
ltd's for $[\pi (z_{m,i,j-1,+}),\pi (z_{m,i,j,+})]$ where $\pi 
=\pi _{\omega _{m,i,j,+}}$. So then we can apply \ref{7.13} to 
$[\pi (z_{m,i,j-1,+}),\pi (z_{m,i,j,+})]$, to obtain a totally 
ordered 
sequence of ltd's, and then by \ref{7.7} we have a corresponding 
totally ordered subset along $[z_{m,i,j-1,+},z_{m,i,j,+}]$. 
We construct 
a new family of pleated surfaces for $[y_{m,i,-},y_{m,i,+}]$ as in 
\ref{8.4}, with the properties of \ref{8.2}. The  
the hypotheses of the lemma in \ref{7.12} holds for $y_{\pm 
}=y_{m,i,\pm }$,  $x_{1}=z_{m,i,j-1,+}$ 
and $x_{2}=z_{m,i,j,+}$.   
Then as in case 1, we can 
find $z\in [z_{m,i,j-1,+},z_{m,i,j,+}]$ and $\zeta $ ltd along a 
segment containing $z$, $\zeta \subset \omega _{m,i,j,+}$ and $w=[f]$ 
a pleated surface whose 
pleating locus includes $\partial \zeta $, such that
$$d_{\alpha _{m,i,-},\zeta }'(y_{m,i,-},z)<{1\over 2}d_{\alpha 
_{m,i,-}}'(y_{m,i,-},y_{m,i,+}),$$
and (\ref{10.3.1}) holds, with $w$ replacing $y_{m,+}$, $z$ replacing 
$z_{m,+}$ and $\zeta $ replacing $\alpha _{m,+}$.
Then we put $\alpha _{m,i+1,+}=\alpha _{m,i,+}\cup \zeta $ and define
$f_{m,i+1,+}=f_{m,i,+}$ on $\alpha _{m,i,+}$ and $=f$ on $S
\setminus \alpha _{m,i,+}$ , which includes $\zeta $. 

\Box

\ssubsection{Proof of \ref{10.1} in the interval bundle 
case.}\label{10.5}

We shall now prove \ref{10.1} by induction. We again restrict to the 
interval bundle case. We shall show that the 
statement 
of \ref{10.1} holds for all $(\alpha ,\ell )$ in the decomposition,
by induction on the topological 
type of $S$. In fact, in order to carry out the induction, we shall 
prove something a bit more general than \ref{10.1}. Fix a 
subsurface $\omega $ of $S$ such that all components of 
$\partial \omega $ are nontrivial and $\vert (\partial \omega 
)_{*}\leq L_{0}$. 
As usual, given a nontrivial nonperipheral loop $\gamma \subset S$, 
$\gamma _{*}$ denotes the geodesic representative in $N$, up to free 
homotopy. We shall show 
that \ref{10.1} holds if $f_{\pm }:S\to N$ are pleated surfaces whose 
pleating locus includes $\partial \omega $, and $\alpha \subset 
\omega $
If $\omega $ is  a simple loop on $S$, then the theorem is trivially 
true. So now suppose inductively that \ref{10.1} holds with $S$ 
replaced by any proper essential subsurface $\omega '$ properly 
contained in $\omega $, for any choice of $f_{\pm }$ and 
$z_{0}=y_{-}$. Now let $N$ 
and $f_{\pm }:S\to N$ be given, where $\partial \omega $ is in the 
bending locus of both $f_{\pm }$ and $\vert (\partial \omega 
)_{*}\vert  \leq L_{0}$. By \ref{7.13}, and \ref{7.7} with $\alpha $ 
replaced by $\omega $, there is at least one 
sequence of 
$(\ell _{i},\alpha 
_{i})$ as in \ref{7.13}, with $\alpha _{i}\subset \omega $ for all 
$i$, to which we can apply \ref{10.2}. As already indicated, the 
proof of \ref{10.1} uses \ref{10.2} as the base for an induction.
So \ref{10.2} gives \ref{10.1} for all $(\alpha ,\ell )$ with 
$(\alpha ,\ell )=(\alpha _{i},\ell _{i})$, $1\leq i\leq R_{0}$ for a 
totally ordered set $\{ (\alpha _{i},\ell _{i}):1\leq i\leq R_{0}\} $ 
as in \ref{7.13}, with $L_{1,1}$ replacing $L_{1,1,1}$. Write 
$B_{i,1}=\{ (\alpha _{i},\ell _{i})\} $. 
Inductively, for some integer $m_{1}$, we are going to define sets 
$B_{i,m}$, 
$1\leq i\leq 
R_{0}(m)$, $1\leq m\leq m_{1}$ and $B_{i,m}'$, $1\leq i\leq 
R_{0}'(m)$, 
$2\leq m\leq m_{1}$,
of pairs $(\alpha 
,\ell )$ in a fixed vertically efficient decomposition of $\omega 
\times [y_{-},y_{+}]$, with the following 
properties.

\begin{description}
    
\item[1.] $B_{i,m}$ and $B_{j,m}'$ are order-closed sets of pairs 
$(\ell ,\beta )$ with $\ell \subset [y_{-},y_{+}]$ and $\beta \subset 
\omega $. 

\item[2.] Let $\alpha _{i,m}$, $\alpha _{i,m}'$ to be the convex hull 
of those $\beta $ with $(\ell ,\beta )\in B_{i,m}$ or $\in B_{i,m}'$. 
For each $i$ and $m$, $B_{i,m}$ is the union of  
 a maximal number of consecutive sets $B_{j,m}'$ for which $\alpha 
 _{j,m}$ are the same. 
and $B_{i,m+1}'$ is the  union of  $B_{i,m}$ and $B_{i+1,m}$, 
together 
with any $(\ell, \beta )$ sandwiched between them. For each $m$, 
all the sets $B_{i,m}$ are disjoint, with 
$B_{i+1,m}$ the next greatest in the ordering after $B_{i,m}$.
 
\end{description}

Inductively, we shall show the following. 

Let $(\ell ,\beta )\in B_{i,m}$ and let $y=[\varphi 
]\in \ell $. Take any collection $\Gamma $ of disjoint simple loops 
in 
$\beta $ with 
$\vert \varphi (\Gamma )\vert \leq L_{0}$. Then there is a 
pleated surface $f$ with pleating locus including 
$\partial \beta \cup \Gamma $ and with $L_{1,m}$-Lipschitz impression 
such that
$$\vert (\partial \beta \cup \Gamma )_{*}\vert \leq L_{1,m},$$
and
\begin{equation}\label{10.5.1}d_{\beta }(y,[f])
    \leq L_{1,m}{\rm{\ or\ }}\vert {\rm{Re}}(\pi _{\beta }([g])-\pi 
_{\beta }(y))\vert \leq  L_{1,m},\end{equation}
depending on whether $\beta $ is a gap or a loop. Similar properties 
hold for $B_{j,m}'$

There are  pleated 
surfaces $f_{i,m,\pm}$ 
whose pleating loci include $\omega $ and $\partial \alpha _{i,m}$ 
and 
 pleated 
surfaces $f_{i+1,m,+}'$ and $f_{i,m,-}'$ 
whose pleating loci include $\omega $ and $\partial \alpha _{i,m}$ 
such that 
\begin{equation}\label{10.5.2}
    \begin{array}{l}
   d([f_{i,m,+}],x(E_{i,m,+}))\leq L_{1,m},\ 
    d([f_{i,m,-}],x(E_{i,m,-}))\leq L_{1,m},\cr  
\vert (\partial\alpha _{i,m})_{*}\vert \leq L_{1,m},\cr \end{array}
\end{equation}
\begin{equation}\label{10.5.3}\begin{array}{l}    
    d([f_{i+1,m,+}'],x(E_{i,m,+}'))\leq 
L_{1,m+1},\ d([f_{i,m,-}'],x(E_{i,m,-}'))\leq L_{1,m+1}\cr 
\vert (\partial\alpha _{i,m}')_{*}\vert \leq L_{1,m+1}.\cr \end{array}
\end{equation}
Here, $E_{i,m,+}$ is the order splitting defined using maximal 
elements 
of $B_{i,m,+}$ and so on. The notation $f_{i+1,m,+}'$ and 
$f_{i,m,-}'$ 
is deliberate, since these are defined using respectively maximal 
elements of $B_{i+1,m}$ and some maximal elements of $B_{i,m}$, and 
minimal elements of $B_{i,m}$ and some minimal elements of 
$B_{i,m,-}$. 

By \ref{10.2}, (\ref{10.5.1}) and (\ref{10.5.2}) are satisfied for 
the $B_{i,1}$. In 
general, 
the properties above suffice to define the $B_{i,m}'$ from the 
$B_{j,m}$, and $B_{i,m+1}$ from the $B_{j,m}'$. In order to keep the 
$B_{i,m}$ disjoint, we may discard some of the sets $B_{j,m}'$, 
because the sets $B_{j,m}'$ overlap by definition. In the case of 
$B_{1,m}'$ we also include those $(\beta ,\ell)$ which are $<(\beta 
_{1},\ell _{1})$ for some $(\beta _{1},\ell _{1})\in B_{i,m}\subset 
B_{1,m}'$. 
We make similar inclusions in 
$B_{k,m}'$ for $k=R_{0}'(m)$. By \ref{7.3}, $\alpha _{j,m}'$ is 
simply the convex 
hull of the $\alpha _{i,m}$ with $B_{i,m}\subset B_{j,m}'$. 
 Eventually we 
reach an $m=m_{1}$ such that there is just one set $B_{1,m_{1}}$. 
This happens 
for $m_{1}$ bounded in terms of the topological type of $\omega $, 
because 
$\alpha _{i,t}$ is properly contained in $\alpha _{j,t+1}$ for 
$B_{i,t}$ contained in $B_{j,t+1}$. Then $(\beta ,\ell )\in 
B_{1,m_{1}}$ for all $(\beta ,\ell )$ in the decomposition of 
$\omega \times [y_{-},y_{+}]$ with $\beta \subset \alpha _{1,m_{1}}$.
 Then the proof is finished if $\alpha 
_{1,m_{1}}'=\omega $ --- and also if it is properly contained in 
$\omega 
$, 
because in the latter case we can apply the inductive hypothesis of 
\ref{10.1} to $\omega \setminus \alpha _{1,m_{1}}'$.

Note that (\ref{10.5.2}) for $B_{j,m+1}$ follows from (\ref{10.5.3}) 
for $f_{i,m}'$ for which $B_{i,m}'\subset B_{j,m+1}$, since we can 
take $f_{j,m,+}=f_{i,m,+}'$ for the largest $i$ with $B_{i,m}\subset 
B_{j,m}$ (so that $B_{i-1,m}'\subset B_{j,m+1}$) and 
$f_{j,m,-}=f_{k,m,-}'$ for the smallest $k$ with $B_{k,m}\subset 
B_{j,m+1}$, $B_{k,m}'\subset 
B_{j,m+1}$. We claim that (\ref{10.5.3}) for $B_{i,m}'$
 follows from (\ref{10.5.1}) for $B_{i,m}$ and $B_{i+1,m}$.  The 
construction of $f_{i+1,m,+}'$ and $f_{i,m,-}'$ are similar, so we 
 consider $f_{i+1,m,+}'$. The first step is to define the pleated 
surfaces 
$f_{i,m,+}$ and $f_{i+1,m,-}$. We can choose these off $\alpha 
_{i+1,m}$ 
and $\alpha _{i,m}$ respectively so that the 
pleating loci have number of intersections bounded in 
terms of the ltd parameter functions, and so that $\partial \alpha 
_{i,m}'$ is in the pleating locus of both and so that they have no 
badly bent annuli. We assume now that this 
has been done. 
Then by \ref{4.4}
$$d([f_{i,m,+}],[f_{i+1,m,-}])\leq L_{1,m}'$$
for $L_{1,m}'$ depending on $L_{1,m}$ and $(\Delta 
_{1},r_{1},s_{1},K_{1})$. It follows that   
$$\vert f_{i+1,m,-}(\partial \alpha _{i,m}')\vert \leq L_{1,m}''$$
for $L_{1,m}''$ depending on $L_{1,m}$, $L_{1,m}'$ and $(\Delta 
_{1},r_{1},s_{1},K_{1})$, and this gives the required bound on $\vert 
(\partial \alpha 
_{i,m}')_{*}\vert $. Then we take $f_{i+1,m,+}'=f_{i+1,m,-}$ off 
$\alpha _{i+1,m}$ and $f_{i+1,m,+}'=f_{i+1,m,+}$ on $\alpha 
_{i+1,m}$. 
 Then $f_{i+1,m,+}'$ has the properties required for 
(\ref{10.5.3}).

So it remains to prove (\ref{10.5.1}) for each 
$B_{i,m}$ and $B_{i,m}'$, by induction. The proof for $B_{i,1}$ is 
given by \ref{10.2}.
As already noted, 
(\ref{10.5.1}) for $B_{i,m}$ and $B_{i+1,m}$ implies (\ref{10.5.2}) 
for $B_{i,m}'$.  So we need to obtain 
(\ref{10.5.1}) for $B_{i,m}'$  
from (\ref{10.5.1}) for $B_{i,m}$ and $B_{i+1,m}$, 
and from (\ref{10.5.3}) for $B_{i,m}'$. Similarly, we shall obtain 
(\ref{10.5.1}) for $B_{j,m+1}$ from (\ref{10.5.1}) and (\ref{10.5.3}) 
for the $B_{i,m}'$ contained in $B_{j,m+1}$. First we consider 
(\ref{10.5.1}) for $B_{i,m}'$ when we already have (\ref{10.5.3}). 
 Now we choose two further 
 pleated surfaces $f_{i,m,+,+}$ and $f_{i,m,-,-}$ as 
follows.  An easy 
induction shows that each 
$B_{i,m}$ contains a 
unique maximal element of the form $(\ell _{n},\alpha _{n})$ for some 
$n\leq R_{0}$. Similarly, $B_{i+1,m}$ has a unique minimal element 
which is 
then $(\ell _{n+1},\alpha _{n+1})$, for the same $n$. Then let 
$C(i,i+1,m)$ denote the convex hull of $\alpha _{n}$ and $\alpha 
_{n+1}$. We can then ensure that the pleating locus of $f_{i,m,+}'$ 
includes $\partial \alpha _{n}$ and that the pleating locus of 
$f_{i+1,m,-}'$ includes $\partial \alpha _{n+1}$. The distance 
between 
the right end of $\ell _{n}$ and the 
left end of $\ell _{n+1}$ in $[y_{-},y_{+}]$ is $\leq \Delta _{0}$.
So we can choose $f_{i,m,+,+}$ to have pleating locus including 
$\partial \alpha _{n}$, $\partial 
C(i,i+1,m)$,  $\partial \alpha '_{i,m}$ and 
$\partial \omega $,  with 
(\ref{10.5.1})
satisfied for $[f]=[f_{i,m,+,+}]$, $\beta =\alpha _{n}$ and $y$ the 
right end of $\ell _{n}$. We similarly choose 
$f_{i+1,m,-,-}$ 
to have pleating locus including 
$\partial \alpha _{n+1}$, $\partial 
C(i,i+1,m)$ (and  $\partial \alpha '_{i,m}$, $\partial 
\omega $)  with (\ref{10.5.1})
satisfied for $[f]=[f_{i+1,m,-,-}]$, $\beta =\alpha _{n+1}$ and $y$ 
the 
left end of $\ell _{n}$. By \ref{4.1}, \ref{4.2}, we can choose the
 maximal multicurves in the 
pleating loci of 
$f_{i,m,+,+}$ and $f_{i,m,-,-}$ to have boundedly many intersections, 
 with bound 
depending on $\Delta _{0}$, 
and with corresponding geodesics of length bounded from $0$,and 
hence, by \ref{4.4},
$$d([f_{i,m,+,+}],[f_{i+1,m,-,-}])\leq 
L_{0}'.$$
Indeed, we can choose $\pi _{S\setminus C}([f_{i,m,+,+}])=\pi 
_{S\setminus C}([f_{i+1,m,-,-}])$ for $C=C(i,i+1,m)$.

For $i=1$, we can choose $f_{1,m,-,-}$ to have pleating locus 
including $\partial \alpha _{1}$ so that 
$$d([f_{1,m,-,-}],[f_{-}])\leq L_{0}'.$$
We can also choose $f_{1,m,-}$ to have pleating locus including 
$\partial \alpha _{1}$. We can make similar conditions on $f_{k,m,+}$ 
and $f_{k,m,+,+}$ for  maximal $k$.

So we now have pleated surfaces in order:
$$f_{i,m,-}',\ f_{i,m,+}',\ f_{i,m,+,+},\ f_{i+1,m,-,-},\ 
f_{i+1,m,-}',\ f_{i+1,m,+}',$$
all of which have $\partial \omega $ in their pleating locus and the 
consecutive pairs have also have as common pleating locus respectively
$$\partial \alpha _{i,m},\ \partial \alpha _{n}\cup \partial 
\alpha _{i,m}',\ \partial 
\alpha _{i,m}'\cup C(i,i+1,m),$$
$$\partial \alpha _{n+1}\cup \partial 
\alpha _{i+1,m}',\ \partial \alpha _{i+1,m}.$$

 Now take any $(\ell ,\beta )\in B_{i,m}'$. We need 
to 
show that (\ref{10.5.1}) holds for this $(\ell ,\beta )$. We may as 
well assume that $(\ell 
,\beta )\notin B_{i,m}\cup B_{i+1,m}$, since otherwise there would 
be nothing to prove. Then we consider the five
successive 
intervals 
$$[[f_{i,m,-}'],[f_{i,m,+}']],\ [[f_{i,m,+}'],[f_{i,m,+,+}]],\ 
[[f_{i,m,+,+}],[f_{i+1,m,-,-}]],$$
$$[[f_{i+1,m,-,-}],[f_{i+1,m,-}']],\ 
[[f_{i+1,m,-}'],[f_{i+1,m,+}']].$$
Then $\ell $ splits into at most $5$ segments, each a bounded 
$d_{\beta }$ distance from a segment of one of these intervals, by 
the results of \ref{7.4} and \ref{7.6} with the usual 
modifications if $\beta $ 
is a loop. The interval $[[f_{i,m,+,+}],[f_{i+1,m,-,-}]]$ is 
bounded, so we 
only need to consider the other intervals, such as 
$[[f_{i,m,+}'],[f_{i,m,+,+}]]$ and $[f_{i+1,m,-}'],[f_{i+1,m,+}']]$. 
On 
each of these we can apply the inductive hypothesis of 
\ref{10.1}, to $\omega \setminus \alpha _{n}$ on 
$[[f_{i,m,+}'],[f_{i,,m,+,+}]]$, and to $\omega \setminus \alpha 
_{i+1,m}$ on $[f_{i+1,m,-}'],[f_{i+1,m,+}']]$ --- assuming as we may 
do 
that $(\ell ,\beta )\notin B_{i+1,m}$. Then the inductive hypothesis 
of \ref{10.1} gives (\ref{10.5.1}) for $(\ell ,\beta )$. 

The proof of (\ref{10.5.1}) for $B_{j,m+1}$ is then immediate, 
because 
$B_{j,m+1}$ is the union of the $B_{i,m}'$ with both $B_{i,m}\subset 
B_{j,m+1}$ and $B_{i+1,m}\subset B_{j,m+1}$. 
 \Box
 
 \ssubsection{Modifications of \ref{10.4}, \ref{10.5} in the 
 case of a compressible end.}\label{10.6}
 
 We now show how to modify the proofs of \ref{10.2} and \ref{10.1} 
for a {\em{compressible}} end $e$. Write $z_{e,0}=[\varphi _{e,0}]$. 
Consider the sequences $\Gamma _{i}(e)$ and $f_{e,i}$ ($0\leq i\leq 
n_{+}(e)$) of \ref{8.6}, with $f_{e,n_{+}(e)}=f_{e,+}$. 
We have seen in \ref{8.6} that $z_{e,0}$ can be chosen (moving a 
bounded distance if necessary) so that 
$$d(z_{e,0},x(z_{e,0}))\geq D_{2},$$
where $D_{2}$ is the constant of \ref{8.8}
 We have also seen that, given a constant $D_{0}$, we can ensure that.
moving $z_{e,0}$ a 
bounded distance depending on $D_{0}$ if necessary, we assume that 
for all $[\varphi ]\in [z_{e,0},y_{e,+}]$, if $\gamma \subset S(e)$ 
is 
nontrivial in $S(e)$ but trivial in $N$, then
$$\vert \varphi (\gamma )\vert \geq D_{0}.$$
We assume that $z_{e,0}$ was chosen in this way in the first place.

We now apply  the proofs of \ref{10.1} and \ref{10.2} in 
\ref{10.3} to \ref{10.5} with $f_{e,+}$, $f_{e,0}$ replacing 
$f_{\pm }$ and $y_{e,+}$, $z_{e,0}$ replacing $y_{\pm }$. 
We consider what other modifications are necessary for such ends, 
in both \ref{10.2} and \ref{10.1}. First we consider \ref{10.2}. 
Let $(\alpha _{i},\ell _{i})$, $1\leq i\leq R_{0}$ be exactly as in 
the 
statement of \ref{10.2}. We drop the index $e$.
Then by 
\ref{8.8}, if $D_{0}$ is sufficiently large given $D_{2}$ and 
$D_{1}$, $\vert f_{j}(\gamma )\vert \geq D_{1}$ whenever 
$T([f_{j}],+,[z_{0},y_{+}])$ (in the notation of \ref{7.9}) does not 
contain $(\alpha _{1},\ell _{1})$. This is automatically true for 
$j=n_{+}$.

So then choose $n_{2}$ so that containment does not happen for $j\geq 
n_{2}$ 
but does happen 
for $j=n_{2}-1$. For ease of notation, renumber so that 
$g_{n_{2}}=g_{1}$. Then there are $\ell _{j}'\subset 
[[f_{0}],y_{+}]$ for $1\leq j\leq R_{0}$ such that $\ell _{j}'$ 
is within a bounded 
$d_{\alpha _{j}}$-distance of $\ell_{j}$. 
Then, using this sequence, we have all the 
properties 
of \ref{8.2}, including \ref{8.2.6}, for $j\geq n_{2}$. We also have
$$\sum _{j=0}^{n_{+}-1}d([f_{j}],[f_{j+1}])\leq 
\kappa _{1}d([f_{1}],y_{+}),$$
$$\delta _{0}d([f_{1}],y_{+})\leq \sum _{j=1}^{R_{0}}\vert \ell 
_{j}'\vert $$
as before. We then 
proceed with the construction of $y_{m,\pm }=[f_{m,\pm }]$, $z_{m,\pm 
}$, $(\alpha 
_{m,\pm },\ell _{m,\pm })$  and $\Gamma _{m,\pm}$ as in \ref{10.3} 
and \ref{10.4}, but with some minor differences. We have 
$z_{0,-}=z_{0}$ and $y_{0,-}=[f_{1}]$, but probably not 
$y_{0,-}=z_{0,-}$. 
Also, so long as $y_{k,-}=y_{0,-}$ for $k\leq m$, we do not try to 
bound $\vert f_{k,-}(\Gamma _{k,-})\vert $ or $d_{\alpha 
_{k,-}}(z_{k,-},y_{k,-})$
The proof of \ref{10.2} is 
exactly 
as before, and so is the proof of \ref{10.1}. Even when we change to 
other sequences of pleated surfaces, as happens in case 2 in 
\ref{10.5}, and again in the proof of \ref{10.1}, we never change the 
definition of the $f_{j}$ on surfaces intersecting $\alpha _{1}$. So 
all 
pleated surfaces $[f]$ that we use have the not-containment 
property.
So the properties of \ref{8.2} are satisfied for all the sequences 
that we use, and the proofs go through as before.  
 
This proves assumption \ref{8.12.1} of Theorem \ref{8.12}.

\ssubsection{Modifications of  \ref{10.4}, \ref{10.5} in the 
 incompressible case.}\label{10.14}

Let $e$ be any incompressible end. We  define $z_{e,0}$, $z_{e',0}$ 
as in 
 \ref{6.17} and $f_{e,0}$ as in \ref{8.5}. Let $\omega (e,e')$ be as 
in \ref{6.17}, and as in \ref{6.17},
write
 $$\beta (e) =S(e)\setminus \cup _{e'}\omega (e,e').$$ 
Take the sequence of mutlicurves and pleated surfaces $\Gamma _{e,0}$ 
and $f_{e,i}$ as in \ref{8.5}. By \ref{8.12}, for $\Sigma '$ as 
defined 
 there and in \ref{6.14}
 $$\vert f_{e,0}(\Sigma ')\vert \leq L_{2},$$
 and hence, enlarging $L_{2}$ if necessary,
 $$d_{\beta (e)}([f_{e,0}],z_{e,0})=d_{\beta 
(e)}([f_{e,0}],z_{e,0}')\leq L_{2}.$$

 Now we consider the proof of \ref{10.1}. 
 By \ref{8.12}, since the hypothesis \ref{8.12.1} has now been 
proved, we have 
 $$\vert (\partial \beta (e))_{*}\vert \leq 
 L_{2},$$
 and we have a pleated surface $f_{-}$ with pleating 
 locus including $\partial \beta (e)$, with pleating locus of length 
 $\leq L_{2}$ restricted to $\beta (e)$, such that
 $$d_{\beta (e)}([f_{-}],z_{e,0})\leq L_{2}.$$
  By the definition of $z_{e,0}'$ in 
 \ref{6.17}, if $(\alpha ,\ell )$ 
is 
 ltd  for $S(e)\times [z_{e,0}',y_{e,+}]$, $\alpha \cap \beta (e)\neq 
 \emptyset $. Take any chain of ltd's as in \ref{7.13} for 
 $S(e)\times [z_{e,0}',y_{e,+}]$. Then \ref{10.2} is proved as in 
\ref{10.3} 
 and \ref{10.4}, using 
 $f_{-}=f_{0,-}$, $f_{e,+}=f_{0,+}$. As in the compressible case, the 
 fact that we do not have a bound on the pleating locus $\Gamma 
 _{0,-}$ of $f_{0,-}$ off $\beta (e)$ does not matter. We 
 will get bounds for the first $k$ with $f_{k,-}\neq f_{0,-}$. Then 
 using this first inductive step, we can prove \ref{10.2} for all 
 $(\alpha ,\ell)\in \cap _{e'}E(e,e',+)$, exactly as before. 
 
 Now, for the remaining $(\alpha ,\ell )$ in the decomposition for 
 $S(e)\times [z_{e,0},y_{e,+}]$, we also need to consider $(\alpha 
 ,\ell )$ in a decomposition for $\omega (e,e')\times 
 [\pi _{\omega }(z_{e',0}',\pi _{\omega }(z_{e,0}')]$ for $\omega 
 =\omega (e,e')$ for varying $e'\neq e$. The technique of proof of 
 \ref{10.1} is the same as before. We  start by obtaining 
 (\ref{10.1.1}) for a totally ordered set of $(\alpha ,\ell)$ with 
 $\alpha \cap \beta \neq \emptyset $ for all $\alpha $ in the 
 chain --- unless $z_{e,0}'=[f_{e,+}]$. Having 
 got this, we can extend to other $(\alpha ',\ell ')$ which are 
 bounded above and below by elements of the totally ordered set. When 
 we do this, we are likely, at some point, to move into the domain of 
 $\omega (e,e')\times 
 [\pi _{\omega }(z_{e',0}'),\pi _{\omega }(z_{e,0}')]$ for some $e'$. 
 Since we now have bounds above and below here (by $f_{e,+}$ and 
 $f_{e',+}$ if not by anything lower down) 
 we are now in a position to start the induction, 
 applying \ref{10.1} to a totally ordered chain for $\omega 
(e,e')\times 
 [\pi _{\omega }(z_{e',0}'),\pi _{\omega }(z_{e,0}')]$, initially 
 using the multicurve and pleated surface sequences of \ref{8.5}.
  The rest of the proof is exactly as in \ref{10.5}.
 
\ssubsection{Where short loops can occur.}\label{10.7}

Now, by the same argument as in \ref{9.2}, we can obtain information 
about where Margulis tubes can occur, and on the geometry of those 
which do occur. Precisely, we have the following. For 
$\gamma \subset S\subset N$ with $\gamma $ nontrivial in $N$ and 
not represented by a parabolic element, we let $\gamma _{*}$ denote 
the closed 
geodesic freely homotopic to $S$, as in Section \ref{3}.

It may be worth pointing out that if $e$ is a {\em{compressible}} 
end, and (\ref{10.1.2}) holds, and $(\alpha ,\ell)$ is ltd for 
$S(e)\times [z_{e,0},y_{e,+}]$ and $\alpha \neq S(e)$, then 
 $\alpha $ is incompressible, by taking $\mu '=\partial \alpha /\vert 
\partial 
 \alpha \vert $. But then there cannot be 
another end, $e'$, compressible or not, with $\alpha \subset S(e')$ 
up 
to homotopy. For if there is such an $e'$, any loop $\gamma '$ on 
$S(e)$ 
which is nontrivial on $S$ but trivial in $N$ must be either 
decomposable into two such loops, one of which is disjoint from 
$\alpha $, or is already disjoint from $\alpha $. Either way, this 
contradicts (\ref{10.1.2}). 

\begin{ulemma} We continue with the notation and hypotheses of 
\ref{10.1}. Fix suitable ltd parameter functions and vertically 
efficient ltd decompositions of $S(e)\times [z_{e,0},y_{e,+}]$. 
  Then there exists 
$\varepsilon _{1}>0$ which depends only on  the topological type of 
$(N_{d,W},\partial N_{d})$, and, if some end $e$ is compressible, on 
the constant 
$c_{0}$ of \ref{10.1}.  Fix a Margulis constant $\varepsilon _{0}$. 
Let $\gamma $ be a nontrivial nonperipheral loop in $N$, and suppose
that $\vert \gamma _{*}\vert <\varepsilon _{1}$ and $T(\gamma
_{*},\varepsilon _{0})$ is in region of $N$ which is disjoint from the
sets $f_{e,+}(S(e))$ but bounded by all of them, and contains a tamely
embedded relative Scott core.  Then $\gamma \subset \partial \alpha
_{j}$ for some $(\alpha _{j},\ell _{j})$ in the vertically efficient
decomposition of $S(e)\times [z_{e,0},y_{e,+}]$ for some end $e$.

\end{ulemma}

Note that any Margulis tube $T(\gamma _{*},\varepsilon _{0})$ is in
such a region of $N$ for some choice of end pleated surfaces
$f_{e,+}$, by choosing them with images $f_{e,+}(S(e))$ in
sufficiently small neighbourhoods of the ends.

\noindent {\em{Proof.}} We use the sequence of pleated surfaces for
$[z_{e,0},y_{e,+}]$ for each end $e$ of \ref{8.4}, and, extending
this, the family of pleated surfaces of \ref{8.12}.  We assume without
loss of generality that the pleated surfaces for different ends have
disjoint images, by taking pleated surfaces in a sufficiently small
neighbourhoods of the geometrically infinite ends.  The definitions
made in \ref{8.3} mean we are using a generalised pleated surface in
$\partial N$ for each geometrically finite end.  We use the same
argument as in \ref{9.2}.  We claim that every point in $N$ which is
in a component of the complement of the images of the boundary pleated
surfaces as described in the statement of the lemma, is in the image
of the homotopy between two successive pleated surfaces in an end
sequence, or in the noninterval part of the core, in the image of the
homotopy equivalence $f:W'\to N$ of Theorem \ref{8.12}.  The
homotopies all have bounded tracks, except for the homotopy between
the last two generalised pleated surfaces in a geometrically finite
end --- but in that case, by \ref{4.5}, the Teichm\" uller distance is
bounded.  If the claim is not true, let $W$ be the relative Scott core
in this component of the complement of the images of the pleated
surfaces.  We can assume an omitted point $w_{0}$ of the homotopy lies in the
interior of $W$.  By \cite{Wal}, we can assume, after composing on the
left with a homotopy which is the identity on $W\subset N$ that the
homotopy equivalence from $W$ to $N$ maps $(W,\partial W)$ to
$(W,\partial W)$.  So we have a homotopy equivalence which maps
$(W,\partial W)$ into $(W\setminus \{ w_{0}\} ,\partial W)$.  But then
$H_{3}(W,\partial W)\cong \mathbb Z$ is isomorphic to
$H_{3}(W\setminus \{w_{0}\} ,\partial W)=0$, which is a contradiction. 
(This can be proved by considering the short exact sequence of chain
complexes $0\to C_{*}(W\setminus \{ w_{0}\} ,\partial W)\to
C_{*}(W,\partial W)\to C_{*}(W,W\setminus \{ w_{0}\} )\to 0$.)  Now for
$\varepsilon _{1}$ sufficiently small, the image of the homotopy
between $f_{e,i}$ and $f_{e,i+1}$ can only intersect a Margulis tube
$T(\gamma_{*},\varepsilon _{1})$ if there is no loop $\zeta $
intersecting $\gamma $ transversely for which $\vert f_{e,i}(\zeta
)\vert <L_{1}$.  Now by \ref{10.1}, in particular by (\ref{10.1.1}),
if $T(\gamma _{*},\varepsilon _{1})$ intersects the homotopy between
$f_{e,i}$ and $f_{e,i}$, then $\gamma $ has no transverse
intersections with the pleating locus of $f_{e,i}$, nor with the
pleating locus of $f_{e,j}$ for $\vert i-j\vert \leq t$, if
$\varepsilon _{1}$ is sufficiently small given $t$.  Now if $t$ is
sufficiently large given the ltd parameter functions, the only loops
which have no transverse intersections with the pleating loci of $2t$
successive $g_{e,j}$ are loops which are $\partial \alpha _{k}$ for
some $\alpha _{k}\times \ell _{k}$ in the decomposition.

Let $f_{W'}:W'\to N$ be the map of \ref{8.12}. By Theorem 
\ref{8.12}, $f_{W'}(W')$ has diameter $\leq L_{3}$, for some constant 
$L_{3}$ which, ultimately, dependds only on the topological type of 
$W'$ and 
the constant $c_{0}$ of \ref{10.1.3}. 
So any Margulis tube $T(\gamma 
_{*},\varepsilon _{1})$ can only intersect $f_{W'}(W')$, in a set of 
diameter $\leq L_{3}$. So then, by the 
first part of the proof, $\gamma \subset \partial \alpha _{j}$ for 
some $\alpha _{j}\times \ell _{j}$ in the decomposition for 
$S(e)\times [z_{e,0},y_{e,+}]$ for some end $e$. In fact, this can be 
true for at most two ends, and for two ends $e$, $e'$ only if $e$ and 
$e'$ are both incompressible and $\omega (e,e')\neq \emptyset $.\Box 

\ssubsection{Bounded distance between corresponding Margulis tubes.} 
\label{10.8}

 Let 
$M$ be a model manifold of topological type 
$S\times [0,1]$.
As in \ref{6.5}, we denote by $T(\gamma _{**})$ the Margulis tube (if 
any) with core loop $\gamma _{**}$ homotopic in $S\times [0,u]$ to 
$\gamma \times \{ t\} $.  In \ref{6.4}, we recalled how to 
parametrise Margulis 
tubes $T(\gamma _{*},\varepsilon _{0})$ by points $w_{*}(\gamma )$ in 
the upper 
half-plane $H^{2}$. In \ref{6.5}, we showed how the boundary of 
a model Margulis tube $\partial T(\gamma _{**})$ also determines a 
point $w_{**}(\gamma )$ of $H^{2}$ up to bounded distance. The 
following lemma gives a 
bound  on the distance between the corresponding points $w_{*}(\gamma 
)$ 
and $w_{**}(\gamma )$. 

We have the following, building on \ref{10.7}.

\begin{ulemma} We continue with the notation and hypotheses of 
\ref{10.1}. Fix a Margulis constant $\varepsilon _{0}$.  
    The following holds for suitable $\varepsilon _{2}>0$ and $L'$ 
  depending on 
$\varepsilon _{0}$, and the topological type of 
$(N_{d,W},\partial N_{d})$, and, if there are compressible ends, the 
constant $c_{0}$ of \ref{10.1}. If one of the following holds: 
\begin{description}
    \item[.]
$\vert \gamma _{*}\vert <\varepsilon 
_{2}$ and  $T(\gamma
_{*},\varepsilon _{0})$ is in region of $N$ which is disjoint from the
sets $f_{e,+}(S(e))$ but bounded by all of them, and contains a tamely
embedded relative Scott core;
\item[.]  $\vert \gamma _{**}\vert
<\varepsilon _{2}$ and $T(\gamma _{**})$ does not intersect $\partial
M$;
\end{description}
then
$$d_{H}(w_{*}(\gamma ),w_{**}(\gamma ))<L',$$
where $d_{H}$ denotes hyperbolic distance in $H^{2}$. \end{ulemma} 

\noindent {\em{Proof.} } By \ref{10.7}, we can assume that 
$\gamma  =\partial \alpha _{e,j}$ for some 
$\alpha _{e,j}\times \ell _{e,j}$ in a vertically efficient 
decomposition 
of $S(e)\times [z_{e,0},y_{e,+}]$  for some end $e$ of $N_{d,W}$

First, we assume that this is true for exactly one end $e$. 
We drop the suffix $e$.
From the construction  of the model, $\vert \gamma 
_{**}\vert $ can also only be small if $\gamma =\partial \alpha _{j}$ 
for some $j$. First, suppose that $\vert \gamma _{*}\vert 
<\varepsilon _{2}$, for a bound on $\varepsilon _{2}$ yet to be 
determined, but with $\varepsilon _{2}<\varepsilon _{1}/2$ for 
$\varepsilon _{1}$ as in \ref{10.7}. Let $A$ be the set of 
$k$ such that $\alpha _{k}$ is a gap with $\partial \alpha 
_{k}=\gamma $. Let $A_{1}$ be the set of $k\in A$ such that $\ell 
_{k}\times \alpha _{k}$ is ltd and $A_{2}=A\setminus A$. Then by the 
definition of the metric in $M$ in \ref{6.5}, the area of $\partial 
T(\gamma _{**})$ is 
boundedly proportional to
\begin{equation}\label{10.8.1} \sum \{ \vert \ell _{k}\vert :k\in 
A_{1}\} +\# (A_{2}).\end{equation}
Note that we can ignore any contribution from $W'$, because we have a 
bound on the diameter of the image of the homotopy equivalence 
$f:W'\to 
N$ of \ref{8.12}.

Now if $k$, $j\in A_{2}$, $\alpha _{k}\cap \alpha _{j}$ contains no 
nontrivial subsurface in the interior, apart from $\gamma $, by the 
vertically 
efficient condition \ref{5.6}. Now  by the 
last part of the statement of \ref{10.7}, the area of $\partial 
T(\gamma 
_{*},\varepsilon _{1})$ is bounded above by the number of different 
 pleated surfaces $f_{n}\vert \alpha _{k}$ for which 
$\partial \alpha 
_{k}$ is in the pleating locus of $f_{n}$, for $k\in A$, which is 
bounded above by a multiple of the sum in \ref{10.8.1}. But the area 
of $\partial T(\gamma _{*},\varepsilon _{1})$ is also bounded below 
by a multiple of the sum in (\ref{10.8.1}), with multiple depending 
on 
$\varepsilon _{1}$. This is because, by the properties of the 
Margulis 
constant, for any constant $L'$ there is 
$k(L')$ such that the following holds. For any 
$x\in \partial T(\gamma _{*},\varepsilon 
_{1})$, there are at most $k(L')$ loops $\zeta $ such that for some 
$n$, $\vert f_{n}(\zeta )\vert \leq L'$ and $f_{n}(\zeta )$ is within 
distance 
$L'$ of $x$. This works even in the case when $S(e)$ is compressible, 
by Lemma \ref{10.12} below.

   So the areas of $\partial T(\gamma _{**})$ and $\partial 
T(\gamma _{*},\varepsilon _{1})$ are boundedly proportional as 
claimed, and the same is true for the areas of 
$\partial T(\gamma _{**})$ and $\partial 
T(\gamma _{*},\varepsilon _{0})$ for a constant depending on 
$\varepsilon _{1}$. The areas determine the imaginary coordinates of 
the corresponding points of $H^{2}$. So we have a bound on the 
ratio of ${\rm{Im}}(w_{*}(\gamma ))$ and ${\rm{Im}}(w_{**}(\gamma ))$
if one of $\vert \gamma 
_{*}\vert $, $\vert \gamma _{**}\vert $ is $<\varepsilon _{2}$.

It remains to bound the distance between ${\rm{Re}}(w_{*}(\gamma ))$ 
and ${\rm{Re}}(w_{**}(\gamma ))$
 by a constant which is  $O({\rm{Im}}(w_{*}(\gamma 
)))+O(1/\varepsilon _{1})$, for $\varepsilon _{1}$.
Now
$${\rm{Re}}(w_{**}(\gamma ))=-n_{\gamma }(z_{2})+n_{\gamma 
}(z_{1})+O(1),$$  
if $[z_{1},z_{2}]$ is the union of all 
the $\ell _{k}$ with $\gamma \subset \partial \alpha _{k}$. (See 
(\ref{6.5.1}).) We recall that $n_{\gamma }$ is defined relative to a 
fixed loop $\zeta \subset S$ which has at least one, and at most two, 
essential intersections with $\gamma $. A different choice of $\zeta 
$ does not change the quantity  $-n_{\gamma 
}(z_{2})+n_{\gamma }(z_{1})$ by more than a bounded additive 
constant, 
independent of $\zeta $. If  the $m$ for which $\gamma $ is in the 
pleating locus of $f_{m}$ include an $m$ such that $f_{m+1}$ and 
$f_{m+2}$ are related by a Dehn twist, rather than an elementary 
move, 
then 
there is at most one 
such $n$, by the properties listed in \ref{8.2}, and there is a 
unique 
loop $\zeta $ intersecting $\gamma $ at most twice, in the pleating 
locus of $f_{m+1}$, and $\tau _{\gamma }^{n}(\zeta )$ is in the 
pleating locus of $f_{m+2}$, where
$$n-{\rm{Re}}(w_{**}(\gamma ))=O(1)+O({\rm{Im}}(w_{**}(\gamma )).$$
If there is no such $m$, we can still find $m$ such that $\gamma $ is 
in the pleating locus of $f_{m}$ and some loop $\zeta \subset S$ with 
at least one, and at most two, essential intersections with $\gamma 
$, 
is in the pleating locus of $f_{m+1}$.  We choose any such $\zeta $, 
and 
can indeed choose such an $m$ so that $f_{m}(S)$ intersects 
$T(\gamma _{*},\varepsilon _{1})$, and then 
$f_{m+1}(\zeta 
)=\zeta _{*}$ comes within a bounded distance of 
$T(\gamma _{*},\varepsilon _{1})$. By \ref{10.7} and \ref{10.1}, 
if we take the 
minimal $m_{1}$ and maximal $m_{2}$ such that $f_{m}(S)$ 
intersects 
$T(\gamma _{*},\varepsilon _{1})$ for $m=m_{1}$, $m_{2}$, then 
$$d_{\beta _{j}}(z_{j},[f_{m_{j}}])\leq L_{1}'$$
for a constant $L_{1}'$ depending only on $(\Delta 
_{1},r_{1},s_{1},K_{1})$, where $\beta _{j}=\alpha _{i(j)}$ for some 
$\alpha 
_{i(j)}$ 
intersecting $\gamma $ transversely with $z_{j}\in \ell _{i(j)}$, 
$j=1$, $2$. Write $n_{j}=n_{\gamma }(z_{j})$, where this 
is defined using $\zeta $. Consider the geodesics  in 
$S(f_{m_{j}})$ homotopic to
$$f_{m_{j}}(\tau _{\gamma }^{n_{j}}(\zeta )).$$
Take the intersection with $T(\gamma _{*},\varepsilon 
_{0})$ of the images under $f_{m_{j}}$ of these geodesics in $N$, and 
homotope, via endpoint preserving isotopy, to two segments $\gamma 
_{3}$, $\gamma _{4}$ in 
$\partial T(\gamma _{*},\varepsilon _{0})$ of length $O(1/\varepsilon 
_{1})$. Here, we are using numbering which parallels the numbering 
$\partial _{k}T(\alpha )$ of \ref{6.5}.
Use the homotopy through the loops $f_{m}(\zeta )$ to give two 
arcs in $\partial T(\gamma _{*},\varepsilon _{0})$ of length 
$O({\rm{Im}}(w_{**}(\gamma ))=O({\rm{Im}}(w_{*}(\gamma ))$ joining 
$\gamma _{1}$ and $\gamma _{2}$. Then $\gamma _{1}\cup \gamma 
_{3}\cup \gamma _{2}\cup \gamma _{4}$ is a loop which is freely 
homotopic 
to $\zeta _{2}*\overline{\zeta _{1}}$, where 
$\zeta _{j}=\tau _{\gamma} ^{n_{j}}(\zeta )$.
This means that
$${\rm{Re}}(w_{*}(\gamma ))={\rm{Re}}(w_{**}(\gamma 
))+O({\rm{Im}}(w_{*}(\gamma )))+O(1/\varepsilon _{1}).$$
This gives the required bound on $d_{H}(w_{*}(\gamma ),w_{**}(\gamma 
))$.

Now we consider the possibility that $\gamma \subset \partial 
\alpha $, $\gamma \subset \partial \alpha '$ for $(\alpha ,\ell )$ 
$(\alpha ',\ell ')$ in the decompositions for two different ends $e$, 
$e'$. This can only happen if $\omega (e,e')\neq \emptyset $, and 
$\alpha $, $\alpha '\subset \omega (e,e')$, and up to bounded 
distance, $(\alpha ,\ell )$ and $(\alpha ',\ell ')$ are sets in the 
decomposition for $\omega \times [\pi _{\omega }(z_{e',0}'),\pi 
_{\omega 
}(z_{e,0}')]$, $\omega =\omega (e,e')$. Then we can argue exactly as 
above, with $\omega \times [\pi _{\omega }(z_{e',0}'),\pi _{\omega 
}(z_{e,0}')]$ replacing $S(e)\times [z_{e,0},y_{e,+}]$. 
\Box

\begin{lemma}\label{10.12}
Let $N_{d,W}$ be as throughout this section. We continue with the 
notation and hypotheses of \ref{10.1}.
Fix a compressible end $e$ of $N_{d,W}$ bounded by $S_{d}(e)$. Let 
$S=S(e)$, 
$z_{0}=z_{e,0}$, 
$y_{+}=y_{e,+}$, $\Gamma _{+}=\Gamma _{e,+}$ be as in \ref{10.1}, in 
particular satisfying  
(\ref{10.1.2}) with respect 
to constants $c_{0}$ and $L_{0}$.
Let $z_{1}=[\varphi _{1}]$, $z_{2}=[\varphi _{2}]\in 
[z_{0},y_{+}]\subset 
{\cal{T}}(S(e))$ with $z_{1}\in 
[z_{0},z_{2}]$. Let $L$ be given.
The following holds for a constant $L'$ depending 
only on the the topological type of $(N_{d,W},\partial N_{d,W})$, on 
the constants 
$c_{0}$, $L_{0}$, on suitable ltd 
parameter functions. Suppose that we have 
fixed a vertically efficient ltd decomposition of $S\times 
[z_{e,0},y_{e,+}]$ 
with respect to the  parameter functions. Let
$\zeta _{1}$, $\zeta _{2}$ be simple loops such that $\vert \varphi 
_{j}(\zeta _{j})\vert \leq L$,  and 
with $\zeta _{j}\subset \alpha _{j}$, $z_{j}\in \ell _{j}$ 
and $z'\in \ell '$ for sets $\alpha _{j}\times \ell _{j}$  
in the ltd decomposition for $j=1$, $2$. Let
$d_{\alpha _{1}}'(z_{1},z_{0})\geq L'$
and let $(\zeta _{1})_{*}=(\zeta _{2})_{*}$. Then $\zeta _{1}=\zeta 
_{2}$.\end{lemma}

\noindent{\em{Proof.}} Suppose that  $(\zeta _{1})_{*}=(\zeta 
_{2})_{*}$ and 
$\zeta _{1}\neq \zeta _{2}$. We 
can assume without loss of generality that $\zeta _{1}$ and $\zeta 
_{{2}}$ 
have 
a common basepoint. So $\zeta _{1}*\overline{\zeta _{2}}$ is a 
nontrivial 
closed loop.  Then  by the Loop Theorem 4.10 of 
\cite{Hem}, we can find a loop 
$\gamma \subset \zeta _{1}\cup \zeta _{2}$ such that $\gamma \subset 
S$ is a simple loop and, regarding $S$ as a submanifold of $N$, 
$\gamma $ bounds an embedded disc in $N$. Now for any loop $\zeta 
_{3}$, 
$$i(\gamma ,\zeta _{3})\leq i(\zeta _{1},\zeta _{3})+i(\zeta 
_{2},\zeta 
_{3}).$$
If $L'$ is large enough then there is at least one ltd $(\alpha 
_{3},\ell _{3})\leq (\alpha _{1},\ell _{1})$ and $z_{3}=[\varphi 
_{3}]\in \ell _{3}$, $\zeta _{3}\subset \alpha _{3}$ with
$$\vert \varphi _{3}(\zeta _{3})\vert \leq L_{0},$$
$$i(\zeta _{1},\zeta _{3})\leq C,$$
$$d_{\alpha _{3}}'(z_{3},z_{0})\geq L''.$$
Here, $L_{0}$ is bounded in terms of topological type, 
$C$ is bounded in terms of the ltd parameter functions, 
and $L''$ can be taken arbitrarily large by choice of $L'$, by 
\ref{7.5}.

Now we claim that (\ref{10.1.2}) implies that for a constant $c_{1}>0$ depending only on $c_{0}$ of (\ref{10.1.2}), 

\begin{equation}\label{10.12.1}
    \begin{array}{l}
    i(\zeta ,\mu )\geq c_{1}\vert \zeta \vert {\rm{\ for\ some\ 
}}\zeta 
    \in \Gamma _{+}\cr
    {\rm{\ if\ }}i(\mu ,\mu ')\leq c_{1}{\rm{\ and\ }}i(\mu 
    ',\gamma )\leq c_{1}.\vert \gamma \vert \cr
    \end{array}\end{equation}

For if (\ref{10.12.1}) fails for a $c_{1}$ depending only on $c_{0}$, we can assume that there are sequences $\zeta _{n}$, $\mu _{n}$, $\mu _{n}'$ and $\gamma _{n}$  such that $\Gamma _{+,n}$ satisfies  the condition (\ref{10.1.2}), but  such that $\mu _{n}\to \mu _{\infty }$ and $\mu _{n}'\to \mu _{\infty }$, and 
$$\lim _{n\to \infty }{\rm{Max}}\{ i(\zeta   _{n}/\vert \zeta _{n}\vert ,\mu _{n}):\zeta _{n}\in \Gamma _{n,+}\} =0$$
and
$$\lim _{n\to \infty }i(\mu _{n},\mu _{n}')=\lim _{n\to \infty }i(\mu _{n}',\gamma _{n}/\vert \gamma _{n}\vert )=0.$$
Then by the Lipschitz properties of $i(,.,)$ (\cite{R3}), we see that (\ref{10.1.2}) fails for $\Gamma _{+,n}$ and $\gamma _{n}$ for sufficiently large $n$, with $\mu =\mu _{\infty }$ and $\mu '=\mu '_{\infty }$, giving
the required contradiction.

Then by \ref{7.5}, again for $C$ depending on the ltd parameter 
functions, and on 
 $L_{0}$, for any $\zeta \in \Gamma _{+}$,
 $$\vert \zeta \vert \geq C^{-1}i(\zeta ,\zeta _{3})\exp d_{\alpha 
 _{3}}'(z_{3},z_{0}),$$
 $$\vert \zeta _{2}\vert \geq C^{-1}i(\zeta _{2},\zeta _{3})\exp 
d_{\alpha 
 _{3}}'(z_{3},z_{0}),$$
 $$\vert \zeta _{3}\vert \geq C^{-1}\exp d_{\alpha 
 _{3}}'(z_{3},z_{0}),$$
We also have, by \ref{2.8},
 $$i(\zeta _{3},\zeta _{2})\leq Cd_{\alpha _{3},\alpha 
 _{2}}'(z_{3},z_{2}),$$
again for $C$ depending on the ltd parameter functions, because  we 
only have a bound on $\vert \varphi _{2}(\partial \alpha _{2})\vert $ 
in terms of these.
Again by \ref{7.5}, if $z_{0}=[\varphi _{0}]$,
 $$\vert \varphi _{0}(\zeta _{2})\vert \geq 
 C^{-1}\vert i(\zeta _{3},\zeta _{2})\exp d_{\alpha 
_{3}}'(z_{3},z_{0}).$$ 
 So we have
$$i(\zeta ,\zeta _{3})\leq 
C\vert \zeta \vert .\exp -d_{\alpha _{3}}'(z_{3},z_{0}),$$
$$i(\gamma ,\zeta _{2})\leq C+
i(\zeta _{3},\zeta _{2})\leq C(\vert \zeta _{3}\vert +\vert \zeta 
_{2}\vert ).
\exp -d_{\alpha _{3}}'(z_{3},z_{0}).$$
This contradicts (\ref{10.12.1}) with $\mu =\zeta _{3}/\vert \zeta 
_{3}\vert $ and $\mu '=\zeta _{2}/\vert \zeta _{2}\vert $,
if $L''$ is 
sufficiently large given $c_{1}$, that is, if $L'$ is sufficiently 
large.
\Box

\ssubsection{The order on Margulis tubes is the same.}\label{10.9}

\ref{10.7} also implies the following, using the partial order of 
\ref{7.3}. 

\begin{ucorollary} We continue with the notation and hypotheses of 
\ref{10.1}. Fix an end $e$ of $N_{d,W}$ and write $S(e)=S$, 
$z_{e,0}=z_{0}$, $y_{e,+}=y_{+}$.
Let $(\gamma ,\ell )<(\gamma '
,\ell ')$ for loops $\gamma $, $\gamma '$ with $\gamma \times \ell $, 
$\gamma '\times \ell '$ in the 
vertically efficient decomposition of $S\times [z_{0},y_{+}]$,  where 
 $\vert 
\gamma _{*}\vert <\varepsilon _{1}$, $\vert (\gamma ')_{*}\vert 
<\varepsilon _{0}$.
Then the Margulis tube 
$T((\gamma ')_{*},\varepsilon _{0})$ can be homotoped to an 
arbitrarily small neighbourhood of $e$, in  
$N_{d,W}\setminus T(\gamma _{*},\varepsilon _{0})$.  Similarly, if
$\gamma \subset {\rm{int}}(\alpha )$, and $(\gamma ,\ell )<(\alpha
,\ell _{j}')$ for $j=1$, $2$ for $\alpha \times \ell _{j}$ in the
decomposition with $\alpha $ a gap and $f_{k}$, $f_{m}$ are pleated
surfaces from the sequences of \ref{8.5} associated to $(\alpha ,\ell
_{1}')$, $(\alpha ,\ell _{2}')$ respectively, with pleating loci
including $\partial \alpha $, then $f_{k}(\alpha )$ can be homotoped
to $f_{m}(\alpha )$ via a homotopy in $N_{d,W}\setminus T(\gamma
_{*},\varepsilon _{0})$ and keeping $\partial \alpha $ in $T((\partial
\alpha )_{*},\varepsilon _{1})$.  A similar statement holds if
$(\alpha ,\ell _{j}')<(\gamma ,\ell )$ for $j=1$, $2$.
    \end{ucorollary} 

\noindent {\em{Proof.}} 
For the first statement, since Margulis 
tubes 
are disjoint, it suffices to find a homotopy avoiding $T(\gamma 
_{*},\varepsilon _{1})$, because we can then compose with a 
homeomorphism which expands out $T(\gamma _{*},\varepsilon _{1})$ to 
$T(\gamma _{*},\varepsilon _{0})$. Then we simply use the homotopy 
defined 
by the sequence $f_{k}$, for $k\geq p$ for some $p$ such that 
$\gamma '$ is in the pleating locus of $f_{p}$. By 
\ref{10.7}, $f_{k}(S)$ does not intersect $T(\gamma 
_{*},\varepsilon _{1})$ for $k>p$.  The statements involving
homotopies restricted to $\alpha $ are proved similarly. \Box

Let $U$ be a component of $N\setminus f(S)$ which is a neighbourhood 
of the end of $N_{d,W}$ which is the unique end of $N_{d,W}$ in 
component of $N_{d}\setminus W$ bounded by $S_{d}$. It is 
natural 
to say that $T(\gamma _{*},\varepsilon _{0})<
T(\gamma '_{*},\varepsilon _{0})$ in $U$ if, as in the lemma 
above, 
$T(\gamma '_{*},\varepsilon _{0})$ can be homotoped to an arbitrarily 
small neighbourhood of some end $e$ in $U$ in the complement of 
$T(\gamma _{*},\varepsilon _{0})$, because this order is clearly 
transitive. 
It is not yet clear that it is antireflexive. 

\ssubsection{First stage in the production of a geometric relative 
Scott core.}\label{10.13}

We already know, from \ref{8.12} that $N$ contains a continuous 
Lipschitz (but not necessarily, so far as we yet know, homeomorphic) 
image of the model relative Scott core in $M$. The following lemma 
gives us more information, about a relative Scott core up to 
homeomorphism, 
such that the complementary neighbourhoods of ends are products. This 
lemma does not say that the relative Scott core produced is 
biLipschitz to the model one, although we shall see later that it is. 
Note that in the case of a geometrically finite end, it is possible 
that every pleated surface in the sequence of the end intersects the 
$D$-neighbourhood of $f_{e,0}(S(e))$

\begin{ulemma} We continue with the notation and hypotheses of 
\ref{10.1}. 
    The following holds for a constant $D$ depending only 
on the topological type of $(N_{d,W},\partial N_{d,W})$, on 
$c_{0}$, $L_{0}$, and the loop sets $\Gamma _{0}(e)$. Fix an end $e$ 
of 
$N_{d,W}$ and $S_{d}=S_{d}(e)$. Then 
there is $S'\subset N_{d,W}$ such that 
$(N_{d,W},S')$ is homeomorphic to $(N_{d,W},S_{d})$, and  the 
component of $N_{d,W}\setminus S'$ which is a neighbourhood of $e$ is 
homeomorphic to $S_{d}(e)\times (0,\infty)$, each component of 
$S'\cap N_{\geq 
\varepsilon _{1}}$ has diameter $\leq D$ and for all $f_{e,j}$ in the 
sequence for $e$, either $f_{e,j}(S(e))$ is contained in the 
neighbourhood of $e$ bounded by $S'$, or $f_{e,j}(S(e))$ intersects 
the $D$-neighbourhood of $S'\cap N_{\geq 
\varepsilon _{1}}$.\end{ulemma}

\noindent{\em{Proof.}} 
Write $S=S(e)$.  Take the 
corresponding sequence $f_{e,j}=f_{j}:S\to N$, $0\leq j\leq 
n=n_{+}(e)$, 
of pleated 
surfaces, given a choice of pleated 
surface $f_{+}=f_{n}$ with image in a small neighbourhood of $e$, 
disjoint from $W$.   For sufficiently large $j$, 
there is a bounded track homotopy between $f_{j}$ and 
$f_{j+1}$ which  is disjoint from both $S_{d}\subset 
\partial _{d}W$ 
and from $f_{0}(S)$.  Now suppose that $j$ is such that the homotopy 
between $f_{k}$ and $f_{k+1}$ is disjoint 
from $f_{0}(S)$ for all $k\geq j$, but $f_{j}(S)$ 
intersects the $D$-neighbourhood of $f_{0}(S)$. This is possible for 
$D$ depending only on the topological type of $S$, by \ref{8.2.6}. 
Then we claim 
that $f_{k}:S\to N\setminus f_{0}(S)$ is injective on 
$\pi 
_{1}$ for $k\geq j$. It suffices to prove this for sufficiently large 
$k$, since all the $f_{k}$, $k\geq j$, are homotopic in 
$N\setminus f_{0}(S)$. Now for all sufficiently large $k$, 
$f_{k}$ and $f_{k+1}$ are homotopic in 
$N\setminus W$ and $f_{k}:S\to N\setminus W$ is injective 
on $\pi _{1}$.  So for sufficiently large $k$, $f_{k}$ is homotopic in
$N\setminus W$ to an embedding $f'$.  It suffices to show that $f'$
is injective on $\pi _{1}$ in $N\setminus f_{0}(S)$.  Suppose not. 
Then we can find an embedded disc $D$ in $N\setminus f_{0}(S)$, with
boundary in $f'(S)$ but otherwise not intersecting $f'(S)$, which is nontrivial in $f'(S)$, and such that $D$ has
homotopically nontrivial intersection  with $f'(\gamma )$, for some
$\gamma \subset S$ which is nontrivial in $N$.  But then, since $f'(\gamma )$ and $f_{0}(\gamma )$ are homotopic in the subset of $N$ bounded by $f'(S)$, $D$ must
intersect $f_{0}(\gamma )$, which is a contradiction.  Then the
hypothesis of the Freedman-Hass-Scott result \cite{F-H-S} is
satisfied, and $f_{k}:S\to N\setminus f_{0}(S)$ is homotopic to an
embedding with image in an arbitrarily small neighbourhood of
$f_{k}(S)$.  Then we can take the image of the embedding to be $S'$,
which is in a small neighbourhood of $f_{j}(S)$ and hence intersects a
$D/2$-neighbourhood of $f_{0}(S)$ for suitable $D$.

Now if some $f_{k}(S(e))$ is separated from $e$ by $S'$ and does not 
intersect the $D$-neighbourhood of $S'\cap N_{\geq \varepsilon 
_{1}}$, 
then we can repeat the construction with $f_{k}$ replacing $f_{0}$, 
and can obtain a larger neighpourhood of $e$, enlarged by a set of 
diameter $D/2$. This construction can only be repeated 
finitely many times, because the sequence $f_{k}$ is finite. \Box  

\ssubsection{Coarse biLipschitz in a ltd piece}\label{10.10}
At this point, we start to use the freedom afforded by different 
choices of ltd parameter functions. Let $(\Delta 
_{1},r_{1},s_{1},K_{1})$ be the parameter functions and flat constant 
used until now, and we fix a vertically efficient partition using 
these. 
We note that, by \ref{6.6}, geometric models defined 
with different choice of ltd parameter functions are coarse Lipschitz 
equivalent. In the following lemma, $\Delta (\nu ,C)$ is likely to be 
much larger than $\Delta _{1}(\nu )$. There is a use of geometric 
limits in this lemma, but only in the context of bounded geometry.

\begin{ulemma} We continue with the notation and hypotheses of 
\ref{10.1}. Let $\varepsilon _{1}$ be as in 
\ref{10.7}. There is a 
function $\Delta (\nu ,C)$, depending on $(\Delta 
_{1},r_{1},s_{1},K_{1})$,
the topological type of $N_{d,W}$, $L_{0}$, $c_{0}$,
such that the following hold. Fix an end $e$ of $N_{d,W}$ and 
$S=S(e)$. Let 
$f_{e,j}=f_{j}:S\to N$ be the corresponding sequence  of pleated 
surfaces 
for an end, as in \ref{8.4},
between pleated surfaces $f_{-}=f_{0}$ and $f_{+}=f_{n}$, and 
let $z_{k}=[\varphi _{k}]\in \cal{T}(S)$ be the point from which the 
pleating locus of $f_{k}$ is constructed, $0\leq k\leq n$. Let 
$\alpha $ be  a gap 
which is long $\nu $-thick and dominant along a segment of length 
$\geq 2\Delta (\nu ,C)$ centred on $z_{i}$, and let 
$d(z_{j},z_{i})\geq \Delta (\nu ,C)$.  

Then $f_{j}(\alpha )\setminus T((\partial \alpha)_{*},\varepsilon 
_{1})$ 
and $f_{i}(\alpha ) \setminus T((\partial \alpha)_{*},\varepsilon 
_{1})$ 
cannot be joined by a path of length $\leq C$ in 
$N\setminus T((\partial \alpha)_{*},\varepsilon _{1})$ 
 If, in addition, $i<j<k$, $z_{i}$, $z_{j}$, $z_{k}$ are all in a 
segment of
 along which $\alpha $ is long, thick and dominant, all distance $\geq
 C$ from the endpoints of this segment, $d(z_{j},z_{i})\geq C$, and 
$C$ is
 sufficiently large, then $f_{j}\vert \alpha $ can be homotoped to
 $f_{k}\vert \alpha $ in the complement of $N\setminus
 (f_{i}(\alpha)\setminus T((\partial \alpha _{j})_{*},\varepsilon
 _{1})$, via a homotopy which maps $\partial \alpha $ into 
$T((\partial
 \alpha _{j})_{*},\varepsilon _{1})$ and with image in $U$.  A 
similar statement
 holds if $k<j<i$.\end{ulemma}

\noindent {\em{Proof.}}  Recall that in \ref{8.4}, the pleating 
locus of $f_{k}$ was chosen to have no short loops, and thus, with 
$\partial \alpha $ not in the pleating locus, which included the 
maximal multicurve $\Gamma _{k}$. However, $\Gamma _{k}$ was 
obtained by modifying a sequence  of maximal multicurves $\Gamma 
_{k}'$ which did 
contain
$\partial \alpha $ and such that, for each $n$, $\#(\Gamma _{n}\cap 
\Gamma _{k}')$ was bounded for some $k$. By \ref{4.3} and \ref{4.4}, 
we can return to 
this sequence for $k$ such that $z_{n}\in \ell $, which we continue 
to call 
$\Gamma _{n}$, and the 
corresponding pleated surface $f_{n}$. But we also keep the same 
first 
pleated surface $f_{0}$ in the sequence.  
The distance between 
$f_{n}(\alpha )$ and $f_{n+1}(\alpha )$ is  bounded in terms of 
$\nu $, using \ref{4.3} and the bound on $\# (\alpha \cap \Gamma 
_{n}\cap \Gamma _{n+1})$ in terms of $\nu $ in \ref{8.4}.  If $e$ is
compressible and $\alpha =S$, and $f_{n}(S)$ is not contained in the set $U$ of
\ref{10.13} then by \ref{10.13}, $f_{n}(S)\cap S'\neq \emptyset $.  So
in this case, after discarding finitely many $n$, $f_{n}(S)\subset U$
for $U$ as in \ref{10.13}.  If $\alpha \neq S$ then, as pointed out in \ref{10.7}, 
$f_{n}\vert \alpha $ is injective on $\pi _{1}$.  Also, $f_{n}\vert
\alpha $ is homotopic to an embedding via a homotopy preserving
$f_{n}\vert \partial \alpha $.  We see this as follows.  For
sufficiently large $p$, $f_{p}\vert S$ is homotopic to an embedding,
by \cite{F-H-S} applied to $N\setminus W$.  Then using the homotopy
through the $f_{m}$ from $F_{n}$ to $f_{p}$, we have a continuous map
$F:\partial \alpha \times [0,1]\to N$ where $F(x,0)=f_{n}(x)$ and
$F(x,1)=f_{p}(x)$ for all $x\in \partial \alpha $.  Replacing $f_{p}$
by $f_{q}$ for an even larger $q$ if necessary, and applying
\cite{Wal}, we can also ensure that $F^{-1}(F(x,j))=F(x,j)$ for all
$x\in \partial \alpha $ and $j=0$, $1$ and that the image of $F$
intersects $f_{p}(\alpha )$ only in $f_{p}(\partial \alpha )$. By the
Generalised Loop Theorem 4.13 of \cite{Hem}, this map is homotopic to an
embedding, keeping boundaries fixed.  So then we have an embedding of
$\alpha $, by combining $F$ with $f_{p}\vert \alpha $, and the
combined map homotopic to $f_{n}\vert \alpha $ via a homotopy which is
constant on $\partial \alpha $.  So then the main result of
\cite{F-H-S} (as extended in Section 7 of \cite{F-H-S}) applies and
$f_{n}\vert \alpha $ is homotopic to an embedding, via a homotopy
fixing $\partial \alpha $, with image in an arbitrarily small
neighbourhood of $f_{p}(\alpha )$.  Note that the same argument
applies to $\alpha \cup S\setminus \alpha $.  So this gives an
alternative proof, in a more general setting, of Otal's proof
\cite{Ot3} that Margulis tubes round sufficiently short loops in $N$
are unknotted. It cannot be considered an easier proof, because we need \ref{10.1} (or something similar) to show that 
$F^{-1}F(\partial \alpha \times \{ 0\})=\alpha \times \{ 0\} $. We need this in order to apply the Generalised Loop Theorem.

Now we use the argument of \ref{9.3}.  We fix $\nu $ and $C$ and subsurface 
$\alpha $, and suppose that $\Delta (\nu ,C)$ does not exist.  We are in the
bounded geometry case, exactly as in \ref{9.3}, except that $\alpha $
replaces $S$ and we have short loops corresponding to the components
of $\partial \alpha $, rather than cusps.  We have bounds on the
pleated surfaces and the pleating loci of the corresponding pleated
surfaces by \ref{10.5} and \ref{10.6}, replacing the use of \ref{9.1}
for \ref{9.3}.  By \ref{10.8}, we have bounds, in terms of $\nu $, on
the Margulis tubes intersected by $f_{i}(\alpha )$ with $\partial
\alpha $ in the pleating locus of $f_{i}$.  If $\Delta (\nu ,C)$ does
not exist, then as in \ref{9.3}, we can take geometric limits of
pieces of hyperbolic manifolds $N_{p}$ containing neighbourhoods of
pleated surfaces $f_{i}^{p}(\alpha )$, bounded by pleated surfaces
which are injective on $\pi _{1}$ and take limits of pleated surfaces
within them, exactly as in \ref{9.3}.  The only difference is that we
restrict the pleated surfaces to $\alpha $ and the Margulis tubes
round the loops $(\partial \alpha )_{*}$, in the different manifolds,
$\to \infty $ as $\Delta \to \infty $.

The statement about homotopy follows simply by using the homotopy
through the pleated surfaces between $f_{j}$ and $f_{k}$, since we
have a bound on the diameter of the homotopy between $f_{n}$ and
$f_{n+1}$.

\Box

\section{Proof of coarse biLipschitz.}\label{11}

In this section we prove the following theorem, which is very close 
to 
\ref{1.3}.

\begin{theorem}\label{11.6} Let $e_{i}$ ($1\leq i\leq n$) be the ends 
of 
    $N_{d,W}$. For any choice of end pleated surfaces 
$f_{e,+}$ satisfying the conditions of \ref{10.1}, and 
$M=M([f_{e_{1},+}],\cdots [f_{e_{1},+}])$, 
there is a map $\Phi :M\to N$, with image containing 
$f_{e_{i},+}(S(e_{i}))$ for all $i\leq n$, which is $\Lambda 
_{2}$-coarse-biLipschitz,  for a constant $\Lambda _{2}$ which 
depends only 
on the topological type of $N_{d,W}$ and the constant $c_{0}$ of
\ref{10.1.2}.\end{theorem}

The 
fact that the map is 
coarse Lipschitz follows from the results in Section \ref{10}. The 
main tools for proving coarse biLipschitz are given in \ref{10.8} 
to \ref{10.10}. I believe that this may constitute an important 
difference from the proofs of \cite{B-C-M} and \cite{Bow2}. The work 
done in \ref{10.8} and \ref{10.10} means that the map $\Phi $, which 
we are about to construct formally, already maps bits of surfaces in 
the ends 
in the correct order.

First we need to get into a position to apply \ref{10.1} and the 
other results of section \ref{10}. We need the following lemma.  
For Theorem \ref{1.1} we consider a sequence with $N_{n}=N$ for all 
$n$. For 
Theorem \ref{1.2} we take each $N_{n}$ to be geometrically finite, 
with 
specified ending lamination data.
 
\begin{lemma}\label{11.1} Let $M_{d}$ be a topological model for the 
horoball 
deletion of a hyperbolic 
$3$ manifold with finitely generated fundamental group with ends 
$e_{i}$, $1\leq i\leq m$. Let $[\mu (e_{1}),\cdots \mu (e_{m})]$ be 
any permissible ending invariant.    Let $N_{n}$ be a sequence of 
hyperbolic 
manifolds of this topological type. Identify the ends of 
$(N_{n})_{d}$ with 
the ends of $M_{d}$, and label them $e_{i}$, $1\leq i\leq m$, 
accordingly. Let the ending invariants be $[\mu (e_{1},n),\cdots \mu 
(e_{m},n)]$ and let $\mu (e_{i},n)\to \mu (e_{i})$ as $n\to \infty $.
Regard $S(e)$ as a subsurface of each of the $N_{n}$, 
bounding a neighbourhood of $e$. Fix an end $e$.  Then we can find 
$\Gamma _{+}(e,n)$,  $z_{e,0}$, $c_{0}$, $L_{0}$, 
 $f_{e,+,n}$ so that the conditions of \ref{10.1} are satisfied 
for all but finitely many $n$, and so that any limit of 
$[f_{e_{i},+,n}]$ (using the topology of \ref{3.6}) is in 
${\cal{T}}(S(e_{i}))$ or 
${\cal{O}}_{a}(S(e_{i},N))$ or ${\cal{W}}(S(e_{i},\Gamma ))$ 
according 
to 
what is true for $\mu (e_{i})$, for the same $\Gamma $ in the last 
case, and  all geodesic lamination components equal to the 
geodesic lamination components of $\mu (e_{i})$.
\end{lemma}

\noindent {\em{Proof.}}

We always choose  $\Gamma _{+}(e,n)$ to include the loops which 
bound the support of the geodesic lamination 
components of $\mu (e)$. If $\mu (e,n)\in {\cal{O}}_{a}(S(e,N))$,
then we choose $f_{e,n}$ with $f(S(e,n))$ far out in the end, and 
with pleating locus $\vert f(\Gamma _{+}(e,n)\vert $ bounded, by 
\ref{3.8}. This gives (\ref{10.1.3}). If $\mu (e,n)\in 
{\cal{T}}(S(e_{i}))$,
 then we use \ref{3.9}. 
instead. In this case, we have freedom in the choice of $\Gamma 
_{e,n}$, with $[f_{e,+,n}]$ taken to be $[f_{2}]$ of \ref{3.9}, if 
$\Gamma _{+}(e,n)$ is the $\Gamma $ of \ref{3.9}. In the interval 
bundle 
case we are then finished.

Now  suppose that $\mu (e)\in 
{\cal{W}}(\Gamma ,S(e))$ for some  $\Gamma \in {\cal{O}}(S,N)$. Then
\ref{10.1.2} is satisfied with any $\zeta \in \Gamma $ 
replacing $\Gamma _{+}(e)$, for some constant $c_{0}'$. If a 
component of 
$\mu (e)$ is in ${\cal{T}}(S(\alpha ))$ for some gap $\alpha $ of 
$S(e)\setminus (\cup \Gamma (e))$, then we choose 
$\alpha \cap \Gamma _{+}(e,n)$ to be constant in $n$, and 
in ${\cal{O}}(\alpha ,N)$, so that (\ref{10.1.2}) is satisfied for 
some $c_{0,\alpha }$, for any $\gamma '$, $\mu $, $\mu '$ as in 
(\ref{10.1.2}) with support in $\alpha $. For $\alpha $ such that the 
corresponding component 
$\mu (e,\alpha )\in {\cal{O}}_{a}(\alpha ,N)$, choose any $\zeta 
_{n}\in \Gamma _{e,n}$
. Then we claim that (\ref{10.1.2}) is satisfied 
for suffciently large $n$, with $\Gamma _{+}(e,n)$ replacing $\Gamma 
_{+}(e)$, for some $c_{0}>0$.
For suppose not. Then there are geodesic laminations $\mu _{n}'$ and 
loops $\gamma _{n}'$ which are trivial in $N$ but nontrivial in 
$S(e)$ 
such that, after taking subsequences,
$$i(\mu _{n}',\gamma _{n}')=0,\ \ \lim _{n\to \infty }i\left( \mu 
_{n}',{\zeta _{n}\over \vert \zeta _{n}\vert }\right) <c_{0}.$$
Let $\mu '$ and $\mu ''$ be limits of $\mu _{n}'$ and $\zeta 
_{n}/\vert \zeta _{n}\vert $ respectively. Then we have $i(\mu 
'',\mu ')=0$. In the notation of \ref{3.9}, $\mu (e,n)=[f_{2}]$ and 
$\vert f_{2}(\zeta _{n})\vert $ is bounded. 
So if $\zeta _{n}\subset \alpha $, we have $\mu ''=\mu 
(e,\alpha )$. Now $i(\zeta  ,\mu ')\geq c_{0}'$ for any $\zeta \in 
\Gamma $. So $\mu '$ has 
support at least $c_{0}'$ in any $\alpha $. It is then impossible to 
have $i(\mu '',\mu ')=0$ unless $\mu '$ is a nonzero multiple of $\mu 
''=\mu 
(e,\alpha )$, giving a contradiction.

\Box

\ssubsection{Construction of Lipschitz $\Phi $.}\label{11.2}
We construct the coarse biLipschitz map $\Phi :M\to N$ as a limit of 
maps $\Phi _{n}$, where each $\Phi 
_{n}$ is defined from  a choice of pleated surfaces $f_{e,n}$ for 
each end $e$ of $N_{d}$. Let $e_{i}$, $1\leq i\leq m$, be the ends 
of $N_{d}$. Then we have
$$\Phi _{n}:M_{n}=M([f_{e_{1},n}],\cdots [f_{e_{m},n}])\to N.$$
Then by \ref{6.13}, $M_{n}$ converges geometrically, with a suitable 
basepoint, to $M=M(\mu (e_{1},\cdots \mu (e_{m}))$.
 
Then, to prove \ref{11.6}, it sufficies to show that the 
$\Phi _{n}$ are coarse biLipschitz with 
respect to a constant $\Lambda _{2}$ which is bounded in terms of the 
topological 
type of $N_{d,W}$, and the constant $c_{0}$ of \ref{11.1} (and 
(\ref{10.1.2})), because then we shall have a convergent 
subsequence to $\Phi :M\to N$ which is coarse biLipschitz with the 
same bound on constant, whose image is all of the convex hull except 
for bounded neighourhoods of the convex hull boundary components 
corresponding to geometrically finite ends. In this subsection, we 
complete the formal construction of the map, which, by the results of 
section \ref{10} and earlier, is Lipschitz with respect to a constant 
$\Lambda _{1}$ depending only on  the 
topological 
type of $N_{d,W}$, and the constant $c_{0}$. From now on we drop the 
index $n$ and write $\Phi $ for $\Phi _{n}$.

 We fix a Margulis 
constant $\varepsilon _{0}$. We fix ltd 
parameter functions $(\Delta _{1},r_{1},s_{1},K_{1})$.  
The geometric model is defined 
using this fixed choice of parameter 
functions and a vertically efficient decomposition of $S(e)\times 
[z_{0}(e),y_{+}(e)]$, although, as we saw in \ref{6.6}, different 
choices 
give the same model up to biLipschitz equivalence.
The geometric model 
is a union of pieces  $M(z_{0}(e),y_{+}(e))$, one for each end $e$
of $N_{d,W}$, a model for the noninterval part $W'$ of $W$, and 
possibly some model Margulis tubes wedged in between.
If $N$ is homeomorphic to the interior of an interval 
bundle, when $W'=\emptyset $ and the two end models are combined in 
one. 

Let $f_{j}=f_{e,j}$ be the sequence of pleated surfaces constructed 
in \ref{8.4}, \ref{8.5} using  the ltd parameter functions 
$(\Delta _{1},r_{1},s_{1},K_{1})$.
This was constructed using a sequence of points $z_{j}=[\varphi 
_{t_{j}}]\in 
[z_{0},y_{+}]$, $0\leq j\leq n_{+}$. We have
$$d([f_{j}],[f_{j+1}])\leq L_{0}$$
and by property \ref{8.2.4}, there is a bound by 
$L_{0}$ on homotopy tracks of a homotopy between $f_{j}$ and 
$f_{j+1}$.
Let $S_{t_{j}}$ be as in \ref{6.2}, that 
is, the hyperbolic surface determined by $z_{j}$. If $[\varphi 
_{t}]\in \ell 
=\ell _{k}$ where $\alpha =\alpha _{k}$ is a gap or loop with 
$\alpha \times \ell $ in the decomposition, write $S_{\alpha ,t}$ 
for the subsurface 
$S_{k,t}$ of $S_{t}$ of \ref{6.2}.  
We can also assume 
that the homeomorphism $\varphi _{t_{k}}:S\to S_{t_{k}}$ is part of 
the family $\varphi _{t}$ with the properties of \ref{6.2}. The 
results of \ref{10.1} give bounds on 
$d_{\alpha }(z_{j},[f_{j}])$, whenever $z_{j}\in \ell $ and $\alpha 
\times \ell $ is a set in the vertically efficient decomposition. 
\ref{10.1} also gives an upper bound on $\vert (\partial \alpha 
)_{*}\vert $, but there is also a better upper bound for short loops
in \ref{10.8}. We recall that $S(f_{j})$ is the hyperbolic surface 
structure 
such that $f_{j}:S(f_{j})\to f_{j}(S)$ is isometric restricted to 
sets 
with totally geodesic image. (we sometimes call this map 
${\rm{Imp}}(f_{j})$.) 
We write $S(f_{j},\alpha )$ for the 
subset of $S(f_{j})$ which is homotopic to $f_{j}(\alpha )$, and is 
a  
complementary component of $\varepsilon (\vert f_{j}(\gamma )\vert 
)$-neighbourhoods 
of the loops $f_{j}(\gamma )$ for $\gamma \subset \partial \alpha $, 
where $\varepsilon (L)$ is the function of \ref{6.2}, taken to be 
$\varepsilon _{0}$ whenever $L<\varepsilon _{0}/2$.

Then for each $j$ and each gap $\alpha $ with $\alpha \times \ell $ 
in 
the decompsition of $S\times [z_{0},y_{+}]$, and for $z_{t_{j}}\in 
\ell _{j}$, we map $\varphi _{t_{j}}^{-1}(S_{\alpha ,t_{j}})$ to 
$f_{j}(S(f_{j},\alpha ))$. Because the results of \ref{10.1}, 
\ref{10.6}, \ref{10.8} imply that $S_{\alpha ,t_{j}}$ and 
$S(f_{j},\alpha )$ are boundedly Lipschitz equivalent, we can make 
this map coarse Lipschitz. For $t\in [t_{j},t_{j+1}]$, if 
$[z_{j},z_{j+1}]\subset \ell $, we map $\varphi _{t}^{-1}(S_{\alpha 
,t})$ to the image of the bounded track homotopy between $f_{j}\vert 
\varphi 
_{t_{j}}^{-1}(S_{\alpha ,t_{j}})$ and $f_{j+1}\vert \varphi 
_{t_{j+1}}^{-1}(S_{\alpha ,t_{j+1}})$. Doing this for all $j$ and all 
gaps $\alpha $ for $\alpha \times \ell $ in the decomposition, the 
map $\Phi $ is defined on all of $M$ except for the model Margulis 
tubes and is coarse Lipschitz, with respect to the hyperbolic metric 
on the 
image. For loops $\gamma $ with model Margulis tubes for which at 
least one of $\vert \gamma _{**}\vert $ or $\vert 
\gamma _{*}\vert <\varepsilon _{2}$, the map $\Phi $ maps $\partial 
T(\gamma _{**})$ to a bounded neighbourhood of $\partial 
T(\gamma _{*},\varepsilon _{0})$ (which is nonempty). The map $\Phi 
$ is then also coarse Lipschitz with respect to the induced metric on 
$\partial T(\gamma _{_*})$, and the metric induced from the 
hyperbolic 
metric on a bounded neighbourhood of $\partial T(\gamma 
_{*},\varepsilon _{0})$. Then we can extend the definition of $\Phi $ 
across the model tubes $T(\gamma _{**})$, so that the map is 
Lipschitz and so that $T(\gamma _{**})$ maps into  a bounded 
neighbourhood of $T(\gamma _{*},\varepsilon _{0})$ whenever 
$\vert \gamma _{**}\vert $ or $\vert \gamma _{*}\vert <\varepsilon 
_{2}$. Defining the map in a Lipschitz way is clear when there is  a 
bound on the 
diameter of $T(\gamma _{**})$, that is, if $\vert \gamma _{**}\vert 
\geq \varepsilon _{2}$. It is also clear if $\vert \gamma _{**}\vert 
<\varepsilon _{2}$ because then, by \ref{10.8}, $\vert \gamma 
_{**}\vert $ and $\vert \gamma _{*}\vert $ are boundedly 
proportional (and more) and since the metric is chosen to that 
$T(\gamma _{**})$ is isometric to a Margulis tube, the coarse 
Lipschitz map 
on $\partial T(\gamma _{**})$ can be continuously extended to be 
coarse
Lipschitz on the $r$-neighbourhood of $\gamma _{**}$, for each $r\leq 
R$ where $T(\gamma _{**})$ is the $R$-neighbourhood of $\gamma _{**}$.

We define $\Phi =f_{W'}$ on $W'$, where $W'$ is the 
non-interval-bundle part of $W$, and $f_{W'}$ is the map of 
\ref{8.12}. This is a Lipschitz  map, with Lipschitz constant 
$\Lambda 
_{1}$ 
depending only on the topological type of $W'$, and on $c_{0}$, by 
\ref{8.12}. Also, by \ref{8.12}, the definition of $f_{W'}$ on 
$\partial W'\setminus N_{d}$ agrees with the definition of $\Phi $ on 
the boundaries of the model end manifolds, where the pleated 
surfaces $f_{e,0}$ are used.

\ssubsection{Proof of coarse biLipschitz: vertical 
length.}\label{11.3}

So it remains to show that $\Phi:M\to N $ is coarse biLipschitz onto 
its image.  Some preliminary work on this has 
already been done, notably in \ref{10.9} and \ref{10.10}. 
We write $\tilde{\Phi}$ for the lift $\tilde{\Phi 
}:\tilde{M}\to H^{3}=\tilde{N}$. If we use 
$d_{\tilde{M}}$ and $d_{H}$ to denote the metric on $\tilde{M}$ and 
on $H^{3}$, then all we actually need to do is to show that
\begin{equation}\label{11.2.1}\liminf _{d_{\tilde{M}}(x,y)\to \infty 
}d_{H}(\tilde{\Phi}(x),\tilde{\Phi}(y))>0,\end{equation}
where we take any $x$, $y$ lifting points $x'$, $y'$ where $\Phi 
(x')$, $\Phi (y')$ are in the region bounded by all the pleated 
surfaces 
$f_{e,+}(S(e))$ (for all ends $e$, for the surfaces $f_{e,+}$ used to 
define $\Phi $. This is to ensure that the 
geodesic segment joining $\tilde{\Phi }(x)$ and $\tilde{\Phi }(y)$ 
is in the image of $\tilde{\Phi }$. (\ref{11.2.1})
is sufficient, because if this is true, then for some $\delta >0$ 
given any points $x$, 
$y$, we can choose $w_{i}$, $0\leq i\leq n$ with $w_{i}$ successive 
 points on the geodesic segment joining $\tilde{\Phi}(x)$ and 
$\tilde{\Phi}(y)$ with
$${\delta \over 2}\leq d_{H}(w_{i},w_{i+1})\leq \delta $$
Then since $w_{i}$ is in the image of $\tilde{\Phi}$, we can choose 
$x_{i}$ with $\tilde{\Phi }(x_{i})=w_{i}$. Then (\ref{11.2.1}) gives 
an upper bound on $d_{\tilde{M}}(x_{i},x_{i+1})$ by some $C$ and hence
$$d_{H}(\tilde{\Phi}(x)\tilde{\Phi}(y))\geq {n\delta \over 2}\geq 
{\delta \over 2C}\sum _{i=0}^{n-1}d_{\tilde{M}}(x_{i},x_{i+1})$$
$$\geq {\delta \over 2C}d_{\tilde{M}}(x,y).$$
We shall actually show the slightly stronger statement than 
(\ref{11.2.1}):
\begin{equation}\label{11.2.2}\lim _{d_{\tilde{M}}(x,y)\to \infty 
}d_{H}(\tilde{\Phi}(x),\tilde{\Phi}(y))=+\infty .\end{equation}

Already, in \ref{10.10}, we switched to a stronger set of ltd 
parameter 
functions, and we continue to do this. All our pleated surfaces, the 
geometric model and $\Phi $, are defined using the ltd parameter 
functions $(\Delta _{1},r_{1},s_{1},K_{1})$. We now use a stronger 
set of ltd parameter functions $(\Delta _{2},r_{2},s_{2},K_{2})$, 
such that 
$\Delta _{2}(\nu )\geq 3\Delta (\nu ,C)$ for all $\nu $, for a 
suitable $\Delta (\nu ,C)$ as in 
\ref{10.10}, so that the last part of \ref{10.10} works,
$r_{2}(\nu )\leq \varepsilon _{2}$ for $\varepsilon 
_{2}$ 
as in \ref{10.8}, $s_{2}(\nu )\leq s_{1}(\nu )$ and $K_{2}\geq 
K_{1}$. 
We also choose 
the 
parameter functions strong enough for all the results of Section 
\ref{6} 
to work. Having made this second choice of ltd parameter functions, 
as 
in \ref{6.6} we can choose a vertically efficient decomposition 
with respect to $(\Delta _{2},r_{2},s_{2},K_{2})$ which is coarser 
than the 
vertically efficient decomposition used for $(\Delta 
_{1},r_{1},s_{1},K_{1})$. Naturally, the reason for choosing these 
new 
parameter functions is that we wish to use 
\ref{10.8} and \ref{10.10}. The biLipschitz constant will depend on 
this second choice of parameter functions, while the Lipschitz 
constant, which has already been found, depends on the first ltd 
parameter functions. 

The idea of the proof of (\ref{11.2.2}) is to show that, for homotopy 
classes of paths $\gamma $ in $M$ whose images under $\Phi $ can be 
homotoped 
keeping endpoints fixed to geodesics of bounded length, first, 
$\gamma $ must have bounded vertical length up to homotopy, and then, 
given a bound on the vertical length, we can bound what we call the 
{\em{horizontal 
length}}. In fact, the horizontal length is simply the length. The 
term 
is only used when the {\em{vertical length}} is bounded. (The terms 
arose in earlier, more simplistic, incomplete versions 
of this part of the proof. But they still have a good intuitive 
meaning.) In the combinatorially bounded geometry Kleinian surface 
case, 
the measure of vertical length up to homotopy is just the difference 
of the vertical (second) coordinates of the endpoints in $S\times 
\mathbb R$. Obviously this needs to be refined in the general case, 
when the measurement of length of the vertical coordinate varies 
depending on which subsurface we are in. We say that a path in $M$ 
has 
{\em{bounded vertical length}} (up to homotopy) if it can be 
homotoped 
into a union of bounded diameter subsets of model Margulis tubes and 
boundedly many sets 
$$M(\alpha ,I)=M(\alpha ,\ell )=\cup \{ \varphi _{t}^{-1}(S_{\alpha 
,t}):t\in I\} .$$
Here, $\varphi _{t}$, $S_{t,\alpha }$ are as in the definition of 
geometric 
model $M(z_{0},y_{+})$ with respect to $(\Delta 
_{2},r_{2},s_{2},K_{2})$ for some end $e$, and $\ell =\{ [\varphi 
_{t}]:t\in I$. The set  $\alpha \times 
\ell $ is either  in the chosen vertically 
efficient 
decomposition of $S(e)\times [z_{0},y_{+}]$ (for $(\Delta 
_{2},r_{2},s_{2},K_{2})$) or is a subset of some $\alpha ,\ell ')$ in 
the decomposition. In either case we also have that $\ell $ has 
$d_{\alpha }$-length  
$\leq L=L(\Delta 
_{2},r_{2},s_{2},K_{2})$ for  suitable $L$.

\ssubsection{Basic topological principles}\label{11.4}

We shall use the following.  Let $N$ be a hyperbolic manifold, and for
$Y\subset N$ we let $B_{\delta }(Y)$ denote the $\delta
$-neighbourhood of $Y$, in the hyperbolic metric.  We let $C$ be the
constant of the Injectivity Radius Lemma in \ref{3.3}.  As in
\ref{10.10}, in the following lemma, we shall use the surface $S'$ of
\ref{10.13}, and the neighbourhood $U$ of $e$ bounded by $S'$, in the
case when $e$ is incompressible.  Let $\varepsilon _{1}$ be as in
\ref{10.7}.  By \ref{10.7} and \ref{10.8} we can assume, and shall do
so, that $\varepsilon _{1}$ is small enough that if $\vert \gamma
_{*}\vert<\varepsilon _{1}$, then a model Margulis tube with core
$\gamma _{**}$ in the model manifold does exist, and $\vert \gamma
_{**}\vert <\varepsilon _{0}$.
\begin{ulemma}

Let $S=S(e)$. Let $h_{i}:S\to N$ ($0\leq i\leq 3$) be pleated 
surfaces 
 which satisfy the conclusion of the Radius of 
Injectivity Lemma \ref{3.3}. Let $\alpha $, $\alpha '$ be subsurfaces 
of $S$.
If $e$ is compressible and $\alpha =S$, let $h_{i}(\alpha )\subset U$ for
$i=0$, $1$, $2$.

Let 
$\vert (\partial \alpha )_{*}\vert <\varepsilon _{1}$.  and $\vert 
(\partial \alpha 
')_{*}\vert <\varepsilon _{1}$. Let $\partial \alpha $ be in the 
pleating locus for $h_{i}$, $i=0$, $1$, $2$, and let $\partial 
\alpha '$ be in the pleating locus of $h_{3}$.

Let $\gamma _{0}$, $\gamma _{2}$ be nontrivial, nonperipheral,
 closed loops in $\alpha $, not homotopic to boundary components, 
such that
 $h_{0}$ and $h_{1}$ are homotopic  via a
 homotopy which maps $\partial \alpha $ into $T((\partial \alpha
 )_{*},\varepsilon _{1})$, and with image in $U$ if $e$ is 
compressible and $\alpha =S$, and
 disjoint from $h_{2}(\gamma _{2})$, and similarly with
 $(h_{0},h_{1},\gamma _{2})$ replaced by $(h_{2},h_{1},\gamma _{0})$. 
 Suppose also that any two of $h_{0}$, $h_{1}$, $h_{2}$ are homotopic  via a
 homotopy which maps $\partial \alpha $ into $T((\partial \alpha
 )_{*},\varepsilon _{1})$, and with image in $U$ if $e$ is 
compressible and $\alpha =S$, and disjoint from $h_{3}(\gamma _{3})$.  Here, 
$\gamma _{3}\subset \overline{\alpha
'}$ is a
nontrivial, nonperipheral loop  which is not homotopic to any component of $\partial \alpha $,
but is allowed to be homotopic to a component of $\partial \alpha
'\setminus \partial \alpha $.

Suppose that
\begin{equation}\label{11.4.1}h_{3}(\alpha ')\cap B_{\delta
}(h_{1}(\alpha ))\setminus T((\partial \alpha )_{*},\varepsilon
_{1})\neq \emptyset
\end{equation}
and 
\begin{equation}\label{11.4.2}(h_{i}(\alpha )\cap 
B_{\delta }(h_{j}(\alpha ))\setminus  
T((\partial \alpha )_{*},\varepsilon _{1})
= \emptyset {\rm{\ for\ }}(i,j)=(0,1),\ 
(1,2),\ (0,2).\end{equation}
Then
\begin{equation}\label{11.4.3}(h_{0}(\alpha )\cup 
h_{2}(\alpha ))\cap 
h_{3}(\alpha '))\setminus  
T((\partial \alpha )_{*},\varepsilon _{1})\neq \emptyset .
\end{equation}

\end{ulemma}

We remark that if $\alpha =\alpha '=S(e)$ we can make a slighty 
simpler 
statement, but the above is 
sufficient for our purposes.

\noindent {\em{Proof.}} 

Let $C$ be the constant of he Radius of Injectivity Lemma \ref{3.3}. 
By our hypothesis on the Radius of Injectivity Lemma  Radius of 
Injectivity Lemma \ref{3.3}, $h_{i}$ maps the 
boundary components  of $(S(h_{i}))_{\geq 
\varepsilon _{1}}$ which are homotopic to $\partial \alpha $ into 
$T((\partial \alpha )_{*},\varepsilon _{1})\setminus 
T((\partial \alpha )_{*},\varepsilon _{1}/C)$. Reparametrising and 
taking a perturbation by a 
distance $\log C$ with support inside $T((\partial \alpha 
)_{*},\varepsilon 
_{1})$, 
 we can assume, for a subsurface $\alpha _{1}$ of 
$\alpha $ which is homotopic to $\alpha $, that $\partial \alpha 
_{1}$ is mapped homeomorphically to 
$\partial T((\partial \alpha )_{*},\varepsilon _{1}/C)$. As 
pointed out in \ref{10.7}, $\alpha $ is incompressible unless $\alpha 
=S(e)$, even if $e$ is compressible, and if $e$ is compressible and
$\alpha =S(e)$, we are assuming that $h_{i}(\alpha )\subset U$.  So in
all cases, as in \ref{10.10}, we can apply the results of \cite{F-H-S}. So we can
perturb to a diffeomorphism on $\alpha _{1}$, with image in an
arbitrarily small neighbourhood of the original.

Then, taking another perturbation, we can assume that the image is 
transverse to $\partial T((\partial \alpha )_{*},\varepsilon 
_{1})$. Then we can reparametrise again and assume that $\partial 
\alpha _{1}$ is contained in $\partial T((\partial \alpha 
)_{*},\varepsilon 
_{1})$. We can perturb so that $h_{i}(\alpha _{1})\subset N\setminus 
T((\partial 
\alpha )_{*},\varepsilon _{1})$, perturbing the intersection of the 
image  with 
$N\setminus 
T((\partial 
\alpha )_{*},\varepsilon _{1})$ by an arbitrarily small amount. 
So altogether we can assume that 
$$h_{i}:(\alpha _{1},\partial \alpha _{1})\to (N\setminus 
T((\partial 
\alpha )_{*},\varepsilon _{1}),\partial T((\alpha )_{*},\varepsilon 
_{1}))$$
 is a homeomorphism, keeping the intersection of the image with 
 $N\setminus T((\partial 
\alpha )_{*},\varepsilon _{1})$ in an arbitrarily small 
 neighbourhood of the original, for $0\leq 
i\leq 2$, and, similarly for a subsurface $\alpha _{1}'$ of $\alpha 
'$ which is 
homotopic to $\alpha '$, we can assume that $h_{3}$ is a 
homeomorphism, with
$$h_{3}:(\alpha _{1}',\partial \alpha _{1}')\to (N\setminus 
T((\partial 
\alpha ')_{*},\varepsilon _{1}),\partial T((\alpha ')_{*},\varepsilon 
_{1}))).$$

We are assuming that the $h_{j}(\alpha _{1})$, $j=0$, $1$, $2$, are
disjoint.  By Waldhausen's theorem \cite{Wal} the closed region
bounded by $h_{j}(\alpha _{1})$ and $h_{k}(\alpha _{1})$ and
$T((\partial \alpha )_{*},\varepsilon _{1})$, for $j$, $k\in \{ 
0,1,2\} $ and
$j\neq k$, is homeomorphic to $\alpha _{1}\times [0,1]$.  Given the
boundary-preserving homotopy between $h_{0}(\alpha _{1})$ and
$h_{1}(\alpha _{1})$ in the complement of $h_{2}(\gamma _{2})$, it
follows that $h_{2}(\alpha _{1})$ is not between $h_{0}(\alpha _{1})$
and $h_{1}(\alpha _{1})$.  Similarly, $h_{0}(\alpha _{1})$ is not
between $h_{1}(\alpha _{1})$ and $h_{2}(\alpha _{1})$.  So the region
bounded by $h_{0}(\alpha _{1})$ and $h_{1}(\alpha _{1})$ and
$T((\partial \alpha )_{*},\varepsilon _{1})$, and the
region bounded by $h_{1}(\alpha _{1})$ and $h_{2}(\alpha _{1})$ and
$T((\partial \alpha )_{*},\varepsilon _{1})$, must
be on opposite sides of $h_{1}(\alpha _{1})$.  So $h_{1}(\alpha _{1})$
is in the region bounded by $h_{0}(\alpha _{1})$ and $h_{2}(\alpha 
_{1})$ and
$T((\partial \alpha )_{*},\varepsilon _{1})$.  But because of the
boundary-preserving homotopy between $h_{0}(\alpha _{1})$ and
$h_{2}(\alpha _{1})$ in the complement of $h_{3}(\gamma _{3})$,
$h_{3}(\alpha _{1}')$ is not contained in the region bounded by 
$h_{0}(\alpha _{1})$ and $h_{2}(\alpha _{1})$ and
$T((\partial \alpha )_{*},\varepsilon _{1})$, and does not intersect 
$T((\partial \alpha )_{*},\varepsilon _{1})$ except in components of
$h_{3}(\partial \alpha _{1}')$ for components of $\partial \alpha
_{1}'$ which are also in $\partial \alpha _{1}$.  So if, as we are
assuming in (\ref{11.4.1}) and (\ref{11.4.2}), $h_{3}(\alpha _{1}')$
does intersect the interior of the region bounded by $h_{0}(\alpha 
_{1})$ and
$h_{2}(\alpha _{1})$ and $T((\partial \alpha )_{*},\varepsilon _{1})$,
then $h_{3}(\alpha _{1}')$ must intersect $h_{0}(\alpha _{1})\cup
h_{2}(\alpha _{1})$.  We are working with an arbitrarily small
perturbation of the original $h_{j}$.  This gives (\ref{11.4.3}), as
required.

\Box

\ssubsection{Cross-sections.}\label{11.13}

We extend the notion of 
partial 
order of \ref{7.3}. A {\em{complete cross-section}} (for $S\times 
[z_{0},y_{+}]$) is 
a set of $n$-
tuples for varying $m$:   
$$b=((\beta _{1},x_{1}),\cdots (\beta _{m},x_{m})),$$
such that the $\beta _{i}$ are disjoint but the union of the 
$\overline{\beta _{i}}$ is $S$, and one of the following holds.
\begin{description}
    \item[.]
$(\beta _{i},\ell _{i})$ is an ltd in the decomposition with 
$x_{i}\in \ell _{i}$.
\item[.] $\beta _{i}=\cup _{j}\beta _{i,j}$, where 
$\beta _{i,j}\times \ell _{i,j}$ is a bounded set in the 
decomposition, $x_{i}\in \cap _{j}\ell _{i,j}$,
$\vert \gamma _{**}\vert <\varepsilon _{2}$ for all $\gamma 
\subset \beta _{i}$, and $\vert \gamma 
_{**}\vert \geq \varepsilon _{2}$
for all $\gamma \subset {\rm{int}}(\beta _{i})$.
\end{description}

Then we define
$$((\beta _{1},x_{1}),\cdots (\beta _{m},x_{m}))<
((\beta _{1}',x_{1}'),\cdots (\beta _{n}',x_{n}'))$$
if the following hold.
\begin{description}
    \item[1.] $x_{i}$ is strictly to the left of $x_{j}'$ in 
$[z_{0},y_{+}]$ whenever $\beta _{i}\cap \beta _{j}'\neq \emptyset $.

\item[2.] If $\beta _{i}=\beta _{j}'$ 
then $\beta _{i}$ is ltd with respect to $(\Delta 
_{2},r_{2},s_{2},K_{2})$  along a segment  
centred at $x_{i}$, and similarly for $(\beta _{j}',x_{j}')$.
\end{description}

A {\em{cross-section for $\alpha \times \ell $}} is similarly 
defined, 
whenever $\alpha \times \ell $ is a union of sets in the vertically 
efficient decomposition.

A cross-section 
$b= ((\beta _{1},[\varphi _{t_{1}}]),\cdots (\beta _{m},[\varphi 
_{t_{m}}]))$ for $\alpha $ defines a surface  $S(b)$ in the model 
manifold 
which contains, in the language of \ref{11.3},
$$\cup _{j=1}^{m}M(\beta _{j},\{ t_{j}\} ),$$ 
where 
these are joined up by annuli in model Margulis tubes 
$T(\gamma _{**})$, for each $\gamma $ which is a component of 
$\partial \beta _{j}\cap \partial \beta _{k}$ for some $j\neq k$. A 
cross-section for $\alpha $ also defines an element of 
$x(b)\in {\cal{T}}(S(\alpha )$ up to bounded distance by
$$\pi _{\beta _{j}}(x(b))=[\varphi _{t_{j}}]$$
if $\beta _{j}$ is a gap, and if $\beta _{i}$ is a loop,
$${\rm{Im}}(\pi _{\beta _{j}}(x(b)))={\rm{Im}}(\pi _{\beta 
_{i}}([\varphi _{t_{j}}]),\ \ {\rm{Re}}(\pi _{\beta 
_{j}}(x(b)))={1\over 
\varepsilon _{0}}$$
for some fixed Margulis constant $\varepsilon _{0}$.

A totally ordered chain 
of cross sections $\{ b_{i}:1\leq i\leq r\} $ is {\em{maximal}} 
if it is not possible to insert any $b$ between any $b_{i}$ and 
$b_{i+1}$ with $b_{i}<b<b_{i+1}$, any $1\leq i<r$. 

\begin{lemma}\label{11.14}
Let  $b=((\beta _{1},x_{1}),\cdots (\beta _{m},x_{m}))$, 
$b'=((\beta _{1}',x_{1}'),\cdots (\beta _{n}',x_{n}'))$ 
be cross-sections for $\omega \times \lambda $ for some  
union of sets $\omega \times \lambda $ in the vertically 
efficient decomposition. Suppose that $b<b'$.
Then each $\gamma \subset \partial \beta 
_{i}$, for each $i$, is 
either equal to some $\gamma '\subset \partial \beta _{j}'$, or is 
contained in 
some $\beta _{j}'$ for $\beta _{j}'$ ltd along a segment in
$[x_{i},x_{j}']$, or $\gamma $ intersects $\partial \beta _{j}'$ 
transversally. Similar properties hold for each $\gamma 
'\subset \partial \beta _{j}'$.\end{lemma} 

\noindent{\em{Proof.}} Suppose not so for $\gamma $. 
Then there must be some $j$ such that 
$\gamma $ is contained in the interior of $\beta _{j}'$ and $\beta 
_{j}'\times \ell _{j}'$ is not ltd in the decomposition for 
$x_{j}'\in 
\ell _{j}'$. Consider the set of $(\alpha ,\ell )$ in the 
decomposition such that $\gamma \cap \alpha \neq \emptyset $, 
$\gamma \notin \partial \alpha $, $(\gamma ,x_{i})<(\alpha ,\ell 
)<(\beta _{j}',x_{j}')$, using the usual ordering of \ref{7.3}. By 
the definition of vertically efficient, in 
particular, property 3 of \ref{5.7} applied to $\gamma $, 
the set of such $(\alpha ,\ell )$ is non-empty: since $\gamma \subset 
{\rm{int}}(\beta _{j}')$, ${\rm{int}}(\beta _{j}')\cap \alpha \neq 
\emptyset $ for all such $\alpha $. Take a maximal such 
$(\alpha ,\ell )$. If some $\gamma '\subset \partial \beta _{j}'$ 
intersects $\alpha $ transversally, then by \ref{7.3} we have 
$\gamma '\cap \gamma \neq \emptyset $. So $\alpha \subset \partial 
\beta _{j}'$. By \ref{7.3}, there is no ltd $(\alpha ',\ell ')$ with 
$(\alpha ,\ell )<(\alpha ',\ell ')$ and with $\ell '$ to the left of 
$x_{j}'$, because that would contradict the maximality of $(\alpha 
,\ell )$. So by property 3 of \ref{5.7} applied to components of 
$\partial \alpha $, $\beta _{j}'\subset \alpha $. So $\beta 
_{j}'=\alpha $.

\Box

\begin{corollary}\label{11.15}
The ordering on cross-sections is transitive. 
\end{corollary}

\noindent{\em{Proof.}} There might, perhaps, be a question about 
condition 2 of 
the ordering. So let $b$, $b'$ be as in \ref{11.14}, with $b<b'$ and 
$b'<b''=((\beta _{1}'',\ell 
_{1}''),\cdots (\beta _{p}'',\ell _{p}''))$ and  $\beta _{i}=\beta 
_{k}''$, 
$\beta 
_{j}'\neq \beta _{i}$, $\beta _{j}'\cap \beta _{i}\neq \emptyset $. 
if $\beta _{j}'$ is ltd we have $\beta _{j}'\subset \beta _{i}$ by 
\ref{7.3}. Then $\beta _{j}'=\beta _{i}$ by \ref{11.14}, which is a 
contradiction.  Now for any ltd $(\alpha ,\ell )$ in the 
decomposition 
with $\ell \subset [x_{i},x_{k}'']$ we have $\alpha \cap \beta 
_{i}=\emptyset $ or $\alpha \subset \beta _{i}=\beta _{k}''$ by 
\ref{7.3}. It follows that no component of $\partial \beta _{j}'$ is 
transverse to $\partial \beta _{i}$, by property 3 of vertically 
efficient in \ref{5.7} applied to $\partial \beta _{j}'$. So $\beta 
_{i}\subset \beta _{j}'$ and $\beta _{j}'\subset \beta _{i}$ by 
\ref{11.13}. 
Thus $\beta _{i}=\beta _{j}'=\beta _{k}'$. Then by \ref{11.12} we 
have $\beta 
_{i}=\beta _{j}'$ unless $\beta _{i}$ is ltd on a segment centred at 
$x_{i}$ for $({1\over 2}\Delta 
_{2},r_{2},s_{2},K_{2})$. \Box

\ssubsection{Bounded sets in the model manifold and 
cross-sections.}\label{11.16}

In order to bound vertical length of preimages under $\Phi $ of 
bounded paths, we shall use the following. We use the language and 
notation of \ref{11.3}, \ref{11.13}. We continue to use the ltd 
parameter functions $(\Delta 
_{2},r_{2},s_{2},K_{2})$, the associated constant $\nu _{2}>0$ and 
$L=L(\Delta _{2},r_{2},s_{2},K_{2},\nu _{2})$. We also fix a model 
end 
manifold $M(z_{0},y_{+})$ in $M$.

\begin{ulemma} 
 
\noindent 1. Given $p_{0}$, there is $n_{0}$ such that the following 
 holds. Suppose that $(\beta ,t)$ and $(\beta ',t')$ are such that 
 $M(\beta \{ t\} )$ and  $M(\beta '\{ t'\} )$ are 
disjoint from model Margulis tubes $T(\gamma _{**})$ with $\vert 
\gamma _{**}\vert <\varepsilon _{2}$. 

\noindent 1a) If $\beta \cap \beta 
'=\emptyset $, suppose that, for any $\alpha $ 
with $\beta \subset \alpha \subset S\setminus \beta '$, 
$M(\beta ,\{ t\} )$ and  $M(\beta ',\{ t'\} )$ are not separated by 
totally ordered chains of $p_{0}$ cross-sections for $\alpha $, 
$S\setminus \alpha $, that is, a chain for $\alpha $ above and below 
$M(\beta ,\{ t\} )$ or similarly for $S\setminus \alpha $ and 
$\beta '$, or by a single chain of complete cross-sections.

\noindent 1b) If $\beta \cap \beta '\neq \emptyset $, 
suppose that $M(\beta ,\{ t\} )$ is 
not separated from $M(\beta ',\{ t'\} )$ by  chains of $p_{0}$ 
cross-sections for $\alpha $ above and below $M(\beta ,\{ t\} )$  for 
any $\alpha \supset \beta $ with $\alpha $ contained in the convex 
hull of $\beta $ and $\beta '$ 
nor by a chain of $p_{0}$ complete cross-sections (that is, 
cross-sections for $S$) if $S$ is the convex hull of $\beta $ and 
$\beta '$.

Then there is a connected set $M_{0}\subset M$ which is a union of 
$\leq n_{0}$ sets which are either sets $M(\alpha _{i}',I_{i})$ 
bounded by $L$, or  subsets of model Margulis 
tubes of diameter $\leq L$, and 
$M_{0}$ contains both  $M(\beta ,\{ t\} )$ and 
$M(\beta ',\{ t'\} )$. 

\noindent 2. Conversely, given $n_{0}$  there exists $p_{0}$, such 
that, 
if $M_{0}\subset M$ contains a totally ordered chain of $\geq p_{0}$ 
cross-sections for $\alpha $, for some $\alpha $, then $M_{0}$ is not 
contained in a connected union of $\leq n_{0}$ submanifolds $M(\omega 
,I)$ bounded by $L$,  and $\leq n_{0}$ subsets of model Margulis 
tubes of diameter $\leq 
L$ \end{ulemma}

\noindent {\em{Proof of 1.}}
Fix $p_{0}$, and suppose $(\beta ,t)$ and $(\beta ',t')$ 
are as in the statement of 1. Let ${\cal{P}}$ be the set of all 
$(\alpha ,\ell )$ with $\alpha \times \ell $ in our fixed vertically 
efficient decomposition of $S\times [z_{0},y_{+}]$ with respect to 
$(\Delta _{2},r_{2}s_{2},K_{2},\nu _{2})$. 

The basis of the proof is to show how to travel from $S(b)$  to 
$M(\beta ,\{ t\} )$, where $b=b(E,\pm )$ is defined using the upper 
or 
lower  boundary of 
a set $E=E(b)\subset {\cal{P}}$. Here, $E$ is closed under $\leq $ or 
$\geq $, and 
 $b$ is a cross-section for some $\alpha $ 
 which is above or below  by $(\beta ,[\varphi _{t}])$,
 and satisfies certain conditions, and similarly with $\beta '$ 
replacing 
 $\beta $.

To get into this situation, if $(\beta ',[\varphi _{t'}])<(\beta 
,[\varphi _{t}])$ and $\beta \neq \beta '$, we take $\alpha $ to be 
the convex hull of $\beta $ 
and $\beta '$. If $\beta '$ is not contained in $\beta $,
then we define $E$ to be the set of 
all $(\omega ,\ell )\leq (\beta ,[\varphi _{t}])$ with 
$\omega $ contained in $\beta '$ or contained in $\alpha \setminus 
\beta '$. If $\beta '\subset \beta $, $\beta '\neq \beta $ then
we interchange the role of $\beta $ and $\beta '$. If $\beta =\beta 
'$ 
and, without loss of generality, $(\beta ',[\varphi _{t'}])\leq 
(\beta ,[\varphi _{t}])$, then we can take 
$$E=\{ (\omega ,\ell )\in {\cal{P}}:(\omega ,\ell )\leq (\beta 
',[\varphi 
_{t'}])\} $$ 
if there is no chain of $p_{0}$ cross-sections between 
$(\beta ',[\varphi _{t'}])$ and $(\beta ,[\varphi _{t}])$. We also 
define $\alpha =\beta $ in this case. If there 
is such a chain, then we take $\alpha =S$, and 
\begin{equation}\label{11.16.1}E=\{ (\omega ,\ell )\in {\cal{P}}:
    (\omega ,\ell )\leq (\beta 
,[\varphi _{t}])\}\cup \{ (\omega ,\ell):\omega \cap \beta =\emptyset 
\}  
\end{equation} 
  If $\beta \cap 
\beta '=\emptyset $, and there is no chain of $p_{0}$ cross-sections 
for $S\setminus \beta '$ below $(\beta ,[\varphi _{t}])$, then we 
also 
define $E$ by (\ref{11.16.1}).
Other cases are similarly dealt with by interchanging $\beta $ and 
$\beta '$, or replacing $\leq $ by $\geq $. In all cases, using 
\ref{7.3}, 
we see that if $(\omega _{i},\ell _{i})$, $i=1$, $2$, are ltd with 
$(\omega 
_{1},\ell _{1})<(\omega _{2},\ell _{2})$ and $(\omega 
_{2},\ell _{2})\in E$, then $(\omega _{1},\ell _{1})\in 
E$ also. We then define $b=b(E,+)$ to be the cross-section associated 
with the upper boundary of $b$.

If $b$ is below $(\beta ,[\varphi _{t}])$, we then define
$$E'=\{ (\omega ,\ell )\in {\cal{P}}:\omega \subset \beta , \} \cup 
\{ (\omega ,\ell )\in 
E:\omega \subset \alpha ,\omega \cap \beta =\emptyset \} .$$
Define $b'=b(E',+)$ to be the cross-section associated to the upper 
boundary of 
$E'$. Note that $b$ has been chosen so that $S(b)$ does 
not have long intersections with model Margulis tubes, 
except in the boundary 
around $(\partial \alpha _{*})$. Then the same is true for $S(b')$, 
because if $(\gamma ,\ell )\in {\cal{P}}$ and $(\gamma ,\ell )\in 
E'$, 
then either $(\gamma ,\ell ')\in 
E$ for some $\ell '$, or $(\gamma ,\ell ')\in 
E'$. Also, if $\alpha 
\setminus \beta \neq \emptyset $, $d_{\alpha \setminus \beta }
(x(E,+),x(E',+))\leq L$. So now, replacing $b$ by $b'$ if necessary, 
assume that any totally ordered chain of cross-sections for $\alpha $
between $b$ and $(\beta ,[\varphi _{t}])$ has $\leq p_{0}$ elements 
and any chain between $b$ and $b'$ for $\alpha _{1}$ properly 
contained in $\alpha $ has $\leq p_{0}$ elements. We continue to 
assume that $b\leq b'$, since the modifications if $b'\leq b$ are 
trivial.

Now let
$$E''=\{(\omega ,\ell )\in E':(\omega ,\ell )\leq (\beta ,[\varphi 
_{t}])\} \cup \{ (\omega ,\ell )\in E:\omega \cap \beta =\emptyset \} 
.$$
Let $b''=x(E'',+)$ be the upper boundary of $E''$. Then by definition 
of $E'$, 
$b\leq b''\leq b'$. Let $b_{0}<\cdots <b_{p}$ be a maximal totally 
ordered chain of cross-sections for $\alpha $ with $b_{0}=b$ and 
$b_{p}=b''$. Since the chain is maximal,  there is a gap or loop 
$\omega _{i}\subset \alpha $ such that $(\omega 
_{i},[x(b_{i}),x(b_{i+1})])$ 
is bounded by $L=L(\Delta 
_{2},r_{2},s_{2},K_{2},\nu _{2})$. Otherwise we would be able to 
insert a 
cross-section strictly between $b_{i}$ and $b_{i+1}$. We are now 
going to 
produce $M_{0}$ by induction on $p$ and on the topological type of 
$\alpha $. If $p=1$, we simply take 
$$M_{0}=M(\omega ,I_{0})\cup M(\beta ,\{ t\} ),$$
where 
$$I_{0}=\{ s:\pi _{\omega }([\varphi _{s}])\in [\pi _{\omega 
}(x(b_{0})),\pi _{\beta }(x(b_{1}))]\} .$$
If $p=2$, then the definitions of $E$, $E''$ ensure that we can 
assume that $\beta _{1}\cap \omega \neq \emptyset $ and $\beta 
_{1}\cap \beta \neq \emptyset $, and we then define 
$$M_{0}=M(\omega ,I_{0})\cup 
M(\beta _{1},I_{1})\cup M(\beta ,\{ t\} ),$$
where $I_{1}$ is defined similarly to $I_{0}$, but with $\beta 
_{1}$ replacing $\omega $ and $b_{1}$, $b_{2}$ replacing $b_{0}$ and 
$b_{1}$. It is not immediately clear that we can proceed similarly in 
the case of $p>2$. But let $E_{1}$ be defined similarly to $E$ in the 
case when $(\beta ',[\varphi _{t'}])<(\beta ,[\varphi _{t'}])$
but using $(\beta _{1},x(b_{1}))$ in place of $(\beta ,[\varphi 
_{t'}])$, and 
then define $E_{2}$, $E_{3}$ similarly to $E'$, $E''$, but  using 
$E_{1}$ in place of $E$. These define cross-sections 
$b_{1,j}=x(E_{j},+)$ 
and $b_{1}'$ with $b_{1,1}\leq b_{1,3}<b_{1,2}\leq b'$. Let $\alpha 
_{1}$ be 
the convex hull 
of $\beta _{1}$ and $\beta $. Thus, $\alpha _{1}\subset \alpha $. A 
maximal chain of cross-sections between $b_{1,1}$ and $b_{1,3}$  for 
$\alpha _{1}$ has $\leq p-1$ elements if $\alpha _{1}=\alpha $, 
because otherwise the previous chain of length $p$ is not maximal. If 
$\alpha _{1}=\alpha $ we can use the inductive hypothesis on $p$. If 
$\alpha _{1}$ is strictly contained in $\alpha $, then, by the 
definition of the $E_{j}$, $d_{\alpha \setminus \alpha 
_{1}}(x(b_{1,1}),x(b_{1,2}))\leq L$ (where $L=L(\Delta 
_{2},r_{2},s_{2},K_{2},\nu 
_{2})$).  So if $\alpha _{1}\neq \alpha $, there is a bounded 
distance 
between the surfaces in the model manifold defined by $b$ and $b'$, 
within the set $M(\alpha \setminus \alpha _{1},I)$ where
$$I=\{ s:\pi _{\alpha \setminus \alpha _{1}}([\varphi _{s}])\in \pi 
_{\alpha \setminus \alpha _{1}}([x(b),x(b')]).$$
So we are then reduced to finding a bounded path from at least one of 
$S(b_{1,1})$ or $S(b_{1,2})$ to $M(\beta ,t)$, that is, in between 
two 
cross-sections for $\alpha _{1}$, and the proof of the 
inductive step is 
completed.

\noindent {\em{Proof of 2.}} If $\alpha $ is a loop and we have a 
chain of $p_{0}$ cross-sections for $\alpha $, then the 
cross-sections 
bound a piece of the model Margulis tube $T(\alpha _{**})$ which can 
be made arbitrarily large by choice of $p_{0}$. Now suppose that 
$\alpha $ 
is a gap. 
A cross-section for $\alpha $ defines a 
surface in the model manifold with boundary components in the set 
$(\partial \alpha )_{**}$. The surfaces corresponding to two 
successive cross-sections $b_{1}$, $b_{2}$ with $b_{1}<b_{2}$
either intersect different sets of model 
Margulis tubes in their interiors, or if neither intersects any
 Margulis tubes in the interior, both have single coordinates 
$(\alpha ,
 x_{1})$ and $(\alpha ,x_{2})$ such that $\alpha $ is long, $\nu 
 $-thick and dominant along $[x_{1},x_{2}]$. So if $p_{0}$ is 
 sufficiently large, the set  containing a chain of $p_{0}$ 
successive 
 cross-sections cannot be contained in $\leq n_{0}$ submanifolds 
 $M(\omega ,I)$ 
 bounded by $L$.
 
 \Box\ssubsection{Bounding vertical length.}\label{11.5}

We are now ready to bound vertical length of preimages under $\Phi $, 
of paths which are bounded in $N$. The idea is that parts of the 
image of $\Phi $ are ordered in the correct way. $\varepsilon _{1}$-
Margulis tubes are ordered the same way as the corresponding model 
Margulis  tubes by \ref{10.9}. Images corresponding to ltd gaps are 
mapped in 
the correct order by \ref{10.10}. We use the 

\begin{utheorem} Given $\Delta _{1}$, there is $n_{0}$
    depending only on $\Delta _{1}$, the topological type of $N$, and 
    the constant $c_{0}$ of \ref{10.1.2}, such that the 
following holds. Let $\gamma $ be a path in the model manifold $M$ 
such that $\Phi (\gamma )$ is homotopic, fixing endpoints to a 
geodesic in $N$ of length 
$\leq \Delta _{1}$. Then $\gamma $ is homotopic, fixing endpoints, 
to a path in a union of $\leq n_{0}$ sets $M(\alpha _{i}',I_{i})$ 
which 
are bounded by $L=L(\Delta _{2},r_{2},s_{2},K_{2},\nu _{2})$. 

\end{utheorem}

\noindent {\em{Proof.}}  
It suffices to look at paths with both endpoints in a single end 
model manifold $M(z_{0}(e),y_{+}(e))$, because then $\Phi (M(S(e),\{ 
0\} ))$ cannot separate $\Phi (x)$ from the associated end for $x\in 
M(z_{0}(e),y_{+}(e))$
sufficiently far, in the model metric, from $M(S(e),\{ 
0\} )$. Let  one end of the path $\gamma $ be in $M(\beta ,\{ t\} 
)\subset 
M(z_{0}(e),y_{+}(e))$, and the other end in $M(\beta ',\{ t'\} )$. 
Suppose that $\gamma $ is not 
homotopic to a path  contained in $\leq n_{0}$ sets. Then we apply 
\ref{11.16} and use the equivalent description of bounded sets in 
terms of chains of cross-sections. 
If $M(\beta ',\{ t'\} )$ is in the model in a different 
end or in the model for $W'$, or $(\beta ,[\varphi _{t}])<(\beta 
',[\varphi _{t'}])$ or $(\beta ',[\varphi _{t'}])<(\beta 
,[\varphi _{t}])$, then by \ref{11.16} there is a set of 
$2p_{0}$ cross-sections $b_{i}$ ($\vert i\vert \leq p_{0}$)
either for some $\alpha \supset \beta $
centred on a cross-section with $(\beta ,[\varphi _{t}])$ as a 
coordinate  or for some 
$\alpha \supset \beta '$
centred on a cross-section with $(\beta ',[\varphi _{t'}])$ as a 
coordinate, and in both cases 
with $S(b_{\pm i})$ separating $M(\beta ,\{ t\} )$ from $M(\beta '\{ 
t'\} )$ , or $p_{0}$ cross-sections for $\alpha $ below or above 
$(\beta ,[\varphi _{t}])$ 
and another $p_{0}$ for $S\setminus \alpha $ above or below 
$(\beta ',[\varphi _{t'}])$, 
such that again the $S(b_{i})$
separate 
$M(\beta ',[\varphi _{t}'])$ from $M(\beta ,[\varphi _{t}])$. 

Let 
$U$ be the neighbourhood of $e$ bounded by $S'$ of \ref{10.13}, if 
$e$ is compressible.

Suppose that $b_{k}$ is a totally ordered chain of cross-sections for 
$\alpha $, for $0\leq k\leq 3$,
$$b_{k}=((\beta _{1,k},x_{1,k}),\cdots (\beta _{n,k},x_{n,k})).$$
Let $S(b_{k})$ be the surface for $b_{k}$ as defined in \ref{11.13}. 
Then up to moving $x_{j,k}$ a bounded $d_{\beta _{j,k}}$ distance, 
$\Phi (S(b_{k}))$ is a pleated surface $h_{b}:\alpha \to N$ with 
$$h_{b}\vert \omega =f_{n(\omega )}\vert \omega .$$ 
Here, 
$\omega \subset \beta _{j,k}$ for some $(j,k)$ --- possibly with  
$\omega 
=\beta _{j,k}$ --- and  $\omega \times \{ x_{j,k}\} \subset \omega 
\times 
\ell $, where $\omega \times \ell $ is one of the sets in the 
vertically efficient decomposition for $(\Delta 
_{2},r_{2},s_{2},K_{2},\nu _{2})$ and $f_{n(\omega )}$ is one of the 
original pleated surfaces $f_{j}$  of \ref{8.4}.  So 
we have good bounds on the pleating locus with respect to 
$(\Delta _{1},r_{1},s_{1},K_{1})$. Then
using \ref{10.9} and \ref{10.10}, the $h_{b_{k}}$ satisfy the 
hypotheses 
of the $h_{k}$ in \ref{11.4}, with $\alpha =\alpha '=S$, apart 
perhaps for 
$h_{k}(S)\subset U$. 
But  by \ref{10.12}, there is a bound on the number of complete 
cross-sections in any totally ordered chain which are not contained 
in 
$U$, in the cases when $U$ is needed, that is, $e$ is compressible.

Now suppose the first of the situations deduced from \ref{11.16} 
occurs, 
that the $b_{i}$,$\vert i\vert \leq p_{0}$ are such that $b_{0}$ 
has $(\beta ,[\varphi _{t}])$ as a coordinate. The other cases are 
dealt with very similarly, so we shall just do this case.
Fix any suitable $\delta >0$. We claim that, by 
\ref{11.4}, there is an integer $p_{1}$ which is bounded in terms 
of $\varepsilon _{2}$ and $(\Delta _{2},r_{2},s_{2},K_{2})$
such that, given any totally 
ordered chain of $4p_{1}+1$ cross-sections for $\alpha $, the pleated 
surfaces for the first and the last cannot intersect the $\delta 
$-neighbourhood of the middle 
one. If $\alpha =S$, we can assume that the number of complete 
cross-sections 
which are not contained in $U$, is $<p_{1}$.  So we now have at 
least 
$p_{1}$ complete cross-sections on either side of the middle one. 
and the corresponding pleated surfaces satisfy the hypotheses of 
\ref{11.5}. Number cross-sections and pleated surfaces
$b_{i}$, $h_{b_{i}}=h_{i}$, $-p_{1}\leq i\leq p_{1}$. 
 Suppose that 
\begin{equation}\label{11.5.1}h_{-p_{1}}(\alpha )\cap B_{\delta }
(h_{0}(\alpha ))\setminus  
T((\partial \alpha )_{*},\varepsilon _{1})\neq 
\emptyset \end{equation}
and fix $i$ 
with $-p_{1}<i<0$. Apply \ref{11.4} with 
$h_{i}$, $h_{0}$, $h_{j}$, $h_{-p_{1}}$ replacing $h_{0}$, $h_{1}$, 
$h_{2}$, $h_{3}$, for varying $j$, $0<j\leq p_{1}$. Then either
\begin{equation}\label{11.5.2}B_{\delta }(h_{i}(\alpha ))\cap 
(h_{0}(\alpha )
    \cup h_{-p_{1}}(\alpha ))\setminus  
T((\partial \alpha )_{*},\varepsilon _{1})\neq \emptyset 
\end{equation}
or
\begin{equation}\label{11.5.3}B_{\delta }(h_{j}(\alpha ))\cap 
(h_{0}(\alpha )\cup 
    h_{-p_{1}}(\alpha ))\setminus  
T((\partial \alpha )_{*},\varepsilon _{1})\neq \emptyset {\rm{\ for\ 
}}0<j\leq p_{1}.\end{equation}
But (\ref{11.5.3}) is impossible for $p_{1}$ sufficiently large, 
because, by the definition of a totally ordered set of 
cross-sections (\ref{11.13}), 
it would imply the existence of $\geq p_{1}$ different 
bounded geodesics  and boundaries of $\varepsilon _{1}$-Margulis 
tubes 
intersecting
a set of diameter $C_{2}$, for $C_{2}$ depending only on 
 the topological type and $\varepsilon _{1}$ . Any two complete 
cross-sections 
 differ by at 
 least one Margulis tube or by at least one bounded loop, and recall 
 that $\varepsilon _{1}$ depends only on $(\Delta 
_{1},r_{1},s_{1},K_{1})$. Similarly, 
(\ref{11.5.2}) cannot hold for all $i$ if $p_{1}$ is sufficiently 
large. So if $p_{1}$ is sufficiently large, 
(\ref{11.5.1}) cannot hold, as required. A similar argument works for 
$h_{p_{1}}(\alpha )$.

Then apply \ref{11.4} again with $h_{3}$ of \ref{11.4} taken to be  
a pleated surface $h'$ with domain $\beta '$. If 
$p_{1}$ is sufficiently large given $p_{2}$ we have, after 
renumberings, an ordered set of surfaces $h_{i}(\alpha )$ 
$-2p_{2}\leq 
i\leq 2p_{2}$ such that 
all the surfaces $h_{i}(\alpha )\setminus T((\partial \alpha 
)_{*},\varepsilon _{1})$ are disjoint. Then by \ref{11.4}, 
$h'(\beta ')$ cannot intersect $h_{i}(\alpha )\setminus T((\partial 
\alpha 
)_{*},\varepsilon _{1})$ for $\vert i\vert \leq p_{2}$. If $p_{2}$ is 
sufficiently large given $\Delta _{1}$ and $\delta $, the $\Delta 
_{1}$ 
neighbourhood of $\Phi (M(\beta ,\{ t\} )$ in $N$ is contained in the 
region 
bounded by $h_{\pm p_{2}}(\alpha )$ and the union of the Margulis 
tubes $T(\zeta _{*},\varepsilon _{1})$ for $\zeta  \subset \partial 
\alpha $. So then if $\gamma _{1}$ is the geodesic  segment homotopic 
to $\Phi (\gamma )$ with endpoints fixed, $\gamma _{1}$ has length 
$> \Delta _{1}$, giving the required contradiction.

\Box

\ssubsection{Bounded diffeomorphism on Margulis tube
boundaries.}\label{11.7}

We now claim that the map $\Phi $ can be chosen to be a bounded
diffeomorphism with bounded inverse between model and actual Margulis
tubes, with respect to the model Riemannian metric and the hyperbolic
metric on $N$.  Since the model Margulis tube is an actual Margulis
tube, and the right one up to bounded distortion, it suffices to make
$\Phi $ into a diffeomorphism with bounded inverse between the
boundaries of the model and actual Margulis tubes.  As usual, the constant
$\Lambda _{1}$ (enlarged from the previous one if necessary) depends,
ultimately, only on the topological type of $N_{d,W}$ and the constant
$c_{0}$.  So fix $\gamma $.  It suffices to consider only $\gamma $
with $\vert \gamma _{**}\vert <\varepsilon _{2}$.  Consider all those
$(\alpha ,\ell)$ in the vertically efficient partition ${\cal{P}}$
with $\gamma \subset \partial \alpha $.  It suffices to show that, for
some $\Delta _{1}$, and given any $M(\alpha ,\ell )$ with $\gamma
\subset \partial \alpha $, there is at least one $(\alpha ',\ell ')$
between distance $\Delta _{1}$ and $2\Delta _{1}$ away along $\partial
T(\gamma _{**})$ which is bounded away in $M$ itself: we can then
apply \ref{11.5} to see that $\Phi (M(\alpha ,\ell ))$ and $\Phi
(M(\alpha ',\ell '))$ are bounded apart, and hence the intersections
with $\Phi (T(\gamma _{**}))$ are bounded apart.  But this is
immediate, because if we fix $M(\alpha ,\ell )$ and some bounded
neighbourhood $M_{1}$ of $M(\alpha ,\ell )$ in the model manifold,
only finitely many disjoint bounded sets $M(\alpha ',\ell ')$ can
intersect $M_{1}$.  So by choosing $\Delta _{1}$ large enough
(assuming that $\varepsilon _{2}$ is sufficiently small, as we may do)
the claim is proved.

\ssubsection{Bounding horizontal length in the model 
ends: the nature of the bounded set.}\label{11.8}

Now let $\gamma $ be a path in $M$ of bounded vertical length with 
lift $\tilde{\gamma }$ to the universal cover. We need to show that 
if $\Phi (\gamma )$ is homotopic, keeping endpoints fixed,  to a 
geodesic 
$\gamma _{1}$ of 
bounded length, then $\gamma $ is homotopic in $M$ to a path of 
bounded 
length. We now know that if $\gamma _{1}$ is bounded, then $\gamma $ 
lies in a bounded subset $M_{1}$  of $M$, up to 
homotopy preserving endpoints, where the bounds on $M_{1}$ depend on 
topological type, the ltd parameter functions (which depend only on 
topological type) and our constant $c_{0}$ (of \ref{10.1.2}). We 
shall 
choose $M_{1}$ in a particular way. First, suppose that $M_{1}$ is 
in an end model manifold homeomorphic to $S\times [0,1]$ or $S\times 
[0,\infty )$ for some finite type surface $S$. Theorem \ref{11.5} 
implies 
that, if $\gamma $ is sufficiently far from the model core, as we 
assume for the moment, it 
is possible to choose $M_{1}$ so that the boundary of $\partial 
M_{1}$ 
is the union of two sets which we call the {\em{vertical boundary}} 
$\partial _{v}M_{1}$ and {\em{horizontal boundary}} $\partial 
_{h}M_{1}$, where these have the following properties. The vertical 
boundary $\partial _{v}M_{1}$ is contained in the union of boundaries 
of model Margulis tubes. The horizontal boundary $\partial _{h}M_{1}$ 
is a union of sets of the form $\varphi _{x}^{-1}(S_{\beta 
,x})=M(\beta ,x)$, and, by \ref{11.16} each set $M(\beta ,x)$ is 
separated from 
$\gamma $ by one or two totally ordered chains of cross-sections for 
$\alpha $, for some $\alpha $. There are three possibilities. A 
single  
totally ordered chain might encase $M(\beta ,x)$ and separate it 
from $\gamma $, or might encase $\gamma $, separating $\gamma $ from 
$M(\beta ,x)$, or might be a totally ordered chain for $S$ which 
separates a component of $\partial _{h}M_{1}$ from $\gamma $. We 
call components of $\partial _{h}M_{1}$ of the first type 
{\em{encased}} and components of the last two types {\em{separated}}

By \ref{11.16}, the length $2p_{1}+1$ of the totally ordered chain or 
chains can be taken as large as we like
be anything we like, by allowing $M_{1}$ to have sufficiently large 
diameter in the model metric (but still bounded). We 
choose $p_{1}=k_{0}p_{0}$, where $p_{0}$ is as in \ref{11.5} and for 
$k_{0}$ depending on the topological type of $S$. As in 
\ref{11.5} we take $p_{0}$ large enough that if a set of $p_{0}$ 
cross-sections for $\beta $ exists, then $(\partial \beta )_{*}\vert 
<\varepsilon _{1}$ and there are model Margulis tubes round all 
components of $\partial \beta $ --- unless, of course, $\partial 
\beta =\emptyset $. Note that there is no a priori claim that the 
regions 
bounded by different chains of cross-sections are disjoint, but 
chains of cross-sections for $\alpha $ and $\alpha '$ can only 
interlink if one of $\alpha $, $\alpha '$ is contained in the other. 
So the number of interlinkings of different chains is bounded by the 
topological type of $S$.    Also by 
\ref{11.5}, if $p_{0}$ is large enough, the images under $\Phi $
of $S(b_{kp_{0},\alpha })$, for distinct integers $k$, are disjoint. 
We 
reindex these cross-sections as $b_{k,\alpha }$, $\vert k\vert\leq 
k_{0}$. 
If $k_{0}$ is large enough, given $k_{1}$, we can reduce some chains 
of cross-sections to $b_{k,\alpha }$ for $k\leq k_{0,\alpha }$ for 
some $k_{0,\alpha }\geq k_{1}$ and $\leq k_{0}$, so that each chain 
is either encased by some other or disjoint from all others. Then we 
can discard any  chains which does not encase $\gamma $
and is encased by another which also does not encase $\gamma $, or is 
disjoint from some chain that does encase $\gamma $.
Similarly if a chain for $\alpha '$ 
encases a chain for $\alpha $ which encases $\gamma $, then we can 
discard the chain for $\alpha '$. So now we have a set of disjoint 
chains $b_{k,\alpha _{i}}$ ($\vert k\vert \leq k_{1}$) for 
$\alpha _{i}$, $2\leq i\leq r$, some $r$ 
bounded in terms of the diameter of $M_{1}$, encasing 
components of $\partial _{h}M_{1}$, and (probably) two chains 
$b_{k,\alpha _{1},1}$ and $b_{k,\alpha _{1},2}$ for 
some $\alpha _{1}$ encasing $\gamma $. If $\alpha _{1}$ does not 
exist, then there must be bounded paths from $\gamma $ to $M(S,\{ 
z_{0}\} )$ or $M(S,\{ y_{+}\} )$, possibly to both. In that case, we 
still take chains $b_{k,\alpha _{1},j}$ of some length $\geq 1$, 
and take $\alpha _{1}=S$.   So $b_{0,\alpha _{1},j}$ exists for 
$j=1$, $2$.
Index so that $S(b_{k,\alpha _{1},j})\subset M_{1}$ for $k\geq 0$, 
$j=1$ $2$.

Now fix some suitable $k_{2}<k_{1}$. Let $M_{1}'$ be the intersection 
of $M_{1}$ with the set  bounded by 
$S(b_{0,\alpha _{1},1})$ and $S(b_{0,\alpha _{1},2})$ and 
model 
Margulis tubes round $\partial \alpha _{1}$. For $2\leq 
j\leq r$, let $M_{i}'$ be the set bounded by $S(b_{\pm k_{2},\alpha 
_{i}})$ and model 
Margulis tubes round $\partial \alpha _{i}$. Thus 
$$M_{i}'\subset M_{1}'{\rm{\ for\ }}2\leq i\leq r.$$
We can now assume without loss of generality that
\begin{equation}\label{11.8.1}M_{1}=M_{1}\cap 
    (M_{1}'\setminus \cup _{i=2}^{r}M_{i}'),\end{equation}
and
\begin{equation}\label{11.8.2}\partial _{h}M_{1}\subset S(b_{0,\alpha 
_{1},1})
    \cup S(b_{0,\alpha 
_{1},2})\cup \cup _{i=2}^{r}S(b_{\pm k_{2},\alpha 
_{i}}).\end{equation}

The 
most obvious way to estimate length of a path  in $M_{1}$ is to take 
some cell 
decomposition of $M_{1}$, and count the number of cell boundaries 
crossed by the 
path. The most obvious way to make a cell decomposition is, to first 
decompose $M_{1}$ into interval  bundles by cutting along horizontal 
surfaces, and then to decompose each interval bundle by annuli. We 
can use annuli whose boundaries have bounded length. Essential 
intersections 
with the horizontal surfaces give a lower bound on length. The annuli 
do not have boundary in $\partial M_{1}$, but we are going to modify 
them to obtain surfaces with boundary in $\partial M_{1}$.

\begin{lemma}\label{11.9} Given $C_{1}$ there is $C_{2}$ such that 
the following 
holds. Let $M_{2}\subset M_{1}$ be homeomorphic to $\alpha 
\times [0,1]$ for some $\alpha \times [0,1]$ and such that $\partial 
\alpha \times [0,1]$ is the intersection of $M_{2}$ with $\partial 
_{v}
M_{1}$. Let $\zeta \subset 
M_{2}$ be a nontrivial  simple closed loop of length $\leq C_{1}$ in 
the model metric. Then there is a surface $S(\zeta )\subset M_{1}$ of 
area and diameter $\leq C_{2}$, in the metric on $S(\zeta )$ induced 
by the model metric, with boundary in $\partial M_{1}$, and such 
that $\zeta $ is contained in $S(\zeta )$ up to free homootpy in 
$M_{1}$. In fact, in the cover of $M_{1}$ determined by $M_{2}$, 
$S(\zeta )$ lifts to an annulus $A(\zeta )$ homotopic to $\zeta $.
\end{lemma}

\noindent {\em{Proof.}}  We start by considering an annulus 
$A_{1}(\zeta )\subset 
M$ which is 
homotopic to $\zeta $ and has one boundary component above $M_{1}$, 
that is, nearer the end, and other boundary component below $M_{1}$. 
For example, let $\zeta _{s}$ be the geodesic on $S_{s}$ which is 
freely homotopic to $\varphi _{s}(\zeta )$, and define 
$$A_{1}(\zeta )=\cup _{s\in 
I}\varphi _{s}^{-1}(\zeta _{s})\times \{ s\} $$
for some suitable interval $I$. 
Now 
we consider $T_{1}(\zeta )=A(\zeta )\cap M_{1}$, and $\partial 
T(\zeta )=A(\zeta )\cap \partial 
M_{1}$. There is no reason why $T(\zeta )$ should have
bounded area or diameter. 
Nevertheless, the intersection of $\partial T(\zeta )$ with 
$\partial _{v}M_{1}$ is bounded, the number of components 
of intersection with $\partial _{h}M_{1}$ is bounded, intersection 
with $M_{2}$ is bounded and 
intersections with $M_{1}$ of  homotopy 
tracks between the boundary components of $\partial A_{1}(\zeta )$ 
are 
bounded. The natural homotopy is given by 
$$(s',z)\mapsto \varphi _{s'}^{-1}\circ \varphi _{s}(z):\varphi 
_{s}^{-1}(\zeta _{s})\to \varphi _{s'}^{-1}(\zeta _{s'}).$$
The definition of the model metric in Section \ref{7} ensures that 
this homotopy is boundedly Lipschitz restricted to $M_{1}$. 
We now change $\partial T(\zeta )$ on horizontal boundary pieces 
which 
are matched by this homotopy in $M_{1}$, and on the homotopy between 
them, so as to reduce length. We do not change $T(\zeta )$ 
 in $M_{2}$, and do not need to, because 
we already have bounds there. The topological type of the surface 
$T(\zeta )$ will not change, but the isotopy class of the embedded 
surface probably will. If $\gamma _{1}$ and $\gamma _{2}$ are 
horizontal pieces matched by the homotopy, with the homotopy entirely 
in $M_{1}$, we simply replace $\gamma 
_{j}$ by a bounded path $\gamma _{j}'$ with the same endpoints,
such that the natural homotopy in $M_{1}$ matches up $\gamma _{1}'$ 
and $\gamma _{2}'$ instead. Now there is a decomposition of the 
horizontal boundary of $T(\zeta )$ into pieces which are paired by 
the 
homotopy, because every time a homotopy track enters $M_{1}$, it does 
so along horizontal boundary and then leaves again. The path $\varphi 
_{t}^{-1}(\zeta )$ has bounded length so need not be involved in any 
changes.  So performing 
this operation a finite number of times we have a surface $S(\zeta )$ 
of bounded area and with  of bounded length in the metric induced 
from the model metric. By construction, 
the components of $\partial S(\zeta )$ in the horizontal boundary are 
homotopically nontrivial. The claim about $A(\zeta )$ follows, 
because $S(\zeta )\cap M_{2}=A_{1}(\zeta )\cap M_{2}$ is homotopic to 
$\zeta $.\Box

\ssubsection{Bounding horizontal length: intersection 
number.}\label{11.10}

Suppose we have fixed a set of loops $\zeta $,  such that the 
corresponding annuli cut the corresponding annuli into cells, and we 
have a bound on the length of the loops $\zeta $ used. Then the 
surfaces $S(\zeta )$, and the horizontal surfaces used, cut $M_{1}$ 
into cells. Let the set of surfaces be $S_{j}$, $1\leq j\leq r$.  
Then 
for a constant $C_{3}>0$ depending only on the ltd parameter functions
$$\vert \gamma \vert \geq C_{3}{\rm{Max}}\{ i(\gamma ,S_{j}):1\leq 
j\leq r\} ,$$
where $i(.,.)$ denotes essential intersection number. The aim now is 
to obtain a parallel lower bound on $\vert \gamma _{1}\vert $, and 
hence an upper bound on $\vert \gamma \vert $ in terms of $\vert 
\gamma _{1}\vert $. For this, we have to notice something about the 
paths $\zeta $, and have to be a little careful about the choices. 

First, we can easily find a surface homotopic to $\Phi (S_{j})$ with 
boundary in $\Phi (\partial M_{1})$, by choosing a maximal multicurve 
in $S_{j}$ which is 
bounded in the model metric and using the pleated surface with this 
int its pleating locus, for 
example. We assume this surface is in fact $\Phi (S_{j})$. If $S_{j}$ 
is a horizontal surface then, as in \ref{11.4}, we can redefine 
$\Phi $ inside Margulis tubes so that $\Phi (S_{j})\subset N_{\geq 
\varepsilon _{1}}$ and $\Phi (\partial S_{j})\subset \partial 
N_{\geq 
\varepsilon _{1}}$. Then if $T_{1}$, $T_{2}$ denotes the Margulis 
tubes intersected by $\partial S_{j}$, $\Phi (\partial S_{j})$, 
$(M,T_{1},S_{j},\gamma )$ is homotopy equivalent to $(N,T_{2},
\Phi (S_{j}),\gamma _{1})$ and so
$$i(\gamma ,S_{j})=i(\gamma _{1},\Phi (S_{j}))$$
for any horizontal surface. It remains to obtain something similar 
for $S(\zeta )$. The key result is the following.

\begin{lemma}\label{11.11} There is a constant $L$ such that the 
following holds.
There is a  submanifold with boundary $N_{1}$ of $N$ and a 
continuous map $\Phi _{1}:M\to N$ homotopic to $\Phi $ such that the 
following hold. $\Phi $ is a homeomorphism in a neighbourhood of $M$ 
with $\Phi _{1}(\partial M_{1})=\partial N_{1}$, $\Phi _{1}(S(\zeta 
))\subset N_{1}$  and $\Phi $ is homotopic to $\Phi _{1}$ under a 
homotopy which is trivial on $M_{1}'$
.\end{lemma}

Note that there is no claim of a bound on the homotopy distance 
between $\Phi $ and $\Phi _{1}$. 

\noindent {\em{Proof.}} 

We use the set-up of \ref{11.8}, in particular, (\ref{11.8.1}) and 
(\ref{11.8.2}).  By \cite{F-H-S}, there is a surface $S'(b_{k,\alpha 
_{i}})$ (or $S'(b_{k,\alpha _{1},j})$) in an arbitrarily small 
neighbourhood of $\Phi (S(b_{k,\alpha _{i}}))$ (or $\Phi 
(S(b_{k,\alpha 
_{1},j}))$ if $i=1$) which is homeomorphic to 
$\alpha _{i}$.
As in the proof of \ref{11.5}, we can also assume that these 
 surfaces intersect 
 $\partial T(\zeta _{*},\varepsilon _{1})$ transversally for $\zeta 
 \subset \partial \alpha _{i}$,
 and only in the homotopy class of the core loop $\zeta $. We take 
 $N_{1}$ to be  the set homotopic to $\Phi (M_{1})$  bounded by 
 the $\varepsilon _{1}$-model Margulis tubes 
 which intersect $\Phi (\partial _{v}M_{1})$ and  the surfaces 
 $S'(b_{\pm k_{2},\alpha 
_{i}})$, $i\geq 2$, and $S'(b_{0,\alpha _{1},j})$, $j=1$, $2$. 
Since $\Phi $ is $\Lambda _{1}$-Lipschitz, we can assume (by choice 
of $p_{0}$) that the image under $\Phi $ of the intersection of 
$M_{1}$ with the region bounded by $S(b_{k,\alpha _{1},j})$ and 
$S(b_{\pm k,\alpha _{i}})$ ($i\geq 2$) is contained in the region 
bounded by $S'(b_{k-1,\alpha _{1},j})$ and 
$S'(b_{\pm (k-1),\alpha _{i}})$ for all $1\leq k\leq k_{1}$. In 
particular, this is true for the image under $S(\zeta )$ of the 
intersection of $S(\zeta )$ with this region. Similar facts hold 
for images under $\Phi $ of interval bundles bounded by 
$S(b_{k,\alpha _{i}})$ and $S(b_{m,\alpha _{i}})$, or by 
$S(b_{k,\alpha _{1},j})$ and $S(b_{m,\alpha _{1},j})$. So we only 
need to 
change the definition of $\Phi $ in the region bounded by 
 $S(b_{\pm 2,\alpha _{1},j})$ and 
$S(b_{\pm (k_{2}\pm 2),\alpha _{i}})$. We could, for example, first 
extend $\Phi (S(\zeta )$ by attaching annuli and put the boundaries 
in 
$S'(b_{-1,\alpha _{1},j})$, $S'(b_{\pm (k_{2}-1),\alpha _{i}})$, and 
then compose with homeomorphisms to compress the interval bundles 
bounded by $S'(b_{\pm 2,\alpha _{1},j})$, $S'(b_{\pm (k_{2}\pm 
2),\alpha _{i}})$
to the interval bundles bounded by $S'(b_{2,\alpha _{1},j})$, 
$S(b_{0,\alpha _{1},j})$ $S(b_{\pm (k_{2}+2),\alpha _{i}})$, 
$S'(b_{\pm (k_{2}),\alpha _{i}})$. The resulting map $\Phi _{1}$ is 
homotopic to $\Phi $ via a homotopy which is trivial outside the 
region bounded by $S(b_{\pm 2,\alpha _{1},j})$ and 
$S(b_{\pm (k_{2}\pm 2),\alpha _{i}})$.
\Box

Now since $\Phi _{1}$ is homotopic to a homeomorphism from $M_{1}$ 
to $N_{1}$,
$$i(\gamma _{1},\Phi (S(\zeta ))=i(\gamma _{1},\Phi _{1}(S(\zeta ))
=i(\Phi _{1}(\gamma ),\Phi _{1}(S(\zeta ))=i(\gamma ,S(\zeta )).$$

Let $M_{2}$ be any interval bundle in $M_{1}$ which contains 
 $\zeta $ and in which $\vert \zeta \vert $ is bounded.
Let $M_{2}'$ be the natural compactification of the 
cover of $M_{1}$ corresponding to $M_{2}$. By this we mean: $M_{1}$ 
is 
a sum of interval bundles, attached along subsurfaces. The universal 
cover of an interval bundle has natural compactification $D\times 
[0,1]$, where $D$ is the closed unit disc. Then (unless $M_{1}=M_{2}$ 
up to homotopy) $M_{2}'$ is an 
infinite tree of $M_{2}$ and  infinitely many copies of 
${\rm{int}}(D/\Gamma _{j})\times [0,1]$, for subgroups $\Gamma _{j}$ 
of $\pi _{1}(M_{2})$ (mostly for $\Gamma _{j}$ the 
trivial group) attached along subsurfaces. Compactify by adding the 
copies of $\partial D$ and points at each end. Then $M_{2}'$ is 
naturally homeomorphic to $M_{2}$.
Let $\gamma _{2}$ be the lift of $\gamma $ to $M_{2}'$. Then $\gamma 
_{2}$ has infinitely many components.
One component is a closed path if $\gamma \subset M$, but otherwise 
all components are paths 
with endpoints in $\partial M_{2}'$. The number of homotopically 
nontrivial paths, that is, not homotopic into the boundary, is the 
number of components of intersection of $\gamma $ with $M_{2}$.
As already 
noted, we have lift of $S(\zeta )$ to $A(\zeta )$ in $M_{2}'$ which 
extends to the boundary, and $(M_{2}',A(\zeta ))$ is homeomorphic to 
$(M_{2},A_{2}(\zeta )$, where $A_{2}(\zeta )$ is an annulus in 
$M_{2}$ 
with boundary components in $\partial M_{2}$. Now clearly we have
$$i(\gamma ,S(\zeta ))\geq i(M_{2}',A(\zeta ))=i(M_{2},A_{2}(\zeta ))
=i(\gamma _{2},\zeta ).$$
For $i(\gamma _{2},\zeta )$, we regard $\gamma _{2}$ as an arc in 
$\partial M_{2}$ and $\zeta $ as being in the same boundary 
component. Each intersection gives a different element of the coset 
space $\pi 
_{1}(M_{2})/<\zeta >$.
Each essential intersection between $\gamma _{1}$ and $S(\zeta )$ 
gives rise to an actual intersection, and 
a homotopy class of path along $\gamma _{1}$ from an initial point to
 another point on $\Phi _{1}(S)\cap \gamma _{1}$, that is, from an 
 initial point to a point on $\Phi (S)$. The intersections
  given by
 cosets in $\pi _{1}(S)/\pi _{1}(S(\zeta ))$   are distinct. So 
 the bound on the length of $\gamma _{1}$ and on the area of 
 $\Phi (S(\zeta ))$ gives a bound on $i(\gamma 
 ,\zeta )$. We have this for all $\zeta $. Since we have already 
 seen that $i(\gamma ,S_{j})$ is bounded for horizontal surfaces, we 
 have the required bound on $\vert \gamma \vert $

\ssubsection{Bounding horizontal length near $W'$.}\label{11.12}

As usual, let $W'$ denote the non-interval bundle part of the core 
$W$. 
Let $\gamma $ be a path near $W$ and $\gamma 
_{1}$ 
the geodesic segment homotmopic to $\Phi (\gamma )$ with the same 
endpoints, as in \ref{11.8}. Choose a model piece $M_{1}$ such that 
$\gamma \subset M_{1}$ and 
$\partial _{h}M_{1}$ has similar properties to those described in 
\ref{11.8}, but containing $W'$ and intersecting the model manifolds 
of several ends. Choose the model piece large enough to include 
anything 
which was not previously included in a model piece $M_{1}$ as in 
\ref{11.8}. So $W'$ may be properly contained in $M_{1}$. 
Again, we assume that $\gamma $ and $\gamma _{1}$ are closed 
loops.  

Now we construct surfaces in $W'$ which decompose $W'$ into 
cells, and with boundary in $\partial W'$.
We look at intersections with the $W_{i}$ of \ref{6.11}. 
Recall that $W=W_{n}$, and $W_{i}$ decomposes into submanifolds 
$W_{i-1}$.
 In \ref{6.11} we constructed a decomposition of $W$ into 
submanifolds 
which were balls or interval bundles. Interval bundles adjacent to 
$\partial W$ have been removed but there might be some interval 
bundles in the interior of $W'$. We use the annuli in the sets 
$\Sigma _{i}$ of \ref{6.14} to extend the decomposition into balls.  
So 
now we have an extended sequence $W_{j}$, with $W_{n}=W'$ (renaming 
$n$ of \ref{6.11}) such that the components of 
$W_{j}$ decompose into the components of $W_{j-1}$, all components 
of $W_{0}$ are balls, all components of $W_{r}$ (for some $r$) are 
either balls or interval bundles, and all components of $\partial 
W_{j-1}\setminus \partial W_{j}$, for $j\leq r$, are annuli.
Now we want to apply the same  method above with 
$M_{1}=W$ and suitable surfaces $S_{i}$ obtained from the components 
of the $\partial W_{j}$. Instead of using components of $\partial 
W_{i}\setminus 
\partial W_{i+1}$, for $i\leq n-2$, we use surfaces with boundary in 
$\partial W'$. To do this, if $S$ is a component of $\partial 
W_{i}\setminus W_{i+1}$ and $i\leq n-2$, we form  a union of two 
surfaces, 
$T_{i+2}(S)$, with boundary in 
$\partial W_{i+2}$, by taking a tubular 
neighbourhood of $S$, taking the boundary components of the tubular 
neighbourhood, adding surfaces in the boundary of tubular 
neighbourhoods of the adjoining components of $\partial 
W_{i+1}\setminus \partial W_{i+2}$. Inductively, if $j<n$ we define 
$T_{j+1}(S)=T_{j+1}(T_{j}(S))$. We 
write $T(S)=T_{n}(S)$, which is a union of $2^{n-i-1}$ surfaces, 
since $S$ is two-sided and all $W_{i}$ have oriented boundary. 

 Now we extend each surface $T(S)$ across the remainder of $M_{1}$, 
using the 
surfaces of \ref{11.9}. We have bounds on the areas of the 
surfaces, as before. If these surfaces do not cut $M_{1}$ into cells, 
then we also use some additional surfaces $S(\zeta )$. We may need to 
use some extra surfaces $S(\zeta )$ for $\zeta =\zeta _{1}\cup \zeta 
_{2}$ of the following form: $\zeta \subset \alpha $ for some 
$M(\alpha ,\ell )\subset M_{1}$, $\zeta _{1}\subset \partial T(S)$ 
for some surface $T(S)$ constructed in 
$W'$ and $\zeta _{2}\subset \cup _{i=1}^{r}\alpha _{i}$, where 
$\alpha 
_{i}\times \ell _{i}\subset S(e)\times 
[z_{e,0},y_{e,+}]$, $(\alpha _{i},\ell _{i})<(\alpha ,\ell )$ and 
$M_{1}\cap M(\alpha _{i},\ell _{i})=\emptyset $. We can choose a set 
of such $\zeta $ such that the corresponding sets $S(\zeta )$, 
together with the sets $T(S)$, cut $M_{1}$ into 
cells.  
\Box

\section{Proofs of Theorems \ref{1.1}, \ref{1.2}, 
\ref{1.3}.}\label{12}
As noted in Section \ref{1}, Theorems \ref{1.1} and \ref{1.2} are 
derived 
from \ref{1.3}. 

\ssubsection{Proof of \ref{1.3}.}\label{12.1}
Item 1 of \ref{1.3} was proved in \ref{6.13}, since 
we only take limits on the end parts of the geometric manifolds  - 
and 
in \ref{6.12} in the combinatorially bounded geometry Kleinian 
surface 
case. The existence of a coarse biLipschitz $\Phi :M\to \overline{N}$ 
is 
entirely proved in sections \ref{10} and \ref{11}, with the required 
bounds on 
the biLipschitz constant, in the case of all 
ends being geometrically infinite. In the case of a geometrically 
finite end, we need to show how to extend the map $\Phi $ to the 
complement of the convex hull. For each 
geometrically finite end $e_{i}$ of $N_{d}$, we choose an 
embedding  $f^{i}:S(e_{i})\to N$ as $f_{3}$ in  
\ref{3.9},and the  pleated surface $f_{e_{i},+}$  as $f_{4}$ of 
\ref{3.9}. By \ref{3.9}, we therefore have bounds on 
$d(\mu (e_{i}),[f_{e_{i},+}])$, and a bound on the hyperbolic length 
of 
homotopy tracks of a homotopy between $f^{i}$ and $f_{e_{i},+}$.
 So now if 
$M=M(\mu (e_{1}),\cdots \mu (e_{n}))$ we have, by the results of 
sections 
\ref{10} and \ref{11}, a coarse biLipschitz map $\Phi :M\to N$ whose 
image  has any surfaces $f^{i}(S(e_{i}))$ in its boundary. The 
surfaces $f^{i}(S(e_{i}))$ bound convex sets containing the convex 
hull of $N$, by 
construction. Perturb $f^{i}$ so that $f^{i}$ is smooth and the 
complementary components are still convex. Then perpendicular 
geodesics pointing outwards from $f^{i}(S(e_{i}))$ hit the boundary 
of the $t$-neighbourhood in exactly one point. This defines a map 
from $f^{i}(S(e_{i}))$ to the $t$-neighbourhood which expands 
distances by $e^{t}$. So we can adjust $\Phi $ to a map onto $N\cup 
\Omega (\Gamma )/\Gamma $ which is coarse biLipschitz with respect to 
the adjusted metric, and maps a set $S(e_{i})\times [u_{i}-1,u_{i}]$ 
in the 
model manifold  to the corresponding component of $N\setminus CH(N)$, 
with sets $\{ x\} \times [u_{i}-1,u_{i}]$ mapping to geodesic rays, 
and mapping $S(e_{i})\times \{ u_{i}-1+u'\} $ to the set distance 
$t(u')$ from 
$f^{i}(S(e_{i})$ for some fixed function $t:[0,1]\to (0,\infty )$. In 
particular, $\Phi $ is coarse biLipschitz with 
respect to the hyperbolic metric on the preimage of the convex subset 
of $N$ bounded by the surfaces $f^{i}(S(e_{i}))$.

\ssubsection{Proof of \ref{1.1}.}\label{12.2}

Let $N_{1}$ and $N_{2}$ be two homeomorphic hyperbolic $3$-manifods 
with finitely generated fundamental groups and the same end 
invariants. Let $M$ be the model manifold for both. 
Let $\Phi _{i}:M\to N\cup \Omega (\Gamma _{i})/ \Gamma _{i}$
be the  map of \ref{1.3}, constructed in \ref{12.1}, with lifts 
$\widetilde{\Phi _{1}}$ and $\widetilde{\Phi _{2}}$. Then 
the set-valued map $\tilde {\Phi }=\widetilde{\Phi _{2}}\circ 
\widetilde{\Phi 
_{1}}^{-1}:H^{3}\to H^{3}$ is coarse biLipschitz with respect to the 
metric which agrees with the hyperbolic metric inside the convex 
hulls, except in bounded (hyperbolic metric) neighbourhoods of the 
convex hull boundaries, and is an adjustment of the hyperbolic metric 
outside. Let $f^{j,i}:S(e_{j})\to N_{i}$ be as in \ref{12.1}, and 
$U_{i}$ be the lift of 
$\cup _{j}f^{j,i}(S(e_{j}))$ to $H^{3}$.
 Then $U_{i}\subset 
H^{3}\setminus \widetilde{CH(N_{i})}$ and (although we do not need 
this) 
$U_{i}$ is a bounded 
hyperbolic distance from a component of $\widetilde{CH(N_{i})}$. 
Taking a 
small perturbation, we can assume that $U_{1}$ and $U_{2}$ have 
smooth boundaries, and let 
$U_{i,t}$ be the surface hyperbolic distance $t$ from $U_{i}$, 
further from $\widetilde{CH(N_{i})}$, as in \ref{12.1}. The component 
of 
$H^{3}\setminus 
U_{i,t}$ containing $\widetilde{CH(N_{i})}$ is convex. So 
perpendicular 
geodesic segments between $U_{i}$ and $U_{i,t}$ only intersect these 
surfaces at the endpoints.  Now we have chosen $\widetilde{\Phi 
_{1}}$ and 
$\widetilde{\Phi _{2}}$ so that $\widetilde{\Phi }$ maps $U_{1}$ to 
$U_{2}$, 
perpendicular segments to perpendicular segments preserving length, 
$U_{1,t}$ to $U_{2,t}$ and $\Omega (\Gamma _{1})$ to $\Omega (\Gamma 
_{2})$.

Now we claim that 
$\tilde{\Phi }$ is coarse biLipschitz with respect to the hyperbolic 
metric. It suffices to show that $\tilde{\Phi }$ and $\widetilde{\Phi 
^{-1}}$ are coarse Lipschitz and we only need to consider the maps 
on $H^{3}\setminus \widetilde{CH(N_{i})}$, $i=1$, $2$. Here, we only 
need 
to consider the norms of the derivatives with respect to the 
hyperbolic metric. We have a splitting of the tangent space at each 
point into the tangent to the geodesic segment and the tangent to the 
$U_{i,t}$ foliation. The derivative of $\tilde{\Phi}$ maps one 
splitting to the other and $\widetilde{\Phi _{i}}$ maps the tangent 
space  to 
$S(e_{j})\times \{ u_{j}-1+u'\} $ to the tangent space to 
$U_{i,t(u')}$, with derivative 
$e^{t(u')}$ times the identity.
Meanwhile, the tangent space to  $\{ x\} \times [u_{j-1},u_{j}]$ is 
mapped to the geodesic tangent space with the same derivative for 
both $\widetilde{\Phi _{1}}$ and $\widetilde{\Phi _{2}}$. So the 
derivative 
of  $\tilde{\Phi }$ is actually the identity with respect to suitable 
framings on the domain and range tangent spaces, and $\tilde{\Phi }$ 
is Lipschitz, as required.

Let $x_{t}$ be a geodesic ray in $H^{3}$ from a basepoint $x_{0}$, 
parametrised by length, and converging to a point $x_{\infty }\in 
\partial H^{3}$. First we claim that $\lim _{t\to \infty } 
\tilde{\Phi}(x_{t})$ exists and is a single point, remembering that 
$\tilde{\Phi}(x_{t})$ is set-valued in general. In fact, 
$\tilde{\Phi}(x_{t})$ has uniformly bounded diameter in the 
hyperbolic metric, because $\tilde{\Phi}$ is coarse biLipschitz. Then 
the limit exists, 
because the coarse biLipschitz properties of $\Phi _{i}$ with respect 
to hyperbolic metric imply that $\tilde{\Phi }(\{ x_{t}:t\in 
[n,n+1]\} )$ has bounded hyperbolic diameter and geometrically small 
Euclidean diameter, because the distance of the set from a fixed 
basepoint is $\geq n/\Lambda $, where $\Lambda $ is the product of 
the 
biLipschitz constants for $\Phi _{1}$ and $\Phi _{2}$. Moreover, this 
also shows that the limit is uniform. Since $\tilde{\Phi}$ is coarse 
biLipschitz with respect to the hyperbolic metric, it is biLipschitz 
with respect to the Euclidean metric on $\partial H^{3}$, using the 
unit ball model for $H^{3}\cup \partial H^{3}$, and hence 
quasi-conformal.
Then we can change the conformal structure on $\Omega (\Gamma _{2})$ 
so that with respect to this new conformal structure, $\tilde{\Phi }$ 
is conformal, because the ending invariants of $N_{1}$, $N_{2}$ are 
the same. So then we can find a quasiconformal map $\psi :\Omega 
(\Gamma _{2})\to \Omega (\Gamma _{2})$ which extends continuously to 
the identity on the complement in $\partial H^{3}$, so that $\psi 
\circ \tilde{\Phi }$ is conformal on $\Omega (\Gamma _{1})$. It is 
equivariant, and hence if not conformal on the limit set, there is a 
nontrivial $\Gamma _{2}$-invariant line-field on the limit set, 
where $\Gamma _{2}$ is the covering group of $N_{2}$. So then by 
Sullivan's theorem \cite{Sull} the extension of  $\psi \circ 
\tilde{\Phi }$ is conformal, and 
$N_{1}$ and $N_{2}$ are isometric.

\ssubsection{Proof of \ref{1.2}.}\label{12.3}

Let $[\mu _{1},\cdots \mu _{r}]$ be any admissible invariant for the 
topological type of a hyperbolic manifiold $N$ with ends $e_{i}$ for 
$N_{d}$, $1\leq i\leq r$. Let 
$[y_{1,n},\cdots y_{r,n}]\in \left( \prod _{i=1}^{r}{\cal 
T}(S(e_{i})\right) /G$ with 
$y_{i,n}\to \mu _{i}$. Then by \ref{6.13}, $M_{n}=M(y_{1,n},\cdots 
y_{r,n})$ 
converges geometrically to $M(\mu _{1},\cdots \mu _{r})$ with a 
suitable choice of basepoint. By \ref{1.3}, if $N_{n}$ denotes the 
hyperbolic manifold homeomorphic to $N$ under an end-preserving 
homeomorphism between $N_{n,d}$ and $N_{d}$ then the map $M_{n}\to 
N_{n}$ is coarse biLipschitz restricted to the preimage of $C(N_{n})$ 
for a biLipschitz constant uniform in $n$. So taking geometric 
limits, 
the $N_{n}$ converge, with suitable basepoint, to a hyperbolic 
manifold 
$N_{\infty}$. The topological type remains the same, because we can 
find $m_{n}\to \infty $ and a sequence of $m_{n}$ successive pleated 
surfaces in each end of $N_{n}$ with hyperbolic distance $\geq 1$ 
between any adjacent pair. Also using these surfaces, we 
have  closed geodesics  arbitrarily 
far out in the end $e_{i}$ converging to $\mu _{i}$ in the limit. So 
$(\mu _{1}, 
\cdots \mu _{r})$ is the invariant of $N_{\infty }$.


\begin{thebibliography} {99}
     
\bibitem{Abi}  Abikoff, W.: The real analytic theory 
of Teichm\"uller Space. Springer Lecture Notes in 
Math. 820 (1980). 

\bibitem{Ag} Agol, I.: Tameness of hyperbolic 3-manifolds, 
arXiv:math.GT/0405568.

\bibitem{A-B}  Ahlfors, L. and Bers, L.: The Riemann 
Mapping Theorem for Variable Metrics, Ann. of 
Math. 72b (1960) 385-404.

\bibitem{Bon2} Bonahon, F.: Cobordism of automorphisms of surfaces,
Annales Scientifiques de
l'ƒcole Normale SupŽrieure SŽr. 4, 16 no. 2 (1983), p. 237-270. 

\bibitem{Bon} Bonahon, F.: Bouts des Vari\'et\'es hyperboliques de 
dimension trois, Ann. of Math. 124 (1986) 71-158.

\bibitem{Bow1} Bowditch, B.: Intersection numbers, and the 
hyperbolicity of the curve complex, to appear in J. reine angew. 
math..

\bibitem{Bow2} Bowditch, B.: Geometric models for hyperbolic 
manifolds. Preprint Southampton 2005.

\bibitem{Bow3} Bowditch, B.: End invariants of hyperbolic 
3-manifolds. Preprint, Southampton 2005.

\bibitem{B-B} Brock, J. and Bromberg, K.: On the density of 
geometrically finite Kleinian groups. Acta Math 192 (2004) 33-93.

\bibitem{B-C-M} Brock, J., Canary, R., Minsky, Y.: 
The classification of Kleinian surface groups II: the 
Ending Lamination Conjecture, arXiv:math.GT/0412006.

\bibitem{Brom} Bromberg, K.: Projective structures wit degenerate 
holonomy and the Bers density conjecture, arXiv:math.GT/0211402.

\bibitem{B-S} Bromberg, K. and Souto, J.: The density conjecture: A 
prehistoric approach, in preparation.

\bibitem{Bus} Buser, P.: Geometry and Spectra of Compact Riemann 
Surfaces. Birkhouser 1992.

\bibitem{C-G} Calegari, D. and Gabai, D.: Shrinkwrapping and the 
taming of hyperbolic 3-manifolds, J. Amer. Math. Soc. 19 (2006) 
385-446.

\bibitem{C-E-G} Canary, R.D., Epstein, D.B.A. and Green, P.: Notes on 
notes of Thurston, Analytical and Geometric Aspects of Hyperbolic 
Space, C.U.P. 1987, L.M.S. Lecture Notes Series no 111, 3-92.

\bibitem{Can} Canary, R.D.: Ends of hyperbolic 3-manifolds, J. Amer. 
Math. Soc., 6 (1993) 1-35.

\bibitem{C-F} Cannon, J.W. and Feustel, C.D.: Essential embedding of annuli and M\"obius bands in 3-manifolds, TAMS 215 (1976), 219-239

\bibitem{E-M} Epstein, D.B.A. and Marden, A.: Convex hulls in 
hyperbolic space, a theorem of Sullivan, and measured pleated 
surfaces, Analytical and Geometric Aspects of hyperbolic space, CUP 
1987, 
London 
Math. Soc. Lecture Notes Series no. 111, 113-254.

\bibitem{Ev} Evans, R.: The ending lamination conjecture for 
hyperbolic 3-manifolds with slender 
end-invariants, to appear in Pac. J. Math.. 

\bibitem{F-L-P} Fathi, A., Laudenbach, F. and Po\' enaru, V. et al. 
Travaux de Thurston sur les surfaces, Ast\' erisque 66-67 (1979).

\bibitem{F-H-S} Freedman, M., Hass, J. and Scott, P.: Least area 
incompressible surfaces in 3-manifolds, Invent. Math. 71 (1983), 
609-642.

\bibitem{H-T} Hatcher, A. and Thurston, W.P.: A presentation of the 
mapping class group of an orientable surface, Topology 19 (1980) 
221-237.

\bibitem{Hem} Hempel, J.: 3-Manifolds. Annals of Math. Studies, 
Princeton University Press, 1976.

\bibitem{K-L-O} Kim, I., Lecuire, C. and Ohshika, K.: Convergence of 
freely decomposable Kleinian groups. Preprint, 2004

\bibitem{Le1} Lecuire, C.: An extension of Masur domain, preprint, 
2004.

\bibitem{Mas}, Masur, H.: Measured Foliations and handlebodies, 
Ergodic 
Theory and Dynamical Systems 6 (1986), 99-116.

\bibitem{M-M1} Masur, H. and Minsky, Y.: Geometry of the Complex of 
Curves I, Hyperbolicity, Invent. Math. 138 (1999) 103-149.

\bibitem{M-M2} Masur, H. and Minsky, Y.: Geometry of the Complex of 
Curves II, Hierarchical Structure, Geom. Funct. Anal. 10 (2000) 
902-974.

\bibitem{McC} McCullough, D.: Compact submanifolds of 3-manifolds 
with 
boundary, Quart. J. Math. Oxford Ser. 37(147) (1986) 299-307.

\bibitem{McC-Mil} McCullough, D. and Miller, A.: Homeomorphisms of 
3-manifolds with Compressible Boundary, Memoirs of the AMS 344 
(1986).

\bibitem{McC-M-S} McCullough, D., Miller, A. and Swarup, G.A.: 
Uniqueness of cores of noncompact 3-manifolds, J. London Math. Soc. 
52 (1985) 548-556.

\bibitem{McM} McMullen, C.: Iteration on Teichm\"uller Space. Invent. 
Math. 99 (1990) 425-454.

\bibitem{Min3} Minsky, Y.: Teichm\" uller geodesics and ends of 
hyperbolic $3$-manifolds, Topology 32 (1993) 624-647.

\bibitem{Min4} Minsky, Y.: Harmonic  maps into hyperbolic 
3-manifolds, 
Trans Amer. Math. Soc. 332 (1992) 607-632.

\bibitem{Min5} Minsky, Y.: Quasi-projections in Teichm\"uller Space, 
J. Reine Angew. Math. 473 (1996), 121-136.

\bibitem{Min0} Minsky, Y.: The classification of punctured 
torus-groups, Annals of Math. 149 (1999) 559-626.

\bibitem{Min1} Minsky, Y.: Bounded Geometry for Kleinian Groups, 
Invent. Math. 146 (2001) 143-192.

\bibitem{Min2} Minsky, Y.: The classification of Kleinian Surface 
Groups, I: 
Models and Bounds, preprint 2002. arXiv:math.GT/0302208.

\bibitem{Mor} Morgan, J.W.: On Thurston's Uniformisation Theorem for 
Three-Dimensional Manifolds, Chapter V, The Smith Conjecture, edited 
by J.W. Morgan and H. Bass, Academic Press 1984.

\bibitem{Nam} Namazi, Hossein: Heegard splittings and hyperbolic 
geometry. Thesis, Stony Brook University, 2005.

\bibitem{Oh3} Ohshika, K.: On limits of quasi-conformal deformations 
of 
Kleinian groups. Math.Z. 201 (1989) 167-176.

\bibitem{Oh1} Ohshika, K.: Limits of geometrically tame Kleinian 
groups, 
Invent. Math. 99 (1990) 185-203.

\bibitem{Oh2} Ohshika, K.: Ending laminations and boundaries for 
deformation spaces of Kleinian groups, J. London Math. Soc. 42 
(1990), 111-121.

\bibitem{Oh4} Ohshika, K.: Realising end invariants by limits of 
minimally parabolic, geometrically finite groups. Preprint 2005. 
arXiv:math.GT/0504546.

\bibitem{Ot1} Otal, J-P.: Courant g\'eodesiques et produits libres, 
Th\` 
ese d\'etat, Universit\'e Paris-Sud, Orsay 1988.

\bibitem{Ot2} Otal, J-P.: Thurston's hyperbolization of Haken 
manifolds, Surveys of Differential Geometry III (1986) Int. Press.

\bibitem{Ot3} Otal, J-P.: Les g\'eod\'esiques ferm\'ees d'une 
vari\'et\'e hyperbolique en tant que noeuds, 95-104, Kleinian Groups 
and Hyperbolic 3-Manifolds, Proceedings of the Warwick workshop,
September 2001, eds Y. Komori, V. Markovic, C. Series, Cambridge U.P.
2003.

\bibitem{Raf} Rafi, K.: Hyperbolic 
3-manifolds and Geodesics in Teichm\"uller Space. 
Thesis, SUNY at Stony Brook 2001.

\bibitem{Raf2} Rafi, K.: A characterization of short curves of a 
Teichm\" uller geodesic. arXiv:math.GT/0404227v2 (2004). 

\bibitem{R1} Rees, M.: Views of Parameter 
Space: Topographer and Resident.    
Ast\'erisque 288 (2003). 

\bibitem{R3} Rees, M.: An alternative approach to the ergodic theory 
of measured foliations on surfaces, Ergod. Th. and Dynam Sys. 1 
(1981) 
461-488.

\bibitem{Sco1} Scott, G.P.: Finitely generated 3-manifold groups are 
finitely presented, J. London Math. Soc. 6 (1973) 437-440.

\bibitem{Sco2} Scott, G.P.: Compact submanifold of 3-manifolds, J. 
London Math. Soc. 7 (1973) 246-250.

\bibitem{Som} Soma, T:  Existence of polygnal wrapping in hyperbolic 
3-manifolds.  Preprint, Tokyo Denki (2005).

\bibitem{Sou} Souto, J.: A note on the tameness of compact 
3-manifolds, Topology 44 (2005).

\bibitem{Sul2} Sullivan, D.P.: A finiteness theorem for cusps. Acta 
Math. 147 (1981) 289-299. 

\bibitem{Sull} Sullivan, D. P.: Quasiconformal homeomorphisms and 
dyanmics  II. Structural Stability implies hyperbolicity for Kleinian 
groups. Acta Math 155 (1985) 243-260.

\bibitem{T} Thurston, W.T.: The Geometry and Topology of 3-manifolds, 
notes, 1979.

\bibitem{T1} Thurston, W.T.: Hyperbolic structures on 3-manifolds, 
I: deformation of acylindrical manifolds. Ann. of Math. 124 (1986) 
203-246.

\bibitem{T2} Thurston, W.T.: Hyperbolic structures on 3-manifolds, 
II: Surface groups 
and 3-manifolds which fiber over the circle. 
arXiv:math.GT/9801045 (1998).

\bibitem{T3} Thurston, W.T.: Hyperbolic structures on 3-manifolds, 
III: Deformations of 3-manifolds with incompressible boundary.
arXiv:math.GT/9801058 (1998).

\bibitem{Wal2} Waldhausen, F.: Eine Verallgemeinering des Schleitfestes. Topology 6 (1967), 501-504.

\bibitem{Wal} Waldhousen, F.: On irreducible 3-manifolds which are 
sufficiently large, Ann. of Math. 87 (1968) 56-88.

\end{thebibliography}
 \end{document}